\documentclass[dvips,12pt,twoside,openright]{report}

\usepackage{color}

\usepackage{hyperref}

\usepackage{amssymb}
\usepackage{amsmath}
\usepackage{amsthm}
\usepackage[english,catalan]{babel}
\usepackage[latin1]{inputenc}
\usepackage{mathrsfs}
\usepackage[all]{xy}

\usepackage{fancyhdr}

\makeatletter
\def\cleardoublepage{\clearpage\if@twoside \ifodd\c@page\else
\hbox{}
\thispagestyle{empty}
\newpage
\if@twocolumn\hbox{}\newpage\fi\fi\fi}
\makeatother

\newcommand{\Aut}{\ensuremath{\mathit{Aut}}}

\newcommand{\Ext}{\ensuremath{\mathit{Ext}}}

\newcommand{\F}{\ensuremath{\mathbf{F}}}

\newcommand{\h}[1]{\-\mbox{-#1}}
\newcommand{\Hom}{\ensuremath{\mathit{Hom}}}
\newcommand{\I}{\ensuremath{\mathbf{I}}}

\newcommand{\kk}{\ensuremath{K}}
\newcommand{\nn}{\ensuremath{\mathbb{N}}}

\newcommand{\OO}{\ensuremath{\mathcal{O}}}

\newcommand{\pp}{\ensuremath{\mathbb{P}}}
\newcommand{\qq}{\ensuremath{\mathbb{Q}}}

\newcommand{\rr}{\ensuremath{\mathbb{R}}}

\newcommand{\syz}{\ensuremath{\mathrm{Syz}}}

\newcommand{\zz}{\ensuremath{\mathbb{Z}}}

\newtheorem{te}{Theorem}[chapter]
\newtheorem*{te*}{Theorem}
\newtheorem*{tecat*}{Teorema}
\newtheorem{p}[te]{Proposition}
\newtheorem*{p*}{Proposition}
\newtheorem*{pcat*}{Proposició}
\newtheorem{co}[te]{Corollary}
\newtheorem*{co*}{Corollary}
\newtheorem*{cocat*}{Coro{\lgem}ari}

\newtheorem{lem}[te]{Lemma}
\newtheorem*{lem*}{Lemma}

\theoremstyle{definition}
\newtheorem{defin}[te]{Definition}
\newtheorem{ex}[te]{Example}
\newtheorem{pr}[te]{Problem}
\newtheorem*{pr*}{Problem}
\newtheorem*{prcat*}{Problema}
\newtheorem{st}[te]{Strategy}
\newtheorem{ob}[te]{Remark}

\newenvironment{dem}{\noindent{\sl Proof. }}{\hfill$\Box$\bigskip}

\newlength\meiabase
\newlength\altura

\newlength\meiabasegrande
\newlength\alturagrande

\newlength\meiabasex
\newlength\alturax
\newlength\meiaalturax

\newlength\meiabasexgrande
\newlength\alturaxgrande

\newlength\meiabasexpeq
\newlength\alturaxpeq

\makeatletter
\renewcommand{\bigoplus}{\operatorname{\mathchoice{\textstyle\bigoplus@}{\bigoplus@}{\bigoplus@}{\bigoplus@}}\slimits@}
\makeatother

\title{Stability and moduli spaces\\ of syzygy bundles}

\author{Pedro Correia Gonçalves Macias Marques}

\makeatletter
\newcommand{\spchapter}{
\if@openright \cleardoublepage \else \clearpage \fi%
\thispagestyle{empty} \enlargethispage{2\baselineskip}%
\global \@topnum \z@ \@afterindentfalse \secdef \@chapter \@schapter}

\renewcommand{\tableofcontents}{%
\if@twocolumn \@restonecoltrue \onecolumn \else \@restonecolfalse \fi
\spchapter*{\contentsname \@mkboth {\MakeUppercase \contentsname }{\MakeUppercase \contentsname}}\@starttoc {toc}\if@restonecol \twocolumn \fi}
\makeatother

\begin{document}

\newpage
\thispagestyle{empty}

\begin{center}

{\large\sc Universitat de Barcelona}

\vspace{\stretch{1}}

{\LARGE
\makeatletter
\@title
\makeatother
}

\bigskip

{\sc Tesi de doctorat}

\end{center}

\vspace{\stretch{1.61803399}}

\begin{flushleft}

{
\makeatletter
\@author
\makeatother\\

\smallskip

Directora de tesi: Rosa María Miró-Roig

\medskip

Programa de doctorat de matemàtiques, curs 2008--2009.

\medskip

Barcelona, setembre de 2009
}
\end{flushleft}

\newpage

\phantom{x}

\vspace{\stretch{1}}

\begin{flushleft}

{\small
}

\end{flushleft}

\cleardoublepage

\thispagestyle{empty}

\phantom{x}

\vspace{\stretch{1}}

\begin{flushright}
\emph{To Maria Inês, Paula and Noémio}
\end{flushright}

\vspace{\stretch{1.61803399}}

\phantom{x}

\cleardoublepage

\selectlanguage{english}

\begin{center}
\textbf{\Large Acknowledgements}
\end{center}

\thispagestyle{empty}

\bigskip

The present work sees the light of day mostly thanks to the close support of my supervisor, Rosa María Miró\h{Roig}. Her contagious enthusiasm, perceptiveness, and tireless dedication gave me the will and strength to put my mind to it and overcome any obstacles I came across.

I wish to thank Laura Costa for always being available to help, both making difficult questions more clear and providing useful suggestions.

I was introduced to algebraic geometry many years ago by Margarida Mendes Lopes, who has the gift of pointing out the most interesting approach to any given problem. I am grateful for her help and advise throughout these years.

I would also like to thank Holger Brenner for sharing his notes on the case of monomials in three variables. The idea behind lemma~\ref{Brenner} is his.

For two years Helena Soares and I shared an office at the University of Barcelona. We had the chance to help each other and keep up our old habit of discussing in detail any problem we did not quite understand. Her dedication helping me from the beginning of my stay until the last proofreading deserves my gratitude. A special word of appreciation is due to her friendship, wisdom, and kindness. \emph{Obrigado}.

A number of colleagues of mine made my work here more interesting and pleasant. They taught me Catalan, helped with the little everyday problems of living in a new place, and together we had many fruitful maths discussions: Elisenda Feliu, Ferran Espuny, Gemma Colomé, Jesús Fernández, José Gil, Martí Lahoz, Mónica Manjarín, \emph{moltes gràcies a tots}. Ferran, Gemma, and José deserve a special mention, for they also helped with the summary in Catalan. For all my questions on how to write in \LaTeX{} correctly, I thank again José, the department's \LaTeX pert. I thank Alberto López for having the patience to make a final proofreading. I have also learnt a lot from colleagues with whom I had classes or sessions on several topics: Alberto Fernández, Abdó Roig, Andratx Bellmunt, Federico Cantero, Florian Heiderich, Joan Pons, Ricardo Zarzuela, Víctor González.

Pedro Tavares is responsible for this publication's cover. I thank him both for it and for the pleasure of listening to his ideas on design and art.

Barcelona proved to be a very kind city for me during these years. I am thankful to all the people who have made it so.

This work had the support of Fundação para a Ciência e a Tecnologia, under grant SFRH/BD/27929/2006, and University of Évora, for which I thank both institutions. I would also like to thank University of Barcelona for the hospitable atmosphere.

\thispagestyle{empty}

\cleardoublepage
\thispagestyle{empty}

\phantom{x}

\vspace{\stretch{1}}

\begin{flushright}
\begin{minipage}{16em}
{\small
Entre a terra e os astros, flor intensa,\\
Nascida do silêncio, a lua cheia\\
Dá vertigens ao mar e azula a areia,\\
E a terra segue-a em êxtases suspensa.

\bigskip

Sophia de Mello Breyner Andresen\\
Lua, \emph{in} Dia do Mar, 1947
}
\end{minipage}
\end{flushright}

\vspace{\stretch{1.61803399}}

\phantom{x}

\newpage

\tableofcontents
\thispagestyle{empty}

%

\pagestyle{fancy}

\newcommand{\starchapter}{Introduction}
\chapter*{\starchapter\markboth{\starchapter}{\starchapter}}
\addcontentsline{toc}{chapter}{\starchapter}

\pagenumbering{roman} 

\begin{flushright}
\begin{minipage}{.61803399\textwidth}
{\footnotesize
\textbf{syzygy} /sizziji/ \emph{noun} (pl.\ syzygies) \textbf{1}~\emph{Astronomy} conjunction or opposition, especially of the moon with the sun. \textbf{2}~a pair of connected or corresponding things.

\textsc{Origin} Greek suzugia, from suzugos `yoked, paired'

\smallskip

Compact Oxford English Dictionary
}
\end{minipage}
\end{flushright}

\bigskip

To determine whether a syzygy bundle on $\pp^N$ is stable, or semistable, is a long\h{standing} problem in algebraic geometry. It is closely related to the problem of finding the Hilbert function and the minimal free resolution of the coordinate ring of the variety defined by a family of general homogeneous polynomials ${f_1,\ldots,f_n}$ in $\kk[X_0,\ldots,X_N]$. This problem goes back at least to the eighties, when Fröberg addresses it in his paper \cite{Fro85}, to find a lower estimate for the Hilbert series of such a ring in terms of the degrees of ${f_1,\ldots,f_n}$.

A syzygy bundle is defined as the kernel of an epimorphism \[\xymatrix{{\bigoplus\limits_{i=1}^n\OO_{\pp^N}(-d_i)}
    \ar[rr]^-{f_1,\ldots,f_n}&&\OO_{\pp^N}\mbox{,}}\]
given by ${(g_1,\ldots,g_n)\mapsto f_1g_1+\cdots+f_ng_n}$, where ${f_1,\ldots,f_n}$ are homogeneous\linebreak[4] polynomials in $\kk[X_0,\ldots,X_N]$ of degrees ${d_1,\ldots,d_n}$, respectively, such that the ideal  ${(f_1,\ldots,f_n)}$ is \mbox{$\mathfrak{m}$-primary}, with ${\mathfrak{m}=(X_0,\ldots,X_N)}$.

In this thesis we consider the case of syzygy bundles defined by general forms ${f_1,\ldots,f_n}$ of the same degree~$d$, and prove their stability and unobstructedness for ${N\ge2}$, except for the case ${(N,d,n)=(2,2,5)}$, where only semistability is guaranteed. To this end, we focus on the case of monomials and derive consequences for general forms from here. The main goal of this work is therefore to give a complete answer to the following problem:
\begin{pr*}[\ref{6.9}]
Does there exist for every~$d$ and every $n\le\tbinom{d+N}{N}$ a family of~$n$ monomials in $\kk\left[X_0,\ldots,X_N\right]$ of degree~$d$ such that their syzygy bundle is semistable?
\end{pr*}
Stability is an Zariski\h{open} propriety, as is semistability. Therefore if this problem has a positive answer, or even better, if we find a family of monomials of a fixed degree~$d$ whose syzygy bundle is stable, then we know that a family of general forms ${f_1,\ldots,f_n}$ of degree~$d$ also gives rise to a stable syzygy bundle.

Problem~\ref{6.9} was presented by Brenner in \cite{Bre08b},
where he gave a geometric interpretation of tight closure in terms of vector bundles. Tight closure is a technique in positive characteristic ring theory introduced by Hochster and Huneke in the late eighties (see \cite{HH88} and \cite{HH90}). It is an operation on ideals containing fields, or submodules of a given module. In the case of parameter ideals, the existence of big Cohen\h{Macaulay} algebras is inherent in studying their tight closures \cite{HH92}. They also give us information on rational singularities \cite{Smi97a}, the Kodaira vanishing theorem \cite{HS97}, and the Briançon\h{Skoda} theorem (see \cite{HH90}, \cite{HH94}, and \cite{Smi97b}). An account of these results can be found in \cite{Hun98}. In 1994 Hochster introduced a generalisation of tight closure, the solid closure, and gave a characterisation of the tight closure of an ideal ${I=(f_1,\ldots,f_n)}$ in a ring~$R$ in terms of local cohomology of forcing algebras
\[\frac{R[T_1,\ldots,T_n]}{(f_1T_1+\cdots+f_nT_n+f)}\mbox{.}\]
This characterisation was described by Brenner in \cite[theorem~2.3]{Bre08b}, and is a consequence of a combination of corollary~2.4, proposition~5.3 and theorems~8.5 and~8.6 in \cite{Hoc94}. Brenner then associated to these forcing algebras a vector bundle, the syzygy bundle of the polynomials ${f_1,\ldots,f_n}$, denoted by ${\syz(f_1,\ldots,f_n)}$, which is the restriction of
\[\mathrm{Spec}\
    \frac{R[T_1,\ldots,T_n]}{(f_1T_1+\cdots+f_nT_n)}\]
to the open set ${U:=D(f_1,\ldots,f_n)}$, of the points that are not common zeros of ${f_1,\ldots,f_n}$ (see \cite[definition 2.12]{Bre08b}). This bundle fits into the exact sequence
\[\xymatrix{0\ar[r]&\syz(f_1,\ldots,f_n)\ar[r]&
    {{\OO_U}^n}\ar[rr]^{f_1,\ldots,f_n}&&\OO_U\ar[r]&0\mbox{,}}\]
called its presenting sequence. For graded rings of dimension two, some open problems on tight closure were solved using these bundles, namely that tight closure equals plus closure if the base field is finite \cite{Bre06b}, and that the Hilbert\h{Kunz} multiplicity is a rational number \cite{Bre06a}.

\bigskip

There is a number of other important reasons for studying vector bundles on algebraic varieties. Among them are the information they provide on the underlying varieties and the moduli spaces they give rise to, which become interesting examples of varieties in high dimensions.

Moduli spaces of vector bundles are schemes that parameterise given families of isomorphism classes of vector bundles. Their existence depends on the way such a family is chosen. For instance the set of all isomorphism classes of vector bundles over a given variety~$X$ is in general too large to be parameterised. Even if we impose restrictions such as a given rank or fixed Chern classes, there is no hope of finding a scheme of finite type parameterising this family of isomorphism classes of vector bundles. However, if we impose an additional restriction, that of stability, there is a natural way of making this parametrisation, as was shown by Maruyama in \cite{Mar76}. Mumford introduced the definition of stability of vector bundles over curves \cite{Mum62} precisely to get round this problem. Takemoto generalised it later to vector bundles over surfaces in \cite{Tak72} and \cite{Tak73}. The notion of stability known today as Mumford\h{Takemoto} stability or \mbox{$\mu$-stability} is a generalisation of this one to torsion\h{free} sheaves over varieties in any dimension.

\bigskip

As we have said before, the present work is dedicated to the study of syzygy bundles over the projective space. They are a particular case of syzygy sheaves, which are defined in a similar manner as the kernel of an epimorphism
\[\xymatrix{{\bigoplus\limits_{i=1}^n\OO_{\pp^N}(-d_i)}
    \ar[rr]^-{f_1,\ldots,f_n}&&\OO_{\pp^N}\mbox{,}}\]
given by ${(g_1,\ldots,g_n)\mapsto f_1g_1+\cdots+f_ng_n}$, where ${f_1,\ldots,f_n}$ are homogeneous\linebreak[4] polynomials in $\kk[X_0,\ldots,X_N]$ of degrees  ${d_1,\ldots,d_n}$, respectively. The difference here is that for the definition of syzygy sheaves the ideal ${(f_1,\ldots,f_n)}$ does not have to be \mbox{$\mathfrak{m}$-primary}. 
Syzygy sheaves are locally free on the open set ${D(f_1,\ldots,f_n)}$, and this is why if we want these sheaves to be vector bundles on the whole projective space, then there should be no common zeros of the polynomials ${f_1,\ldots,f_n}$, which is the same as saying that the ideal ${(f_1,\ldots,f_n)}$ should be \mbox{$\mathfrak{m}$-primary}. This last condition implies that ${n\ge N+1}$, since an \mbox{$\mathfrak{m}$-primary} ideal cannot have a lower number of generators than the number of variables \cite[theorem~11.14]{AM69}.


To see a very simple and well\h{known} example of a syzygy bundle, consider the cotangent bundle $\Omega_{\pp^N}^1$, for it is the kernel of
\[\xymatrix{\OO_{\pp^N}(-1)^{N+1}\ar[rr]^-{X_0,\ldots,X_N}&&
    \OO_{\pp^N}\mbox{.}}\]

Until now, very little was known on stability of syzygy bundles. The cotangent bundle is known to be stable. In characteristic zero, if instead of linear forms we take ${N+1}$ homogeneous polynomials of a given degree~$d$, we also know that the corresponding syzygy bundle is stable, due to a result by \mbox{Bohnhorst} and Spindler \cite{BS92}. Still in the context of characteristic zero, another syzygy bundle known to be stable is the one obtained from a family of $\tbinom{d+N}{N}$ linearly independent homogeneous polynomials of degree~$d$. Flenner proved its semistability in 1984 \cite{Fle84}, and Ballico proved in 1992 its stability \cite{Bal92}. In 2008 Hein proved the semistability of a syzygy sheaf defined by a family of $n$ general homogeneous polynomials, with ${2\le n\le d(N+1)}$, in any characteristic \cite[appendix~A]{Bre08a}.

Also in 2008 more progress was made, when Brenner presented the following sufficient condition for a syzygy bundle defined by monomials to be stable, or semistable (see \cite{Bre08a} or \cite{Bre08b}).
\begin{co*}[\ref{cor6.6}]
Let $f_i$, with $i\in I$, denote monomials in ${\kk[X_0,\ldots,X_N]}$ of degrees~$d_i$, such that the ideal ${(f_i,i\in I)}$ is \mbox{$\mathfrak{m}$-primary}. Suppose that, for every subset ${J\subseteq I}$, with ${|J|\ge2}$, the inequality
\[\frac{d_J-\sum_{i\in J}d_i}{|J|-1}\le\frac{-\sum_{i\in I}d_i}{|I|-1}\]
holds, where $d_J$ is the degree of the highest common factor of the subfamily $\{f_i:i\in J\}$. Then the syzygy bundle ${\syz(f_i,i\in I)}$ is semistable (and stable if strict inequality holds for ${J\subset I}$).
\end{co*}
The inequality to be satisfied in this corollary depends only on the cardinalities of the sets involved and on degrees of monomials. This reduces the problem of deciding if a given syzygy bundle defined by monomials is stable or semistable to a finite number of calculations. Of course, since the number of subsets of a family increases exponentially with its cardinality, the number of calculations to make can become quite cumbersome. Therefore any strategy on the choice of families of monomials that decreases this number is helpful.

Having presented this result, Brenner proposed problem~\ref{6.9}, mentioned before. Chapters \ref{t2} and \ref{t3} of this thesis are dedicated to answering this question for ${N\ge2}$. We will see that the answer is affirmative in all cases. In fact, the results obtained here are stronger. With the exception of the case ${(N,d,n)=(2,5,5)}$, we have found in each case a family of $n$ monomials such that their syzygy bundle is stable.

To get an idea of what has to be done when solving this problem, first note that since only monomials of the same degree~$d$ are considered, the inequality in corollary~\ref{cor6.6} becomes simply
\[\frac{d_J-kd}{k-1}\le\frac{-nd}{n-1}\mbox{,}\]
where ${n:=|I|}$ and ${k:=|J|}$. If ${d_J=0}$, this inequality holds, since we always have ${k\le n}$ and the sequence  $\left(a_{d,j}\right)_{j\ge2}$, with $a_{d,j}:=-\tfrac{jd}{j-1}$, is monotonically increasing. Therefore if we wish to find a set $I$ satisfying this inequality, we have to guarantee that, for each subset~${J\subseteq I}$, the degree~$d_J$ of the greatest common divisor of its elements is low enough.

As we will see in more detail in chapter~\ref{t2}, the set of monic monomials in three variables of a given degree can be represented in suitable manner in the shape of a triangle as in this following example, for degree three.
{\small%
\[
\xymatrix@C=\meiabasex@R=\meiaalturax@!0{
&&&{X_2}^3\\
\\
&&{X_0}{X_2}^2&&{X_1}{X_2}^2\\
\\
&{X_0}^2X_2&&{X_0}X_1X_2&&{X_1}^2X_2\\
\\
{X_0}^3&&{X_0}^2X_1&&{X_0}{X_1}^2&&{X_1}^3\\
}
\]
}%
In this way, the closer any two monomials are, the higher the degree of their greatest common divisor is. From here, we see that to find a family of monomials whose syzygy bundle is stable, monomials should be chosen to be sufficiently spread apart in this representation.

Different approaches were adopted to achieve this. For the lowest values of~$n$ (up to $18$), an individual solution was found for each case. Then for ${18<n\le d+2}$, we looked for the largest triangular number~$T$ not greater than~$n$, and chose $T$ monomials arranged in a position as evenly spread as possible, putting the remaining $n-T$ monomials in appropriate places on the sides of the triangle. For ${d+2<n\le3d}$, all monomials were chosen from the sides of the triangle, and form ${3d+1}$ on, its interior was filled, starting from the layers that are closer to the sides.

In four or more variables, analogous representations are possible, with the shape of a tetrahedron or a hypertetrahedron, respectively. These cases were solved in a simpler fashion. Informally this problem becomes a question of finding enough space to accommodate monomials in these triangles, tetrahedrons and hypertetrahedrons. The fact that in higher dimensions we have more space to work in should not come as a complete surprise.

For the first instances of~$n$, the solution was to consider a family of monomials in ${\kk[X_0,\ldots,X_{N-1}]}$ whose syzygy bundle over $\pp^{N-1}$ is stable and add the monomial ${X_N}^d$. This corresponds to choosing suitable monomials in one face of the hypertetrahedron, and adding the opposite vertex. Then, in analogy with the case of three variables, the faces of the hypertetrahedron are filled. Finally, following an idea form Brenner, the interior of the hypertetrahedron is filled with a family of monomials, obtained from a set of monomials in lower degree whose syzygy bundle is stable.

In this way, we were able to give a complete answer to Brenner's problem. Furthermore, the fact that stability is an open condition allowed us to conclude that a syzygy bundle given by a family of~$n$ general homogeneous polynomials in ${\kk[X_0,\ldots,X_N]}$ of a given degree~$d$ is stable (except for the case ${(N,d,n)=(2,2,5)}$, where only semistability is guaranteed).

\bigskip

Obtaining stability of a syzygy bundle associated to a family of general polynomials of a given degree provides information on the moduli space corresponding to its isomorphism class. Very little is known about moduli spaces of stable vector bundles in general. Having their existence established, one is bound to ask what is there to be said about its structure, both locally and globally. In particular, it would be interesting to know whether it is connected, irreducible, rational or smooth, and what it looks like as a topological space.

In what moduli spaces of syzygy bundles are concerned, we asserted that a stable syzygy bundle is unobstructed and its isomorphism class belongs to a generically smooth irreducible component of the moduli space, and gave the dimension of this component. We also determined whether the closure of the stratum of the moduli space corresponding to stable syzygy bundles is the whole irreducible component it belongs to, and determined its codimension in case it is not.

\bigskip

Let us describe in more detail the structure of this thesis, and highlight the main results.

\smallskip

In \textbf{chapter~\ref{t1}} we provide definitions for the objects studied in the present work, and some tools which will be used later on. We start by defining syzygy bundles and syzygy sheaves in section~\ref{t1syzygy}. In section~\ref{t1stability}, we define stability of coherent sheaves and discuss some of its properties. We then give an account on the contributions to study stability of syzygy bundles done so far. In section~\ref{t1segment}, we solve problem~\ref{6.9} for the case of monomials in two variables. Finally, in section~\ref{t1moduli}, we present a formal definition of moduli space and a result giving sufficient conditions for it to be nonsingular at a given point.

\smallskip

In \textbf{chapter~\ref{t2}} an answer to problem~\ref{6.9} is presented for ${N=2}$. Stability is guaranteed in all cases but ${(N,d,n)=(2,2,5)}$, in which only semistability is obtained. The chapter is divided into different sections, according to different values of~$n$. Recall that for ${N=2}$, we have ${3\le n\le\tbinom{d+2}{2}}$. For the first cases, with ${3\le n\le18}$, as we have said before, an individual answer is given in section~\ref{t2first}. As we have seen, monic monomials in three variables of a given degree can be sketched in a triangle, in a suitable way for the purposes of this problem.  For this reason, the key to solve cases ${18<n\le d+2}$, still in section~\ref{t2first} is to choose a good position for $T$ monomials in this triangle, where $T$ is the greatest triangular number not greater than~$n$. Cases ${d+2<n\le3d}$ are solved in section~\ref{t2sides} by taking a triangle with one complete side and filling the other two sides. The last cases, in section~\ref{t2interior}, for ${3d<n\le\tbinom{d+2}{2}}$, are done by taking a triangle with all sides and filling its interior in a suitable way. Finally, in section~\ref{t2main}, summing up all results in this chapter, the following theorem is presented.

\begin{te*}[\ref{main2}]
Let~$d$ and~$n$ be integers such that $3\le n\le\tbinom{d+2}{2}$, and $(d,n)\ne(2,5)$. Then there is a family of~$n$ monomials in $K\left[X_0,X_1,X_2\right]$ of degree~$d$ such that their syzygy bundle is stable. If $(d,n)=(2,5)$, there is a family of~$n$ monomials such that their syzygy bundle is semistable.
\end{te*}

\smallskip

In \textbf{chapter~\ref{t3}} an answer to problem~\ref{6.9} is presented for ${N\ge3}$. In each case, there is a family of monomials whose corresponding syzygy bundle is stable.

In general monic monomials in $\kk[X_0,\ldots,X_N]$ of a given degree~$d$ can be represented in a hypertetrahedron, in an analogous manner to the triangles we see in case ${N=2}$. We shall call the \emph{$i$th face} of this hypertetrahedron the set of monic monomials where the variable $X_i$ does not occur.

This chapter is divided into different sections, as the previous one, according to different values of~$n$. Recall that we have ${N+1\le n\le\tbinom{d+N}{N}}$. For the first cases, in section~\ref{t3first}, with
\[{N+1\le n\le\tbinom{d+N-1}{N-1}+1}\mbox{,}\]
each family of ${n-1}$ monomials in $\kk[X_0,\ldots,X_{N-1}]$ whose syzygy bundle over $\pp^{N-1}$ is stable generates a family of~$n$ monomials in $\kk[X_0,\ldots,X_N]$ whose syzygy bundle over $\pp^N$ is also stable. Cases \[{\tbinom{d+N-1}{N-1}+1<n\le\tbinom{d+N}{N}-\tbinom{d-1}{N}}\]
are solved in section~\ref{t3faces} by taking the $N$th face and the vertex ${X_N}^d$ of the hypertetrahedron, and adding monomials in the remaining faces. Taking the set of all the hypertetrahedron's faces and adding the monomials in its interior which are closest to the vertexes gives us a solution to cases \[{\tbinom{d+N}{N}-\tbinom{d-1}{N}<n\le\tbinom{d+N}{N}
    -\tbinom{d-1}{N}+N+1}\mbox{,}\]
treated in the beginning of section~\ref{t3interior}. The last cases, with
\[{\tbinom{d+N}{N}-\tbinom{d-1}{N}+N+1<n\le\tbinom{d+N}{N}}
    \mbox{,}\]
also in section~\ref{t3interior}, are solved by taking a family of monomials of lower degree whose syzygy bundle is stable, multiplying them by ${X_0\cdots X_N}$, and adding all monomials in the faces of the hypertetrahedron. This is a generalisation of a lemma by Brenner, made for the case ${N=2}$ in his notes \cite{Bre}, which he kindly shared. Finally, gathering all information from both this chapter and the previous one, the following theorem, which is the main one in this work, is stated in section~\ref{t3main}.
\begin{te*}[\ref{main3}]
Let~$N$, $d$ and~$n$ be integers such that ${N\ge2}$, ${(N,d,n)\ne(2,2,5)}$, and ${N+1\le n\le\tbinom{d+N}{N}}$. Then there is a family of~$n$ monomials in $\kk\left[X_0,\ldots,X_N\right]$ of degree~$d$ such that their syzygy bundle is stable.

For ${(N,d,n)=(2,2,5)}$, there are~$5$ monomials of degree~$2$ in $\kk\left[X_0,X_1,X_2\right]$ such that their syzygy bundle is semistable.
\end{te*}

We conclude this chapter stating two open questions that arise naturally, now that problem~\ref{6.9} is solved. The first one is the generalisation of this problem one gets if we consider homogeneous monomials of arbitrary degree. The second one is to know whether syzygy bundles of higher order are stable, or semistable. This last question involves the open problem of finding the minimal free resolutions of general forms.

\smallskip

In \textbf{chapter~\ref{t4}} we present three results on syzygy bundles and their moduli space. The first one is a consequence of the fact that the conditions of stability and semistability are open in the moduli space.

\begin{te*}[\ref{geral4}]
Let $N\ge2$, $d\ge1$ and $N+1\le n\le\tbinom{d+N}{N}$, with ${(N,d,n)\ne(2,2,5)}$. If ${f_1,\ldots,f_n}$ are general forms of degree~$d$ in $\kk[X_0,\ldots,X_N]$, such that the ideal ${(f_1,\ldots,f_n)}$ is \mbox{$\mathfrak{m}$-primary}, then the syzygy bundle ${\syz(f_1,\ldots,f_n)}$ is stable.

If ${(N,d,n)=(2,2,5)}$, and ${f_1,\ldots,f_5}$ are general quadratic forms such that the ideal ${(f_1,\ldots,f_5)}$ is \mbox{$\mathfrak{m}$-primary}, then the syzygy bundle ${\syz(f_1,\ldots,f_5)}$ is semistable.
\end{te*}

Let $N$, $d$ and $n$ be integers such that ${N\ge2}$, ${d\ge1}$, ${N+1\le n\le\tbinom{d+N}{N}}$ and ${(N,d,n)\ne(2,2,5)}$. We denote by $M_{N,d,n}$ the moduli space of rank ${n-1}$ stable vector bundles on $\pp^N$ with Chern classes ${c_i=\tbinom{n}{i}(-d)^i}$, for ${1\le i\le N}$, and by $S_{N,d,n}$ the stratum of this moduli space corresponding to syzygy bundles.

The second result in the chapter gives us the dimension of the stratum $S_{N,d,n}$, in terms of $N$, $d$ and $n$.

\begin{p*}[\ref{syzlocus}]
Fix integers $N$, $d$ and $n$ such that ${N\ge2}$, ${d\ge1}$,
\[{N+1\le n\le \tbinom{d+N}{N}}\mbox{,}\]
and ${(N,d,n)\ne (2,2,5)}$. Then
\[\dim S_{N,d,n}=n\tbinom{d+N}{N}-n^2.\]
\end{p*}

The third result asserts that the irreducible component of the moduli space corresponding to the isomorphism class of a stable syzygy bundle $E_{N,d,n}$ over $\pp^N$, with ${N\ge2}$, is generically smooth, and gives its dimension.

\begin{te*}[\ref{moduli}]
Let $N$, $d$ and $n$ be integers such that ${N\ge2}$, ${d\ge1}$,
\[{N+1\le n\le\tbinom{d+N}{N}}\mbox{,}\]
and ${(N,d,n)\ne(2,2,5)}$. Then
\begin{enumerate}
  \item the syzygy bundle $E_{N,d,n}$ is unobstructed and its isomorphism class belongs to a generically smooth irreducible component of the moduli space~$M_{N,d,n}$, of dimension ${n\tbinom{d+N}{N}-n^2}$, if $N\ge3$, and ${n\tbinom{d+2}{2}+n\tbinom{d-1}{2}-n^2}$, if ${N=2}$;
  \item if $N\ge 3$, then the closure of the syzygy locus $S_{N,d,n}$ is an irreducible component of $M_{N,d,n}$; if ${N=2}$, the closure of $S_{N,d,n}$ has codimension $n\tbinom{d-1}{2}$ in $M_{N,d,n}$.
\end{enumerate}
\end{te*}

\renewcommand{\starchapter}{Notation and conventions}
\chapter*{\starchapter\markboth{\starchapter}{\starchapter}}
\addcontentsline{toc}{chapter}{\starchapter}

Throughout this thesis $\kk$ will be an algebraically closed field of arbitrary characteristic. Given a vector space~$V$ over~$\kk$, $\pp(V)$ will denote the projective space over~$V$, the set of one dimensional subspaces of~$V$, regarded as an algebraic variety over~$\kk$. In case ${V=\kk^{N+1}}$, the projective space will be denoted as~$\pp^N$.

Given a scheme~$(X,\OO_X)$, and a sheaf~$E$ over~$X$, we will denote its dual $\mathcal{H}om_{\OO_X}(E,\OO_X)$ by $E^\vee$, and, for any $a\in\zz$, the tensor product ${E\otimes\OO_X(a)}$ will be denoted by ${E(a)}$. The $i$th cohomology group of $E$ will be denoted by $H^i(X,E)$, or just $H^i(E)$, if no confusion arises. Its dimension as a vector space over~$\kk$ will be denoted by $h^i(X,E)$, or simply $h^i(E)$.

Unless otherwise stated, $\OO$ will stand for $\OO_{\pp^N}$, the sheaf of regular functions on~$\pp^N$, either for general~$N$, or for a particular case, which will be clear from context.

If~$E$ is a rank~$r$ vector bundle over an \mbox{$N$-dimensional} variety~$X$, its Chern classes will be denoted by $c_1(E),\ldots,c_s(E)$, where $s=\min(r,N)$. In case $X=\pp^N$, we will identify these Chern classes with integers.

We denote by
\[{\mathfrak{m}:=(X_0,\ldots,X_N)}\]
the irrelevant maximal ideal of $\kk[X_0,\ldots,X_N]$.

Recall that
\begin{align*}
  \Hom\big(\OO(-d)^n,\OO\big)&\cong
      \bigoplus_{i=1}^n\Hom\big(\OO(-d),\OO\big)\\
  &\cong\bigoplus_{i=1}^nH^0\big(\OO(-d)\big)\\
  &\cong\bigoplus_{i=1}^n\kk^{\binom{d+N}{N}}\\
  &\cong\kk^{n\binom{d+N}{N}}.
\end{align*}
When we say that a \emph{general} morphism $\xymatrix@1{{\varphi:\bigoplus_{i=1}^n\OO(-d)}\ar[r]&\OO}$ has a certain property, we mean that there is an open dense subset~$U$ of $\Hom\big(\OO(-d)^n,\OO\big)$ such that all members of~$U$ have this property. Analogously, when we say that ${f_1,\ldots,f_n}$ are general homogeneous polynomials in ${R:=\kk[X_0,\ldots,X_N]}$ of degrees ${d_1,\ldots,d_n}$, respectively, we mean that the sequence ${(f_1,\ldots,f_n)}$ belongs to a suitable open dense subset of ${R_{d_1}\times\cdots\times R_{d_1}}$.

In what follows, for any ${x\in\rr}$, $\lceil x\rceil$ will represent the least integer greater or equal to~$x$. Given $a,b\in\nn$, we shall assume ${\tbinom{a}{b}=0}$, whenever $a<b$.

\chapter{Preliminaries}
\label{t1}

\pagenumbering{arabic} 

This chapter provides definitions for the objects studied in the present work, and some tools which will be used later on. We start by defining syzygy bundles and syzygy sheaves in section~\ref{t1syzygy}. In section~\ref{t1stability} we define stability of coherent sheaves and discuss some of its properties. We then present the main problem addressed in this thesis, on syzygy bundles of polynomials, and give an account on the contributions done to solve it so far. In section~\ref{t1segment} we solve this problem for the case of monomials in two variables. Finally, in section~\ref{t1moduli} we present a formal definition of moduli space and a result giving sufficient conditions for it to be nonsingular at a given point.

\bigskip

Let $\kk$ be an algebraically closed field of arbitrary characteristic, and $\pp^N$ be the $N$th dimensional projective space over~$\kk$.

If~$X$ is an algebraic variety of dimension~$N$ and~$E$ is a coherent sheaf of \mbox{$\OO_X$-modules} over~$X$, we say that~$E$ is \emph{torsion free} if, for any ${x\in X}$, the stalk~$E_x$ is a torsion\h{free} \mbox{$\OO_{X,x}$-module}, i.e.\ given any non\h{zero} elements $a\in\OO_{X,x}$ and $v\in E_x$, we always have ${av\ne0}$.

The \emph{Euler characteristic} of $E$ is
\[\chi(E):=\sum_{i=0}^N(-1)^ih^i(E).\]

Let~$E$ be a torsion\h{free} sheaf. We say that~$E$ is \emph{normalised} if ${-r<c_1(E)\le0}$. We set ${E_{\mathrm{norm}}:=E(k_E)}$, where $k_E$ is the unique integer such that $E(k_E)$ is normalised.

Let~$E$ be a vector bundle over~$\pp^N$ of rank~$r$. We say that~$E$ is \emph{homogeneous} if for every projective transformation ${t:\pp^N\to\pp^N}$, we have ${t^*E\cong E}$.

\section{Syzygy bundles and syzygy sheaves}
\label{t1syzygy}

The purpose of this dissertation is to study stability of syzygy bundles and their moduli spaces. Let us start by giving the definitions of syzygy bundle and syzygy sheaf on~$\pp^N$.

\begin{defin}
Let ${f_1,\ldots,f_n}$ be homogeneous polynomials in ${\kk[X_0,\ldots,X_N]}$ of degrees $d_1,\ldots,$ $d_n$, respectively, such that the ideal ${(f_1,\ldots,f_n)}$ is \mbox{$\mathfrak{m}$-primary}. We say that the kernel of the morphism of vector bundles over~$\pp^N$
\begin{equation}\label{syzmorph}
\xymatrix{{\bigoplus\limits_{i=1}^n\OO(-d_i)}
    \ar[rr]^-{f_1,\ldots,f_n}&&\OO}\mbox{,}
\end{equation}
defined by ${(g_1,\ldots,g_n)\mapsto f_1g_1+\cdots+f_ng_n}$, is the \emph{syzygy bundle} associated to the family ${f_1,\ldots,f_n}$. We denote this bundle by $\syz(f_1,\ldots,f_n)$ and call the short exact sequence
\[\xymatrix{0\ar[r]&\syz(f_1,\ldots,f_n)\ar[r]&
    {\bigoplus\limits_{i=1}^n\OO(-d_i)}
    \ar[r]&\OO\ar[r]&0}\]
its \emph{presenting sequence}.

If we drop the condition that the ideal generated by ${f_1,\ldots,f_n}$ is \mbox{$\mathfrak{m}$-primary}, the kernel of morphism (\ref{syzmorph}) is still a torsion\h{free} sheaf, but it is locally free only on $D(f_1,\ldots,f_n)$, the open subset of $\pp^N$ of points that are not common zeros of ${f_1,\ldots,f_n}$. We call it the \emph{syzygy sheaf} associated to the family ${f_1,\ldots,f_n}$.
\end{defin}

\begin{ex}
A very simple example of a syzygy bundle is $\Omega_{\pp^N}^1$, since its presenting sequence is the Euler exact sequence
\[\xymatrix{0\ar[r]&\Omega_{\pp^N}^1\ar[r]&\OO(-1)^{N+1}\ar[r]&
    \OO\ar[r]&0}\]
(see for instance \cite[chapter~II, theorem~8.13]{Har77}).
\end{ex}

\begin{ex}
If we let
\[n=h^0\big(\OO(d)\big)=\tbinom{d+N}{N}\mbox{,}\]
and take for instance the family of all monomials of a fixed degree~$d$, we get the syzygy bundle ${\syz\left(\big\{ {X_0}^{i_0} \cdots {X_N}^{i_N}: i_0+\cdots+i_N=d \big\}\right)}$. The relevance of this bundle was pointed out in several articles. In 1988 Green made use of it to present a second proof of the explicit Noether\h{Lefshetz} theorem \cite{Gre88}. In 1991 Ballico used this bundle to assert some results on generation and diameter of the Hartshorne\h{Rao} module of a curve in~$\pp^N$ \cite{Bal91}, and in 1992 he proved its stability \cite{Bal92}. Migliore and Miró\h{Roig} had already stressed the importance of this syzygy bundle when studying \mbox{$k$-Buchbaum} curves in~$\pp^3$ \cite{MMR90}, and in 1994 Miró\h{Roig} characterised it cohomologically \cite{MR94}.
\end{ex}

\begin{ex}
Choose polynomials ${{X_0}^2+X_1X_2}$ and ${X_0+X_1}$ in $\kk[X_0,X_1,X_2]$. Since their common zeros are ${P:=[0:0:1]}$ and ${Q:=[1:-1:1]}$, their syzygy sheaf is torsion\h{free} on~$\pp^2$, and locally free on ${\pp^2\setminus\{P,Q\}}$.
\end{ex}

\bigskip

From a syzygy bundle's presenting sequence, we get immediately that its rank is ${n-1}$. We also see that its Chern polynomial satisfies
\[\xymatrix{c_t\left({\bigoplus\limits_{i=1}^n\OO(-d_i)}\right)=
    c_t\big(\syz(f_1,\ldots,f_n)\big)\cdot c_t(\OO)\mbox{,}}\]
and therefore
\[c_t\big(\syz(f_1,\ldots,f_n)\big)=\prod_{i=1}^n(1-d_it)\in
    \tfrac{\zz[t]}{(t^{N+1})}.\]
In particular, its first Chern class is
\[c_1\big(\syz(f_1,\ldots,f_n)\big)=-(d_1+\cdots+d_n).\]

Chern classes can be generalised to torsion\h{free} coherent sheaves (we can find a description of this for instance in \cite[sections~1 and~2]{Har80}). As Brenner computed in \cite[lemma~2.1]{Bre08a}, the first Chern class of a syzygy sheaf is
\[c_1\big(\syz(f_1,\ldots,f_n)\big)=d'-(d_1+\cdots+d_n)\mbox{,}\]
where~$d'$ is the degree of the greatest common divisor~$f'$ of the family ${f_1,\ldots,f_n}$. 

In the present work, we study stability of syzygy bundles, in the case when the degrees ${d_1,\ldots,d_n}$ are all equal.

\section{Stability of vector bundles}
\label{t1stability}

In this section we present two different definitions of stability and semistability, one from Mumford and Takemoto, which is the one adopted here, and another from Gieseker and Maruyama. We then state the problem addressed in this work, and give an account of the contributions made to solve it until now.

\bigskip

If~$E$ is a torsion\h{free} coherent sheaf over~$\pp^N$ of rank~$r$, we define the \emph{slope} of~$E$ as
\[\mu(E):=\frac{c_1(E)}{r}\mbox{,}\]
where ${c_1(E)}$ is the first Chern class of~$E$, defined as the degree of~$E$, i.e.\ ${c_1(E):=\deg\ E}$, with
\[\deg\ E:=\deg(\det E)=
    \deg\big((\Lambda^rE)^{\vee\vee}\big).\]

\begin{defin}
We say that a torsion\h{free} sheaf $E$ over~$\pp^N$ is \emph{semistable}, in the sense of Mumford\h{Takemoto}, if for every non\h{zero} coherent proper subsheaf $F\subset E$,
\[\mu(F)\le\mu(E).\]
We say that a semistable sheaf~$E$ is \emph{stable} if equality never occurs. If $E$ is semistable, but not stable, we say that it is \emph{strictly semistable}.
\end{defin}

Note that if $r$ and $c_1(E)$ are coprime, the notions of stability and semistability coincide.

This definition of stability is sometimes referred to as \mbox{$\mu$-stability}. Mumford introduced it in his paper \cite{Mum62} for vector bundles on curves, and Takemoto generalised it for vector bundles on surfaces in \cite{Tak72} and \cite{Tak73}. The notion of \mbox{$\mu$-stability} (and \mbox{$\mu$-semistability}) arose in Mumford and Takemoto's work in connection with the following problem:
\begin{quote}
Given a scheme~$X$, is there a scheme of finite type such that its closed points are in bijective correspondence with isomorphism classes of vector bundles with fixed Chern classes?
\end{quote}
Unfortunately the set of isomorphism classes of vector bundles with fixed Chern classes on a scheme~$X$ is too big to be endowed with a natural algebraic structure. The following well\h{known} example illustrates this fact.


\begin{ex}
Consider the set ${\left\{\OO(n)\oplus\OO(-n)\right\}_{n\in\nn}}$, of vector bundles over~$\pp^1$. Since the rank of the vector bundles in this family is two, and the first Chern class is zero, these topological invariants are fixed, and yet this is a infinite numerable set of points, which cannot form an algebraic variety.
\end{ex}

The natural class of vector bundles which admits a natural algebraic structure comes from Mumford's geometric invariant theory \cite{Mum65} (or \cite{MFK94}, for its most recent version). They correspond to stable vector bundles. In fact, if we impose semistability, we get round this problem, as Maruyama established in \cite{Mar77}, because the family of semistable sheaves is bounded \cite[theorem~3.3.1]{HL97}.

A different notion of stability from the one stated above was introduced by Gieseker in \cite{Gie77} and Maruyama in \cite{Mar76} and \cite{Mar77}. We denote, for each torsion\h{free} coherent sheaf~$E$ over~$\pp^N$ of rank~$r$, and each integer ${k\in\zz}$,
\[p_E(k):=\frac{\chi\big(E(k)\big)}{r}.\]

\begin{defin}
We say that $E$ is \emph{GM\h{semistable}} if, for all non\h{zero} coherent proper subsheaves ${F\subset E}$, we have
\[p_F(k)\le p_E(k)\mbox{,}\]
for all sufficiently large integers ${k\in\zz}$. It is said to be \emph{GM\h{stable}} if strict inequality holds for all sufficiently large integers ${k\in\zz}$.
\end{defin}

A connection between these two notions of stability is the following result, which can be found in \cite[chapter~II, lemma~1.2.12]{OSS80}.

\begin{lem}
Stable torsion\h{free} coherent sheaves over~$\pp^N$ are also GM\h{stable}. GM\h{semistable} sheaves over~$\pp^N$ are also semistable.
\end{lem}


The notion of stability (in its different versions) arises naturally also from a gauge theoretical point of view, for, as pointed out in \cite[section~2]{MR07}, there is a deep relation between stability of vector bundles and the existence of Hermite\h{Einstein} metrics. This relation is known as the Kobayashi\h{Hitchin} correspondence and was established by works in Narasimhan\h{Seshadri} \cite{NS65}, Donaldson (\cite{Don85} and \cite{Don87}), and Uhlenbeck\h{Yau} (\cite{UY86} and \cite{UY89}).

A useful result on stability (and semistability) of vector bundles is the following theorem, which can be found in \cite[chapter~II, theorem~1.2.2]{OSS80}.
\begin{te}\label{OSS1.2.2}
Let~$E$ be a torsion\h{free} sheaf over~$\pp^N$. The following statements are equivalent:
\begin{enumerate}
  \item \label{OSS1.2.2i} $E$ is stable (respectively semistable).
  \item \label{OSS1.2.2ii} For every non\h{zero} coherent proper subsheaf ${F\subset E}$ whose quotient $\tfrac{E}{F}$ is torsion\h{free}, ${\mu(F)<\mu(E)}$ (respectively ${\mu(F)\le\mu(E)}$).
  \item \label{OSS1.2.2iii} For every non\h{zero} torsion\h{free} quotient
      \[\xymatrix{E\ar[r]&Q\ar[r]&0}\]
      of~$E$, with ${\mathrm{rk}\ Q<\mathrm{rk}\ E}$,
      ${\mu(E)<\mu(Q)}$ (respectively ${\mu(E)\le\mu(Q)}$).
\end{enumerate}
\end{te}

One easy consequence of this theorem is that direct sums cannot be stable. To see this, observe that if $F$ is a summand of $E$, then $F$ is simultaneously a proper subsheaf and a quotient. If $E$ were stable, by definition, we would have ${\mu(F)<\mu(E)}$, but by number~\ref{OSS1.2.2iii} of the theorem above, we would also have ${\mu(F)>\mu(E)}$.

Another useful result is the following lemma, which can be found in \cite[chapter~II, lemma~1.2.4]{OSS80}, and will be frequently used in arguments to follow.
\begin{lem}\label{OSS1.2.4}
Let~$E$ and~$F$ be torsion\h{free} sheaves over~$\pp^N$.
\begin{enumerate}
  \item \label{OSS1.2.4i} Line bundles are stable.
  \item \label{OSS1.2.4ii} If $E$ and $F$ are both semistable, then their direct sum ${E\oplus F}$ is semistable if and only if ${\mu(E)=\mu(F)}$.
  \item \label{OSS1.2.4iii} $E$ is (semi)stable if and only if~$E^\vee$ is (semi)stable.
  \item \label{OSS1.2.4iv} $E$ is (semi)stable if and only if, for any $a\in\zz$, $E(a)$ is (semi)stable.
\end{enumerate}
\end{lem}

The first statement in this lemma is a direct consequence of the definition of stability. Line bundles have no non\h{zero} proper sub\h{bundles}. For torsion\h{free} sheaves of low rank, working with the definition may be feasible, but as the rank grows, the number of possible subsheaves becomes too large to cope with. It is therefore a very important problem to find criteria to check a given vector bundle's stability. In 1984 Hoppe presented the following cohomological criterion \cite[lemma~2.6]{Hop84}.

\begin{lem}\label{Hop2.6}
Let~$E$ be a normalised vector bundle on~$\pp^N$ of rank~$r$. If for every ${q\in\{1,\ldots,r-1\}}$,
\[H^0\big((\Lambda^qE)_{\mathrm{norm}}\big)=0\mbox{,}\]
then $E$ is stable.
\end{lem}

With little change in Hoppe's proof, one can derive the following lemma, which for this reason is sometimes presented as part of the above (see for instance \cite{MR07}).

\begin{lem}\label{Hop2.6adapt}
Let~$E$ be a normalised vector bundle on~$\pp^N$ of rank~$r$. Then $E$ is semistable if and only if for every ${q\in\{1,\ldots,r-1\}}$,
\[H^0\big((\Lambda^qE)_{\mathrm{norm}}(-1)\big)=0.\]
\end{lem}

The conditions in lemma~\ref{Hop2.6} are not necessary. Indeed, let~$N$ be odd and consider a nullcorrelation bundle~$E$ on~$\pp^N$, i.e.\ a vector bundle~$E$ fitting into an exact sequence
\[\xymatrix{0\ar[r]&\OO(-1)\ar[r]&\Omega_{\pp^N}^1(1)
    \ar[r]&E\ar[r]&0\mbox{.}}\]
Then~$E$ is a stable bundle of rank~$N-1$, but ${H^0\big((\Lambda^qE)_{\mathrm{norm}}\big)\ne0}.$
In fact, $(\Lambda^qE)_{\mathrm{norm}}$ admits~$\OO$ as a direct summand (see \cite[lemma~1.10]{AO94}).

\bigskip

The main problem addressed in this thesis is the following.

\begin{pr}\label{mainprob}
Given a family ${f_1,\ldots,f_n}$ of homogeneous polynomials of degree~$d$ in ${\kk[X_0,\ldots,X_N]}$, such that the ideal ${(f_1,\ldots,f_n)}$ is \mbox{$\mathfrak{m}$-primary}, is the syzygy bundle $\syz(f_1,\ldots,f_n)$ stable (or, at least, semistable)?
\end{pr}

So far, very few contributions to solve this problem exist. We list them in the following examples, briefly describing the techniques used in each one.

Note that since ${(f_1,\ldots,f_n)}$ is an \mbox{$\mathfrak{m}$-primary} ideal, by \cite[theorem~11.14]{AM69}, its number of generators is at least ${N+1}$. On the other hand, if the family of polynomials ${f_1,\ldots,f_n}$ were linearly dependent, we would have a syzygy of degree zero for these polynomials, since there would be ${\alpha_1,\ldots,\alpha_n\in\kk}$ such that ${\alpha_1f_1+\cdots+\alpha_nf_n=0}$. Therefore $\OO(-d)$ would be a sub\h{bundle} of $\syz(f_1,\ldots,f_n)$, and we would have
\[\mu\big(\OO(-d)\big)=-1>-d\cdot\tfrac{n}{n-1}=
    \mu\big(\syz(f_1,\ldots,f_n)\big)\mbox{,}\]
which would mean that this bundle is not semistable. Since ${\tbinom{d+N}{N}}$ is the number of all monomials of degree~$d$ in ${\kk[X_0,\ldots,X_N]}$, and therefore the dimension of the linear space of degree~$d$ polynomials in these variables, we will always assume
\[N+1\le n\le\tbinom{d+N}{N}.\]

In case ${N=1}$, by the Birkhoff\h{Grothendieck} theorem (see \cite[chapter~I, theorem~2.1.1]{OSS80}, or \cite{Gro57} and \cite{Bir13} for the original sources), any vector bundle over~$\pp^1$ is isomorphic to a sum of line bundles. Therefore, by number~\ref{OSS1.2.4i} of lemma~\ref{OSS1.2.4}, syzygy bundles are stable if and only if they have rank~$1$, since direct sums cannot be stable. For syzygy bundles of higher rank, all we can ask is whether they are semistable or not.

\begin{ex}
For ${d=1}$, the only case is ${n=N+1}$. In this case, the syzygy bundle is isomorphic to $\Omega_{\pp^N}^1$, since its presenting sequence is the Euler exact sequence
\[\xymatrix{0\ar[r]&\Omega_{\pp^N}^1\ar[r]&\OO(-1)^{N+1}\ar[r]&
    \OO\ar[r]&0}\]
(see for instance \cite[chapter~II, theorem~8.13]{Har77}). This vector bundle is known to be stable. One proof of this is in  \cite[chapter~II, theorem~1.3.2]{OSS80}, where the stability of the tangent bundle $T_{\pp^N}$ is stated. Since the tangent bundle is the dual of $\Omega_{\pp^N}^1$, by number \ref{OSS1.2.4iii} in lemma~\ref{OSS1.2.4}, we get thus the stability of $\Omega_{\pp^N}^1$.
\end{ex}

\begin{ex}
For ${N\ge2}$, ${n=N+1}$, and arbitrary~$d$, stability of syzygy bundles is guaranteed as a particular case of theorem~2.7 in Bohnhorst and\linebreak[4] Spindler's paper~\cite{BS92}. This theorem gives a criterion for stability of vector bundles~$E$ of rank~$N$ fitting into a short exact sequence
\[\xymatrix{0\ar[r]&
    {\bigoplus\limits_{i=1}^k\OO(a_i)}\ar[r]&
    {\bigoplus\limits_{j=1}^{N+k}\OO(b_j)}
    \ar[r]&E\ar[r]&0\mbox{,}}\]
with ${a_1\ge\cdots\ge a_k}$ and ${b_1\ge\cdots\ge b_{N+k}}$. It turns out that $E$ is stable if and only if
\[b_1<\mu(E)=\tfrac{1}{N}\left[\sum_{j=1}^{N+k}b_j
    -\sum_{i=1}^ka_i\right].\]
Since the dual of a syzygy bundle fits into a short exact sequence
\[\xymatrix{0\ar[r]&\OO\ar[r]&\OO(d)^{N+1}
    \ar[r]&\syz(f_1,\ldots,f_{N+1})^\vee\ar[r]&0\mbox{,}}\]
this bundle is stable, and therefore $\syz(f_1,\ldots,f_{N+1})$ is also stable, again by number \ref{OSS1.2.4iii} in lemma~\ref{OSS1.2.4}.
\end{ex}

\begin{ex}
In 1984 Flenner proved the semistability of the syzygy bundle associated to a family of~${n=\tbinom{d+N}{N}}$ homogeneous polynomials of degree~$d$, in the case of a field of characteristic zero \cite[corollary~2.2]{Fle84}. For ${d>1}$ and ${N\ge2}$, Ballico proved in 1992 the stability of the same syzygy bundle \cite[theorem~0.2]{Bal92}. In 1995 Paoletti presented a different proof of the same result, and asserted the stability of the exterior powers ${\Lambda^q\syz(f_1,\ldots,f_n)}$, for ${q\in\{1,\ldots,N-1\}}$, and the semistability of the homogeneous vector sub\h{b}un\-dles of $\syz(f_1,\ldots,f_n)$ \cite{Pao95}.

In positive characteristic, and also for ${n=\tbinom{d+N}{N}}$, Langer \cite{Lan09} has given sufficient conditions for a syzygy bundle~$E$ to be strongly semistable, i.e.\ to be a semistable vector bundle such that all its Frobenius pullbacks $\left(F^k\right)^*E$, where $\xymatrix{F:\pp^N\ar[r]&\pp^N}$ is the Frobenius map, are semistable. Trivedi improved these results in \cite{Tri}, and recently Mehta found a general proof for any characteristic \cite[appendix]{Lan}.
\end{ex}

\begin{ex}
In 2008 Hein proved the semistability of a syzygy sheaf defined by a family of $n$ general homogeneous polynomials of degree~$d$, with ${2\le n\le d(N+1)}$, in any characteristic \cite[appendix~A]{Bre08a}.
\end{ex}

Also in 2008, Brenner obtained the following necessary numerical conditions for a syzygy sheaf to be semistable, both in his paper \cite[proposition~2.2]{Bre08a} and in his notes for a Winter School on Commutative Algebra and Applications, held in Barcelona in 2006 \cite{Bre08b}.

\begin{p}
Let $f_i\ne0$, with $i\in I$ and ${|I|\ge2}$, denote homogeneous elements in the polynomial ring ${\kk[X_0,\ldots,X_N]}$ of degrees $d_i$. For every subset ${J\subseteq I}$, denote by $d_J$ the degree of the highest common factor of the subfamily $\{f_i:i\in J\}$. Suppose that the syzygy sheaf ${\syz(f_i,i\in I)}$ is semistable. Then for every ${J\subseteq I}$, with ${|J|\ge2}$, we have the numerical condition
\[\frac{d_J-\sum_{i\in J}d_i}{|J|-1}\le\frac{d_I-\sum_{i\in I}d_i}{|I|-1}.\]
If ${\syz(f_i,i\in I)}$ is stable, the strict inequality holds for ${J\subset I}$.
\end{p}

These conditions become sufficient, if we consider monomials instead of polynomials, as Brenner showed in the mentioned works, inspired by Hoppe's ideas \cite{Hop84}, and using results from Klyachko's papers \cite{Kly90}, \cite{Kly98} and \cite{Kly02} on toric bundles (see \cite[corollary~6.4]{Bre08a} or \cite[theorem~6.6]{Bre08b}). In fact, Brenner assumes characteristic zero, but, as was pointed out by Perling \cite{Per04}, Klyachko's results are valid in any characteristic.

\begin{co}\label{cor6.6}
Let $f_i$, with $i\in I$, denote monomials in ${\kk[X_0,\ldots,X_N]}$ of degrees~$d_i$, such that the ideal ${(f_i,i\in I)}$ is \mbox{$\mathfrak{m}$-primary}. Suppose that, for every subset ${J\subseteq I}$, with ${|J|\ge2}$, the inequality
\[\frac{d_J-\sum_{i\in J}d_i}{|J|-1}\le\frac{-\sum_{i\in I}d_i}{|I|-1}\]
holds, where $d_J$ is the degree of the highest common factor of the subfamily $\{f_i:i\in J\}$. Then the syzygy bundle ${\syz(f_i,i\in I)}$ is semistable (and stable if strict inequality holds for ${J\subset I}$).
\end{co}

\begin{ex}
Applying the criterion in this corollary to the family
\[\big\{{X_0}^4,\,{X_1}^4,\,{X_2}^4,\,{X_0}^2X_1X_2\big\}\]
of monomials in ${\kk[X_0,X_1,X_2]}$, we see that it generates an \mbox{$\mathfrak{m}$-primary} ideal, and that the only subsets admitting a greatest common divisor of positive degree are $\big\{{X_0}^4,\,{X_0}^2X_1X_2\big\}$, $\big\{{X_1}^4,\,{X_0}^2X_1X_2\big\}$ and $\big\{{X_2}^4,\,{X_0}^2X_1X_2\big\}$. They all have cardinality two, and the first admits a greatest common divisor of degree two, whereas the others only admit a linear greatest common divisor. Since ${2-(4+4)=-6<-\tfrac{16}{3}}$, the syzygy bundle corresponding to this family is stable.
\end{ex}

\begin{ex}
If we change one monomial in the previous example and consider the family
\[\big\{{X_0}^4,\,{X_1}^4,\,{X_2}^4,\,{X_0}^3X_1\big\}\mbox{,}\]
we see that the subset $\big\{{X_0}^4,\,{X_0}^3X_1\big\}$ admits a greatest common divisor of degree three, and since ${3-(4+4)=-5>-\tfrac{16}{3}}$, the syzygy bundle associated to this family is not even semistable.
\end{ex}

\begin{ex}
Now if we consider the family
\[\big\{{X_0}^2,\,{X_1}^2,\,{X_2}^2,\,X_0X_1,
    \,X_0X_2\big\}\mbox{,}\]
we can check that all subsets strictly satisfy the inequality in corollary~\ref{cor6.6}, except for $\big\{{X_0}^2,\,X_0X_1,\,X_0X_2\big\}$. This subset has three monomials admitting a linear greatest common divisor. Since ${\tfrac{1-(2+2+2)}{2}=-\tfrac{(2+2+2+2+2)}{4}}$, the corresponding syzygy bundle is strictly semistable. We can observe that this is the only family of five quadratic monomials in three variables such that the ideal generated by them is \mbox{$m$-primary}, up to change of variables.
\end{ex}

Since stability is an open property in the moduli space, a positive answer to problem~\ref{mainprob} can be given for a family of general homogeneous polynomials, if an example is presented, for each $N$, $d$ and $n$. This last corollary gives us a tool to find such examples, in the case of monomials. Therefore a general answer to problem~\ref{mainprob} is given if the following problem (presented by Brenner in his notes~\cite{Bre08b}) is solved.

\begin{pr}\label{6.9}
Does there exist for every~$d$ and every $n\le\tbinom{d+N}{N}$ a family of~$n$ monomials in $\kk\left[X_0,\ldots,X_N\right]$ of degree~$d$ such that their syzygy bundle is semistable?
\end{pr}

Since monomials in~$I$ can be taken to be all monic, we can work only with these and in what follows \emph{monomial} will always mean monic monomial.

\begin{ob}\label{adj}
If we use the notation $a_{d,j}:=-\tfrac{jd}{j-1}$, the inequality in the previous corollary can be written
\[{\tfrac{d_J}{k-1}+a_{d,k}\le a_{d,n}}\mbox{,}\]
where ${k=|J|}$. The fact that once~$d$ is fixed, the sequence $\left(a_{d,j}\right)_{j\ge2}$ is monotonically increasing will be useful in many arguments.
\end{ob}

The inequality above is equivalent to
\begin{equation}\label{6.6}
    \left(d-d_J\right)n+d_J-dk\ge0\mbox{,}
\end{equation}
and this is the version we will be using mostly.

\section{Monomials in two variables}
\label{t1segment}

In this section a solution to problem~\ref{6.9} is presented for the case ${N=1}$, of monomials in two variables.

In case ${N=1}$, we get stability if ${n=2}$, since the syzygy bundles is a line bundle. For ${n\ge3}$, as we have seen above, by the Birkhoff\h{Grothendieck}  theorem, all vector bundles are sum of line bundles, and therefore cannot be stable. If ${n=d+1}$, the syzygy bundle $\syz(\I_{1,d,n})$, where
\[\I_{1,d,n}:=\big\{{X_0}^d,\,{X_0}^{d-1}X_1,\ldots,{X_1}^d\big\}\mbox{,}\]
is strictly semistable. Indeed, if~$g$ is the greatest common divisor of monomials in a subset ${J\subseteq\I_{1,d,n}}$, all monomials in~$J$ are of the form $gh$, with~$h$ a monomial of degree ${d-d_J}$, where~$d_J$ is the degree of~$g$. There are ${d-d_J+1}$ monomials of degree ${d-d_J}$, so
\[k:=|J|\le d-d_J+1.\]
Now
\[(d-d_J)n+d_J-dk\ge(d-d_J)(d+1)+d_J-d(d-d_J+1)=0.\]
Therefore inequality~(\ref{6.6}) holds. In fact, if we choose
\[J:=\big\{{X_0}^d,\,{X_0}^{d-1}X_1,
    \ldots,{X_0}^{d-d_J}{X_1}^{d_J}\big\}\mbox{,}\]
we get equality, which means that the syzygy bundle is not stable.

In general, for the remaining values of~$n$, i.e.\ $3\le n\le d$, if ${f_1,\ldots,f_n}$ is a family of homogeneous polynomials in $\kk[X_0,X_1]$ such that the ideal ${(f_1,\ldots,f_n)}$ is \mbox{$\mathfrak{m}$-primary}, their syzygy bundle ${\syz(f_1,\ldots,f_n)}$ has rank ${n-1}$ and first Chern class
\[c_1\big(\syz(f_1,\ldots,f_n)\big)=-dn.\]
By the Birkhoff\h{Grothendieck} theorem, there are integers ${a_1,\ldots,a_{n-1}}$ such that
\[\syz(f_1,\ldots,f_n)\cong
    \bigoplus\limits_{i=1}^{n-1}\OO(a_i).\]
Now by number~\ref{OSS1.2.4ii} of lemma~\ref{OSS1.2.4}, this syzygy bundle is semistable if and only if
\[\mu\big(\OO(a_1)\big)=\cdots=\mu\big(\OO(a_{n-1})\big)\mbox{,}\]
and this happens if and only if ${a_1=\cdots=a_{n-1}}$. Therefore
\[(n-1)a_1=a_1+\cdots+a_{n-1}=
    c_1\big(\syz(f_1,\ldots,f_n)\big)=-dn\mbox{,}\]
and since $n$ and ${n-1}$ are coprime, $d$ is a multiple of ${n-1}$. We have thus found a necessary condition for such a syzygy bundle to be semistable.

Now for the converse, suppose $d$ is a multiple of ${n-1}$, say ${d=(n-1)e}$, with ${e\in\zz}$, and consider the set
\[\I_{1,d,n}:=\big\{{X_0}^d,\,{X_0}^{(n-2)e}{X_1}^e,
    \,{X_0}^{(n-3)e}{X_1}^{2e},\ldots,\,{X_0}^e{X_1}^{(n-2)e},
    \,{X_1}^d\big\}.\]
To check inequality~(\ref{6.6}), consider a subset $J\subseteq\I_{1,d,n}$, with~$k$ elements, ${k\ge2}$. The degree $d_J$ of the greatest common divisor of the monomials in~$J$ is at most ${(n-k)e}$, and we get
\begin{align*}
(d-d_J)n+d_J-dk&=(n-k)d-(n-1)d_J\\
    &\ge(n-k)d-(n-1)(n-k)e\\
    &=(n-k)(n-1)e-(n-1)(n-k)e=0.
\end{align*}
Therefore, $\syz(\I_{1,d,n})$ is a strictly semistable bundle.

Since in cases ${n=2}$ and ${n=d+1}$, $d$ is also a multiple of ${n-1}$, we have proved the following result.

\begin{te}\label{main1}
Let~$d$ and~$n$ be integers such that $2\le n\le d+1$ and~$d$ is a multiple of ${n-1}$. Then there is a family of~$n$ monomials in $K\left[X_0,X_1\right]$ of degree~$d$ such that their syzygy bundle is semistable. It is stable for ${n=2}$, and strictly semistable otherwise. If~$d$ is not a multiple of ${n-1}$, there is no such family.

Moreover, if~$d$ is not a multiple of ${n-1}$, and ${f_1,\ldots,f_n}$ is any family of homogeneous polynomials in $\kk[X_0,X_1]$ of degree~$d$ such that the ideal ${(f_1,\ldots,f_n)}$ is \mbox{$\mathfrak{m}$-primary}, their syzygy bundle ${\syz(f_1,\ldots,f_n)}$ is not semistable.
\end{te}

This result gives a complete answer to problem~\ref{6.9} in case ${N=1}$.
Cases ${N=2}$ and ${N\ge3}$ are the object of chapters~\ref{t2} and~\ref{t3}, respectively.

\section{Moduli spaces of stable vector bundles}
\label{t1moduli}

Here is a formal definition of a moduli space, as described in \cite[chapter~II, 4.1]{OSS80}, adapted to the environment of schemes over~$\kk$.

Let $S$ be a scheme over~$\kk$ and ${r,c_1,\ldots,c_s\in\zz}$ be such that ${r\ge1}$ and ${s=\min(r,N)}$. We say that a family $E$ of stable vector bundles of rank~$r$ over $\pp^N$ with Chern classes ${c_1,\ldots,c_s}$ is \emph{parameterised} by~$S$ if~$E$ is a rank~$r$ vector bundle over ${S\times\pp^N}$ such that for all ${s\in S}$ the bundle
\[E(s):=E_{|\{s\}\times\pp^N}\]
over ${\pp^N\cong\{s\}\times\pp^N}$ is stable and, for each ${i\in\{1,\ldots,s\}}$, its $i$th Chern class is $c_i\big(E(s)\big)=c_i$.

Let $\xymatrix{p:S\times\pp^N\ar[r]&S}$ be the projection onto the first factor. Two families $E$ and $E'$ parameterised by~$S$ are \emph{equivalent} if there is a line bundle~$L$ over~$S$ such that ${E'\cong E\otimes p^*L}$. We denote by ${\mathcal{F}_{r,c_1,\ldots,c_s}(S)}$ the set of equivalence classes of families of stable bundles of rank~$r$ over~$\pp^N$ with Chern classes ${c_1,\ldots,c_s}$ and parameterised by~$S$. This defines a contravariant functor
\[\xymatrix{\mathcal{F}_{r,c_1,\ldots,c_s}:\mathcal{S}\mathit{ch}_\kk\ar[r]&
    \mathcal{S}\mathit{ets}}\]
from the category of schemes over~$\kk$ to the category of sets.

\begin{defin}
A \emph{fine moduli space} for stable rank~$r$ vector bundles over~$\pp^N$ with given Chern classes ${c_1,\ldots,c_s\in\zz}$ is a scheme~$M$ together with a vector bundle~$U$ over ${M\times\pp^N}$ such that the contravariant functor $\mathcal{F}_{r,c_1,\ldots,c_s}$ is represented by $(M,U)$.
\end{defin}

This is the nicest definition of a moduli space for vector bundles. If it exists, it is unique up to isomorphism. Unfortunately in many cases no such space exists, and we have to weaken the conditions to find a scheme still able to parameterise a given family. This leads to the following definition.

\begin{defin}
A scheme~$M$ is a \emph{coarse moduli space} for $\mathcal{F}_{r,c_1,\ldots,c_s}$ if the following conditions are satisfied:
\begin{enumerate}
  \item there is a natural transformation of contravariant functors
      \[\xymatrix{\mathcal{F}_{r,c_1,\ldots,c_s}\ar[r]&\mathrm{Hom}(-,M)}\]
      which is bijective for any (reduced) point $x_0$;
  \item for every scheme~$N$ and every natural transformation
      \[\xymatrix{\mathcal{F}_{r,c_1,\ldots,c_s}\ar[r]&\mathrm{Hom}(-,N)}\]
      there is a unique morphism $\xymatrix{f:M\ar[r]&N}$ for
      which the diagram
      \[\xymatrix{\mathcal{F}_{r,c_1,\ldots,c_s}\ar[rr]\ar[rd]&&
            \mathrm{Hom}(-,M)\ar[ld]^{f_*}\\
        &\mathrm{Hom}(-,N)}\]
      commutes.
\end{enumerate}
\end{defin}

Maruyama showed in \cite{Mar77} that GM\h{stable} rank~$r$ vector bundles over~$\pp^N$ with given Chern classes ${c_1,\ldots,c_s\in\zz}$ can be parameterised by a coarse moduli space ${M(r,c_1,\ldots,c_s)}$. Since stable bundles are GM\h{stable}, the existence of a moduli space parameterising stable bundles is guaranteed. However, very little is known about this space. In most cases there is no information on its dimension, its number of irreducible components, or whether it is non\h{singular} in a generic point.

The following proposition gives some information about the local structure of the moduli space.

\begin{p}\label{Har78,4.1}
\cite[proposition~4.1]{Har78}
Let $E$ be a stable bundle on a nonsingular projective variety~$X$. Then $H^1(X,\mathcal{E}nd\ E)$ is naturally isomorphic to the Zariski tangent space of the moduli scheme~$M$ at the point corresponding to~$E$. If ${H^2(X,\mathcal{E}nd\ E)=0}$, then~$M$ is nonsingular at that point and its dimension is equal to ${\dim H^1(X,\mathcal{E}nd\ E)}$.
\end{p}

\begin{defin}
We say that a stable vector bundle~$E$ over a scheme~$X$ is \emph{unobstructed} if its moduli space~$M$ is smooth at the point corresponding to~$E$.
\end{defin}

\begin{ob}\label{obstruction}
Since the obstruction space of the local ring $\OO_{M,[E]}$ is a subspace of ${H^2(X,\mathcal{E}nd\ E)}$, when this group is zero, $E$ is unobstructed.
\end{ob}

In chapter~\ref{t4} we will study the irreducible component of the moduli space containing the points corresponding to stable syzygy bundles.

\chapter{Monomials in three variables}
\label{t2}

In this chapter an answer to problem~\ref{6.9} is presented for ${N=2}$. Stability is guaranteed in all cases but ${(N,d,n)=(2,2,5)}$, in which only semistability is obtained.

The chapter is divided into different sections, according to different values of~$n$. Recall that for ${N=2}$, we have ${3\le n\le\tbinom{d+2}{2}}$. For the first cases, with ${3\le n\le18}$, an individual answer is given in section~\ref{t2first}. As explained below, monomials in three variables can be sketched in a triangle, in a suitable way for the purposes of this problem.  For this reason, the key to solve cases ${18<n\le d+2}$, still in section~\ref{t2first} is to choose a good position for $T$ monomials in this triangle, where $T$ is the greatest triangular number not greater than~$n$. Cases ${d+2<n\le3d}$ are solved in section~\ref{t2sides} by taking a triangle with one complete side and filling the other two sides. The last cases, in section~\ref{t2interior}, for ${3d<n\le\tbinom{d+2}{2}}$, are done by taking a triangle with all sides and filling its interior in a suitable way. Finally, in section~\ref{t2main}, summing up all results in this chapter, the main theorem is presented.

\bigskip

\paragraph{Sketching monomials in a triangle.} Monomials in $K\left[X_0,X_1,X_2\right]$ of a given degree can be sketched in a triangle as in the following examples for degrees two, three and five.

\newpage

{\small%
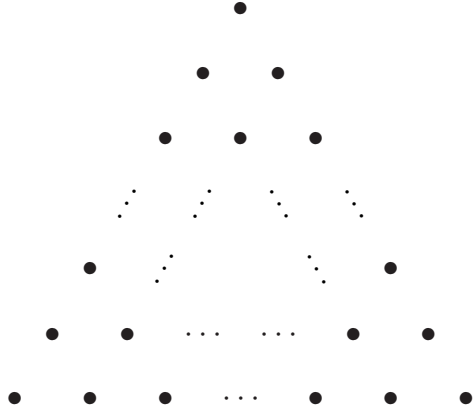
\begin{align*}
\\
\xymatrix@C=\meiabasex@R=\meiaalturax@!0{
\\
&&&{X_2}^2\\
\\
&&{X_0}X_2&&{X_1}X_2\\
\\
&{X_0}^2&&{X_0}X_1&&{X_1}^2\\
\\
}
&&\phantom{mmm}&&
\xymatrix@C=\meiabasex@R=\meiaalturax@!0{
&&&{X_2}^3\\
\\
&&{X_0}{X_2}^2&&{X_1}{X_2}^2\\
\\
&{X_0}^2X_2&&{X_0}X_1X_2&&{X_1}^2X_2\\
\\
{X_0}^3&&{X_0}^2X_1&&{X_0}{X_1}^2&&{X_1}^3\\
}
\\
\end{align*}
\begin{align*}
\\
\xymatrix@C=\meiabasexgrande@R=\alturaxgrande@!0{
&&&&&{X_2}^5\\
&&&&{X_0}{X_2}^4&&{X_1}{X_2}^4\\
&&&{X_0}^2{X_2}^3&&{X_0}X_1{X_2}^3&&{X_1}^2{X_2}^3\\
&&{X_0}^3{X_2}^2&&{X_0}^2X_1{X_2}^2&&{X_0}{X_1}^2{X_2}^2&&
    {X_1}^3{X_2}^2\\
&{X_0}^4X_2&&{X_0}^3X_1X_2&&{X_0}^2{X_1}^2X_2&&
    {X_0}{X_1}^3X_2&&{X_1}^4X_2\\
{X_0}^5&&{X_0}^4X_1&&{X_0}^3{X_1}^2&&{X_0}^2{X_1}^3&&
    {X_0}{X_1}^4&&{X_1}^5\\
}
\\
\end{align*}
}
The general case, for degree~$d$, is sketched in figure~\ref{trix}.

\begin{figure}[p]
{\small%
\[
\xymatrix@C=\meiabasexgrande@R=\alturaxgrande@!0{
&&&&&&{X_2}^d\\
&&&&&{X_0}{X_2}^{d-1}&&{X_1}{X_2}^{d-1}\\
&&&&{X_0}^2{X_2}^{d-2}&&{X_0}{X_1}{X_2}^{d-2}&&
    {X_1}^2{X_2}^{d-2}\\
&&&*[@!60]{\cdots}&&*[@!60]{\cdots}&&
    *[@!-60]{\cdots}&&*[@!-60]{\cdots}\\
&&{{X_0}^{d-2}{X_2}^2}&&*[@!60]{\cdots}&&&&
    *[@!-60]{\cdots}&&{X_1}^{d-2}{X_2}^2\\
&{X_0}^{d-1}X_2&&{X_0}^{d-2}X_1X_2&&{\cdots}&&{\cdots}&&
    {X_0}{X_1}^{d-2}X_2&&{X_1}^{d-1}X_2\\
{X_0}^d&&{X_0}^{d-1}X_1&&{X_0}^{d-2}{X_1}^2&&
    {\cdots}&&{X_0}^2{X_1}^{d-2}&&{X_0}{X_1}^{d-1}&&{X_1}^d\\
}
\]
}%
\caption{monic monomials in $K\left[X_0,X_1,X_2\right]$ of degree~$d$.}\label{trix}
\end{figure}

For the sake of simplicity and economy of space, we can sketch the triangle in figure~\ref{trix} as

\[
\xymatrix@C=\meiabasegrande@R=\alturagrande@!0{
&&&&&&{\bullet}\\
&&&&&{\bullet}&&{\bullet}\\
&&&&{\bullet}&&{\bullet}&&
    {\bullet}\\
&&&*[@!60]{\cdots}&&*[@!60]{\cdots}&&
    *[@!-60]{\cdots}&&*[@!-60]{\cdots}\\
&&{\bullet}&&*[@!60]{\cdots}&&&&
    *[@!-60]{\cdots}&&{\bullet}\\
&{\bullet}&&{\bullet}&&{\cdots}&&{\cdots}&&
    {\bullet}&&{\bullet}\\
{\bullet}&&{\bullet}&&{\bullet}&&
    {\cdots}&&{\bullet}&&{\bullet}&&{\bullet}\\
}
\]

Once arranged in this manner, the closer two monomials are, the higher the degree of their greatest common divisor is. We can see this more precisely if we consider the following graph.

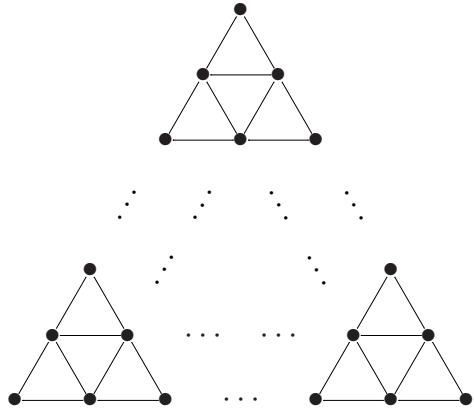
\[
\xymatrix@C=\meiabasegrande@R=\alturagrande@!0@M=0pt{
&&&&&&{\bullet}\ar@{-}[dl]\ar@{-}[dr]\\
&&&&&{\bullet}\ar@{-}[dl]\ar@{-}[dr]\ar@{-}[rr]&&
    {\bullet}\ar@{-}[dl]\ar@{-}[dr]\\
&&&&{\bullet}\ar@{-}[rr]&&{\bullet}\ar@{-}[rr]&&
    {\bullet}\\
&&&*[@!60]{\cdots}&&*[@!60]{\cdots}&&
    *[@!-60]{\cdots}&&*[@!-60]{\cdots}\\
&&{\bullet}\ar@{-}[dl]\ar@{-}[dr]&&*[@!60]{\cdots}&&&&
    *[@!-60]{\cdots}&&
    {\bullet}\ar@{-}[dl]\ar@{-}[dr]\\
&{\bullet}\ar@{-}[dl]\ar@{-}[dr]\ar@{-}[rr]&&
    {\bullet}\ar@{-}[dl]\ar@{-}[dr]&&{\cdots}&&{\cdots}&&
    {\bullet}\ar@{-}[dl]\ar@{-}[dr]\ar@{-}[rr]&&
    {\bullet}\ar@{-}[dl]\ar@{-}[dr]\\
{\bullet}\ar@{-}[rr]&&{\bullet}\ar@{-}[rr]&&{\bullet}&&
    {\cdots}&&{\bullet}\ar@{-}[rr]&&{\bullet}\ar@{-}[rr]&&
    {\bullet}\\
}
\]
Two monomials which are next to each other will admit a greatest common divisor of degree ${d-1}$, and in general if~$r$ is the distance between two monomials (measured as the number of edges in a minimal path connecting them), they will admit a greatest common divisor of degree ${d-r}$.

Since this degree is crucial in inequality~(\ref{6.6}), these triangles will help understanding whether a possible set of monomials gives rise to a stable (or semistable) syzygy bundle.

Now, for instance, the triangle
{\small%
\[
\xymatrix@C=\meiabasex@R=\alturax@!0{
&&{X_0}^2{X_2}^3\\
&{X_0}^3{X_2}^2&&{X_0}^2X_1{X_2}^2\\
}
\]
}%
represents the subset of degree~$5$ monomials that are multiples of ${X_0}^2{X_2}^2$ and the triangle
{\small%
\[
\xymatrix@C=\meiabasexgrande@R=\alturaxgrande@!0{
&&{X_0}{X_2}^4\\
&{X_0}^2{X_2}^3&&{X_0}X_1{X_2}^2\\
{X_0}^3{X_2}^2&&{X_0}^2X_1{X_2}^2&&{X_0}{X_1}^2{X_2}^2\\
}
\]
}%
represents the subset of degree~$5$ monomials that are multiples of $X_0{X_2}^2$. As we will see in proposition~\ref{4to18}, the set
\[\I_{2,5,5}:=\big\{{X_0}^5,\,{X_1}^5,\,{X_2}^5,
    \,{X_0}^2{X_1}^2X_2,\,{X_1}^2{X_2}^3\big\}\]
satisfies inequality~(\ref{6.6}). We can sketch~$\I_{2,5,5}$ in the triangle
\[
\xymatrix@C=\meiabasexpeq@R=\alturaxpeq@!0{
&&&&&{X_2}^5\\
&&&&{\circ}&&{\circ}\\
&&&{\circ}&&{\circ}&&{X_1}^2{X_2}^3\\
&&{\circ}&&{\circ}&&{\circ}&&
    {\circ}\\
&{\circ}&&{\circ}&&{X_0}^2{X_1}^2X_2&&
    {\circ}&&{\circ}\\
{X_0}^5&&{\circ}&&{\circ}&&{\circ}&&
    {\circ}&&{X_1}^5\\
}
\]
where the circles~($\circ$) represent the monic monomials of degree~$5$ that are absent from~$\I_{2,5,5}$. In a more economic fashion, this triangle becomes
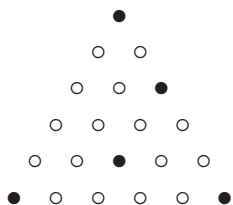
\[
\xymatrix@C=\meiabase@R=\altura@!0{
&&&&&{\bullet}\\
&&&&{\circ}&&{\circ}\\
&&&{\circ}&&{\circ}&&{\bullet}\\
&&{\circ}&&{\circ}&&{\circ}&&
    {\circ}\\
&{\circ}&&{\circ}&&{\bullet}&&
    {\circ}&&{\circ}\\
{\bullet}&&{\circ}&&{\circ}&&{\circ}&&
    {\circ}&&{\bullet}\\
}
\]
and from here we easily see that in order to check inequality~(\ref{6.6}), the most relevant subsets of~$\I_{2,5,5}$ are
\[\big\{{X_0}^2{X_1}^2X_2,\,{X_1}^2{X_2}^3\big\}
    \qquad\mbox{and}\qquad
    \big\{{X_1}^5,\,{X_0}^2{X_1}^2X_2,
    \,{X_1}^2{X_2}^3\big\}\mbox{,}\]
since they are, in cardinalities two and three respectively, the ones admitting a greatest common divisor of highest degree.

\bigskip

In each case, both in this chapter and the next one, except when stated otherwise, we will adopt the following strategy:

\begin{st}\label{X_0}
For each given~$d$ and~$n$, we choose a set of~$n$ monomials~$\I_{N,d,n}$ such that for ${0<d_J<d}$, no monomial of degree~$d_J$ divides a greater number of monomials in~$\I_{N,d,n}$ than~${X_0}^{d_J}$.
\end{st}

This will make proofs simpler, since for each instance of the degree~$d_J$, we will have to check inequality~(\ref{6.6}) only for the subset of multiples of~${X_0}^{d_J}$ in~$\I_{N,d,n}$.

\section{First cases}\label{t2first}

In this section cases ${3\le n\le d+2}$ are solved. For each of the first ones, with ${3\le n\le18}$, an individual answer is given, by choosing a set of monomials sufficiently spread apart in the triangle. For the remaining ones, with ${18<n\le d+2}$, a general answer is presented, keeping the principle of spreading most of the monomials as evenly as possible in the triangle. To this end, we choose the highest triangular number~$T$ not greater than~$n$, and lay~$T$ monomials in a triangular net, placing the others on the sides of the triangle.

\bigskip

The first case in this section (${n=3}$) is already known. As we have seen in the previous chapter, a syzygy bundle~$\syz(f_1,f_2,f_3)$ defined by three monomials $f_1,f_2,f_3$ of degree~$d$ in $\kk\left[X_0,X_1,X_2\right]$ such that the ideal $(f_1,f_2,f_3)$ is \mbox{$\mathfrak{m}$-primary} admits a presenting sequence
\[\xymatrix{0\ar[r]&\syz(f_1,f_2,f_3)\ar[r]&\OO(-d)^3
    \ar[rr]^-{f_1,f_2,f_3}&&\OO\ar[r]&0\mbox{.}}\]
Therefore its dual fits in a sequence
\[\xymatrix{0\ar[r]&\OO\ar[r]&\OO(d)^3
    \ar[r]&\syz(f_1,f_2,f_3)^\vee\ar[r]&0\mbox{,}}\]
and hence is stable, by theorem~2.7 in~\cite{BS92}. Since a vector bundle is stable if its dual is stable (see lemma~\ref{OSS1.2.4}), $\syz(f_1,f_2,f_3)$ is stable. The proof presented here is a simple example of an application of Brenner's criterion (corollary~\ref{cor6.6}).

\begin{p}\label{caseminimaln}
The syzygy bundle $\mathrm{Syz}\big({X_0}^d,\,{X_1}^d,\,{X_2}^d\big)$ is stable in~$\pp^2$.
\end{p}
\begin{dem}
If we set $\I_{2,d,3}:=\big\{{X_0}^d,\,{X_1}^d,\,{X_2}^d\big\}$, any subset~$J$ of~$\I_{2,d,3}$ with more than one element will have the constant polynomial~$1$ as greater common divisor of its monomials. Therefore $d_J=0$ and since the sequence $\left(a_{d,j}\right)_{j\ge2}$ is monotonically increasing (see remark~\ref{adj}), inequality~(\ref{6.6}) strictly holds.
\end{dem}

\begin{p}\label{4to18}
For any integers~$d$ and~$n$ such that $4\le n\le 18$ and ${n\le d+2}$, there is a set~$\I_{2,d,n}$ of~$n$ monomials in $K\left[X_0,X_1,X_2\right]$ of degree~$d$ such that the corresponding syzygy bundle is stable.
\end{p}
\begin{dem}
We shall give a particular solution for each value of~$n$. In each case, the monomials
${X_0}^d$, ${X_1}^d$, and ${X_2}^d$ will belong to the set~$\I_{2,d,n}$, and therefore the ideal generated by this set is \mbox{$\mathfrak{m}$-primary}.

Let $e_0$, $e_1$ and $e_2$ be integers such that
\[e_0+e_1+e_2=d\mbox{,}\qquad e_0\ge e_1\ge e_2\mbox{,}\qquad
    \mbox{and}\qquad e_0-e_2\le1.\]
In particular, ${e_0=\left\lceil\tfrac{d}{3}\right\rceil}$.

\bigskip\noindent{\sc Case ${n=4}$.}\\
In case ${n=4}$ (${d\ge2}$), consider the set
\[\I_{2,d,4}:=\big\{{X_0}^d,\,{X_1}^d,\,{X_2}^d,
    \,{X_0}^{e_0}{X_1}^{e_1}{X_2}^{e_2}\big\}.\]

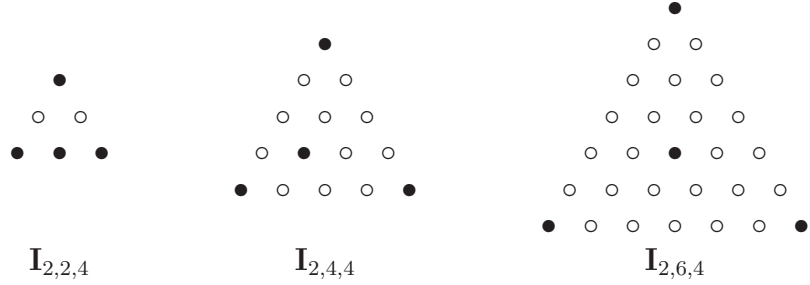
\begin{align*}
\xymatrix@C=\meiabase@R=\altura@!0{
\\
\\
&&{\bullet}\\
&{\circ}&&{\circ}\\
{\bullet}&&{\bullet}&&{\bullet}\\
\\
\\
&&\I_{2,2,4}
}
&&
\xymatrix@C=\meiabase@R=\altura@!0{
\\
&&&&{\bullet}\\
&&&{\circ}&&{\circ}\\
&&{\circ}&&{\circ}&&{\circ}\\
&{\circ}&&{\bullet}&&{\circ}&&{\circ}\\
{\bullet}&&{\circ}&&{\circ}&&
    {\circ}&&{\bullet}\\
\\
&&&&\I_{2,4,4}
}
&&
\xymatrix@C=\meiabase@R=\altura@!0{
&&&&&&{\bullet}\\
&&&&&{\circ}&&{\circ}\\
&&&&{\circ}&&{\circ}&&{\circ}\\
&&&{\circ}&&{\circ}&&{\circ}&&{\circ}\\
&&{\circ}&&{\circ}&&{\bullet}&&
    {\circ}&&{\circ}\\
&{\circ}&&{\circ}&&{\circ}&&
    {\circ}&&{\circ}&&{\circ}\\
{\bullet}&&{\circ}&&{\circ}&&{\circ}&&{\circ}&&
    {\circ}&&{\bullet}\\
&&&&&&\I_{2,6,4}
}
\end{align*}

To verify inequality~(\ref{6.6}), since~$e_0$ is the greatest exponent in the last monomial, it is enough to consider the case
\[J:=\big\{{X_0}^d,\,{X_0}^{e_0}{X_1}^{e_1}{X_2}^{e_2}\big\}.\]
In this case, if~$k$ is the cardinality of~$J$, we have ${k=2}$ and
\[(d-d_J)n+d_J-dk=4(d-e_0)+e_0-2d=2d-3e_0.\]
If ${d=2}$, we get ${e_0=1}$ and ${2d-3e_0=1>0}$. If ${d>2}$, we get
\[2d-3e_0\ge d-2>0.\]
In any case, inequality~(\ref{6.6}) is strictly satisfied.

\bigskip\noindent{\sc Case ${n=5}$.}\\
In case ${n=5}$ (${d\ge3}$), we will make an exception to strategy~\ref{X_0}. Let ${i:=\left\lceil\tfrac{d}{2}\right\rceil}$ and consider the set
\[\I_{2,d,5}:=\big\{{X_0}^d,\,{X_1}^d,\,{X_2}^d,
    \,{X_0}^{e_0}{X_1}^{e_1}{X_2}^{e_2},
    \,{X_1}^{d-i}{X_2}^i\big\}.\]

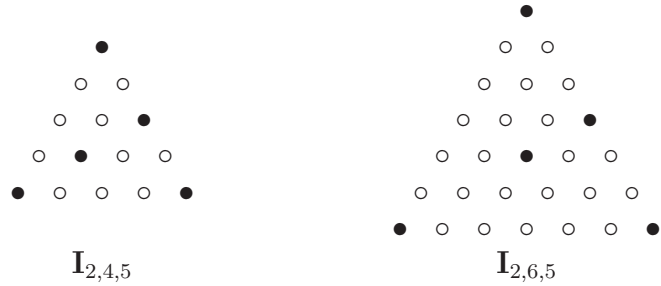
\begin{align*}
\xymatrix@C=\meiabase@R=\altura@!0{
\\
&&&&{\bullet}\\
&&&{\circ}&&{\circ}\\
&&{\circ}&&{\circ}&&{\bullet}\\
&{\circ}&&{\bullet}&&{\circ}&&{\circ}\\
{\bullet}&&{\circ}&&{\circ}&&
    {\circ}&&{\bullet}\\
\\
&&&&\I_{2,4,5}
}
&&
\xymatrix@C=\meiabase@R=\altura@!0{
&&&&&&{\bullet}\\
&&&&&{\circ}&&{\circ}\\
&&&&{\circ}&&{\circ}&&{\circ}\\
&&&{\circ}&&{\circ}&&{\circ}&&{\bullet}\\
&&{\circ}&&{\circ}&&{\bullet}&&
    {\circ}&&{\circ}\\
&{\circ}&&{\circ}&&{\circ}&&
    {\circ}&&{\circ}&&{\circ}\\
{\bullet}&&{\circ}&&{\circ}&&{\circ}&&{\circ}&&
    {\circ}&&{\bullet}\\
&&&&&&\I_{2,6,5}
}
\end{align*}

Note that since ${\tfrac{d-2}{3}\le e_2\le\tfrac{d}{3}}$ and ${\tfrac{d}{2}\le i\le\tfrac{d+1}{2}}$, we always have ${i<e_2}$. To compare ${d-i}$ with $e_1$, note that if ${d=3}$ or ${d=4}$, we get ${i=2}$ and ${e_1=1}$, therefore ${d-i\ge1=e_1}$; if ${d\ge5}$, since ${\tfrac{d-1}{3}\le e_1\le\tfrac{d+1}{3}}$, we get ${d-i-e_1\ge\tfrac{1}{6}(d-5)}$, therefore we have again ${d-i\ge e_1}$.

To verify inequality~(\ref{6.6}), let us check by possible cardinalities of subsets~$J$ of~$\I_{2,d,5}$. In cardinality two, let us observe that the greatest common divisor of ${X_0}^{e_0}{X_1}^{e_1}{X_2}^{e_2}$ and ${X_1}^{d-i}{X_2}^i$ is ${X_1}^{e_1}{X_2}^{e_2}$, since we always have ${d-i\ge e_1}$ and ${i\ge e_2}$, for ${d\ge3}$. Now if ${d=3}$ or ${d=4}$, we get ${i=2}$ and ${e_1=e_2=1}$, therefore ${e_1+e_2-i=0}$; if ${d=5}$, we get ${i=3}$, ${e_1=2}$ and ${e_2=1}$, therefore ${e_1+e_2-i=0}$; if ${d=6}$, we get ${i=3}$ and ${e_1=e_2=2}$, therefore ${e_1+e_2-i=1}$; if ${d\ge7}$, we get
\[e_1+e_2-i=d-e_0-i\ge\tfrac{1}{6}(d-7).\]
Therefore the two closest monomials are ${X_0}^{e_0}{X_1}^{e_1}{X_2}^{e_2}$ and ${X_1}^{d-i}{X_2}^i$. Therefore we get
\begin{align*}
(d-d_J)n+d_J-dk&=5(d-e_1-e_2)+e_1+e_2-2d\\
    &=5e_0+e_1+e_2-2d=4e_0-d\ge\tfrac{4d}{3}-d=\tfrac{d}{3}>0.
\end{align*}
In cardinality three, the subset admitting a greatest common divisor of highest degree is
\[J:=\big\{{X_1}^d,\,{X_0}^{e_0}{X_1}^{e_1}{X_2}^{e_2},
    \,{X_1}^{d-i}{X_2}^i\big\}.\]
Its greatest common divisor is ${X_1}^{e_1}$ and we get
\begin{align*}
(d-d_J)n+d_J-dk&=5(d-e_1)+e_1-3d\\
    &=2d-4e_1\ge2d-\tfrac{4(d+1)}{3}=\tfrac{2(d-2)}{3}>0.
\end{align*}
In cardinality greater than three, the greatest common divisor is always of degree zero therefore we have nothing to check.

\bigskip\noindent{\sc Case ${n=6}$.}\\
In case ${n=6}$ (${d\ge4}$), consider the set
\[\I_{2,d,6}:=\big\{{X_0}^d,\,{X_1}^d,\,{X_2}^d,\,{X_0}^{e_0}{X_1}^{d-e_0},
    \,{X_0}^{d-e_0}{X_2}^{e_0},\,{X_1}^{e_0}{X_2}^{d-e_0}\big\}.\]

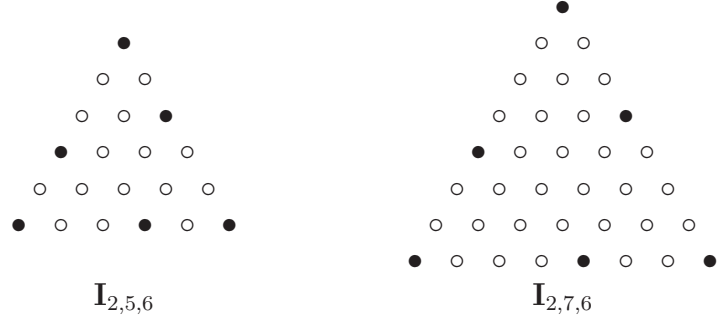
\begin{align*}
\xymatrix@C=\meiabase@R=\altura@!0{
\\
&&&&&{\bullet}\\
&&&&{\circ}&&{\circ}\\
&&&{\circ}&&{\circ}&&{\bullet}\\
&&{\bullet}&&{\circ}&&{\circ}&&{\circ}\\
&{\circ}&&{\circ}&&{\circ}&&{\circ}&&{\circ}\\
{\bullet}&&{\circ}&&{\circ}&&{\bullet}&&{\circ}&&{\bullet}\\
\\
&&&&&\I_{2,5,6}
}
&&
\xymatrix@C=\meiabase@R=\altura@!0{
&&&&&&&{\bullet}\\
&&&&&&{\circ}&&{\circ}\\
&&&&&{\circ}&&{\circ}&&{\circ}\\
&&&&{\circ}&&{\circ}&&{\circ}&&{\bullet}\\
&&&{\bullet}&&{\circ}&&{\circ}&&
    {\circ}&&{\circ}\\
&&{\circ}&&{\circ}&&{\circ}&&
    {\circ}&&{\circ}&&{\circ}\\
&{\circ}&&{\circ}&&{\circ}&&{\circ}&&{\circ}&&
    {\circ}&&{\circ}\\
{\bullet}&&{\circ}&&{\circ}&&{\circ}&&{\bullet}&&{\circ}&&
    {\circ}&&{\bullet}\\
&&&&&&&\I_{2,7,6}
}
\end{align*}

If ${0<d_J\le e_0}$, the multiples of~${X_0}^{d_J}$ in~$\I_{2,d,6}$ are the monomials in the set
\[J:=\big\{{X_0}^d,\,{X_0}^{e_0}{X_1}^{d-e_0},
    \,{X_0}^{d-e_0}{X_2}^{e_0}\big\}.\]
Therefore we have ${k=3}$ and
\begin{align*}
(d-d_J)n+d_J-dk&\ge6(d-d_J)+d_J-3d\\
    &=3d-5d_J\ge3d-5e_0\ge\tfrac{2(2d-5)}{3}>0.
\end{align*}

If ${e_0<d_J\le d-e_0}$, the multiples of~${X_0}^{d_J}$ in~$\I_{2,d,6}$ are the monomials in the set
\[J:=\big\{{X_0}^d,\,{X_0}^{d-e_0}{X_2}^{e_0}\big\}.\]
Therefore we have ${k=2}$ and
\begin{align*}
(d-d_J)n+d_J-dk&=6(d-d_J)+d_J-2d=4d-5d_J\ge
    5e_0-d\ge\tfrac{2}{3}d>0.
\end{align*}
Therefore inequality~(\ref{6.6}) is strictly satisfied.

If ${d-e_0<d_J<d}$, the only multiple of~${X_0}^{d_J}$ in~$\I_{2,d,6}$ is ${X_0}^d$ and we have nothing to check.

\bigskip\noindent{\sc Case ${n=7}$.}\\
In case ${n=7}$ (${d\ge5}$), we will again make an exception to strategy~\ref{X_0}. Consider the set
\begin{align*}
\I_{2,d,7}:=\big\{&{X_0}^d,\,{X_1}^d,\,{X_2}^d,
    \,{X_0}^{e_0}{X_1}^{e_1}{X_2}^{e_2},\\
    &{X_0}^{e_0}{X_1}^{d-e_0},\,{X_0}^{d-e_0}{X_2}^{e_0},
    \,{X_1}^{e_0}{X_2}^{d-e_0}\big\}.
\end{align*}

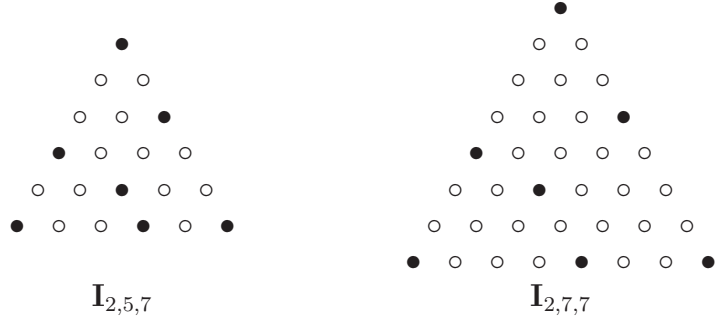
\begin{align*}
\xymatrix@C=\meiabase@R=\altura@!0{
\\
&&&&&{\bullet}\\
&&&&{\circ}&&{\circ}\\
&&&{\circ}&&{\circ}&&{\bullet}\\
&&{\bullet}&&{\circ}&&{\circ}&&{\circ}\\
&{\circ}&&{\circ}&&{\bullet}&&{\circ}&&{\circ}\\
{\bullet}&&{\circ}&&{\circ}&&{\bullet}&&{\circ}&&{\bullet}\\
\\
&&&&&\I_{2,5,7}
}
&&
\xymatrix@C=\meiabase@R=\altura@!0{
&&&&&&&{\bullet}\\
&&&&&&{\circ}&&{\circ}\\
&&&&&{\circ}&&{\circ}&&{\circ}\\
&&&&{\circ}&&{\circ}&&{\circ}&&{\bullet}\\
&&&{\bullet}&&{\circ}&&{\circ}&&
    {\circ}&&{\circ}\\
&&{\circ}&&{\circ}&&{\bullet}&&
    {\circ}&&{\circ}&&{\circ}\\
&{\circ}&&{\circ}&&{\circ}&&{\circ}&&{\circ}&&
    {\circ}&&{\circ}\\
{\bullet}&&{\circ}&&{\circ}&&{\circ}&&{\bullet}&&{\circ}&&
    {\circ}&&{\bullet}\\
&&&&&&&\I_{2,7,7}
}
\end{align*}

In fact, the exception to strategy~\ref{X_0} occurs only for ${e_0<d_J\le e_0+e_1}$, in which case ${X_0}^{e_0}{X_1}^{e_1}$ divides monomials ${X_0}^{e_0}{X_1}^{e_1}{X_2}^{e_2}$ and ${X_0}^{e_0}{X_1}^{d-e_0}$ and no monomial of degree $d_J$ divides a greater number of monomials in~$\I_{2,d,7}$.

If ${0<d_J\le e_0}$, the multiples of~${X_0}^{d_J}$ in~$\I_{2,d,7}$ are the monomials in the set
\[J:=\big\{{X_0}^d,\,{X_0}^{e_0}{X_1}^{e_1}{X_2}^{e_2},
    \,{X_0}^{e_0}{X_1}^{d-e_0},\,{X_0}^{d-e_0}{X_2}^{e_0}\big\}.\]
Therefore we have ${k=4}$ and
\begin{align*}
(d-d_J)n+d_J-dk&\ge7(d-d_J)+d_J-4d\\
    &=3d-6d_J\ge3d-6e_0\ge d-4>0.
\end{align*}

If ${e_0<d_J\le e_0+e_1}$, the multiples of ${X_0}^{e_0}{X_1}^{e_1}$ in~$\I_{2,d,7}$ are the monomials in the set
\[J:=\big\{{X_0}^{e_0}{X_1}^{e_1}{X_2}^{e_2},\,
    {X_0}^{e_0}{X_1}^{d-e_0}\big\}.\]
Therefore we have ${k=2}$ and
\begin{align*}
(d-d_J)n+d_J-dk&=7(d-d_J)+d_J-2d=5d-6d_J\\
    &\ge5d-6(e_0+e_1)=6e_2-d\ge d-4>0.
\end{align*}
Therefore inequality~(\ref{6.6}) is strictly satisfied.

If ${e_0+e_1<d_J<d}$, the only multiple of~${X_0}^{d_J}$ in~$\I_{2,d,7}$ is ${X_0}^d$ and we have nothing to check.

\bigskip\noindent{\sc Case ${n=8}$.}\\
In case ${n=8}$ (${d\ge6}$), we will again make an exception to strategy~\ref{X_0}. Consider the set
\begin{align*}
\I_{2,d,8}:=\big\{&{X_0}^d,\,{X_1}^d,\,{X_2}^d,
        \,{X_0}^{e_0}{X_1}^{e_1}{X_2}^{e_2},
        \,{X_0}^{d-e_2}{X_1}^{e_2},\\
    &{X_0}^{e_2}{X_2}^{d-e_2},
        \,{X_1}^{d-e_2}{X_2}^{e_2},
        \,{X_1}^{e_0}{X_2}^{d-e_0}\big\}.
\end{align*}

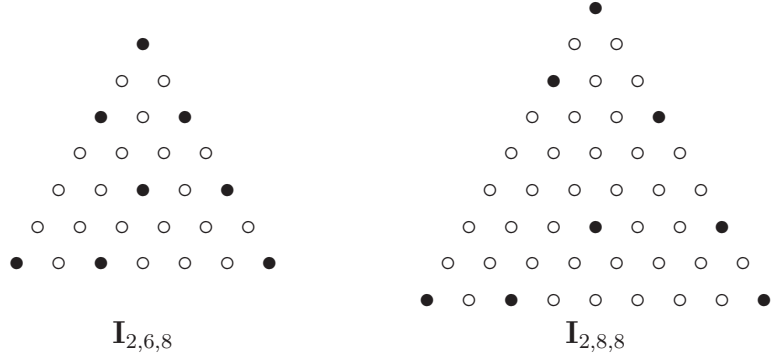
\begin{align*}
\xymatrix@C=\meiabase@R=\altura@!0{
\\
&&&&&&{\bullet}\\
&&&&&{\circ}&&{\circ}\\
&&&&{\bullet}&&{\circ}&&{\bullet}\\
&&&{\circ}&&{\circ}&&{\circ}&&{\circ}\\
&&{\circ}&&{\circ}&&{\bullet}&&{\circ}&&{\bullet}\\
&{\circ}&&{\circ}&&{\circ}&&{\circ}&&{\circ}&&{\circ}\\
{\bullet}&&{\circ}&&{\bullet}&&{\circ}&&{\circ}&&{\circ}&&{\bullet}\\
\\
&&&&&&\I_{2,6,8}
}
&&
\xymatrix@C=\meiabase@R=\altura@!0{
&&&&&&&&{\bullet}\\
&&&&&&&{\circ}&&{\circ}\\
&&&&&&{\bullet}&&{\circ}&&{\circ}\\
&&&&&{\circ}&&{\circ}&&{\circ}&&{\bullet}\\
&&&&{\circ}&&{\circ}&&{\circ}&&{\circ}&&{\circ}\\
&&&{\circ}&&{\circ}&&{\circ}&&{\circ}&&{\circ}&&{\circ}\\
&&{\circ}&&{\circ}&&{\circ}&&{\bullet}&&{\circ}&&{\circ}&&
    {\bullet}\\
&{\circ}&&{\circ}&&{\circ}&&{\circ}&&{\circ}&&{\circ}&&
    {\circ}&&{\circ}\\
{\bullet}&&{\circ}&&{\bullet}&&{\circ}&&{\circ}&&{\circ}&&
    {\circ}&&{\circ}&&{\bullet}\\
&&&&&&&&\I_{2,8,8}
}
\end{align*}

If ${0<d_J\le e_2}$, the most relevant subsets of~$\I_{2,d,8}$ are
\begin{align*}
\big\{&{X_1}^d,\,{X_0}^{e_0}{X_1}^{e_1}{X_2}^{e_2},
    \,{X_0}^{d-e_2}{X_1}^{e_2},\,{X_1}^{d-e_2}{X_2}^{e_2},
    \,{X_1}^{e_0}{X_2}^{d-e_0}\big\}\\
\intertext{and}
\big\{&{X_2}^d,\,{X_0}^{e_0}{X_1}^{e_1}{X_2}^{e_2},
    \,{X_0}^{e_2}{X_2}^{d-e_2},\,{X_1}^{d-e_2}{X_2}^{e_2},
    \,{X_1}^{e_0}{X_2}^{d-e_0}\big\}.
\end{align*}
Therefore we have ${k=5}$ and
\begin{align*}
(d-d_J)n+d_J-dk&\ge8(d-d_J)+d_J-5d\\
    &=3d-7d_J\ge3d-7e_2\ge\tfrac{2}{3}d>0.
\end{align*}

If ${e_2<d_J\le e_1}$, the most relevant set is
\[J:=\big\{{X_1}^d,\,{X_0}^{e_0}{X_1}^{e_1}{X_2}^{e_2},
    \,{X_1}^{d-e_2}{X_2}^{e_2},
    \,{X_1}^{e_0}{X_2}^{d-e_0}\big\}.\]
Therefore we have ${k=4}$ and
\begin{align*}
(d-d_J)n+d_J-dk&=8(d-d_J)+d_J-4d=4d-7d_J\\
    &\ge4d-7e_1=\tfrac{1}{3}(5d-7)>0.
\end{align*}

If ${e_1<d_J\le e_1+e_2}$, the most relevant set is
\[J:=\big\{{X_1}^d,\,{X_1}^{d-e_2}{X_2}^{e_2},
    \,{X_1}^{e_0}{X_2}^{d-e_0}\big\}.\]
Therefore we have ${k=3}$ and
\begin{align*}
(d-d_J)n+d_J-dk&=8(d-d_J)+d_J-3d=5d-7d_J\\
    &\ge5d-7(e_1+e_2)\ge\tfrac{1}{3}d>0.
\end{align*}

If ${e_1+e_2<d_J\le e_0+e_1}$, the most relevant set is
\[J:=\big\{{X_0}^d,\,{X_0}^{d-e_2}{X_1}^{e_2}\big\}.\]
Therefore we have ${k=2}$ and
\begin{align*}
(d-d_J)n+d_J-dk&=8(d-d_J)+d_J-2d=6d-7d_J\\
    &\ge6d-7(e_0+e_1)\ge\tfrac{2}{3}(2d-7)>0.
\end{align*}
Therefore inequality~(\ref{6.6}) is strictly satisfied.

If ${e_0+e_1<d_J<d}$, there is no subset of~$\I_{2,d,8}$ with cardinality greater than one admitting a greatest common divisor of degree $d_J$, and we have nothing to check.

\bigskip\noindent{\sc Case ${n=9}$.}\\
In case ${n=9}$ (${d\ge7}$), we shall look at two cases separately: ${d=8}$ and ${d\ne8}$. For the former, we will once more make an exception to strategy~\ref{X_0}. Consider the set
\begin{align*}
\I_{2,8,9}:=\big\{&{X_0}^8,\,{X_1}^8,\,{X_2}^8,\,{X_0}^3{X_1}^3{X_2}^2,
        \,{X_0}^6{X_1}^2,\\
    &{X_0}^2{X_2}^6,\,{X_0}^5{X_2}^3,
        \,{X_1}^6{X_2}^2,\,{X_1}^3{X_2}^5\big\}.
\end{align*}

\[
\xymatrix@C=\meiabase@R=\altura@!0{
&&&&&&&&{\bullet}\\
&&&&&&&{\circ}&&{\circ}\\
&&&&&&{\bullet}&&{\circ}&&{\circ}\\
&&&&&{\circ}&&{\circ}&&{\circ}&&{\bullet}\\
&&&&{\circ}&&{\circ}&&{\circ}&&{\circ}&&{\circ}\\
&&&{\bullet}&&{\circ}&&{\circ}&&{\circ}&&{\circ}&&{\circ}\\
&&{\circ}&&{\circ}&&{\circ}&&{\bullet}&&{\circ}&&{\circ}&&
    {\bullet}\\
&{\circ}&&{\circ}&&{\circ}&&{\circ}&&{\circ}&&{\circ}&&
    {\circ}&&{\circ}\\
{\bullet}&&{\circ}&&{\bullet}&&{\circ}&&{\circ}&&{\circ}&&
    {\circ}&&{\circ}&&{\bullet}\\
&&&&&&&&\I_{2,8,9}}
\]

The largest subset of~$\I_{2,8,9}$ admitting a greatest common divisor of degree at most~$2$ is
\[J:=\big\{{X_2}^8,\,{X_0}^3{X_1}^3{X_2}^2,\,{X_0}^2{X_2}^6,
    \,{X_0}^5{X_2}^3,\,{X_1}^6{X_2}^2,
    \,{X_1}^3{X_2}^5\big\}.\]
Therefore we have ${k=6}$ and
\begin{align*}
(d-d_J)n+d_J-dk&=9(8-d_J)+d_J-8\cdot6=24-8d_J\ge8>0.
\end{align*}
If ${d_J=3}$, the largest subsets of~$\I_{2,8,9}$ admitting a greatest common divisor of degree~$d_J$ are
\begin{align*}
\big\{&{X_0}^8,\,{X_0}^3{X_1}^3{X_2}^2,
        \,{X_0}^6{X_1}^2,\,{X_0}^5{X_2}^3\big\}\mbox{,}\\
\big\{&{X_1}^8,\,{X_0}^3{X_1}^3{X_2}^2,
        \,{X_1}^6{X_2}^2,\,{X_1}^3{X_2}^5\big\}\mbox{,}\\
\intertext{and}
\big\{&{X_2}^8,\,{X_0}^2{X_2}^6,\,{X_0}^5{X_2}^3,
        \,{X_1}^3{X_2}^5\big\}.
\end{align*}
In any case we have ${k=4}$ and
\begin{align*}
(d-d_J)n+d_J-dk&=9(8-3)+3-8\cdot4=16>0.
\end{align*}

There are no subsets of~$\I_{2,8,9}$ admitting a greatest common divisor of degree~$4$, and the largest subsets of~$\I_{2,8,9}$ admitting a greatest common divisor of degree~$5$ are
\begin{align*}
\big\{&{X_0}^8,\,{X_0}^6{X_1}^2,\,{X_0}^5{X_2}^3\big\}\mbox{,}\\
\big\{&{X_0}^3{X_1}^3{X_2}^2,\,{X_1}^6{X_2}^2,
    \,{X_1}^3{X_2}^5\big\}\mbox{,}\\
\intertext{and}
\big\{&{X_2}^8,\,{X_0}^2{X_2}^6,\,{X_1}^3{X_2}^5\big\}.
\end{align*}
In any case we have ${k=3}$ and
\begin{align*}
(d-d_J)n+d_J-dk&=9(8-5)+5-8\cdot3=8>0.
\end{align*}

If ${d_J=6}$, the largest subsets of~$\I_{2,8,9}$ admitting a greatest common divisor of degree~$d_J$ all have two elements, therefore we have ${k=2}$ and
\begin{align*}
(d-d_J)n+d_J-dk&=9(8-6)+6-8\cdot2=8>0.
\end{align*}

If ${d_J=7}$, all subsets of~$\I_{2,8,9}$ admitting a greatest common divisor of degree~$d_J$ have cardinality one, therefore we have nothing to check.

\bigskip

Let ${d\ne8}$ and let us write ${d=3m+t}$, with ${0\le t<3}$, and for each ${l\in\{1,2\}}$, let ${i_l:=lm+\min(l,t)}$. Since ${d\ge7}$, we get ${m\ge2}$. Consider the set
\begin{align*}
\I_{2,d,9}:=\big\{&{X_0}^d,\,{X_1}^d,\,{X_2}^d,
        \,{X_0}^{i_1}{X_1}^{d-i_1},
        \,{X_0}^{i_2}{X_1}^{d-i_2},\\
    &{X_0}^{d-i_1}{X_2}^{i_1},\,{X_0}^{d-i_2}{X_2}^{i_2},
        \,{X_1}^{i_1}{X_2}^{d-i_1},
        \,{X_1}^{i_2}{X_2}^{d-i_2}\big\}.
\end{align*}

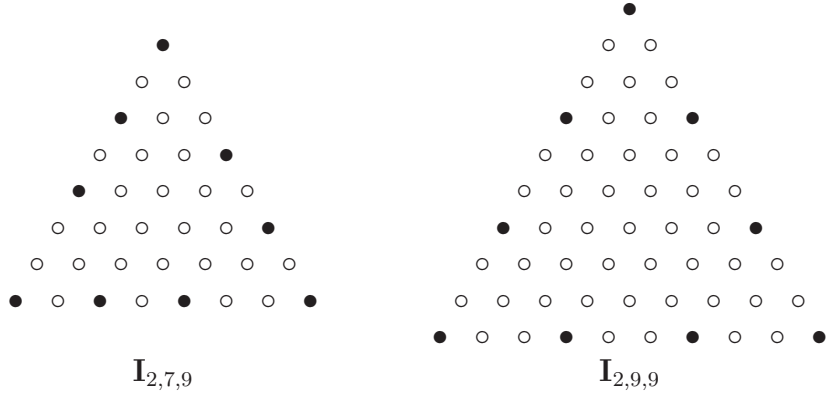
\begin{align*}
\xymatrix@C=\meiabase@R=\altura@!0{
\\
&&&&&&&{\bullet}\\
&&&&&&{\circ}&&{\circ}\\
&&&&&{\bullet}&&{\circ}&&{\circ}\\
&&&&{\circ}&&{\circ}&&{\circ}&&{\bullet}\\
&&&{\bullet}&&{\circ}&&{\circ}&&{\circ}&&{\circ}\\
&&{\circ}&&{\circ}&&{\circ}&&{\circ}&&{\circ}&&{\bullet}\\
&{\circ}&&{\circ}&&{\circ}&&{\circ}&&{\circ}&&{\circ}&&
    {\circ}\\
{\bullet}&&{\circ}&&{\bullet}&&{\circ}&&{\bullet}&&{\circ}&&
    {\circ}&&{\bullet}\\
\\
&&&&&&&\I_{2,7,9}
}
&&
\xymatrix@C=\meiabase@R=\altura@!0{
&&&&&&&&&{\bullet}\\
&&&&&&&&{\circ}&&{\circ}\\
&&&&&&&{\circ}&&{\circ}&&{\circ}\\
&&&&&&{\bullet}&&{\circ}&&{\circ}&&{\bullet}\\
&&&&&{\circ}&&{\circ}&&{\circ}&&{\circ}&&{\circ}\\
&&&&{\circ}&&{\circ}&&{\circ}&&{\circ}&&{\circ}&&{\circ}\\
&&&{\bullet}&&{\circ}&&{\circ}&&{\circ}&&{\circ}&&{\circ}&&
    {\bullet}\\
&&{\circ}&&{\circ}&&{\circ}&&{\circ}&&{\circ}&&{\circ}&&
    {\circ}&&{\circ}\\
&{\circ}&&{\circ}&&{\circ}&&{\circ}&&{\circ}&&{\circ}&&
    {\circ}&&{\circ}&&{\circ}\\
{\bullet}&&{\circ}&&{\circ}&&{\bullet}&&{\circ}&&{\circ}&&
    {\bullet}&&{\circ}&&{\circ}&&{\bullet}\\
&&&&&&&&&\I_{2,9,9}
}
\end{align*}

If ${0<d_J\le i_1}$, the multiples of~${X_0}^{d_J}$ in~$\I_{2,d,9}$ are the monomials in the set
\begin{align*}
J:=\big\{&{X_0}^d,\,{X_0}^{i_1}{X_1}^{d-i_1},
    \,{X_0}^{i_2}{X_1}^{d-i_2},\,{X_0}^{d-i_1}{X_2}^{i_1},
    \,{X_0}^{d-i_2}{X_2}^{i_2}\big\}.
\end{align*}
Therefore we have ${k=5}$ and
\begin{align*}
(d-d_J)n+d_J-dk&\ge9(d-d_J)+d_J-5d=4d-8d_J\\
    &\ge4d-8i_1=4(3m+t)-8\big(m+\min(1,t)\big)\\
    &=4m+4t-8\min(1,t)\ge4m-4>0.
\end{align*}

If ${i_1<d_J\le i_2}$, since ${d-i_2=m\le i_1}$, the multiples of~${X_0}^{d_J}$ in~$\I_{2,d,9}$ are among the monomials in the set
\begin{align*}
J:=\big\{&{X_0}^d,\,{X_0}^{i_2}{X_1}^{d-i_2},
    \,{X_0}^{d-i_1}{X_2}^{i_1}\big\}.
\end{align*}
Therefore we have ${k\le3}$ and
\begin{align*}
(d-d_J)n+d_J-dk&\ge9(d-d_J)+d_J-3d=6d-8d_J\\
    &\ge6d-8i_2\ge2m-2t.
\end{align*}
If ${d=7=3\cdot2+1}$, we get ${m=2}$ and ${t=1}$ therefore ${2m-2t=2>0}$. If ${d\ge9}$, we get ${m\ge3}$ and ${2m-2t\ge2>0}$.

If ${i_2<d_J<d}$, the only multiple of~${X_0}^{d_J}$ in~$\I_{2,d,9}$ is ${X_0}^d$ and we have nothing to check.

\bigskip\noindent{\sc Case ${n=10}$.}\\
In case ${n=10}$ (${d\ge8}$), we shall look at two cases separately: ${d=9}$ and ${d\neq9}$.

Suppose ${d=9}$ and consider the set
\begin{align*}
\I_{2,9,10}:=\big\{&{X_0}^9,\,{X_1}^9,\,{X_2}^9,\,{X_0}^3{X_1}^3{X_2}^3,
        \,{X_0}^6{X_1}^3,\,{X_0}^3{X_1}^6,\\
    &{X_0}^6{X_2}^3,\,{X_0}^3{X_2}^6,\,{X_1}^6{X_2}^3,
        \,{X_1}^3{X_2}^6\big\}.
\end{align*}

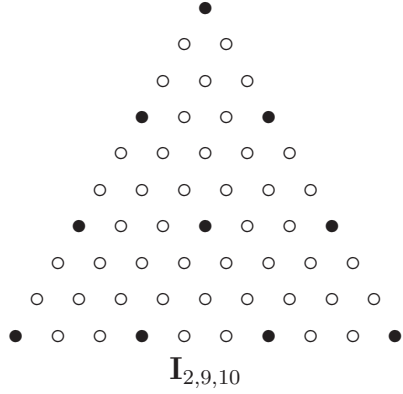
\[
\xymatrix@C=\meiabase@R=\altura@!0{
&&&&&&&&&{\bullet}\\
&&&&&&&&{\circ}&&{\circ}\\
&&&&&&&{\circ}&&{\circ}&&{\circ}\\
&&&&&&{\bullet}&&{\circ}&&{\circ}&&{\bullet}\\
&&&&&{\circ}&&{\circ}&&{\circ}&&{\circ}&&{\circ}\\
&&&&{\circ}&&{\circ}&&{\circ}&&{\circ}&&{\circ}&&{\circ}\\
&&&{\bullet}&&{\circ}&&{\circ}&&{\bullet}&&{\circ}&&{\circ}&&
    {\bullet}\\
&&{\circ}&&{\circ}&&{\circ}&&{\circ}&&{\circ}&&{\circ}&&
    {\circ}&&{\circ}\\
&{\circ}&&{\circ}&&{\circ}&&{\circ}&&{\circ}&&{\circ}&&
    {\circ}&&{\circ}&&{\circ}\\
{\bullet}&&{\circ}&&{\circ}&&{\bullet}&&{\circ}&&{\circ}&&
    {\bullet}&&{\circ}&&{\circ}&&{\bullet}\\
&&&&&&&&&\I_{2,9,10}}
\]

If ${0<d_J\le 3}$, the multiples of~${X_0}^{d_J}$ in~$\I_{2,9,10}$ are the monomials in the set
\begin{align*}
J:=\big\{&{X_0}^9,\,{X_0}^3{X_1}^3{X_2}^3,
    \,{X_0}^6{X_1}^3,\,{X_0}^3{X_1}^6,
    \,{X_0}^6{X_2}^3,\,{X_0}^3{X_2}^6\big\}.
\end{align*}
Therefore we have ${k=6}$ and
\begin{align*}
(d-d_J)n+d_J-dk&\ge10(9-d_J)+d_J-9\cdot6\ge9>0.
\end{align*}

If ${3<d_J\le6}$, the multiples of~${X_0}^{d_J}$ in~$\I_{2,9,10}$ are the monomials in the set
\begin{align*}
J:=\big\{&{X_0}^9,\,{X_0}^6{X_1}^3,\,{X_0}^6{X_2}^3\big\}.
\end{align*}
Therefore we have ${k=3}$ and
\begin{align*}
(d-d_J)n+d_J-dk&=10(9-d_J)+d_J-9\cdot3\ge9>0.
\end{align*}
Therefore inequality~(\ref{6.6}) is strictly satisfied.

If ${6<d_J<9}$, the only multiple of~${X_0}^{d_J}$ in~$\I_{2,9,10}$ is ${X_0}^d$ and we have nothing to check.

\bigskip

If ${d\neq9}$, let ${d=5m+t}$, where ${0\le t<5}$, and for each ${l\in\{1,2,3,4\}}$, let ${i_l:=lm+\min(l,t)}$. Since ${d\ge8}$, we get ${m\ge1}$.

Consider the set
\begin{align*}
\I_{2,d,10}:=\big\{&{X_0}^d,\,{X_1}^d,\,{X_2}^d,
        \,{X_0}^{i_2}{X_1}^{i_1}{X_2}^{d-i_1-i_2},
        \,{X_0}^{i_4}{X_1}^{d-i_4},
        \,{X_0}^{i_2}{X_1}^{d-i_2},\\
    &{X_0}^{i_3}{X_2}^{d-i_3},
        \,{X_0}^{i_1}{X_2}^{d-i_1},
        \,{X_1}^{i_2}{X_2}^{d-i_2},
        \,{X_1}^{i_4}{X_2}^{d-i_4}\big\}.
\end{align*}

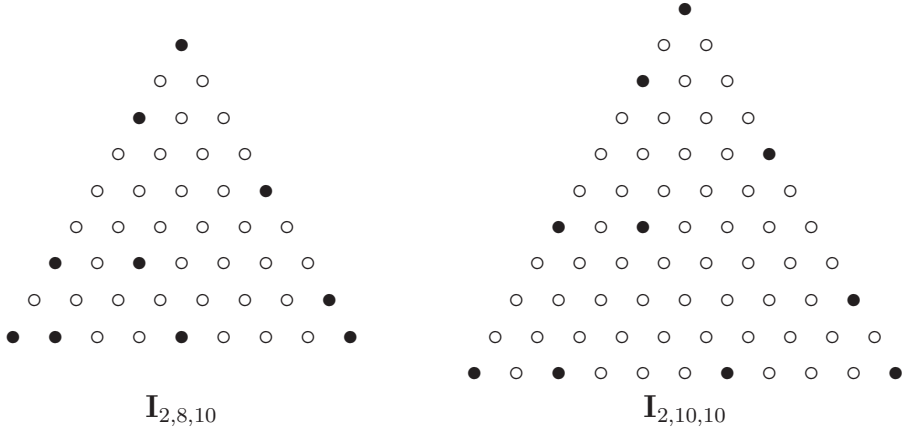
\begin{align*}
\xymatrix@C=\meiabase@R=\altura@!0{
\\
&&&&&&&&{\bullet}\\
&&&&&&&{\circ}&&{\circ}\\
&&&&&&{\bullet}&&{\circ}&&{\circ}\\
&&&&&{\circ}&&{\circ}&&{\circ}&&{\circ}\\
&&&&{\circ}&&{\circ}&&{\circ}&&{\circ}&&{\bullet}\\
&&&{\circ}&&{\circ}&&{\circ}&&{\circ}&&{\circ}&&{\circ}\\
&&{\bullet}&&{\circ}&&{\bullet}&&{\circ}&&{\circ}&&{\circ}&&
    {\circ}\\
&{\circ}&&{\circ}&&{\circ}&&{\circ}&&{\circ}&&{\circ}&&
    {\circ}&&{\bullet}\\
{\bullet}&&{\bullet}&&{\circ}&&{\circ}&&{\bullet}&&{\circ}&&
    {\circ}&&{\circ}&&{\bullet}&&\\
\\
&&&&&&&&\I_{2,8,10}
}
&&
\xymatrix@C=\meiabase@R=\altura@!0{
&&&&&&&&&&{\bullet}\\
&&&&&&&&&{\circ}&&{\circ}\\
&&&&&&&&{\bullet}&&{\circ}&&{\circ}\\
&&&&&&&{\circ}&&{\circ}&&{\circ}&&{\circ}\\
&&&&&&{\circ}&&{\circ}&&{\circ}&&{\circ}&&{\bullet}\\
&&&&&{\circ}&&{\circ}&&{\circ}&&{\circ}&&{\circ}&&{\circ}\\
&&&&{\bullet}&&{\circ}&&{\bullet}&&{\circ}&&{\circ}&&
    {\circ}&&{\circ}\\
&&&{\circ}&&{\circ}&&{\circ}&&{\circ}&&{\circ}&&{\circ}&&
    {\circ}&&{\circ}\\
&&{\circ}&&{\circ}&&{\circ}&&{\circ}&&{\circ}&&{\circ}&&
    {\circ}&&{\circ}&&{\bullet}\\
&{\circ}&&{\circ}&&{\circ}&&{\circ}&&{\circ}&&{\circ}&&
    {\circ}&&{\circ}&&{\circ}&&{\circ}\\
{\bullet}&&{\circ}&&{\bullet}&&{\circ}&&{\circ}&&{\circ}&&
    {\bullet}&&{\circ}&&{\circ}&&{\circ}&&{\bullet}\\
&&&&&&&&&&\I_{2,10,10}
}
\end{align*}

If ${0<d_J\le i_1}$, the multiples of~${X_0}^{d_J}$ in~$\I_{2,d,10}$ are the monomials in the set
\begin{align*}
J:=\big\{&{X_0}^d,\,{X_0}^{i_2}{X_1}^{i_1}{X_2}^{d-i_1-i_2},
        \,{X_0}^{i_4}{X_1}^{d-i_4},
        \,{X_0}^{i_2}{X_1}^{d-i_2},\\
    &{X_0}^{i_3}{X_2}^{d-i_3},
        \,{X_0}^{i_1}{X_2}^{d-i_1}\big\}.
\end{align*}
Therefore we have ${k=6}$ and
\begin{align*}
(d-d_J)n+d_J-dk&\ge10(d-d_J)+d_J-6d=4d-9d_J\ge4d-9i_1\\
    &\ge11m+4t-9\min(1,t)\ge11m-5>0.
\end{align*}

If ${i_1<d_J\le i_2}$, the multiples of~${X_0}^{d_J}$ in~$\I_{2,d,10}$ are the monomials in the set
\begin{align*}
J:=\big\{&{X_0}^d,\,{X_0}^{i_2}{X_1}^{i_1}{X_2}^{d-i_1-i_2},
        \,{X_0}^{i_4}{X_1}^{d-i_4},
        \,{X_0}^{i_2}{X_1}^{d-i_2},
        \,{X_0}^{i_3}{X_2}^{d-i_3}\big\}.
\end{align*}
Therefore we have ${k=5}$ and
\begin{align*}
(d-d_J)n+d_J-dk&=10(d-d_J)+d_J-5d=5d-9d_J\ge5d-9i_2.
\end{align*}
If ${d=8=5\cdot1+3}$, we get ${i_2=4}$ and ${5d-9i_2=4>0}$. If ${d\ge10}$, we get  ${m\ge2}$ and
\[5d-9i_2\ge7m+5t-9\min(2,t)\ge7m-8>0.\]

If ${i_2<d_J\le i_3}$, the multiples of~${X_0}^{d_J}$ in~$\I_{2,d,10}$ are the monomials in the set
\begin{align*}
J:=\big\{&{X_0}^d,\,{X_0}^{i_4}{X_1}^{d-i_4},
        \,{X_0}^{i_3}{X_2}^{d-i_3}\big\}.
\end{align*}
Therefore we have ${k=3}$ and
\begin{align*}
(d-d_J)n+d_J-dk&=10(d-d_J)+d_J-3d=7d-9d_J\ge7d-9i_3\\
    &=8m+7t-9\min(3,t)\ge8m-6>0.
\end{align*}

If ${i_3<d_J\le i_4}$, the multiples of~${X_0}^{d_J}$ in~$\I_{2,d,10}$ are the monomials in the set
\begin{align*}
J:=\big\{&{X_0}^d,\,{X_0}^{i_4}{X_1}^{d-i_4}\big\}.
\end{align*}
Therefore we have ${k=2}$ and
\begin{align*}
(d-d_J)n+d_J-dk&=10(d-d_J)+d_J-2d=8d-9d_J\ge8d-9i_4.
\end{align*}
If ${d=8=5\cdot1+3}$, we get ${i_4=7}$ and ${5d-9i_4=1>0}$. If ${d\ge10}$, we get  ${m\ge2}$ and ${8d-9i_4=4m-t\ge4m-4>0}$.

Therefore inequality~(\ref{6.6}) is strictly satisfied.

If ${i_4<d_J<d}$, the only multiple of~${X_0}^{d_J}$ in~$\I_{2,d,10}$ is ${X_0}^d$ and we have nothing to check.

\bigskip\noindent{\sc Case ${n=11}$.}\\
In case ${n=11}$ (${d\ge9}$), we shall look at two cases separately: ${d=12}$ and ${d\neq12}$.

Suppose ${d=12}$ and consider the set
\begin{align*}
\I_{2,12,11}:=\big\{&{X_0}^{12},\,{X_1}^{12},\,{X_2}^{12},
        \,{X_0}^9{X_1}^3,\,{X_0}^6{X_1}^6,\,{X_0}^3{X_1}^9,\\
    &{X_0}^9{X_2}^3,\,{X_0}^6{X_2}^6,\,{X_0}^3{X_2}^9,
        \,{X_1}^9{X_2}^3,\,{X_1}^6{X_2}^6\big\}.
\end{align*}

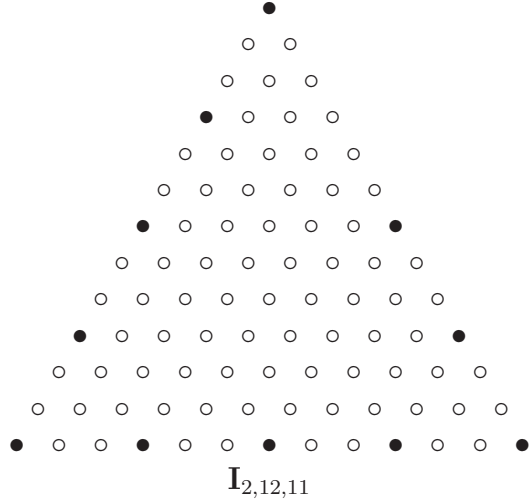
\[
\xymatrix@C=\meiabase@R=\altura@!0{
&&&&&&&&&&&&{\bullet}\\
&&&&&&&&&&&{\circ}&&{\circ}\\
&&&&&&&&&&{\circ}&&{\circ}&&{\circ}\\
&&&&&&&&&{\bullet}&&{\circ}&&{\circ}&&{\circ}\\
&&&&&&&&{\circ}&&{\circ}&&{\circ}&&{\circ}&&{\circ}\\
&&&&&&&{\circ}&&{\circ}&&{\circ}&&{\circ}&&{\circ}&&{\circ}\\
&&&&&&{\bullet}&&{\circ}&&{\circ}&&{\circ}&&{\circ}&&
    {\circ}&&{\bullet}\\
&&&&&{\circ}&&{\circ}&&{\circ}&&{\circ}&&{\circ}&&{\circ}&&
    {\circ}&&{\circ}\\
&&&&{\circ}&&{\circ}&&{\circ}&&{\circ}&&{\circ}&&{\circ}&&
    {\circ}&&{\circ}&&{\circ}\\
&&&{\bullet}&&{\circ}&&{\circ}&&{\circ}&&{\circ}&&{\circ}&&
    {\circ}&&{\circ}&&{\circ}&&{\bullet}\\
&&{\circ}&&{\circ}&&{\circ}&&{\circ}&&{\circ}&&{\circ}&&
    {\circ}&&{\circ}&&{\circ}&&{\circ}&&{\circ}\\
&{\circ}&&{\circ}&&{\circ}&&{\circ}&&{\circ}&&{\circ}&&
    {\circ}&&{\circ}&&{\circ}&&{\circ}&&{\circ}&&
    {\circ}\\
{\bullet}&&{\circ}&&{\circ}&&{\bullet}&&{\circ}&&{\circ}&&
    {\bullet}&&{\circ}&&{\circ}&&{\bullet}&&{\circ}&&
    {\circ}&&{\bullet}\\
&&&&&&&&&&&&\I_{2,12,11}
}
\]

If ${0<d_J\le 3}$, the multiples of~${X_0}^{d_J}$ in~$\I_{2,12,11}$ are the monomials in the set
\begin{align*}
J:=\big\{&{X_0}^{12},\,{X_0}^9{X_1}^3,\,{X_0}^6{X_1}^6,
        \,{X_0}^3{X_1}^9,\,{X_0}^9{X_2}^3,\,{X_0}^6{X_2}^6,
        \,{X_0}^3{X_2}^9\big\}.
\end{align*}
Therefore we have ${k=7}$ and
\begin{align*}
(d-d_J)n+d_J-dk&\ge11(12-d_J)+d_J-12\cdot7\ge6>0.
\end{align*}

If ${3<d_J\le6}$, the multiples of~${X_0}^{d_J}$ in~$\I_{2,12,11}$ are the monomials in the set
\begin{align*}
J:=\big\{&{X_0}^{12},\,{X_0}^9{X_1}^3,\,{X_0}^6{X_1}^6,
        \,{X_0}^9{X_2}^3,\,{X_0}^6{X_2}^6\big\}.
\end{align*}
Therefore we have ${k=5}$ and
\begin{align*}
(d-d_J)n+d_J-dk&=11(12-d_J)+d_J-12\cdot5\ge12>0.
\end{align*}

If ${6<d_J\le9}$, the multiples of~${X_0}^{d_J}$ in~$\I_{2,12,11}$ are the monomials in the set
\begin{align*}
J:=\big\{&{X_0}^{12},\,{X_0}^9{X_1}^3,\,{X_0}^9{X_2}^3\big\}.
\end{align*}
Therefore we have ${k=3}$ and
\begin{align*}
(d-d_J)n+d_J-dk&=11(12-d_J)+d_J-12\cdot3\ge6>0.
\end{align*}
Therefore inequality~(\ref{6.6}) is strictly satisfied.

If ${9<d_J<12}$, the only multiple of~${X_0}^{d_J}$ in~$\I_{2,12,11}$ is ${X_0}^d$ and we have nothing to check.

\bigskip

If ${d\neq12}$, let ${d=5m+t}$, where ${0\le t<5}$, and for each ${l\in\{1,2,3,4\}}$, let ${i_l:=lm+\min(l,t)}$. Since ${d\ge9}$, we get ${m\ge1}$.

Consider the set
\begin{align*}
\I_{2,d,11}:=\big\{&{X_0}^d,\,{X_1}^d,\,{X_2}^d,
        \,{X_0}^{i_2}{X_1}^{i_1}{X_2}^{d-i_1-i_2},
        \,{X_0}^{i_4}{X_1}^{d-i_4},
        \,{X_0}^{i_3}{X_1}^{d-i_3},\\
    &{X_0}^{i_2}{X_1}^{d-i_2},\,{X_0}^{i_3}{X_2}^{d-i_3},
        \,{X_0}^{i_1}{X_2}^{d-i_1},
        \,{X_1}^{i_2}{X_2}^{d-i_2},
        \,{X_1}^{i_4}{X_2}^{d-i_4}\big\}.
\end{align*}

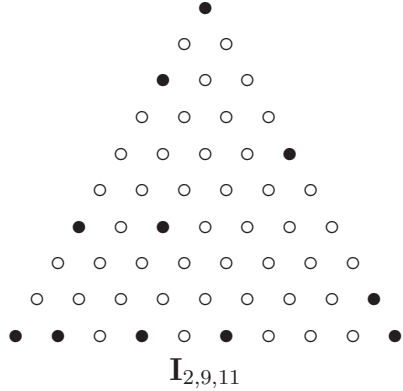
\[
\xymatrix@C=\meiabase@R=\altura@!0{
&&&&&&&&&{\bullet}\\
&&&&&&&&{\circ}&&{\circ}\\
&&&&&&&{\bullet}&&{\circ}&&{\circ}\\
&&&&&&{\circ}&&{\circ}&&{\circ}&&{\circ}\\
&&&&&{\circ}&&{\circ}&&{\circ}&&{\circ}&&{\bullet}\\
&&&&{\circ}&&{\circ}&&{\circ}&&{\circ}&&{\circ}&&{\circ}\\
&&&{\bullet}&&{\circ}&&{\bullet}&&{\circ}&&{\circ}&&{\circ}&&
    {\circ}\\
&&{\circ}&&{\circ}&&{\circ}&&{\circ}&&{\circ}&&{\circ}&&
    {\circ}&&{\circ}\\
&{\circ}&&{\circ}&&{\circ}&&{\circ}&&{\circ}&&{\circ}&&
    {\circ}&&{\circ}&&{\bullet}\\
{\bullet}&&{\bullet}&&{\circ}&&{\bullet}&&{\circ}&&{\bullet}&&
    {\circ}&&{\circ}&&{\circ}&&{\bullet}\\
&&&&&&&&&\I_{2,9,11}
}
\]

If ${0<d_J\le i_1}$, the multiples of~${X_0}^{d_J}$ in~$\I_{2,d,11}$ are the monomials in the set
\begin{align*}
J:=\big\{&{X_0}^d,\,{X_0}^{i_2}{X_1}^{i_1}{X_2}^{d-i_1-i_2},
        \,{X_0}^{i_4}{X_1}^{d-i_4},
        \,{X_0}^{i_3}{X_1}^{d-i_3},\\
    &{X_0}^{i_2}{X_1}^{d-i_2},
        \,{X_0}^{i_3}{X_2}^{d-i_3},
        \,{X_0}^{i_1}{X_2}^{d-i_1}\big\}.
\end{align*}
Therefore we have ${k=7}$ and
\begin{align*}
(d-d_J)n+d_J-dk&\ge11(d-d_J)+d_J-7d=4d-10d_J\ge4d-10i_1\\
    &\ge10m+4t-10\min(1,t)\ge10m-6>0.
\end{align*}

If ${i_1<d_J\le i_2}$, the multiples of~${X_0}^{d_J}$ in~$\I_{2,d,11}$ are the monomials in the set
\begin{align*}
J:=\big\{&{X_0}^d,\,{X_0}^{i_2}{X_1}^{i_1}{X_2}^{d-i_1-i_2},
        \,{X_0}^{i_4}{X_1}^{d-i_4},
        \,{X_0}^{i_3}{X_1}^{d-i_3},\\
    &{X_0}^{i_2}{X_1}^{d-i_2},
        \,{X_0}^{i_3}{X_2}^{d-i_3}\big\}.
\end{align*}
Therefore we have ${k=6}$ and
\begin{align*}
(d-d_J)n+d_J-dk&=11(d-d_J)+d_J-6d=5d-10d_J\ge5d-10i_2.
\end{align*}
If ${d=9=5\cdot1+4}$, we get ${i_2=4}$ and ${5d-10i_2=5>0}$. If ${d=10=5\cdot2+0}$, we get ${i_2=4}$ and hence ${5d-10i_2=10>0}$. If ${d=11=5\cdot2+1}$, we get ${i_2=5}$ and ${5d-10i_2=5>0}$. If ${d=13=5\cdot2+3}$, we get ${i_2=6}$ and ${5d-10i_2=5>0}$. If ${d=14=5\cdot2+4}$, we get ${i_2=6}$ and ${5d-10i_2=10>0}$. If ${d\ge15}$, we get  ${m\ge3}$ and
\[5d-10i_2\ge5m+5t-10\min(2,t)\ge5m-10>0.\]

If ${i_2<d_J\le i_3}$, the multiples of~${X_0}^{d_J}$ in~$\I_{2,d,11}$ are the monomials in the set
\begin{align*}
J:=\big\{&{X_0}^d,\,{X_0}^{i_4}{X_1}^{d-i_4},
        \,{X_0}^{i_3}{X_1}^{d-i_3},
        \,{X_0}^{i_3}{X_2}^{d-i_3}\big\}.
\end{align*}
Therefore we have ${k=4}$ and
\begin{align*}
(d-d_J)n+d_J-dk&=11(d-d_J)+d_J-4d=7d-10d_J\ge7d-10i_3.
\end{align*}
If ${d=9=5\cdot1+4}$, we get ${i_3=6}$ and ${7d-10i_3=3>0}$. If ${d\ge10}$, we get  ${m\ge2}$ and
\[7d-10i_3\ge5m+7t-10\min(3,t)\ge5m-9>0.\]

If ${i_3<d_J\le i_4}$, the multiples of~${X_0}^{d_J}$ in~$\I_{2,d,11}$ are the monomials in the set
\begin{align*}
J:=\big\{&{X_0}^d,\,{X_0}^{i_4}{X_1}^{d-i_4}\big\}.
\end{align*}
Therefore we have ${k=2}$ and
\begin{align*}
(d-d_J)n+d_J-dk&=11(d-d_J)+d_J-2d=9d-10d_J\ge9d-10i_4\\
    &\ge5m-t\ge5m-4>0.
\end{align*}
Therefore inequality~(\ref{6.6}) is strictly satisfied.

If ${i_4<d_J<d}$, the only multiple of~${X_0}^{d_J}$ in~$\I_{2,d,11}$ is ${X_0}^d$ and we have nothing to check.

\bigskip\noindent{\sc Case ${n=12}$.}\\
In case ${n=12}$ (${d\ge10}$), we shall look at two cases separately: ${d=11}$ and ${d\ne11}$. For the former, we will again make an exception to strategy~\ref{X_0}. Consider the set
\begin{align*}
\I_{2,11,12}:=\big\{&{X_0}^{11},\,{X_1}^{11},\,{X_2}^{11},
        \,{X_0}^8{X_1}^3,\,{X_0}^8{X_2}^3,
        \,{X_0}^5{X_1}^2{X_2}^4,\,{X_0}^4{X_1}^4{X_2}^3,\\
    &{X_0}^3{X_1}^8,\,{X_0}^3{X_2}^8,\,{X_0}^2{X_1}^5{X_2}^4,\,
        \,{X_1}^8{X_2}^3,\,{X_1}^3{X_2}^8\big\}.
\end{align*}

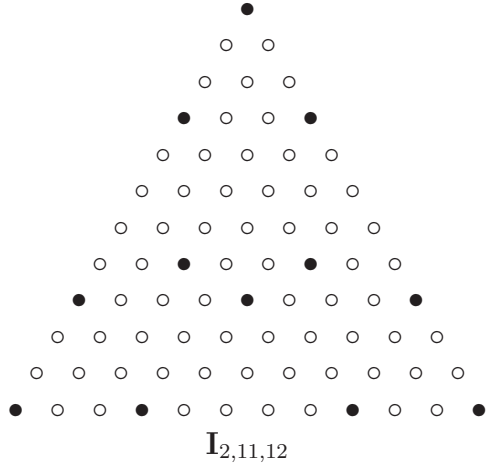
\[
\xymatrix@C=\meiabase@R=\altura@!0{
&&&&&&&&&&&{\bullet}\\
&&&&&&&&&&{\circ}&&{\circ}\\
&&&&&&&&&{\circ}&&{\circ}&&{\circ}\\
&&&&&&&&{\bullet}&&{\circ}&&{\circ}&&{\bullet}\\
&&&&&&&{\circ}&&{\circ}&&{\circ}&&{\circ}&&{\circ}\\
&&&&&&{\circ}&&{\circ}&&{\circ}&&{\circ}&&{\circ}&&{\circ}\\
&&&&&{\circ}&&{\circ}&&{\circ}&&{\circ}&&{\circ}&&{\circ}&&
    {\circ}\\
&&&&{\circ}&&{\circ}&&{\bullet}&&{\circ}&&{\circ}&&{\bullet}&&
    {\circ}&&{\circ}\\
&&&{\bullet}&&{\circ}&&{\circ}&&{\circ}&&{\bullet}&&{\circ}&&
    {\circ}&&{\circ}&&{\bullet}\\
&&{\circ}&&{\circ}&&{\circ}&&{\circ}&&{\circ}&&{\circ}&&
    {\circ}&&{\circ}&&{\circ}&&{\circ}\\
&{\circ}&&{\circ}&&{\circ}&&{\circ}&&{\circ}&&{\circ}&&
    {\circ}&&{\circ}&&{\circ}&&{\circ}&&{\circ}\\
{\bullet}&&{\circ}&&{\circ}&&{\bullet}&&{\circ}&&{\circ}&&
    {\circ}&&{\circ}&&{\bullet}&&{\circ}&&{\circ}&&{\bullet}\\
&&&&&&&&&&&\I_{2,11,12}}
\]

In fact, the exception to strategy~\ref{X_0} occurs only for ${d_J=9}$.

If ${0<d_J\le2}$, the multiples of~${X_0}^{d_J}$ in~$\I_{2,11,12}$ are the monomials in the set
\begin{align*}
J:=\big\{&{X_0}^{11},\,{X_0}^8{X_1}^3,\,{X_0}^8{X_2}^3,
        \,{X_0}^5{X_1}^2{X_2}^4,\,{X_0}^4{X_1}^4{X_2}^3,\\
    &{X_0}^3{X_1}^8,\,{X_0}^3{X_2}^8,
        \,{X_0}^2{X_1}^5{X_2}^4\big\}.
\end{align*}
Therefore we have ${k=8}$ and
\begin{align*}
(d-d_J)n+d_J-dk&\ge12(11-d_J)+d_J-11\cdot8\ge22>0.
\end{align*}

If ${d_J=3}$, the multiples of~${X_0}^{d_J}$ in~$\I_{2,11,12}$ are the monomials in the set
\begin{align*}
J:=\big\{&{X_0}^{11},\,{X_0}^8{X_1}^3,\,{X_0}^8{X_2}^3,
        \,{X_0}^5{X_1}^2{X_2}^4,\,{X_0}^4{X_1}^4{X_2}^3,
        \,{X_0}^3{X_1}^8,\,{X_0}^3{X_2}^8\big\}.
\end{align*}
Therefore we have ${k=7}$ and
\begin{align*}
(d-d_J)n+d_J-dk&\ge12(11-3)+3-11\cdot7=22>0.
\end{align*}

If ${d_J=4}$, the multiples of~${X_0}^{d_J}$ in~$\I_{2,11,12}$ are the monomials in the set
\begin{align*}
J:=\big\{&{X_0}^{11},\,{X_0}^8{X_1}^3,\,{X_0}^8{X_2}^3,
        \,{X_0}^5{X_1}^2{X_2}^4,
        \,{X_0}^4{X_1}^4{X_2}^3\big\}.
\end{align*}
Therefore we have ${k=5}$ and
\begin{align*}
(d-d_J)n+d_J-dk&\ge12(11-4)+4-11\cdot5=33>0.
\end{align*}

If ${d_J=5}$, the multiples of~${X_0}^{d_J}$ in~$\I_{2,11,12}$ are the monomials in the set
\begin{align*}
J:=\big\{&{X_0}^{11},\,{X_0}^8{X_1}^3,\,{X_0}^8{X_2}^3,
        \,{X_0}^5{X_1}^2{X_2}^4\big\}.
\end{align*}
Therefore we have ${k=4}$ and
\begin{align*}
(d-d_J)n+d_J-dk&\ge12(11-5)+5-11\cdot4=33>0.
\end{align*}

If ${5<d_J\le8}$, the multiples of~${X_0}^{d_J}$ in~$\I_{2,11,12}$ are the monomials in the set
\begin{align*}
J:=\big\{&{X_0}^{11},\,{X_0}^8{X_1}^3,\,{X_0}^8{X_2}^3\big\}.
\end{align*}
Therefore we have ${k=3}$ and
\begin{align*}
(d-d_J)n+d_J-dk&\ge12(11-d_J)+d_J-11\cdot3\ge11>0.
\end{align*}

If ${d_J=9}$, the largest subsets of~$\I_{2,11,12}$ admitting a greatest common divisor of degree~${d_J}$ are
\[
\big\{{X_0}^5{X_1}^2{X_2}^4,\,{X_0}^4{X_1}^4{X_2}^3\big\}
    \qquad\mbox{and}\qquad
\big\{{X_0}^4{X_1}^4{X_2}^3,\,{X_0}^2{X_1}^5{X_2}^4\big\}.
\]
Therefore we have ${k=2}$ and
\begin{align*}
(d-d_J)n+d_J-dk&=12(11-9)+9-11\cdot2=11>0.
\end{align*}
Therefore inequality~(\ref{6.6}) is strictly satisfied.

If ${d_J=10}$, the only multiple of~${X_0}^{d_J}$ in~$\I_{2,11,12}$ is ${X_0}^d$ and we have nothing to check.

\bigskip

If ${d\ne11}$, let ${d=4m+t}$, where ${0\le t<4}$, and for each ${l\in\{1,2,3\}}$, let\linebreak[4] ${i_l:=lm+\min(l,t)}$. Since ${d\ge10}$, we get ${m\ge2}$.

Consider the set
\begin{align*}
\I_{2,d,12}:=\big\{&{X_0}^d,\,{X_1}^d,\,{X_2}^d,
        \,{X_0}^{i_3}{X_1}^{d-i_3},
        \,{X_0}^{i_2}{X_1}^{d-i_2},
        \,{X_0}^{i_1}{X_1}^{d-i_1},
        \,{X_0}^{i_3}{X_2}^{d-i_3},\\
    &{X_0}^{i_2}{X_2}^{d-i_2},\,{X_0}^{i_1}{X_2}^{d-i_1},
        \,{X_1}^{i_1}{X_2}^{d-i_1},
        \,{X_1}^{i_2}{X_2}^{d-i_2},
        \,{X_1}^{i_3}{X_2}^{d-i_3}\big\}.
\end{align*}

\[
\xymatrix@C=\meiabase@R=\altura@!0{
&&&&&&&&&&{\bullet}\\
&&&&&&&&&{\circ}&&{\circ}\\
&&&&&&&&{\circ}&&{\circ}&&{\circ}\\
&&&&&&&{\bullet}&&{\circ}&&{\circ}&&{\bullet}\\
&&&&&&{\circ}&&{\circ}&&{\circ}&&{\circ}&&{\circ}\\
&&&&&{\circ}&&{\circ}&&{\circ}&&{\circ}&&{\circ}&&{\circ}\\
&&&&{\bullet}&&{\circ}&&{\circ}&&{\circ}&&{\circ}&&{\circ}&&
    {\bullet}\\
&&&{\circ}&&{\circ}&&{\circ}&&{\circ}&&{\circ}&&{\circ}&&
    {\circ}&&{\circ}\\
&&{\bullet}&&{\circ}&&{\circ}&&{\circ}&&{\circ}&&{\circ}&&
    {\circ}&&{\circ}&&{\bullet}\\
&{\circ}&&{\circ}&&{\circ}&&{\circ}&&{\circ}&&{\circ}&&
    {\circ}&&{\circ}&&{\circ}&&{\circ}\\
{\bullet}&&{\circ}&&{\bullet}&&{\circ}&&{\bullet}&&{\circ}&&
    {\circ}&&{\bullet}&&{\circ}&&{\circ}&&{\bullet}\\
&&&&&&&&&&\I_{2,10,12}
}
\]

If ${0<d_J\le i_1}$, the multiples of~${X_0}^{d_J}$ in~$\I_{2,d,12}$ are the monomials in the set
\begin{align*}
J:=\big\{&{X_0}^d,\,{X_0}^{i_3}{X_1}^{d-i_3},
        \,{X_0}^{i_2}{X_1}^{d-i_3},
        \,{X_0}^{i_1}{X_1}^{d-i_1},\\
    &{X_0}^{i_3}{X_2}^{d-i_3},\,{X_0}^{i_2}{X_2}^{d-i_2},
        \,{X_0}^{i_1}{X_2}^{d-i_1}\big\}.
\end{align*}
Therefore we have ${k=7}$ and
\begin{align*}
(d-d_J)n+d_J-dk&\ge12(d-d_J)+d_J-7d=5d-11d_J
        \ge5d-11i_1\\
    &\ge9m+5t-11\min(1,t)\ge9m-6>0.
\end{align*}

If ${i_1<d_J\le i_2}$, the multiples of~${X_0}^{d_J}$ in~$\I_{2,d,12}$ are the monomials in the set
\begin{align*}
J:=\big\{&{X_0}^d,\,{X_0}^{i_3}{X_1}^{d-i_3},
        \,{X_0}^{i_2}{X_1}^{d-i_3},
        \,{X_0}^{i_3}{X_2}^{d-i_3},
        \,{X_0}^{i_2}{X_2}^{d-i_2}\big\}.
\end{align*}
Therefore we have ${k=5}$ and
\begin{align*}
(d-d_J)n+d_J-dk&=12(d-d_J)+d_J-5d=7d-11d_J\ge7d-11i_2\\
        &\ge6m+7t-11\min(2,t)\ge6m-8>0.
\end{align*}

If ${i_2<d_J\le i_3}$, the multiples of~${X_0}^{d_J}$ in~$\I_{2,d,12}$ are the monomials in the set
\begin{align*}
J:=\big\{&{X_0}^d,\,{X_0}^{i_3}{X_1}^{d-i_3},
        \,{X_0}^{i_3}{X_2}^{d-i_3}\big\}.
\end{align*}
Therefore we have ${k=3}$ and
\begin{align*}
(d-d_J)n+d_J-dk&=12(d-d_J)+d_J-3d=9d-11d_J\ge9d-11i_3.
\end{align*}
If ${d=10=4\cdot2+2}$, we get ${i_3=8}$ and ${9d-11i_3=2>0}$. If ${d\ge12}$, we get  ${m\ge3}$ and
\[9d-11i_3\ge3m-2t\ge3m-6>0.\]
Therefore inequality~(\ref{6.6}) is strictly satisfied.

If ${i_3<d_J<d}$, the only multiple of~${X_0}^{d_J}$ in~$\I_{2,d,12}$ is ${X_0}^d$ and we have nothing to check.

\bigskip\noindent{\sc Case ${n=13}$.}\\
In case ${n=13}$ (${d\ge11}$), let us write ${d=4m+t}$, where ${0\le t<4}$, and for each ${l\in\{1,2,3\}}$, let ${i_l:=lm+\min(l,t)}$. Since ${d\ge11}$, we get ${m\ge2}$.

Consider the set
\begin{align*}
\I_{2,d,13}:=\big\{&{X_0}^d,\,{X_1}^d,\,{X_2}^d,
        \,{X_0}^{i_2}{X_1}^{d-i_3}{X_2}^{i_3-i_2},\\
    &{X_0}^{i_3}{X_1}^{d-i_3},\,{X_0}^{i_2}{X_1}^{d-i_2},
        \,{X_0}^{i_1}{X_1}^{d-i_1},\\
    &{X_0}^{i_3}{X_2}^{d-i_3},\,{X_0}^{i_2}{X_2}^{d-i_2},
        \,{X_0}^{i_1}{X_2}^{d-i_1},\\
    &{X_1}^{i_1}{X_2}^{d-i_1},
        \,{X_1}^{i_2}{X_2}^{d-i_2},
        \,{X_1}^{i_3}{X_2}^{d-i_3}\big\}.
\end{align*}

\[
\xymatrix@C=\meiabase@R=\altura@!0{
&&&&&&&&&&&{\bullet}\\
&&&&&&&&&&{\circ}&&{\circ}\\
&&&&&&&&&{\circ}&&{\circ}&&{\circ}\\
&&&&&&&&{\bullet}&&{\circ}&&{\circ}&&{\bullet}\\
&&&&&&&{\circ}&&{\circ}&&{\circ}&&{\circ}&&{\circ}\\
&&&&&&{\circ}&&{\circ}&&{\circ}&&{\circ}&&{\circ}&&{\circ}\\
&&&&&{\bullet}&&{\circ}&&{\circ}&&{\circ}&&{\circ}&&{\circ}&&
    {\bullet}\\
&&&&{\circ}&&{\circ}&&{\circ}&&{\circ}&&{\circ}&&{\circ}&&
    {\circ}&&{\circ}\\
&&&{\circ}&&{\circ}&&{\bullet}&&{\circ}&&{\circ}&&{\circ}&&
    {\circ}&&{\circ}&&{\circ}\\
&&{\bullet}&&{\circ}&&{\circ}&&{\circ}&&{\circ}&&{\circ}&&
    {\circ}&&{\circ}&&{\circ}&&{\bullet}\\
&{\circ}&&{\circ}&&{\circ}&&{\circ}&&{\circ}&&{\circ}&&
    {\circ}&&{\circ}&&{\circ}&&{\circ}&&{\circ}\\
{\bullet}&&{\circ}&&{\bullet}&&{\circ}&&{\circ}&&{\bullet}&&{\circ}&&
    {\circ}&&{\bullet}&&{\circ}&&{\circ}&&{\bullet}\\
&&&&&&&&&&&\I_{2,d,13}
}
\]

If ${0<d_J\le i_1}$, the multiples of~${X_0}^{d_J}$ in~$\I_{2,d,13}$ are the monomials in the set
\begin{align*}
J:=\big\{&{X_0}^d,\,{X_0}^{i_2}{X_1}^{i_1}{X_2}^{d-i_1-i_2},\\
    &{X_0}^{i_3}{X_1}^{d-i_3},\,{X_0}^{i_2}{X_1}^{d-i_3},
        \,{X_0}^{i_1}{X_1}^{d-i_1},\\
    &{X_0}^{i_3}{X_2}^{d-i_3},\,{X_0}^{i_2}{X_2}^{d-i_2},
        \,{X_0}^{i_1}{X_2}^{d-i_1}\big\}.
\end{align*}
Therefore we have ${k=8}$ and
\begin{align*}
(d-d_J)n+d_J-dk&\ge13(d-d_J)+d_J-8d=5d-12d_J
        \ge5d-12i_1\\
    &\ge8m+5t-12\min(1,t)\ge8m-7>0.
\end{align*}

If ${i_1<d_J\le i_2}$, the multiples of~${X_0}^{d_J}$ in~$\I_{2,d,13}$ are the monomials in the set
\begin{align*}
J:=\big\{&{X_0}^d,\,{X_0}^{i_2}{X_1}^{i_1}{X_2}^{d-i_1-i_2},
        \,{X_0}^{i_3}{X_1}^{d-i_3},
        \,{X_0}^{i_2}{X_1}^{d-i_3},\\
    &{X_0}^{i_3}{X_2}^{d-i_3},\,{X_0}^{i_2}{X_2}^{d-i_2}\big\}.
\end{align*}
Therefore we have ${k=6}$ and
\begin{align*}
(d-d_J)n+d_J-dk&=13(d-d_J)+d_J-6d=7d-12d_J\ge7d-12i_2.
\end{align*}
If ${d=11=4\cdot2+3}$, we get ${i_2=6}$ and ${7d-12i_2=5>0}$. If ${d\ge12}$, we get  ${m\ge3}$ and
\[7d-12i_2\ge4m+7t-12\min(2,t)\ge4m-10>0.\]

If ${i_2<d_J\le i_3}$, the multiples of~${X_0}^{d_J}$ in~$\I_{2,d,13}$ are the monomials in the set
\begin{align*}
J:=\big\{&{X_0}^d,\,{X_0}^{i_3}{X_1}^{d-i_3},
    \,{X_0}^{i_3}{X_2}^{d-i_3}\big\}.
\end{align*}
Therefore we have ${k=3}$ and
\begin{align*}
(d-d_J)n+d_J-dk&=13(d-d_J)+d_J-3d=10d-12d_J\ge10d-12i_3\\
    &=4m-2t\ge4m-6>0.
\end{align*}
Therefore inequality~(\ref{6.6}) is strictly satisfied.

If ${i_3<d_J<d}$, the only multiple of~${X_0}^{d_J}$ in~$\I_{2,d,13}$ is ${X_0}^d$ and we have nothing to check.

\bigskip\noindent{\sc Case ${n=14}$.}\\
In case ${n=14}$ (${d\ge12}$), let us write ${d=4m+t}$, where ${0\le t<4}$, and for each ${l\in\{1,2,3\}}$, let ${i_l:=lm+\min(l,t)}$. Since ${d\ge12}$, we get ${m\ge3}$.

Consider the set
\begin{align*}
\I_{2,d,14}:=\big\{&{X_0}^d,\,{X_1}^d,\,{X_2}^d,
        \,{X_0}^{i_2}{X_1}^{d-i_3}{X_2}^{i_3-i_2},
        \,{X_0}^{i_1}{X_1}^{d-i_2}{X_2}^{i_2-i_1},\\
    &{X_0}^{i_3}{X_1}^{d-i_3},\,{X_0}^{i_2}{X_1}^{d-i_2},
        \,{X_0}^{i_1}{X_1}^{d-i_1},\\
    &{X_0}^{i_3}{X_2}^{d-i_3},\,{X_0}^{i_2}{X_2}^{d-i_2},
        \,{X_0}^{i_1}{X_2}^{d-i_1},\\
    &{X_1}^{i_1}{X_2}^{d-i_1},
        \,{X_1}^{i_2}{X_2}^{d-i_2},
        \,{X_1}^{i_3}{X_2}^{d-i_3}\big\}.
\end{align*}

\[
\xymatrix@C=\meiabase@R=\altura@!0{
&&&&&&&&&&&&{\bullet}\\
&&&&&&&&&&&{\circ}&&{\circ}\\
&&&&&&&&&&{\circ}&&{\circ}&&{\circ}\\
&&&&&&&&&{\bullet}&&{\circ}&&{\circ}&&{\bullet}\\
&&&&&&&&{\circ}&&{\circ}&&{\circ}&&{\circ}&&{\circ}\\
&&&&&&&{\circ}&&{\circ}&&{\circ}&&{\circ}&&{\circ}&&{\circ}\\
&&&&&&{\bullet}&&{\circ}&&{\circ}&&{\circ}&&{\circ}&&{\circ}&&
    {\bullet}\\
&&&&&{\circ}&&{\circ}&&{\circ}&&{\circ}&&{\circ}&&{\circ}&&
    {\circ}&&{\circ}\\
&&&&{\circ}&&{\circ}&&{\circ}&&{\circ}&&{\circ}&&{\circ}&&
    {\circ}&&{\circ}&&{\circ}\\
&&&{\bullet}&&{\circ}&&{\circ}&&{\bullet}&&{\circ}&&{\circ}&&
    {\bullet}&&{\circ}&&{\circ}&&{\bullet}\\
&&{\circ}&&{\circ}&&{\circ}&&{\circ}&&{\circ}&&{\circ}&&
    {\circ}&&{\circ}&&{\circ}&&{\circ}&&{\circ}\\
&{\circ}&&{\circ}&&{\circ}&&{\circ}&&{\circ}&&{\circ}&&
    {\circ}&&{\circ}&&{\circ}&&{\circ}&&{\circ}&&{\circ}\\
{\bullet}&&{\circ}&&{\circ}&&{\bullet}&&{\circ}&&{\circ}&&
    {\bullet}&&{\circ}&&{\circ}&&{\bullet}&&{\circ}&&
    {\circ}&&{\bullet}\\
&&&&&&&&&&&&\I_{2,d,14}
}
\]

If ${0<d_J\le i_1}$, the multiples of~${X_0}^{d_J}$ in~$\I_{2,d,14}$ are the monomials in the set
\begin{align*}
J:=\big\{&{X_0}^d,\,{X_0}^{i_2}{X_1}^{i_1}{X_2}^{d-i_1-i_2},
        \,{X_0}^{i_1}{X_1}^{i_2}{X_2}^{d-i_1-i_2},\\
    &{X_0}^{i_3}{X_1}^{d-i_3},\,{X_0}^{i_2}{X_1}^{d-i_3},
        \,{X_0}^{i_1}{X_1}^{d-i_1},\\
    &{X_0}^{i_3}{X_2}^{d-i_3},\,{X_0}^{i_2}{X_2}^{d-i_2},
        \,{X_0}^{i_1}{X_2}^{d-i_1}\big\}.
\end{align*}
Therefore we have ${k=9}$ and
\begin{align*}
(d-d_J)n+d_J-dk&\ge14(d-d_J)+d_J-9d=5d-13d_J
        \ge5d-13i_1\\
    &\ge7m+5t-13\min(1,t)\ge7m-8>0.
\end{align*}

If ${i_1<d_J\le i_2}$, the multiples of~${X_0}^{d_J}$ in~$\I_{2,d,14}$ are the monomials in the set
\begin{align*}
J:=\big\{&{X_0}^d,\,{X_0}^{i_2}{X_1}^{i_1}{X_2}^{d-i_1-i_2},
    \,{X_0}^{i_3}{X_1}^{d-i_3},\,{X_0}^{i_2}{X_1}^{d-i_3},\\
    &{X_0}^{i_3}{X_2}^{d-i_3},\,{X_0}^{i_2}{X_2}^{d-i_2}\big\}.
\end{align*}
Therefore we have ${k=6}$ and
\begin{align*}
(d-d_J)n+d_J-dk&=14(d-d_J)+d_J-6d=8d-13d_J\ge8d-13i_2\\
    &\ge6m+8t-13\min(2,t)\ge6m-10>0.
\end{align*}

If ${i_2<d_J\le i_3}$, the multiples of~${X_0}^{d_J}$ in~$\I_{2,d,14}$ are the monomials in the set
\begin{align*}
J:=\big\{&{X_0}^d,\,{X_0}^{i_3}{X_1}^{d-i_3},
    \,{X_0}^{i_3}{X_2}^{d-i_3}\big\}.
\end{align*}
Therefore we have ${k=3}$ and
\begin{align*}
(d-d_J)n+d_J-dk&=14(d-d_J)+d_J-3d=11d-13d_J\ge11d-13i_3\\
    &=5m-2t\ge5m-6>0.
\end{align*}
Therefore inequality~(\ref{6.6}) is strictly satisfied.

If ${i_3<d_J<d}$, the only multiple of~${X_0}^{d_J}$ in~$\I_{2,d,14}$ is ${X_0}^d$ and we have nothing to check.

\bigskip\noindent{\sc Case ${n=15}$.}\\
In case ${n=15}$ (${d\ge13}$), let us write ${d=4m+t}$, where ${0\le t<4}$, and for each ${l\in\{1,2,3\}}$, let ${i_l:=lm+\min(l,t)}$. Since ${d\ge13}$, we get ${m\ge3}$.

Consider the set
\begin{align*}
\I_{2,d,15}:=\big\{&{X_0}^d,\,{X_1}^d,\,{X_2}^d,\\
    &{X_0}^{i_2}{X_1}^{d-i_3}{X_2}^{i_3-i_2},
        \,{X_0}^{i_1}{X_1}^{d-i_2}{X_2}^{i_2-i_1},
        \,{X_0}^{i_1}{X_1}^{d-i_3}{X_2}^{i_3-i_1},\\
    &{X_0}^{i_3}{X_1}^{d-i_3},\,{X_0}^{i_2}{X_1}^{d-i_2},
        \,{X_0}^{i_1}{X_1}^{d-i_1},\\
    &{X_0}^{i_3}{X_2}^{d-i_3},\,{X_0}^{i_2}{X_2}^{d-i_2},
        \,{X_0}^{i_1}{X_2}^{d-i_1},\\
    &{X_1}^{i_1}{X_2}^{d-i_1},
        \,{X_1}^{i_2}{X_2}^{d-i_2},
        \,{X_1}^{i_3}{X_2}^{d-i_3}\big\}.
\end{align*}

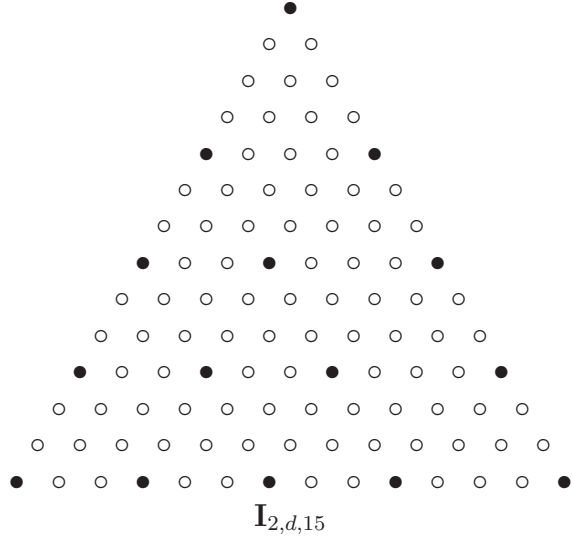
\[
\xymatrix@C=\meiabase@R=\altura@!0{
&&&&&&&&&&&&&{\bullet}\\
&&&&&&&&&&&&{\circ}&&{\circ}\\
&&&&&&&&&&&{\circ}&&{\circ}&&{\circ}\\
&&&&&&&&&&{\circ}&&{\circ}&&{\circ}&&{\circ}\\
&&&&&&&&&{\bullet}&&{\circ}&&{\circ}&&{\circ}&&{\bullet}\\
&&&&&&&&{\circ}&&{\circ}&&{\circ}&&{\circ}&&{\circ}&&{\circ}\\
&&&&&&&{\circ}&&{\circ}&&{\circ}&&{\circ}&&{\circ}&&{\circ}&&
    {\circ}\\
&&&&&&{\bullet}&&{\circ}&&{\circ}&&{\bullet}&&{\circ}&&{\circ}&&
    {\circ}&&{\bullet}\\
&&&&&{\circ}&&{\circ}&&{\circ}&&{\circ}&&{\circ}&&{\circ}&&
    {\circ}&&{\circ}&&{\circ}\\
&&&&{\circ}&&{\circ}&&{\circ}&&{\circ}&&{\circ}&&{\circ}&&
    {\circ}&&{\circ}&&{\circ}&&{\circ}\\
&&&{\bullet}&&{\circ}&&{\circ}&&{\bullet}&&{\circ}&&{\circ}&&
    {\bullet}&&{\circ}&&{\circ}&&{\circ}&&{\bullet}\\
&&{\circ}&&{\circ}&&{\circ}&&{\circ}&&{\circ}&&{\circ}&&
    {\circ}&&{\circ}&&{\circ}&&{\circ}&&{\circ}&&{\circ}\\
&{\circ}&&{\circ}&&{\circ}&&{\circ}&&{\circ}&&{\circ}&&
    {\circ}&&{\circ}&&{\circ}&&{\circ}&&{\circ}&&{\circ}&&
    {\circ}\\
{\bullet}&&{\circ}&&{\circ}&&{\bullet}&&{\circ}&&{\circ}&&
    {\bullet}&&{\circ}&&{\circ}&&{\bullet}&&{\circ}&&
    {\circ}&&{\circ}&&{\bullet}\\
&&&&&&&&&&&&&\I_{2,d,15}
}
\]

If ${0<d_J\le i_1}$, the multiples of~${X_0}^{d_J}$ in~$\I_{2,d,15}$ are the monomials in the set
\begin{align*}
J:=\big\{&{X_0}^d,\,{X_0}^{i_2}{X_1}^{i_1}{X_2}^{d-i_1-i_2},
        \,{X_0}^{i_1}{X_1}^{i_2}{X_2}^{d-i_1-i_2},
        \,{X_0}^{i_1}{X_1}^{i_1}{X_2}^{d-2i_1},\\
    &{X_0}^{i_3}{X_1}^{d-i_3},\,{X_0}^{i_2}{X_1}^{d-i_3},
        \,{X_0}^{i_1}{X_1}^{d-i_1},\\
    &{X_0}^{i_3}{X_2}^{d-i_3},\,{X_0}^{i_2}{X_2}^{d-i_2},
        \,{X_0}^{i_1}{X_2}^{d-i_1}\big\}.
\end{align*}
Therefore we have ${k=10}$ and
\begin{align*}
(d-d_J)n+d_J-dk&\ge15(d-d_J)+d_J-10d=5d-14d_J
        \ge5d-14i_1\\
    &\ge6m+5t-14\min(1,t)\ge6m-9>0.
\end{align*}

If ${i_1<d_J\le i_2}$, the multiples of~${X_0}^{d_J}$ in~$\I_{2,d,15}$ are the monomials in the set
\begin{align*}
J:=\big\{&{X_0}^d,\,{X_0}^{i_2}{X_1}^{i_1}{X_2}^{d-i_1-i_2},
        \,{X_0}^{i_3}{X_1}^{d-i_3},\\
    &{X_0}^{i_2}{X_1}^{d-i_3},
        \,{X_0}^{i_3}{X_2}^{d-i_3},
        \,{X_0}^{i_2}{X_2}^{d-i_2}\big\}.
\end{align*}
Therefore we have ${k=6}$ and
\begin{align*}
(d-d_J)n+d_J-dk&=15(d-d_J)+d_J-6d=9d-14d_J\ge9d-14i_2\\
    &\ge8m+9t-14\min(2,t)\ge8m-10>0.
\end{align*}

If ${i_2<d_J\le i_3}$, the multiples of~${X_0}^{d_J}$ in~$\I_{2,d,15}$ are the monomials in the set
\begin{align*}
J:=\big\{&{X_0}^d,\,{X_0}^{i_3}{X_1}^{d-i_3},
    \,{X_0}^{i_3}{X_2}^{d-i_3}\big\}.
\end{align*}
Therefore we have ${k=3}$ and
\begin{align*}
(d-d_J)n+d_J-dk&=15(d-d_J)+d_J-3d=12d-14d_J\ge12d-14i_3\\
    &=6m-2t\ge6m-6>0.
\end{align*}
Therefore inequality~(\ref{6.6}) is strictly satisfied.

If ${i_3<d_J<d}$, the only multiple of~${X_0}^{d_J}$ in~$\I_{2,d,15}$ is ${X_0}^d$ and we have nothing to check.

\bigskip\noindent{\sc Case ${n=16}$.}\\
In case ${n=16}$ (${d\ge14}$), let us write ${d=5m+t}$, where ${0\le t<5}$, and for each ${l\in\{1,2,3,4\}}$, let ${i_l:=lm+\min(l,t)}$. Since ${d\ge14}$, we get ${m\ge2}$.

Consider the set
\begin{align*}
\I_{2,d,16}:=\big\{&{X_0}^d,\,{X_1}^d,\,{X_2}^d,
        \,{X_0}^{i_2}{X_1}^{d-i_3}{X_2}^{i_3-i_2},
        \,{X_0}^{i_4}{X_1}^{d-i_4},
        \,{X_0}^{i_3}{X_1}^{d-i_3},\\
    &{X_0}^{i_2}{X_1}^{d-i_2},
        \,{X_0}^{i_1}{X_1}^{d-i_1},
        \,{X_0}^{i_4}{X_2}^{d-i_4},
        \,{X_0}^{i_3}{X_2}^{d-i_3},
        \,{X_0}^{i_2}{X_2}^{d-i_2},\\
    &{X_0}^{i_1}{X_2}^{d-i_1},
        \,{X_1}^{i_1}{X_2}^{d-i_1},
        \,{X_1}^{i_2}{X_2}^{d-i_2},
        \,{X_1}^{i_3}{X_2}^{d-i_3},
        \,{X_1}^{i_4}{X_2}^{d-i_4}\big\}.
\end{align*}

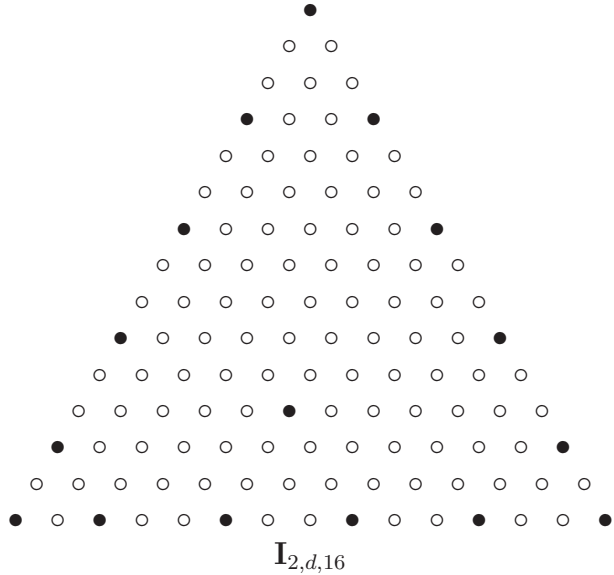
\[
\xymatrix@C=\meiabase@R=\altura@!0{
&&&&&&&&&&&&&&{\bullet}\\
&&&&&&&&&&&&&{\circ}&&{\circ}\\
&&&&&&&&&&&&{\circ}&&{\circ}&&{\circ}\\
&&&&&&&&&&&{\bullet}&&{\circ}&&{\circ}&&{\bullet}\\
&&&&&&&&&&{\circ}&&{\circ}&&{\circ}&&{\circ}&&{\circ}\\
&&&&&&&&&{\circ}&&{\circ}&&{\circ}&&{\circ}&&{\circ}&&
    {\circ}\\
&&&&&&&&{\bullet}&&{\circ}&&{\circ}&&{\circ}&&{\circ}&&
    {\circ}&&{\bullet}\\
&&&&&&&{\circ}&&{\circ}&&{\circ}&&{\circ}&&{\circ}&&
    {\circ}&&{\circ}&&{\circ}\\
&&&&&&{\circ}&&{\circ}&&{\circ}&&{\circ}&&{\circ}&&
    {\circ}&&{\circ}&&{\circ}&&{\circ}\\
&&&&&{\bullet}&&{\circ}&&{\circ}&&{\circ}&&{\circ}&&
    {\circ}&&{\circ}&&{\circ}&&{\circ}&&{\bullet}\\
&&&&{\circ}&&{\circ}&&{\circ}&&{\circ}&&{\circ}&&
    {\circ}&&{\circ}&&{\circ}&&{\circ}&&{\circ}&&{\circ}\\
&&&{\circ}&&{\circ}&&{\circ}&&{\circ}&&{\circ}&&
    {\bullet}&&{\circ}&&{\circ}&&{\circ}&&{\circ}&&{\circ}&&
    {\circ}\\
&&{\bullet}&&{\circ}&&{\circ}&&{\circ}&&{\circ}&&{\circ}&&
    {\circ}&&{\circ}&&{\circ}&&{\circ}&&{\circ}&&
    {\circ}&&{\bullet}\\
&{\circ}&&{\circ}&&{\circ}&&{\circ}&&{\circ}&&{\circ}&&{\circ}&&
    {\circ}&&{\circ}&&{\circ}&&{\circ}&&{\circ}&&{\circ}&&
    {\circ}\\
{\bullet}&&{\circ}&&{\bullet}&&{\circ}&&{\circ}&&
    {\bullet}&&{\circ}&&{\circ}&&{\bullet}&&{\circ}&&
    {\circ}&&{\bullet}&&{\circ}&&{\circ}&&{\bullet}\\
&&&&&&&&&&&&&&\I_{2,d,16}
}
\]

If ${0<d_J\le i_1}$, the multiples of~${X_0}^{d_J}$ in~$\I_{2,d,16}$ are the monomials in the set
\begin{align*}
J:=\big\{&{X_0}^d,\,{X_0}^{i_2}{X_1}^{i_2}{X_2}^{d-2i_2},
        \,{X_0}^{i_4}{X_1}^{d-i_4},\,{X_0}^{i_3}{X_1}^{d-i_3},
        \,{X_0}^{i_2}{X_1}^{d-i_2},\\
    &{X_0}^{i_1}{X_1}^{d-i_1},\,{X_0}^{i_4}{X_2}^{d-i_4},
        \,{X_0}^{i_3}{X_2}^{d-i_3},\,{X_0}^{i_2}{X_2}^{d-i_2},
        \,{X_0}^{i_1}{X_2}^{d-i_1}\big\}.
\end{align*}
Therefore we have ${k=10}$ and
\begin{align*}
(d-d_J)n+d_J-dk&\ge16(d-d_J)+d_J-10d=6d-15d_J\ge6d-15i_1\\
    &\ge15m+6t-15\min(1,t)\ge15m-9>0.
\end{align*}

If ${i_1<d_J\le i_2}$, the multiples of~${X_0}^{d_J}$ in~$\I_{2,d,16}$ are the monomials in the set
\begin{align*}
J:=\big\{&{X_0}^d,\,{X_0}^{i_2}{X_1}^{i_2}{X_2}^{d-2i_2},
        \,{X_0}^{i_4}{X_1}^{d-i_4},\,{X_0}^{i_3}{X_1}^{d-i_3},
        \,{X_0}^{i_2}{X_1}^{d-i_2},\\
    &{X_0}^{i_4}{X_2}^{d-i_4},\,{X_0}^{i_3}{X_2}^{d-i_3},
        \,{X_0}^{i_2}{X_2}^{d-i_2}\big\}.
\end{align*}
Therefore we have ${k=8}$ and
\begin{align*}
(d-d_J)n+d_J-dk&=16(d-d_J)+d_J-8d=8d-15d_J\ge8d-15i_2\\
    &\ge10m+8t-15\min(2,t)\ge10m-14>0.
\end{align*}

If ${i_2<d_J\le i_3}$, the multiples of~${X_0}^{d_J}$ in~$\I_{2,d,16}$ are the monomials in the set
\begin{align*}
J:=\big\{&{X_0}^d,\,{X_0}^{i_4}{X_1}^{d-i_4},
        \,{X_0}^{i_3}{X_1}^{d-i_3},\,{X_0}^{i_4}{X_2}^{d-i_4},
        \,{X_0}^{i_3}{X_2}^{d-i_3}\big\}.
\end{align*}
Therefore we have ${k=5}$ and
\begin{align*}
(d-d_J)n+d_J-dk&=16(d-d_J)+d_J-5d=11d-15d_J\ge11d-15i_3\\
    &\ge10m+11t-15\min(3,t)\ge10m-12>0.
\end{align*}

If ${i_3<d_J\le i_4}$, the multiples of~${X_0}^{d_J}$ in~$\I_{2,d,16}$ are the monomials in the set
\begin{align*}
J:=\big\{&{X_0}^d,\,{X_0}^{i_4}{X_1}^{d-i_4},
        \,{X_0}^{i_4}{X_2}^{d-i_4}\big\}.
\end{align*}
Therefore we have ${k=3}$ and
\begin{align*}
(d-d_J)n+d_J-dk&=16(d-d_J)+d_J-3d=13d-15d_J\ge13d-15i_4\\
    &\ge5m-2t\ge5m-8>0.
\end{align*}
Therefore inequality~(\ref{6.6}) is strictly satisfied.

If ${i_4<d_J<d}$, the only multiple of~${X_0}^{d_J}$ in~$\I_{2,d,16}$ is ${X_0}^d$ and we have nothing to check.

\bigskip\noindent{\sc Case ${n=17}$.}\\
In case ${n=17}$ (${d\ge15}$), let us write ${d=5m+t}$, where ${0\le t<5}$, and for each ${l\in\{1,2,3,4\}}$, let ${i_l:=lm+\min(l,t)}$. Since ${d\ge15}$, we get ${m\ge3}$.

Consider the set
\begin{align*}
\I_{2,d,17}:=\big\{&{X_0}^d,\,{X_1}^d,\,{X_2}^d,
        \,{X_0}^{i_2}{X_1}^{d-i_3}{X_2}^{i_3-i_2},
        \,{X_0}^{i_2}{X_1}^{d-i_4}{X_2}^{i_4-i_2},\\
    &{X_0}^{i_4}{X_1}^{d-i_4},\,{X_0}^{i_3}{X_1}^{d-i_3},
        \,{X_0}^{i_2}{X_1}^{d-i_2},
        \,{X_0}^{i_1}{X_1}^{d-i_1},\\
    &{X_0}^{i_4}{X_2}^{d-i_4},\,{X_0}^{i_3}{X_2}^{d-i_3},
        \,{X_0}^{i_2}{X_2}^{d-i_2},
        \,{X_0}^{i_1}{X_2}^{d-i_1},\\
    &{X_1}^{i_1}{X_2}^{d-i_1},\,{X_1}^{i_2}{X_2}^{d-i_2},
        \,{X_1}^{i_3}{X_2}^{d-i_3},
        \,{X_1}^{i_4}{X_2}^{d-i_4}\big\}.
\end{align*}

If ${0<d_J\le i_1}$, the multiples of~${X_0}^{d_J}$ in~$\I_{2,d,17}$ are the monomials in the set
\begin{align*}
J:=\big\{&{X_0}^d,\,{X_0}^{i_2}{X_1}^{i_2}{X_2}^{d-2i_2},
        \,{X_0}^{i_2}{X_1}^{i_1}{X_2}^{d-i_1-i_2},\\
    &{X_0}^{i_4}{X_1}^{d-i_4},\,{X_0}^{i_3}{X_1}^{d-i_3},
        \,{X_0}^{i_2}{X_1}^{d-i_2},
        \,{X_0}^{i_1}{X_1}^{d-i_1},\\
    &{X_0}^{i_4}{X_2}^{d-i_4},\,{X_0}^{i_3}{X_2}^{d-i_3},
        \,{X_0}^{i_2}{X_2}^{d-i_2},
        \,{X_0}^{i_1}{X_2}^{d-i_1}\big\}.
\end{align*}
Therefore we have ${k=11}$ and
\begin{align*}
(d-d_J)n+d_J-dk&\ge17(d-d_J)+d_J-11d=6d-16d_J\ge6d-16i_1\\
    &\ge14m+6t-16\min(1,t)\ge14m-10>0.
\end{align*}

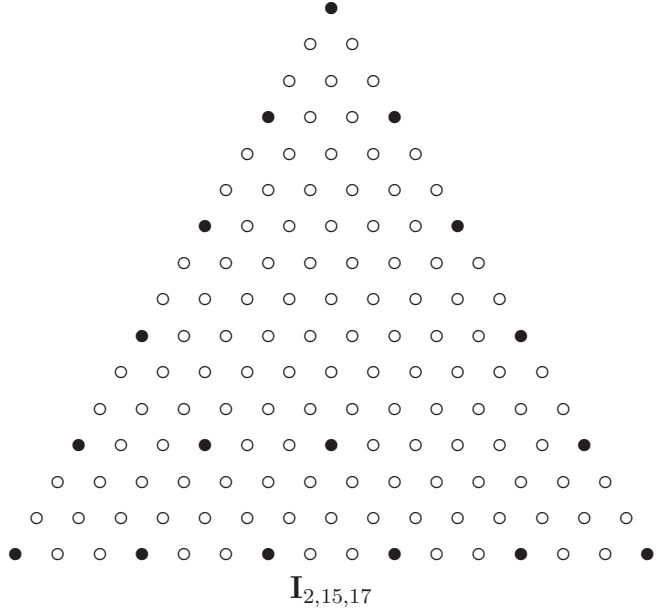
\[
\xymatrix@C=\meiabase@R=\altura@!0{
&&&&&&&&&&&&&&&{\bullet}\\
&&&&&&&&&&&&&&{\circ}&&{\circ}\\
&&&&&&&&&&&&&{\circ}&&{\circ}&&{\circ}\\
&&&&&&&&&&&&{\bullet}&&{\circ}&&{\circ}&&{\bullet}\\
&&&&&&&&&&&{\circ}&&{\circ}&&{\circ}&&{\circ}&&{\circ}\\
&&&&&&&&&&{\circ}&&{\circ}&&{\circ}&&{\circ}&&{\circ}&&
    {\circ}\\
&&&&&&&&&{\bullet}&&{\circ}&&{\circ}&&{\circ}&&{\circ}&&
    {\circ}&&{\bullet}\\
&&&&&&&&{\circ}&&{\circ}&&{\circ}&&{\circ}&&{\circ}&&
    {\circ}&&{\circ}&&{\circ}\\
&&&&&&&{\circ}&&{\circ}&&{\circ}&&{\circ}&&{\circ}&&
    {\circ}&&{\circ}&&{\circ}&&{\circ}\\
&&&&&&{\bullet}&&{\circ}&&{\circ}&&{\circ}&&{\circ}&&
    {\circ}&&{\circ}&&{\circ}&&{\circ}&&{\bullet}\\
&&&&&{\circ}&&{\circ}&&{\circ}&&{\circ}&&{\circ}&&
    {\circ}&&{\circ}&&{\circ}&&{\circ}&&{\circ}&&{\circ}\\
&&&&{\circ}&&{\circ}&&{\circ}&&{\circ}&&{\circ}&&
    {\circ}&&{\circ}&&{\circ}&&{\circ}&&{\circ}&&{\circ}&&
    {\circ}\\
&&&{\bullet}&&{\circ}&&{\circ}&&{\bullet}&&{\circ}&&{\circ}&&
    {\bullet}&&{\circ}&&{\circ}&&{\circ}&&{\circ}&&
    {\circ}&&{\bullet}\\
&&{\circ}&&{\circ}&&{\circ}&&{\circ}&&{\circ}&&{\circ}&&{\circ}&&
    {\circ}&&{\circ}&&{\circ}&&{\circ}&&{\circ}&&{\circ}&&
    {\circ}\\
&{\circ}&&{\circ}&&{\circ}&&{\circ}&&{\circ}&&
    {\circ}&&{\circ}&&{\circ}&&{\circ}&&{\circ}&&
    {\circ}&&{\circ}&&{\circ}&&{\circ}&&{\circ}\\
{\bullet}&&{\circ}&&{\circ}&&{\bullet}&&{\circ}&&{\circ}&&
    {\bullet}&&{\circ}&&{\circ}&&{\bullet}&&{\circ}&&
    {\circ}&&{\bullet}&&{\circ}&&{\circ}&&{\bullet}\\
&&&&&&&&&&&&&&&\I_{2,15,17}
}
\]

If ${i_1<d_J\le i_2}$, the multiples of~${X_0}^{d_J}$ in~$\I_{2,d,17}$ are the monomials in the set
\begin{align*}
J:=\big\{&{X_0}^d,\,{X_0}^{i_2}{X_1}^{i_2}{X_2}^{d-2i_2},
        \,{X_0}^{i_2}{X_1}^{i_1}{X_2}^{d-i_1-i_2},\\
    &{X_0}^{i_4}{X_1}^{d-i_4},\,{X_0}^{i_3}{X_1}^{d-i_3},
        \,{X_0}^{i_2}{X_1}^{d-i_2},\\
    &{X_0}^{i_4}{X_2}^{d-i_4},\,{X_0}^{i_3}{X_2}^{d-i_3},
        \,{X_0}^{i_2}{X_2}^{d-i_2}\big\}.
\end{align*}
Therefore we have ${k=9}$ and
\begin{align*}
(d-d_J)n+d_J-dk&=17(d-d_J)+d_J-9d=8d-16d_J\ge8d-16i_2\\
    &\ge8m+8t-16\min(2,t)\ge8m-16>0.
\end{align*}

If ${i_2<d_J\le i_3}$, the multiples of~${X_0}^{d_J}$ in~$\I_{2,d,17}$ are the monomials in the set
\begin{align*}
J:=\big\{&{X_0}^d,\,{X_0}^{i_4}{X_1}^{d-i_4},
        \,{X_0}^{i_3}{X_1}^{d-i_3},\,{X_0}^{i_4}{X_2}^{d-i_4},
        \,{X_0}^{i_3}{X_2}^{d-i_3}\big\}.
\end{align*}
Therefore we have ${k=5}$ and
\begin{align*}
(d-d_J)n+d_J-dk&=17(d-d_J)+d_J-5d=12d-16d_J\ge12d-16i_3\\
    &\ge12m+12t-16\min(3,t)\ge12m-12>0.
\end{align*}

If ${i_3<d_J\le i_4}$, the multiples of~${X_0}^{d_J}$ in~$\I_{2,d,17}$ are the monomials in the set
\begin{align*}
J:=\big\{&{X_0}^d,\,{X_0}^{i_4}{X_1}^{d-i_4},
        \,{X_0}^{i_4}{X_2}^{d-i_4}\big\}.
\end{align*}
Therefore we have ${k=3}$ and
\begin{align*}
(d-d_J)n+d_J-dk&=17(d-d_J)+d_J-3d=14d-16d_J\ge14d-16i_4\\
    &\ge6m-2t\ge6m-8>0.
\end{align*}
Therefore inequality~(\ref{6.6}) is strictly satisfied.

If ${i_4<d_J<d}$, the only multiple of~${X_0}^{d_J}$ in~$\I_{2,d,17}$ is ${X_0}^d$ and we have nothing to check.

\bigskip\noindent{\sc Case ${n=18}$.}\\
In case ${n=18}$ (${d\ge16}$), let us write ${d=5m+t}$, where ${0\le t<5}$, and for each ${l\in\{1,2,3,4\}}$, let ${i_l:=lm+\min(l,t)}$. Since ${d\ge16}$, we get ${m\ge3}$.

Consider the set
\begin{align*}
\I_{2,d,18}:=\big\{&{X_0}^d,\,{X_1}^d,\,{X_2}^d,\\
    &{X_0}^{i_2}{X_1}^{d-i_3}{X_2}^{i_3-i_2},
        \,{X_0}^{i_2}{X_1}^{d-i_4}{X_2}^{i_4-i_2},
        \,{X_0}^{i_1}{X_1}^{d-i_3}{X_2}^{i_3-i_1},\\
    &{X_0}^{i_4}{X_1}^{d-i_4},\,{X_0}^{i_3}{X_1}^{d-i_3},
        \,{X_0}^{i_2}{X_1}^{d-i_2},
        \,{X_0}^{i_1}{X_1}^{d-i_1},\\
    &{X_0}^{i_4}{X_2}^{d-i_4},\,{X_0}^{i_3}{X_2}^{d-i_3},
        \,{X_0}^{i_2}{X_2}^{d-i_2},
        \,{X_0}^{i_1}{X_2}^{d-i_1},\\
    &{X_1}^{i_1}{X_2}^{d-i_1},\,{X_1}^{i_2}{X_2}^{d-i_2},
        \,{X_1}^{i_3}{X_2}^{d-i_3},
        \,{X_1}^{i_4}{X_2}^{d-i_4}\big\}.
\end{align*}

If ${0<d_J\le i_1}$, the multiples of~${X_0}^{d_J}$ in~$\I_{2,d,18}$ are the monomials in the set
\begin{align*}
J:=\big\{&{X_0}^d,{X_0}^{i_2}{X_1}^{d-i_3}{X_2}^{i_3-i_2},
        \,{X_0}^{i_2}{X_1}^{d-i_4}{X_2}^{i_4-i_2},
        \,{X_0}^{i_1}{X_1}^{d-i_3}{X_2}^{i_3-i_1},\\
    &{X_0}^{i_4}{X_1}^{d-i_4},\,{X_0}^{i_3}{X_1}^{d-i_3},
        \,{X_0}^{i_2}{X_1}^{d-i_2},
        \,{X_0}^{i_1}{X_1}^{d-i_1},\\
    &{X_0}^{i_4}{X_2}^{d-i_4},\,{X_0}^{i_3}{X_2}^{d-i_3},
        \,{X_0}^{i_2}{X_2}^{d-i_2},
        \,{X_0}^{i_1}{X_2}^{d-i_1}\big\}.
\end{align*}
Therefore we have ${k=12}$ and
\begin{align*}
(d-d_J)n+d_J-dk&\ge18(d-d_J)+d_J-12d=6d-17d_J\ge6d-17i_1\\
    &\ge13m+6t-17\min(1,t)\ge13m-11>0.
\end{align*}

If ${i_1<d_J\le i_2}$, the multiples of~${X_0}^{d_J}$ in~$\I_{2,d,18}$ are the monomials in the set
\begin{align*}
J:=\big\{&{X_0}^d,\,{X_0}^{i_2}{X_1}^{d-i_3}{X_2}^{i_3-i_2},
        \,{X_0}^{i_2}{X_1}^{d-i_4}{X_2}^{i_4-i_2},\\
    &{X_0}^{i_4}{X_1}^{d-i_4},\,{X_0}^{i_3}{X_1}^{d-i_3},
        \,{X_0}^{i_2}{X_1}^{d-i_2},\\
    &{X_0}^{i_4}{X_2}^{d-i_4},\,{X_0}^{i_3}{X_2}^{d-i_3},
        \,{X_0}^{i_2}{X_2}^{d-i_2}\big\}.
\end{align*}

\[
\xymatrix@C=\meiabase@R=\altura@!0{
&&&&&&&&&&&&&&&&{\bullet}\\
&&&&&&&&&&&&&&&{\circ}&&{\circ}\\
&&&&&&&&&&&&&&{\circ}&&{\circ}&&{\circ}\\
&&&&&&&&&&&&&{\circ}&&{\circ}&&{\circ}&&{\circ}\\
&&&&&&&&&&&&{\bullet}&&{\circ}&&{\circ}&&{\circ}&&{\bullet}\\
&&&&&&&&&&&{\circ}&&{\circ}&&{\circ}&&{\circ}&&{\circ}&&
    {\circ}\\
&&&&&&&&&&{\circ}&&{\circ}&&{\circ}&&{\circ}&&{\circ}&&
    {\circ}&&{\circ}\\
&&&&&&&&&{\bullet}&&{\circ}&&{\circ}&&{\circ}&&{\circ}&&
    {\circ}&&{\circ}&&{\bullet}\\
&&&&&&&&{\circ}&&{\circ}&&{\circ}&&{\circ}&&{\circ}&&
    {\circ}&&{\circ}&&{\circ}&&{\circ}\\
&&&&&&&{\circ}&&{\circ}&&{\circ}&&{\circ}&&{\circ}&&
    {\circ}&&{\circ}&&{\circ}&&{\circ}&&{\circ}\\
&&&&&&{\bullet}&&{\circ}&&{\circ}&&{\bullet}&&{\circ}&&
    {\circ}&&{\bullet}&&{\circ}&&{\circ}&&{\circ}&&{\bullet}\\
&&&&&{\circ}&&{\circ}&&{\circ}&&{\circ}&&{\circ}&&
    {\circ}&&{\circ}&&{\circ}&&{\circ}&&{\circ}&&{\circ}&&
    {\circ}\\
&&&&{\circ}&&{\circ}&&{\circ}&&{\circ}&&{\circ}&&
    {\circ}&&{\circ}&&{\circ}&&{\circ}&&{\circ}&&{\circ}&&
    {\circ}&&{\circ}\\
&&&{\bullet}&&{\circ}&&{\circ}&&{\circ}&&{\circ}&&{\circ}&&
    {\bullet}&&{\circ}&&{\circ}&&{\circ}&&{\circ}&&
    {\circ}&&{\circ}&&{\bullet}\\
&&{\circ}&&{\circ}&&{\circ}&&{\circ}&&{\circ}&&{\circ}&&{\circ}&&
    {\circ}&&{\circ}&&{\circ}&&{\circ}&&{\circ}&&{\circ}&&
    {\circ}&&{\circ}\\
&{\circ}&&{\circ}&&{\circ}&&{\circ}&&{\circ}&&{\circ}&&
    {\circ}&&{\circ}&&{\circ}&&{\circ}&&{\circ}&&
    {\circ}&&{\circ}&&{\circ}&&{\circ}&&{\circ}\\
{\bullet}&&{\circ}&&{\circ}&&{\bullet}&&{\circ}&&{\circ}&&
    {\bullet}&&{\circ}&&{\circ}&&{\bullet}&&{\circ}&&
    {\circ}&&{\bullet}&&{\circ}&&{\circ}&&{\circ}&&{\bullet}\\
&&&&&&&&&&&&&&&&\I_{2,16,18}
}
\]
Therefore we have ${k=9}$ and
\begin{align*}
(d-d_J)n+d_J-dk&=18(d-d_J)+d_J-9d=9d-17d_J\ge9d-17i_2\\
    &\ge11m+9t-17\min(2,t)\ge11m-16>0.
\end{align*}

If ${i_2<d_J\le i_3}$, the multiples of~${X_0}^{d_J}$ in~$\I_{2,d,18}$ are the monomials in the set
\begin{align*}
J:=\big\{&{X_0}^d,\,{X_0}^{i_4}{X_1}^{d-i_4},
        \,{X_0}^{i_3}{X_1}^{d-i_3},\,{X_0}^{i_4}{X_2}^{d-i_4},
        \,{X_0}^{i_3}{X_2}^{d-i_3}\big\}.
\end{align*}
Therefore we have ${k=5}$ and
\begin{align*}
(d-d_J)n+d_J-dk&=18(d-d_J)+d_J-5d=13d-17d_J\ge13d-17i_3\\
    &\ge14m+13t-17\min(3,t)\ge14m-12>0.
\end{align*}

If ${i_3<d_J\le i_4}$, the multiples of~${X_0}^{d_J}$ in~$\I_{2,d,18}$ are the monomials in the set
\begin{align*}
J:=\big\{&{X_0}^d,\,{X_0}^{i_4}{X_1}^{d-i_4},
        \,{X_0}^{i_4}{X_2}^{d-i_4}\big\}.
\end{align*}
Therefore we have ${k=3}$ and
\begin{align*}
(d-d_J)n+d_J-dk&=18(d-d_J)+d_J-3d=15d-17d_J\ge15d-17i_4\\
    &\ge7m-2t\ge7m-8>0.
\end{align*}
Therefore inequality~(\ref{6.6}) is strictly satisfied.

If ${i_4<d_J<d}$, the only multiple of~${X_0}^{d_J}$ in~$\I_{2,d,18}$ is ${X_0}^d$ and we have nothing to check.

We conclude thus that in all cases the syzygy bundle corresponding to the family $\I_{2,d,n}$ is stable.
\end{dem}

\begin{p}\label{triangularnumbers}
For any integers $d$ and~$n$ such that $18<n\le d+2$, there is a set $\I_{2,d,n}$ of~$n$ monomials in $K\left[X_0,X_1,X_2\right]$ of degree~$d$ such that the corresponding syzygy bundle is stable.
\end{p}
\begin{dem}
For each natural~$j$, let $T_j:=\tfrac{j(j+1)}{2}=\tbinom{j+1}{2}$ be the $j$th triangular number. Choose $j$ such that
\[T_{j+2}\le n<T_{j+3}.\]
Since ${18<n}$, we have ${j\ge3}$. Let us write ${n=T_{j+2}+r}$, with ${0\le r\le j+2}$. Since ${n\le d+2}$, we get ${T_{j+2}\le d+2}$ and therefore ${2d-j^2-5j-2\ge0}$.

As before, following strategy~\ref{X_0}, for each~$d$ and~$n$, we choose a suitable set of monomials $\I_{2,d,n}$ such that for ${0<d_J<d}$, no monomial of degree~$d_J$ divides a greater number of monomials in~$\I_{2,d,n}$ than~${X_0}^{d_J}$. To this end, let us write ${d=m(j+1)+t}$, where ${0\le t<j+1}$. To choose the first~$T_{j+2}$ monomials in such a way that they are as evenly spread as possible, we can define for each ${l\in\{1,\ldots,j\}}$, the integers
\[i_l=lm+\min(l,t).\]
Now these numbers satisfy ${0<i_1<\cdots<i_j<d}$,
\[{d-i_j\le i_j-i_{j-1}\le\cdots\le i_2-i_1\le i_1}\mbox{,}\]
and ${i_1-(d-i_j)\le1}$.

Note that since ${2d\ge j^2+5j+2}$, we get ${d\ge3(j+1)+1}$ and therefore ${m\ge3}$. Let
\[e:=\left\lceil\tfrac{m}{2}\right\rceil.\]

Consider the set
\begin{align*}
\I_{2,d,n}':=\big\{&{X_0}^d,\,{X_0}^{i_j}{X_1}^{d-i_j},
        \,{X_0}^{i_j}{X_2}^{d-i_j},\\
    &{X_0}^{i_{j-1}}{X_1}^{d-i_{j-1}},\,
        {X_0}^{i_{j-1}}{X_1}^{d-i_j}{X_2}^{i_j-i_{j-1}},\,
        {X_0}^{i_{j-1}}{X_2}^{d-i_{j-1}},\\
    &{X_0}^{i_{j-2}}{X_1}^{d-i_{j-2}},\,
        {X_0}^{i_{j-2}}{X_1}^{d-i_{j-1}}
        {X_2}^{i_{j-1}-i_{j-2}},\\
    &\qquad {X_0}^{i_{j-2}}{X_1}^{d-i_j}{X_2}^{i_j-i_{j-2}},\,
        {X_0}^{i_{j-2}}{X_2}^{d-i_{j-2}},\\
    &\ldots\\
    &{X_0}^{i_1}{X_1}^{d-i_1},\,
        {X_0}^{i_1}{X_1}^{d-i_2}{X_2}^{i_2-i_1},\\
    &\qquad {X_0}^{i_1}{X_1}^{d-i_3}{X_2}^{i_3-i_1},
        \ldots,{X_0}^{i_1}{X_1}^{d-i_j}{X_2}^{i_j-i_1},\,
        {X_0}^{i_1}{X_2}^{d-i_1},\\
    &{X_1}^d,\,{X_1}^{i_j}{X_2}^{d-i_j},\ldots,
        {X_1}^{i_1}{X_2}^{d-i_1},\,{X_2}^d\big\},
\end{align*}
The~$T_{j+2}$ monomials in~$\I_{2,18,n}'$ are illustrated in the following picture. In case ${d=18}$, we have ${19\le n\le20}$, therefore, ${j=3}$ and ${T_{j+2}=15}$.

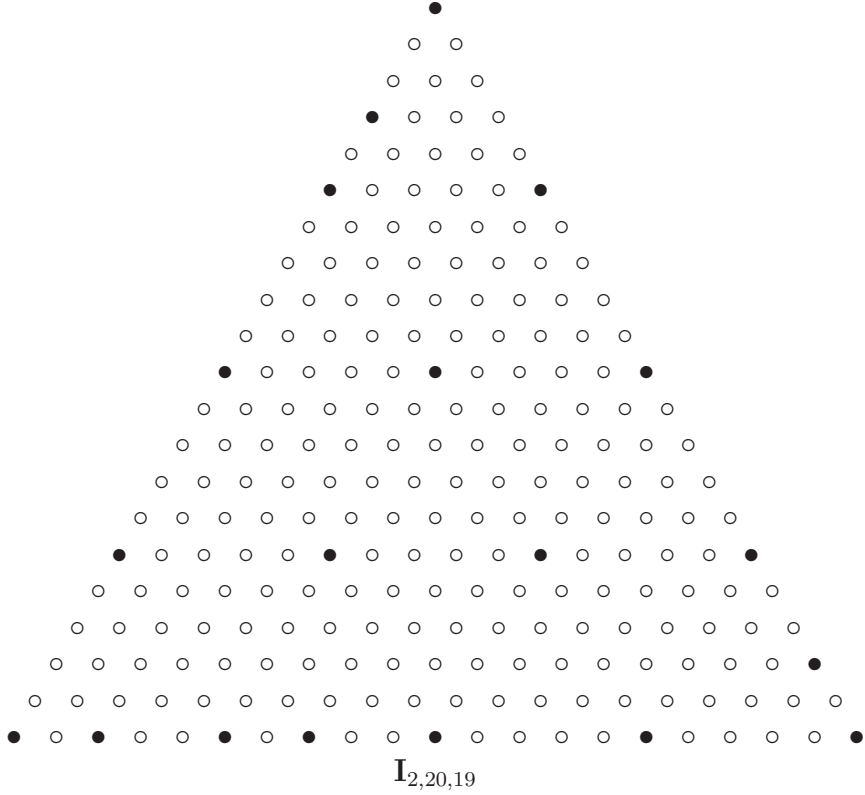
\[
\xymatrix@C=\meiabase@R=\altura@!0{
&&&&&&&&&&&&&&&&&&{\bullet}\\
&&&&&&&&&&&&&&&&&{\circ}&&{\circ}\\
&&&&&&&&&&&&&&&&{\circ}&&{\circ}&&{\circ}\\
&&&&&&&&&&&&&&&{\circ}&&{\circ}&&{\circ}&&{\circ}\\
&&&&&&&&&&&&&&{\circ}&&{\circ}&&{\circ}&&{\circ}&&{\circ}\\
&&&&&&&&&&&&&{\bullet}&&{\circ}&&{\circ}&&{\circ}&&{\circ}&&
    {\bullet}\\
&&&&&&&&&&&&{\circ}&&{\circ}&&{\circ}&&{\circ}&&{\circ}&&
    {\circ}&&{\circ}\\
&&&&&&&&&&&{\circ}&&{\circ}&&{\circ}&&{\circ}&&{\circ}&&
    {\circ}&&{\circ}&&{\circ}\\
&&&&&&&&&&{\circ}&&{\circ}&&{\circ}&&{\circ}&&{\circ}&&
    {\circ}&&{\circ}&&{\circ}&&{\circ}\\
&&&&&&&&&{\circ}&&{\circ}&&{\circ}&&{\circ}&&{\bullet}&&
    {\circ}&&{\circ}&&{\circ}&&{\circ}&&{\circ}\\
&&&&&&&&{\bullet}&&{\circ}&&{\circ}&&{\circ}&&{\circ}&&
    {\circ}&&{\circ}&&{\circ}&&{\circ}&&{\circ}&&
    {\bullet}\\
&&&&&&&{\circ}&&{\circ}&&{\circ}&&{\circ}&&{\circ}&&
    {\circ}&&{\circ}&&{\circ}&&{\circ}&&{\circ}&&{\circ}&&
    {\circ}\\
&&&&&&{\circ}&&{\circ}&&{\circ}&&{\circ}&&{\circ}&&{\circ}&&
    {\circ}&&{\circ}&&{\circ}&&{\circ}&&{\circ}&&{\circ}&&
    {\circ}\\
&&&&&{\circ}&&{\circ}&&{\circ}&&{\circ}&&{\circ}&&{\circ}&&
    {\circ}&&{\circ}&&{\bullet}&&{\circ}&&{\circ}&&{\circ}&&
    {\circ}&&{\circ}\\
&&&&{\bullet}&&{\circ}&&{\circ}&&{\circ}&&{\bullet}&&
    {\circ}&&{\circ}&&{\circ}&&{\circ}&&{\circ}&&{\circ}&&
    {\circ}&&{\circ}&&{\circ}&&{\bullet}\\
&&&{\circ}&&{\circ}&&{\circ}&&{\circ}&&{\circ}&&{\circ}&&
    {\circ}&&{\circ}&&{\circ}&&{\circ}&&{\circ}&&{\circ}&&
    {\circ}&&{\circ}&&{\circ}&&{\circ}\\
&&{\circ}&&{\circ}&&{\circ}&&{\circ}&&{\circ}&&{\circ}&&
    {\circ}&&{\circ}&&{\circ}&&{\circ}&&{\circ}&&{\circ}&&
    {\circ}&&{\circ}&&{\circ}&&{\circ}&&{\circ}\\
&{\circ}&&{\circ}&&{\circ}&&{\circ}&&{\circ}&&{\circ}&&
    {\circ}&&{\circ}&&{\circ}&&{\circ}&&{\circ}&&{\circ}&&
    {\circ}&&{\circ}&&{\circ}&&{\circ}&&{\circ}&&{\circ}\\
{\bullet}&&{\circ}&&{\circ}&&{\circ}&&{\bullet}&&{\circ}&&
    {\circ}&&{\circ}&&{\bullet}&&{\circ}&&{\circ}&&{\circ}&&
    {\circ}&&{\bullet}&&{\circ}&&{\circ}&&{\circ}&&{\circ}&&
    {\bullet}\\
&&&&&&&&&&&&&&&&&&\I_{2,18,n}'
}
\]

Now consider the sequence
\begin{align*}
\big(&{X_0}^{i_j+e}{X_1}^{d-i_j-e},\,
        {X_0}^e{X_2}^{d-e},\,
        {X_1}^{i_j+e}{X_2}^{d-i_j-e},\\
    &{X_0}^{i_{j-1}+e}{X_1}^{d-i_{j-1}-e},\,
        {X_0}^{i_1+e}{X_2}^{d-i_1-e},\,
        {X_1}^{i_{j-1}+e}{X_2}^{d-i_{j-1}-e},\\
    &\ldots,\\
    &{X_0}^{i_{j-q}+e}{X_1}^{d-i_{j-q}-e},\,
        {X_0}^{i_q+e}{X_2}^{d-i_q-e},\,
        {X_1}^{i_{j-q}+e}{X_2}^{d-i_{j-q}-e}\big)\mbox{,}
\end{align*}
where ${q:=\left\lceil\tfrac{j-1}{3}\right\rceil}$. Let~$\I_{2,d,n}''$ be the set of the first~$r$ monomials in this sequence and let ${\I_{2,d,n}:=\I_{2,d,n}'\cup \I_{2,d,n}''}$. Since $\I_{2,d,n}'$ has $T_{j+2}$ monomials, the number of monomials in~$\I_{2,d,n}$ is~$n$.

Here is a picture of~$\I_{2,20,19}$. In this case (${n=19}$ and ${d=20}$), we get ${j=3}$ and ${d=5(j+1)}$, therefore ${m=5}$, ${t=0}$ and ${e=3}$.

\[
\xymatrix@C=\meiabase@R=\altura@!0{
&&&&&&&&&&&&&&&&&&&&{\bullet}\\
&&&&&&&&&&&&&&&&&&&{\circ}&&{\circ}\\
&&&&&&&&&&&&&&&&&&{\circ}&&{\circ}&&{\circ}\\
&&&&&&&&&&&&&&&&&{\bullet}&&{\circ}&&{\circ}&&{\circ}\\
&&&&&&&&&&&&&&&&{\circ}&&{\circ}&&{\circ}&&{\circ}&&{\circ}\\
&&&&&&&&&&&&&&&{\bullet}&&{\circ}&&{\circ}&&{\circ}&&{\circ}&&
    {\bullet}\\
&&&&&&&&&&&&&&{\circ}&&{\circ}&&{\circ}&&{\circ}&&{\circ}&&
    {\circ}&&{\circ}\\
&&&&&&&&&&&&&{\circ}&&{\circ}&&{\circ}&&{\circ}&&{\circ}&&
    {\circ}&&{\circ}&&{\circ}\\
&&&&&&&&&&&&{\circ}&&{\circ}&&{\circ}&&{\circ}&&{\circ}&&
    {\circ}&&{\circ}&&{\circ}&&{\circ}\\
&&&&&&&&&&&{\circ}&&{\circ}&&{\circ}&&{\circ}&&{\circ}&&
    {\circ}&&{\circ}&&{\circ}&&{\circ}&&{\circ}\\
&&&&&&&&&&{\bullet}&&{\circ}&&{\circ}&&{\circ}&&{\circ}&&
    {\bullet}&&{\circ}&&{\circ}&&{\circ}&&{\circ}&&{\bullet}\\
&&&&&&&&&{\circ}&&{\circ}&&{\circ}&&{\circ}&&{\circ}&&
    {\circ}&&{\circ}&&{\circ}&&{\circ}&&{\circ}&&{\circ}&&
    {\circ}\\
&&&&&&&&{\circ}&&{\circ}&&{\circ}&&{\circ}&&{\circ}&&{\circ}&&
    {\circ}&&{\circ}&&{\circ}&&{\circ}&&{\circ}&&{\circ}&&
    {\circ}\\
&&&&&&&{\circ}&&{\circ}&&{\circ}&&{\circ}&&{\circ}&&{\circ}&&
    {\circ}&&{\circ}&&{\circ}&&{\circ}&&{\circ}&&{\circ}&&
    {\circ}&&{\circ}\\
&&&&&&{\circ}&&{\circ}&&{\circ}&&{\circ}&&{\circ}&&{\circ}&&
    {\circ}&&{\circ}&&{\circ}&&{\circ}&&{\circ}&&{\circ}&&
    {\circ}&&{\circ}&&{\circ}\\
&&&&&{\bullet}&&{\circ}&&{\circ}&&{\circ}&&{\circ}&&
    {\bullet}&&{\circ}&&{\circ}&&{\circ}&&{\circ}&&{\bullet}&&
    {\circ}&&{\circ}&&{\circ}&&{\circ}&&{\bullet}\\
&&&&{\circ}&&{\circ}&&{\circ}&&{\circ}&&{\circ}&&{\circ}&&
    {\circ}&&{\circ}&&{\circ}&&{\circ}&&{\circ}&&{\circ}&&
    {\circ}&&{\circ}&&{\circ}&&{\circ}&&{\circ}\\
&&&{\circ}&&{\circ}&&{\circ}&&{\circ}&&{\circ}&&{\circ}&&
    {\circ}&&{\circ}&&{\circ}&&{\circ}&&{\circ}&&{\circ}&&
    {\circ}&&{\circ}&&{\circ}&&{\circ}&&{\circ}&&{\circ}\\
&&{\circ}&&{\circ}&&{\circ}&&{\circ}&&{\circ}&&{\circ}&&
    {\circ}&&{\circ}&&{\circ}&&{\circ}&&{\circ}&&{\circ}&&
    {\circ}&&{\circ}&&{\circ}&&{\circ}&&{\circ}&&{\circ}&&
    {\bullet}\\
&{\circ}&&{\circ}&&{\circ}&&{\circ}&&{\circ}&&{\circ}&&
    {\circ}&&{\circ}&&{\circ}&&{\circ}&&{\circ}&&{\circ}&&
    {\circ}&&{\circ}&&{\circ}&&{\circ}&&{\circ}&&{\circ}&&
    {\circ}&&{\circ}\\
{\bullet}&&{\circ}&&{\bullet}&&{\circ}&&{\circ}&&{\bullet}&&
    {\circ}&&{\bullet}&&{\circ}&&{\circ}&&{\bullet}&&{\circ}&&
    {\circ}&&{\circ}&&{\circ}&&{\bullet}&&{\circ}&&{\circ}&&
    {\circ}&&{\circ}&&{\bullet}\\
&&&&&&&&&&&&&&&&&&&&\I_{2,20,19}
}
\]

Since we are adopting strategy~\ref{X_0}, for ${0<d_J<d}$, no monomial of degree~$d_J$ divides a greater number of monomials in~$\I_{2,d,n}$ than~${X_0}^{d_J}$.

For ${1\le l\le j}$, let $J_l$ be the set of monomials in~$\I_{2,d,n}'$ that are multiples of~${X_0}^{i_l}$. We have
\begin{align*}
J_l:=\big\{&{X_0}^d,\,{X_0}^{i_j}{X_1}^{d-i_j},
        \,{X_0}^{i_j}{X_2}^{d-i_j},\\
    &{X_0}^{i_{j-1}}{X_1}^{d-i_{j-1}},
        \,{X_0}^{i_{j-1}}{X_1}^{d-i_j}{X_2}^{i_j-i_{j-1}},
        \,{X_0}^{i_{j-1}}{X_2}^{d-i_{j-1}},\\
    &{X_0}^{i_{j-2}}{X_1}^{d-i_{j-2}},\,
        {X_0}^{i_{j-2}}{X_1}^{d-i_{j-1}}{X_2}^{i_{j-1}-i_{j-2}},
        \,{X_0}^{i_{j-2}}{X_1}^{d-i_j}{X_2}^{i_j-i_{j-2}},\\
    &\qquad {X_0}^{i_{j-2}}{X_2}^{d-i_{j-2}},\\
    &\ldots\\
    &{X_0}^{i_l}{X_1}^{d-i_l},\,
        {X_0}^{i_l}{X_1}^{d-i_{l+1}}{X_2}^{i_{l+1}-i_l},\\
    &\qquad {X_0}^{i_l}{X_1}^{i_{l+2}}{X_2}^{i_{l+2}-i_l},
        \ldots,{X_0}^{i_l}{X_1}^{d-i_j}{X_2}^{i_j-i_l},\,
        {X_0}^{i_l}{X_2}^{d-i_l}\big\}
\end{align*}
and $\left|J_l\right|=T_{j-l+2}$.

We distinguish two cases: ${r=0}$ and ${r\ne0}$.

\bigskip
\noindent{\sc Case 1.}\\
Suppose ${r=0}$ and therefore ${n=T_{j+2}}$.

If ${0<d_J\le i_1}$, the multiples of~${X_0}^{d_J}$ in~$\I_{2,d,n}$ are among the monomials in the set $J_1$. Therefore if~$k$ is the number of multiples of~${X_0}^{d_J}$, we have ${k=T_{j+1}}$ and
\begin{align*}
(d-d_J)n+d_J-dk&=(d-d_J)T_{j+2}+d_J-dT_{j+1}\\
    &=d(j+2)-d_J(T_{j+2}-1)\\
    &\ge \big(m(j+1)+t\big)(j+2)-i_1(T_{j+2}-1)\\
    &=\big(m(j+1)+t\big)(j+2)-(m+\min(1,t))(T_{j+2}-1).
\end{align*}
This last expression takes the following forms, depending on the different values of~$t$:
\begin{align*}
&\mbox{- for $t=0$, }&&\tfrac{m}{2}j(j+1)\mbox{;}\\
&\mbox{- for $t>0$, }&&\tfrac{(j+2)}{2}
        \big(2t+(m-1)(j-1)-4 \big)+m+1.
\end{align*}
These expressions are positive in both cases, since ${j\ge3}$ and ${m\ge3}$. So inequality~(\ref{6.6}) is strictly satisfied.

If ${i_l<d_J\le i_{l+1}}$, for ${1\le l\le j-1}$, the multiples of~${X_0}^{d_J}$ in~$\I_{2,d,n}$ are the monomials in the set
$J_{l+1}$. Therefore we have ${k=T_{j+1-l}}$ and
\begin{align*}
(d-d_J)n+d_J-dk&=(d-d_J)T_{j+2}+d_J-dT_{j+1-l}\\
    &=d(T_{j+2}-T_{j+1-l})-d_J(T_{j+2}-1)\\
    &\ge d(T_{j+2}-T_{j+1-l})-i_{l+1}(T_{j+2}-1).
\end{align*}
This last expression takes the following forms, depending on the different values of~$t$:
\begin{align*}
\mbox{- for }&t\leq l+1,\\
    &\tfrac{1}{2}(m-1)lj(j-l) +\tfrac{1}{2}(m-1)j(j-l)
        + \tfrac{1}{2}(m-2)l(j-l) \\
    &\qquad +  \tfrac{1}{2}(m-2)(j-l)
    + \tfrac{1}{2}(l-1)^2(j-l) + (l-1)(j-l)\\
    &\qquad + \tfrac{1}{2}\left((j-l)^2 +3(j-l)\right)
        (l+1-t)>0;\displaybreak[0]\\
\mbox{- for }&t>l+1,\\
    &\tfrac{1}{2}(m-1)lj(j-l) + \tfrac{1}{2}(m-1)j(j-l) +
        \tfrac{1}{2}(m-2)l(j-l) \\
    &\qquad +  \tfrac{1}{2}(m-3)(j-l)+ \tfrac{1}{2}l^2(j-l)\\
    &\qquad +\tfrac{1}{2}\left(2l(j-l) +l^2 +2j +3l +4 \right)
        (t-l-1)>0.
\end{align*}
Therefore inequality~(\ref{6.6}) is strictly satisfied.

If ${i_j<d_J<d}$, the only multiple of~${X_0}^{d_J}$ in~$\I_{2,d,n}$ is ${X_0}^d$ and there is nothing to check.

Therefore all possible values of $d_J$ are verified, and hence the syzygy bundle associated to~$\I_{2,d,n}$ is stable.

\bigskip
\noindent{\sc Case~2.}\\
Suppose ${r>0}$ and therefore ${n>T_{j+2}}$.

In this case, since ${n=T_{j+2}+r\le d+2}$, with ${0<r\le j+2}$, we must have ${d+2\ge T_{j+2}+1}$. From this, and in case ${j>3}$, we get
\[2d\ge j^2+5j+4\ge9j+4\ge8(j+1).\]
In case ${j=3}$, since ${d\ge17}$, we have ${d\ge4(j+1)+1}$. In any case, ${m\ge4}$.

We distinguish three subcases: ${r=3s+1}$, ${r=3s+2}$ and ${r=3s}$, with ${s\ge0}$ in the first two, and ${s\ge1}$ in the last one.

\bigskip
\noindent{\sc Subcase 2.1:} ${r=3s+1}$, with ${s\ge0}$.\\
If ${0<d_J\le e}$, the multiples of~${X_0}^{d_J}$ in~$\I_{2,d,n}$ are the monomials in the set
\begin{align*}
J:=J_1\cup\big\{
    &{X_0}^{i_j+e}{X_1}^{d-i_j-e},\ldots,
        {X_0}^{i_{j-s}+e}{X_1}^{d-i_{j-s}-e},\\
    &{X_0}^e{X_2}^{d-e},\ldots,
        {X_0}^{i_{s-1}+e}{X_2}^{d-i_{s-1}-e}\big\}.
\end{align*}
Therefore if~$k$ is the number of multiples of~${X_0}^{d_J}$, we have ${k=T_{j+1}+2s+1}$ and
\begin{align*}
(d-d_J)n+d_J-dk&=(d-d_J)(T_{j+2}+3s+1)+d_J-d(T_{j+1}+2s+1)\\
    &=d(j+2+s)-d_J(T_{j+2}+3s)\\
    &\ge (m(j+1)+t)(j+2+s)-e(T_{j+2}+3s)\\
    &\ge (m(j+1)+t)(j+2+s)-\tfrac{m+1}{2}(T_{j+2}+3s)\\
    &=\tfrac{1}{4}(3m-1)j^2+\tfrac{1}{4}(7m-5)(j-2)+4(m-1)\\
        &\qquad+\tfrac{1}{2}(2m(j-2)+3(m-1))s+t(j+2+s)>0.
\end{align*}

If ${e<d_J\le i_1}$, the multiples of~${X_0}^{d_J}$ in~$\I_{2,d,n}$ are the monomials in the set
\begin{align*}
J:=J_1\cup\big\{
    &{X_0}^{i_j+e}{X_1}^{d-i_j-e},\ldots,
        {X_0}^{i_{j-s}+e}{X_1}^{d-i_{j-s}-e},\\
    &{X_0}^{i_1+e}{X_2}^{d-i_1-e},\ldots,
        {X_0}^{i_{s-1}+e}{X_2}^{d-i_{s-1}-e}\big\}.
\end{align*}
Therefore we have ${k=T_{j+1}+\max(2s,1)}$ and
\begin{align*}
(d-d_J)n+d_J-dk&=(d-d_J)(T_{j+2}+3s+1)+d_J-d(T_{j+1}+
        \max(2s,1))\\
    &=d(j+3+3s-\max(2s,1))-d_J(T_{j+2}+3s)\\
    &\ge d(j+3+3s-\max(2s,1))-i_1(T_{j+2}+3s)\\
    &=(m(j+1)+t)(j+3+3s-\max(2s,1))\\
        &\qquad-(m+\min(1,t))(T_{j+2}+3s).
\end{align*}
This last expression takes the following forms, depending on the different values of~$s$ and~$t$:
\begin{align*}
&\mbox{- for $s=t=0$, }&&\tfrac{m}{2}(j+2)(j-1)\mbox{;}\\
&\mbox{- for $s=0$ and $t>0$, }&&\tfrac{1}{2}
        \big((m-1)(j-2)^2+(5m-7)(j-2)+4(m-3)\big)\\
    &&&\qquad +(t-1)(j+2)\mbox{;}\\
&\mbox{- for $s>0$ and $t=0$, }&&m\left(\tfrac{1}{2}
        j(j+3)+(j-2)s\right)\mbox{;}\\
&\mbox{- for $s>0$ and $t>0$, }&&\tfrac{1}{2}
        \big((m-1)j^2+(3m-5)(j-2)\big)+3(m-3)\\
    &&&\qquad +(j+3+s)(t-1)+m(j-2)s+j+4-2s.
\end{align*}
These expressions are positive in all cases, since ${m\ge4}$, ${j\ge3}$ and ${s\le\tfrac{j+1}{3}}$.

Therefore inequality~(\ref{6.6}) is strictly satisfied.

If ${i_l<d_J\le i_l+e}$, for ${1\le l\le j-2}$, the multiples of~${X_0}^{d_J}$ in~$\I_{2,d,n}$ are the monomials in the set
\begin{align*}
J:=J_{l+1}\cup\big\{
    &{X_0}^{i_j+e}{X_1}^{d-i_j-e},\ldots,
        {X_0}^{i_a+e}{X_1}^{d-i_a-e},\\
    &{X_0}^{i_l+e}{X_2}^{d-i_l-e},\ldots,
        {X_0}^{i_{s-1}+e}{X_2}^{d-i_{s-1}-e}\big\},
\end{align*}
where ${a=\max(j-s,l)}$ and the second line is understood to be empty if ${s\le l}$. Therefore we have ${k=T_{j+1-l}+\min(s+1,j+1-l)+\max(s-l,0)}$ and
\begin{align*}
(d-d_J)n+d_J-dk&=(d-d_J)(T_{j+2}+3s+1)+d_J\\
    &\qquad -d\big(T_{j+1-l}+1+\min(s,j-l)+\max(s-l,0)\big)\\
    &=d\big(T_{j+2}-T_{j+1-l}+3s
        -\min(s,j-l)-\max(s-l,0)\big)\\
    &\qquad -d_J(T_{j+2}+3s)\\
    &\ge d\big(T_{j+2}-T_{j+1-l}+3s-\min(s,j-l)
        -\max(s-l,0)\big)\\
    &\qquad -\left(i_l+\tfrac{m+1}{2}\right)(T_{j+2}+3s).
\end{align*}
We can rewrite this last expression in the following forms, depending on the different values of~$j$, $l$, $s$ and~$t$, so that they become sums of non\h{negative} numbers (and at least one of these is strictly positive):
\begin{align*}
\mbox{- for }&\mbox{$s\leq j-l$, $s\leq l$, and $t\leq l$,}\\
    &\tfrac{1}{2}(m-2)lj(j-l-1) + \tfrac{3}{4}(m-1)j(j-l) +
        \tfrac{1}{4}(m-4)l(j-l) \\
    &\qquad +\tfrac{3}{4}(m-4)l^2+\tfrac{7}{4}(m-3)(j-l)+
        \tfrac{1}{4}(m-4)l+\tfrac{1}{2}m\\
    &\qquad +\tfrac{1}{2}(l-1)(j-l-2)^2 +(l-1)^2(j-l-2)
        +j(j-l-2) \\
    &\qquad +\tfrac{3}{4}l(j-l-2)+
        \tfrac{15}{4}l(l-1)+3(j-l) +\tfrac{5}{2} \\
    &\qquad + \tfrac{1}{2}(4m(j-l)+m-3)s +l(m+1)(j-l-s)\\
    &\qquad + \tfrac{1}{2}(l-t)\big((j-l)^2+3(j-l)+2+2s\big)
        >0\mbox{;}\displaybreak[0]\\
\mbox{- for }&\mbox{$s\leq j-l$, $s\leq l$, and $t>l$,}\\
    &\tfrac{1}{2}(m-4)jl(j-l-1) +\tfrac{3}{4}(m-1)j(j-l)+
        \tfrac{1}{4}(m-4)lj +\tfrac{1}{2}ml^2\\
    &\qquad +\tfrac{7}{4}(m-1)(j-l)+\tfrac{1}{4}ml+
        \tfrac{1}{2}(m-4)+\tfrac{3}{2}lj(j-l-2)\\
    &\qquad +\tfrac{1}{2}(l-1)^2(j-l)+
        \tfrac{1}{2}j(j-l-2)+\tfrac{1}{4}l(j-l-2)+j\\
    &\qquad  +\tfrac{7}{4}l(l-1)+
        \tfrac{1}{2}+
    \tfrac{1}{2}(t-l)\big((2j-l)l+2j+3l+4+4s\big)\\
    &\qquad +\tfrac{1}{2}(4m(j-l)+m-3)s+
        (m+1)l(j-l-s)>0\mbox{;}\displaybreak[0]\\
\mbox{- for }&\mbox{$s\leq j-l$, $s>l$, and $t\leq l$,}\\
    &\tfrac{1}{2}(m-1)lj(j-l) +\tfrac{3}{4}(m-1)j(j-l)+
        \tfrac{7}{4}(m-3)l(j-l) +\tfrac{1}{4}ml^2\\
    &\qquad +\tfrac{1}{4}(m-2)(3j+4l+2) +
        \tfrac{1}{2}l^2(j-l)+\tfrac{1}{2}j^2 +l(j-l) \\
    &\qquad +\tfrac{1}{4}l(l-1)+\tfrac{1}{4}(j-l-2) +
        \tfrac{1}{2}(2m(j-2l)+m-3)(s-l)\\
    &\qquad + \big((m+2)l+m\big)(j-l-s) \\
    &\qquad + \tfrac{1}{2}\left((j-l)^2+3(j-2l)+l+2+4s\right)
        (l-t)>0\mbox{;}\displaybreak[0]\\
\mbox{- for }&\mbox{$s\leq j-l$, $s>l$, and $t>l$,}\\
    &\tfrac{1}{2}(m-1)lj(j-l) +\tfrac{3}{4}(m-1)j(j-l)+
        \tfrac{7}{4}(m-3)l(j-l) +\tfrac{1}{4}ml^2 \\
    &\qquad +\tfrac{1}{4}(m-1)(3j+4l+2) +
        \tfrac{1}{2}l^2(j-l)+\tfrac{1}{2}j(j-2l-1) \\
    &\qquad +\tfrac{1}{2}(2m(j-2l)+m-3)(s-l)+\tfrac{7}{4}(l-1)+
        \tfrac{3}{4}\\
    &\qquad +\big((m+2)l+m\big)(j-l-s)+2l(j-2l)+\tfrac{13}{4}l(l-1) \\
    &\qquad +\tfrac{1}{2}\left(2(l+1)(j-l)+l^2+7l+4+2s\right)
        (t-l)>0\mbox{;}\displaybreak[0]\\
\mbox{- for }&\mbox{$s>j-l$, $s\leq l$, and $t\leq l$,}\\
    &\tfrac{1}{2}(m-4)lj(j-l-2) +\tfrac{11}{4}(m-2)j(j-l)+
        \tfrac{5}{4}(m-1)l(2l-j) \\
    &\qquad +\tfrac{9}{4}m(j-l) +\tfrac{1}{4}ml +
        \tfrac{1}{2}m +\tfrac{3}{2}(l-1)(j-l-2)^2  + \tfrac{1}{4}l +
        \tfrac{15}{2}\\
    &\qquad +2(l-1)^2(j-l-2) +\tfrac{27}{4}j(j-l-2) +
        \tfrac{1}{2}l(j-l-2)\\
    &\qquad +(l-1)^2+\tfrac{11}{4}(j-l)+
    \tfrac{1}{2}(s-j+l)( 6m(j-l) +  3m - 3 )\\
    &\qquad + \tfrac{1}{2}(l-t)\left( (j-l)^2 +
        5(j-l) + 2 \right)>0\mbox{;}\displaybreak[0]\\
\mbox{- for }&\mbox{$s>j-l$, $s\leq l$, $t>l$,}\\
    &\tfrac{1}{2}(m-4)lj(j-l-2) + \tfrac{11}{4}(m-2)j(j-l) +
        \tfrac{5}{4}(m-1)l(2l-j) \\
    &\qquad +\tfrac{9}{4}m(j-l) +\tfrac{1}{4}ml +
        \tfrac{1}{2}m+\tfrac{3}{2}(l-1)(j-l-2)^2 \\
    &\qquad + 2(l-1)^2(j-l-2) + \tfrac{27}{4}j(j-l-2) +
        \tfrac{1}{2}l(j-l-2) + (l-1)^2\\
    &\qquad +\tfrac{11}{4}(j-l)+\tfrac{1}{4}l +
        \tfrac{15}{2}+\tfrac{1}{2}(s-j+l)(6m(j-l)+ 3m-3)\\
    &\qquad +\tfrac{1}{2}(t-l)
        \left(2l(j-l)+l^2+5l+4+6s\right)>0.
\end{align*}
Since ${s\le\tfrac{j+1}{3}}$, if ${l<s}$, we get ${j-l>\tfrac{2j-1}{3}\ge\tfrac{j+1}{3}\ge s}$. Therefore all possible cases are checked, and inequality~(\ref{6.6}) is strictly satisfied.

If ${i_l+e<d_J\le i_{l+1}}$, for ${1\le l\le j-2}$, the multiples of~${X_0}^{d_J}$ in~$\I_{2,d,n}$ are the monomials in the set
\begin{align*}
J:=J_{l+1}\cup\big\{
    &{X_0}^{i_j+e}{X_1}^{d-i_j-e},\ldots,
        {X_0}^{i_a+e}{X_1}^{d-i_a-e},\\
    &{X_0}^{i_{l+1}+e}{X_2}^{d-i_{l+1}-e},\ldots,
        {X_0}^{i_{s-1}+e}{X_2}^{d-i_{s-1}-e}\big\}\mbox{,}
\end{align*}
where ${a=\max(j-s,l+1)}$ and the second line is understood to be empty if ${s\le l+1}$. Therefore we have ${k=T_{j+1-l}+\min(s+1,j-l)+\max(s-l-1,0)}$ and
\begin{align*}
(d-d_J)n+d_J-dk&=(d-d_J)(T_{j+2}+3s+1)+d_J\\
    &\qquad -d\big(T_{j+1-l}+\min(s+1,j-l)+
        \max(s-l-1,0)\big)\\
    &=d\big(T_{j+2}-T_{j+1-l}+3s+1-\min(s+1,j-l)\\
    &\qquad -\max(s-l-1,0)\big)-d_J(T_{j+2}+3s)\\
    &\ge d\big(T_{j+2}-T_{j+1-l}+3s+1-\min(s+1,j-l)\\
    &\qquad -\max(s-l-1,0)\big)-i_{l+1}(T_{j+2}+3s).
\end{align*}
As before, we can rewrite this last expression in the following forms, depending on the different values of~$j$, $l$, $s$ and~$t$:
\begin{align*}
\mbox{- for }&\mbox{$s+1\le j-l$, $s\le l+1$, and $t\le l+1$,}\\
    &\tfrac{1}{2}(m-4)lj(j-l) +\tfrac{1}{2}(m-4)(j-l)(j-l-1)+
        (m-4)(j-l-1)\\
    &\qquad +\tfrac{3}{2}lj(j-l-2)+ \tfrac{1}{2}l^2(j-l)+
        \tfrac{3}{2}(j-l)(j-l-2)+\tfrac{1}{2}l(j-l)\\
    &\qquad + \tfrac{5}{2}(j-l-2)+ 3l^2+1+m(2(j-l-1)+1)s\\
    &\qquad +(ml+t)(j-l-1-s)+\tfrac{1}{2}
        \left((j-l)^2+5(j-l)\right)(l+1-t)>0\mbox{;}\displaybreak[0]\\
\mbox{- for }&\mbox{$s+1\le j-l$, $s\le l+1$, and $t>l+1$,}\\
    &\tfrac{1}{2}(m-4)lj(j-l) +\tfrac{1}{2}(m-4)(j-l)(j-l-1)+
        (m-4)(j-l-1)\\
    &\qquad  + \tfrac{3}{2}l(j-l-2)^2+ 2l(l-1)(j-l) +
        \tfrac{3}{2}(j-l)(j-l-2)\\
    &\qquad  + \tfrac{7}{2}l(j-l-2)+ \tfrac{1}{2}(j-l-2)
        + 3(l-1)+ 2\\
    &\qquad  +(2m(j-l-1)+m+2t)s+ (ml+3l+3)(j-l-1-s)\\
    &\qquad + \tfrac{1}{2}\left( 2l(j-l) + l^2 + 2(j-l)
        + 5l + 4 \right)(t-l-1)>0\mbox{;}\displaybreak[0]\\
\mbox{- for }&\mbox{$s+1\le j-l$, $s>l+1$, and $t\le l+1$,}\\
    &\tfrac{1}{2}(m-1)lj(j-l) +
        \tfrac{1}{2}(m-4)\left(j^2+l^2\right)+
        \tfrac{5}{2}(m-1)(j-l-1) + ml  \\
    &\qquad + \tfrac{1}{2}m+\tfrac{1}{2}(l-1)^2(j-l) +
        \tfrac{3}{2}j(j-l-1)+ 4l^2+ \tfrac{7}{2}(l-1) + 3\\
    &\qquad + m(j-l-2)(s-l-1) + 2(ml+t)(j-l-1-s)\\
    &\qquad + \tfrac{1}{2}\left( (j-l-2)^2
        + 11(j-l-2) + 14 \right)(l+1-t)
        + lt>0\mbox{;}\displaybreak[0]\\
\mbox{- for }&\mbox{$s+1\le j-l$, $s>l+1$, and $t>l+1$,}\\
    &\tfrac{1}{2}(m-4)lj(j-l) +
        \tfrac{1}{2}(m-1)\left(j^2+l^2\right)+
        \tfrac{5}{2}(m-2)(j-l-1) +ml\\
    &\qquad +\tfrac{1}{2}m +\tfrac{3}{2}lj(j-l-2) +
        \tfrac{1}{2}l(l-1)(j-l)+ 6l^2 +\tfrac{1}{2}j +
        \tfrac{11}{2}(l-1) +\tfrac{9}{2}\\
    &\qquad + ( m(j-l-2) + t )(s-l-1)+
        \big((2m+3)l+3\big)(j-l-1-s)\\
    &\qquad + \tfrac{1}{2}\left( 2l(j-l) + l^2 + 2j
        + 7l + 8 \right)(t-l-1)>0\mbox{;}\displaybreak[0]\\
\mbox{- for }&\mbox{$s+1>j-l$, $s\le l+1$, and $t\le l+1$,}\\
    &\tfrac{1}{2}(m-4)lj(j-l) + \tfrac{5}{2}(m-1)(j-l)^2+
        \tfrac{7}{2}(m-1)(j-l) + \tfrac{3}{2}lj(j-l-2)\\
    &\qquad + \tfrac{1}{2}l^2(j-l) + 2(j-l)^2
        +\tfrac{1}{2}lj +\tfrac{5}{2}l^2 + j-l \\
    &\qquad + 3m(j-l)(s-j+l-1)\\
    &\qquad + \tfrac{1}{2}\left( (j-l)^2 + 5(j-l) \right)
        (l+1-t)>0\mbox{;}\displaybreak[0]\\
\mbox{- for }&\mbox{$s+1>j-l$, $s\le l+1$, and $t>l+1$,}\\
    &\tfrac{1}{2}(m-4)lj(j-l) + \tfrac{5}{2}(m-1)(j-l)^2+
        \tfrac{7}{2}(m-1)(j-l)  \\
    &\qquad + \tfrac{3}{2}lj(j-l-2)+ \tfrac{1}{2}l^2(j-l) +
        2(j-l)^2 +\tfrac{1}{2}lj\\
    &\qquad + \tfrac{5}{2}l^2 + j-l + 3m(j-l)(s-j+l-1)\\
    &\qquad + \tfrac{1}{2}\left( 2l(j-l) + l^2
        + 5l + 6 + 6s \right)(t-l-1)>0.
\end{align*}
Since ${s\le\tfrac{j+1}{3}}$, if ${l+1<s}$, we get ${j-l-1>\tfrac{2j-1}{3}\ge\tfrac{j+1}{3}\ge s}$. Therefore all possible cases are checked, and inequality~(\ref{6.6}) is strictly satisfied.

If ${i_{j-1}<d_J\le i_{j-1}+e}$, the multiples of~${X_0}^{d_J}$ in~$\I_{2,d,n}$ are the monomials in the set
\begin{align*}
J:=J_j\cup\big\{
    &{X_0}^{i_j+e}{X_1}^{d-i_j-e},\,
        {X_0}^{i_{j-a}+e}{X_1}^{d-i_{j-a}-e}\big\}\mbox{,}
\end{align*}
where ${a=\min(s,1)}$. Therefore we have ${k=T_2+\min(s+1,2)}$ and
\begin{align*}
(d-d_J)n+d_J-dk&=(d-d_J)(T_{j+2}+3s+1)+d_J-d\big(T_2+\min(s+1,2)\big)\\
    &=d\big(T_{j+2}-T_2+3s-\min(s,1)\big)-d_J(T_{j+2}+3s)\\
    &\ge d\big(T_{j+2}-3+3s-\min(s,1)\big)\\
    &\qquad -\left(i_{j-1}+\tfrac{m+1}{2}\right)(T_{j+2}+3s).
\end{align*}
This last expression takes the following forms, depending on the different values of~$j$, $s$ and~$t$:
\begin{align*}
\mbox{- for }&\mbox{$s\le1$ and $t\le j-1$,}\\
    &\tfrac{3}{4}(m-2)j(j-1) + \tfrac{1}{2}(m-4)j +
        \tfrac{9}{2}(m-1) + \tfrac{5}{4}(j-3)^2 +
        \tfrac{11}{4}(j-3) + 4 \\
    &\qquad +\big((m+1)(j-3)+2\big)(1-s) +
        \tfrac{1}{2}(m-3)s +
        (s+3)(j-1-t)>0\mbox{;}\displaybreak[0]\\
\mbox{- for }&\mbox{$s\le1$ and $t=j$,}\\
    &\tfrac{3}{4}(m-1)j(j-1) + \tfrac{1}{2}mj +
        \tfrac{9}{2}(m-1)+ (j-3)^2 + \tfrac{5}{2}j\\
    &\qquad +(m+1)(j-3)(1-s) +
        \tfrac{1}{2}(m-3)s>0\mbox{;}\displaybreak[0]\\
\mbox{- for }&\mbox{$s>1$ and $t\le j-1$,}\\
    &\tfrac{3}{4}(m-2)j(j-1)+\tfrac{1}{2}(m-4)j+5(m-1)+
        \tfrac{5}{4}(j-3)^2 + \tfrac{11}{4}(j-3) + 3\\
    &\qquad +\tfrac{1}{2}( 9(m-1) + 6 )(s-1)+
        4(j-1-t)>0\mbox{;}\displaybreak[0]\\
\mbox{- for }&\mbox{$s>1$ and $t=j$,}\\
    &\tfrac{3}{4}(m-1)j(j-1) + \tfrac{1}{2}mj + 5(m-1)+
        (j-3)^2 + \tfrac{5}{2}(j-1) + \tfrac{3}{2}\\
    &\qquad +\tfrac{1}{2}(9m+3)(s-1)>0.
\end{align*}
Therefore inequality~(\ref{6.6}) is strictly satisfied.

If ${i_{j-1}+e<d_J\le i_j}$, the multiples of~${X_0}^{d_J}$ in~$\I_{2,d,n}$ are the monomials in the set
\[
J:=J_j\cup\big\{{X_0}^{i_j+e}{X_1}^{d-i_j-e}\big\}.
\]
Therefore we have ${k=T_2+1=4}$ and
\begin{align*}
(d-d_J)n+d_J-dk&=(d-d_J)(T_{j+2}+3s+1)+d_J-4d\\
    &=d(T_{j+2}+3s-3)-d_J(T_{j+2}+3s)\\
    &\ge d(T_{j+2}+3s-3)-i_j(T_{j+2}+3s)\\
    &=(m(j+1)+t)(T_{j+2}+3s-3)-(mj+t)(T_{j+2}+3s)\\
    &=\tfrac{1}{2}\big((m-3)j(j-1) +3j(j-3)\big) +3(j-t)
        +3ms>0.
\end{align*}

If ${i_j<d_J\le i_j+e}$, the multiples of~${X_0}^{d_J}$ in~$\I_{2,d,n}$ are the monomials in the set
\[J:=\big\{{X_0}^d,\,{X_0}^{i_j+e}{X_1}^{d-i_j-e}\big\}.\]
Therefore we have ${k=2}$ and
\begin{align*}
(d-d_J)n+d_J-dk&=(d-d_J)(T_{j+2}+3s+1)+d_J-2d\\
    &=d(T_{j+2}+3s-1)-d_J(T_{j+2}+3s)\\
    &\ge d(T_{j+2}+3s-1)-(i_j+e)(T_{j+2}+3s)\\
    &=(m(j+1)+t)(T_{j+2}+3s-1)\\
    &\qquad -\left(mj+t+\tfrac{m+1}{2}\right)(T_{j+2}+3s)\\
    &=\tfrac{1}{2}\big((m-1)j^2 +3(m-1)j +3(m-1)s\big)
        +(j-t)\\
    &\qquad +2(m-1)+1>0.
\end{align*}

If ${i_j+e<d_J<d}$, the only multiple of~${X_0}^{d_J}$ in~$\I_{2,d,n}$ is ${X_0}^d$ and we have nothing to check.

Therefore all possible values of $d_J$ are verified, and hence the syzygy bundle associated to~$\I_{2,d,n}$ is stable.

\bigskip
\noindent{\sc Subcase 2.2:} $r=3s+2$, with ${s\ge0}$.\\
The difference between this case and the previous one is that we are adding the monomial ${X_0}^{i_s+e}{X_2}^{d-i_s-e}$ to~$\I_{2,d,n}$. Therefore we should only worry with the cases ${0<d_J\le i_s+e}$, since for degrees greater than ${i_s+e}$ the set~$J$ of multiples of~${X_0}^{d_J}$ has the same number of elements as in the corresponding set of the previous case, whereas the set $\I_{2,d,n}$ has one more element. Given the fact that the sequence $\left(a_{d,j}\right)_{j\ge2}$ is monotonically increasing (see remark~\ref{adj}), inequality~(\ref{6.6}) is strictly satisfied.

If ${0<d_J\le e}$, the multiples of~${X_0}^{d_J}$ in~$\I_{2,d,n}$ are among the monomials in the set
\begin{align*}
J:=J_1\cup\big\{
    &{X_0}^{i_j+e}{X_1}^{d-i_j-e},\ldots,
        {X_0}^{i_{j-s}+e}{X_1}^{d-i_{j-s}-e},\\
    &{X_0}^e{X_2}^{d-e},\ldots,
        {X_0}^{i_s+e}{X_2}^{d-i_s-e}\big\}.
\end{align*}
Therefore if~$k$ is the number of multiples of~${X_0}^{d_J}$, we have ${k=T_{j+1}+2s+2}$ and
\begin{align*}
(d-d_J)n+d_J-dk&=(d-d_J)(T_{j+2}+3s+2)+d_J-d(T_{j+1}+2s+2)\\
    &=d(j+2+s)-d_J(T_{j+2}+3s+1)\\
    &\ge (m(j+1)+t)(j+2+s)-e(T_{j+2}+3s+1)\\
    &\ge (m(j+1)+t)(j+2+s)-\tfrac{m+1}{2}(T_{j+2}+3s+1)\\
    &=\tfrac{1}{4}(3m-1)j^2+\tfrac{1}{4}(7m-5)(j-2)+\tfrac{7}{2}(m-2)+
        \tfrac{5}{2}\\
    &\qquad+\tfrac{1}{2}(2m(j-2)+3(m-1))s+t(j+2+s)>0.
\end{align*}

If ${e<d_J\le i_1}$, the multiples of~${X_0}^{d_J}$ in~$\I_{2,d,n}$ are the monomials in the set
\begin{align*}
J:=J_1\cup\big\{
    &{X_0}^{i_j+e}{X_1}^{d-i_j-e},\ldots,
        {X_0}^{i_{j-s}+e}{X_1}^{d-i_{j-s}-e},\\
    &{X_0}^{i_1+e}{X_2}^{d-i_1-e},\ldots,
        {X_0}^{i_s+e}{X_2}^{d-i_s-e}\big\}.
\end{align*}
Therefore we have ${k=T_{j+1}+2s+1}$ and
\begin{align*}
(d-d_J)n+d_J-dk&=(d-d_J)(T_{j+2}+3s+2)+d_J-d(T_{j+1}+2s+1)\\
    &=d(j+3+s)-d_J(T_{j+2}+3s+1)\\
    &\ge d(j+3+s)-i_1(T_{j+2}+3s+1)\\
    &=\big(m(j+1)+t\big)(j+3+s)\\
        &\qquad-\big(m+\min(1,t)\big)(T_{j+2}+3s+1).
\end{align*}
This last expression takes the following forms, depending on the different values of~$t$:
\begin{align*}
&\mbox{- for $t=0$, }&&\tfrac{m}{2}(j+2)(j-1)+ m(j-2)s +mj\mbox{;}\\
&\mbox{- for $t>0$, }&&\tfrac{1}{2}
        \big((m-1)(j-2)^2+7(m-1)(j-2)+8(m-3)+12\big)\\
    &&&\qquad +(m(j-3)+m-3+t)s+(t-1)(j+3).
\end{align*}
These expressions are both positive, so inequality~(\ref{6.6}) is strictly satisfied.

If ${i_l<d_J\le i_l+e}$, for ${1\le l\le s}$, we get ${j-s\ge\tfrac{2}{3}j>l}$, since ${3s+2\le j+2}$. Therefore the multiples of~${X_0}^{d_J}$ in~$\I_{2,d,n}$ are the monomials in the set
\begin{align*}
J:=J_{l+1}\cup\big\{
    &{X_0}^{i_j+e}{X_1}^{d-i_j-e},\ldots,
        {X_0}^{i_{j-s}+e}{X_1}^{d-i_{j-s}-e},\\
    &{X_0}^{i_l+e}{X_2}^{d-i_l-e},\ldots,
        {X_0}^{i_s+e}{X_2}^{d-i_s-e}\big\}.
\end{align*}
Therefore we have ${k=T_{j+1-l}+2s+2-l}$ and
\begin{align*}
(d-d_J)n+d_J-dk&=(d-d_J)(T_{j+2}+3s+2)+d_J\\
    &\qquad -d(T_{j+1-l}+2s+2-l)\\
    &=d(T_{j+2}-T_{j+1-l}+s+l)-d_J(T_{j+2}+3s+1)\\
    &\ge d(T_{j+2}-T_{j+1-l}+s+l)-
        \left(i_l+\tfrac{m+1}{2}\right)(T_{j+2}+3s+1).
\end{align*}
This last expression takes the following forms, depending on the different values of~$t$:
\begin{align*}
\mbox{- for }&t\le l\mbox{,}\\
    &\tfrac{1}{2}(m-1)lj(j-l) +\tfrac{3}{4}(m-1)j(j-l)+
        \tfrac{7}{4}(m-3)l(j-l) +\tfrac{1}{4}ml^2\\
    &\qquad +\tfrac{3}{4}(m-4)j +
        \tfrac{1}{2}l^2(j-l)+\tfrac{1}{2}j(j-2l) +2l(j-2l)
        +\tfrac{13}{4}l(l-1)\\
    &\qquad +\tfrac{7}{4}(j-l-2) +
        \tfrac{3}{2}(l+1)+\tfrac{1}{2}(2m(j-2l)+m-3)(s-l)\\
    &\qquad + \big((m+2)l+m\big)(j-l-s) \\
    &\qquad + \tfrac{1}{2}\left((j-l)^2+3(j-2l)+l+4+4s\right)
        (l-t)>0\mbox{;}\displaybreak[0]\\
\mbox{- for }&t>l\mbox{,}\\
    &\tfrac{1}{2}(m-1)lj(j-l) +\tfrac{3}{4}(m-1)j(j-l)+
        \tfrac{7}{4}(m-3)l(j-l) +\tfrac{1}{4}ml^2 \\
    &\qquad +\tfrac{3}{4}(m-1)j +
        \tfrac{1}{2}l^2(j-l-2)+\tfrac{1}{2}j(j-2l-1) \\
    &\qquad +2l(j-2l-1)+\tfrac{17}{4}l(l-1) +
        \tfrac{11}{4}(l-1)+\tfrac{3}{4}\\
    &\qquad +\tfrac{1}{2}(2m(j-2l)+m-3)(s-l)\\
    &\qquad +\big((m+2)l+m\big)(j-l-s)\\
    &\qquad +\tfrac{1}{2}\left(2(l+1)(j-l)+l^2+7l+4+2s\right)
        (t-l)>0.
\end{align*}

If ${i_l+e<d_J\le i_{l+1}}$, for ${1\le l\le s-1}$, the multiples of~${X_0}^{d_J}$ in~$\I_{2,d,n}$ are the monomials in the set
\begin{align*}
J:=J_{l+1}\cup\big\{
    &{X_0}^{i_j+e}{X_1}^{d-i_j-e},\ldots,
        {X_0}^{i_{j-s}+e}{X_1}^{d-i_{j-s}-e},\\
    &{X_0}^{i_{l+1}+e}{X_2}^{d-i_{l+1}-e},\ldots,
        {X_0}^{i_s+e}{X_2}^{d-i_s-e}\big\}.
\end{align*}
Therefore we have ${k=T_{j+1-l}+2s+1-l}$ and
\begin{align*}
(d-d_J)n+d_J-dk&=(d-d_J)(T_{j+2}+3s+2)+d_J\\
    &\qquad -d(T_{j+1-l}+2s+1-l)\\
    &=d(T_{j+2}-T_{j+1-l}+s+1+l)-d_J(T_{j+2}+3s+1)\\
    &\ge d(T_{j+2}-T_{j+1-l}+s+1+l)-i_{l+1}(T_{j+2}+3s+1).
\end{align*}
This last expression takes the following forms, depending on the different values of~$t$:
\begin{align*}
\mbox{- for }&t\le l+1\mbox{,}\\
    &\tfrac{1}{2}(m-1)lj(j-l) +
        \tfrac{1}{2}(m-4)\left(j^2+l^2\right)+
        \tfrac{5}{2}(m-1)(j-l-2) + 2m  \\
    &\qquad +\tfrac{1}{2}(l-1)^2(j-l) +
        \tfrac{3}{2}j(j-l-1)+ 4l(l-1)+\tfrac{13}{2}(l-1) +
        \tfrac{5}{2}\\
    &\qquad + m(j-l-2)(s-l-1) + 2(ml+t)(j-l-1-s)\\
    &\qquad + \tfrac{1}{2}\left( (j-l-2)^2
        + 11(j-l-2) + 16 \right)(l+1-t)
        + lt>0\mbox{;}\displaybreak[0]\\
\mbox{- for }&t>l+1\mbox{,}\\
    &\tfrac{1}{2}(m-4)lj(j-l) +
        \tfrac{1}{2}(m-1)\left(j^2+l^2\right)+
        \tfrac{5}{2}(m-2)(j-l-2) +2m\\
    &\qquad +\tfrac{3}{2}lj(j-l-2) +
        \tfrac{1}{2}l(l-1)(j-l)+ 6l(l-1) +\tfrac{1}{2}j +
        \tfrac{21}{2}(l-1) +\tfrac{7}{2}\\
    &\qquad + ( m(j-l-2) + t )(s-l-1)+
        \big((2m+3)l+3\big)(j-l-1-s)\\
    &\qquad + \tfrac{1}{2}\left( 2l(j-l) + l^2 + 2j
        + 7l + 8 \right)(t-l-1)>0.
\end{align*}
Therefore inequality~(\ref{6.6}) is strictly satisfied.

\bigskip
\noindent{\sc Subcase 2.3:} ${r=3s}$, with ${s\ge1}$.\\
The difference between this case and the previous one is that we are adding the monomial ${X_1}^{i_s+e}{X_2}^{d-i_s-e}$ to~$\I_{2,d,n}$. Since this is not a multiple of~${X_0}^{d_J}$, the set~$J$ of multiples of~${X_0}^{d_J}$ has the same number of elements as in the corresponding sets of the previous case, whereas the set $\I_{2,d,n}$ has one more element. Given the fact that the sequence $\left(a_{d,j}\right)_{j\ge2}$ is monotonically increasing (see remark~\ref{adj}), inequality~(\ref{6.6}) is strictly satisfied.

\bigskip

We can conclude that stability is guaranteed in all cases.
\end{dem}

\section{Sides of the triangle}\label{t2sides}

In this section cases ${d+2<n\le3d}$ are solved. We start with a triangle with one complete side \[\big\{{X_0}^d,\,{X_0}^{d-1}X_1,\ldots,X_0{X_1}^{d-1},
    \,{X_1}^d\big\}\]
and the opposite vertex~${X_2}^d$
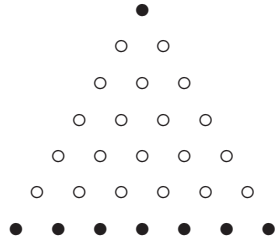
\[
\xymatrix@C=\meiabase@R=\altura@!0{
&&&&&&{\bullet}\\
&&&&&{\circ}&&{\circ}\\
&&&&{\circ}&&{\circ}&&{\circ}\\
&&&{\circ}&&{\circ}&&{\circ}&&{\circ}\\
&&{\circ}&&{\circ}&&{\circ}&&{\circ}&&{\circ}\\
&{\circ}&&{\circ}&&{\circ}&&{\circ}&&{\circ}&&{\circ}\\
{\bullet}&&{\bullet}&&{\bullet}&&{\bullet}&&{\bullet}&&{\bullet}&&{\bullet}\\
}
\]
and add monomials along another side, until we have all but the last monomial.
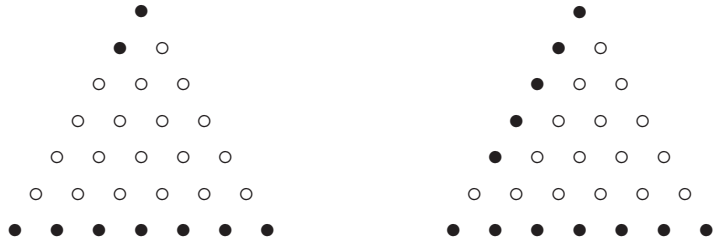
\begin{align*}
\xymatrix@C=\meiabase@R=\altura@!0{
&&&&&&{\bullet}\\
&&&&&{\bullet}&&{\circ}\\
&&&&{\circ}&&{\circ}&&{\circ}\\
&&&{\circ}&&{\circ}&&{\circ}&&{\circ}\\
&&{\circ}&&{\circ}&&{\circ}&&{\circ}&&{\circ}\\
&{\circ}&&{\circ}&&{\circ}&&{\circ}&&{\circ}&&{\circ}\\
{\bullet}&&{\bullet}&&{\bullet}&&{\bullet}&&{\bullet}&&{\bullet}&&{\bullet}\\
}
&&
\xymatrix@C=\meiabase@R=\altura@!0{
&&&&&&{\bullet}\\
&&&&&{\bullet}&&{\circ}\\
&&&&{\bullet}&&{\circ}&&{\circ}\\
&&&{\bullet}&&{\circ}&&{\circ}&&{\circ}\\
&&{\bullet}&&{\circ}&&{\circ}&&{\circ}&&{\circ}\\
&{\circ}&&{\circ}&&{\circ}&&{\circ}&&{\circ}&&{\circ}\\
{\bullet}&&{\bullet}&&{\bullet}&&{\bullet}&&{\bullet}&&{\bullet}&&{\bullet}\\
}
\end{align*}
We then start with the other side, and add the monomial~${X_1}^{d-2}{X_2}^2$, followed by~${X_0}^{d-1}X_2$
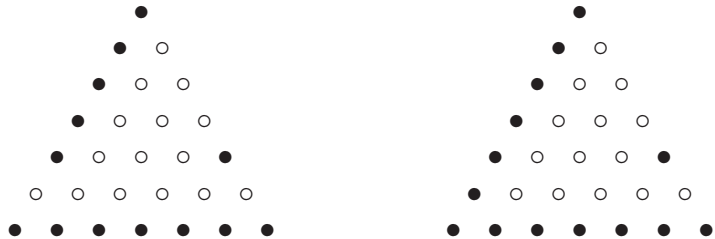
\begin{align*}
\xymatrix@C=\meiabase@R=\altura@!0{
&&&&&&{\bullet}\\
&&&&&{\bullet}&&{\circ}\\
&&&&{\bullet}&&{\circ}&&{\circ}\\
&&&{\bullet}&&{\circ}&&{\circ}&&{\circ}\\
&&{\bullet}&&{\circ}&&{\circ}&&{\circ}&&{\bullet}\\
&{\circ}&&{\circ}&&{\circ}&&{\circ}&&{\circ}&&{\circ}\\
{\bullet}&&{\bullet}&&{\bullet}&&{\bullet}&&{\bullet}&&{\bullet}&&{\bullet}\\
}
&&
\xymatrix@C=\meiabase@R=\altura@!0{
&&&&&&{\bullet}\\
&&&&&{\bullet}&&{\circ}\\
&&&&{\bullet}&&{\circ}&&{\circ}\\
&&&{\bullet}&&{\circ}&&{\circ}&&{\circ}\\
&&{\bullet}&&{\circ}&&{\circ}&&{\circ}&&{\bullet}\\
&{\bullet}&&{\circ}&&{\circ}&&{\circ}&&{\circ}&&{\circ}\\
{\bullet}&&{\bullet}&&{\bullet}&&{\bullet}&&{\bullet}&&{\bullet}&&{\bullet}\\
}
\end{align*}
and end by filling the rest of the remaining side.
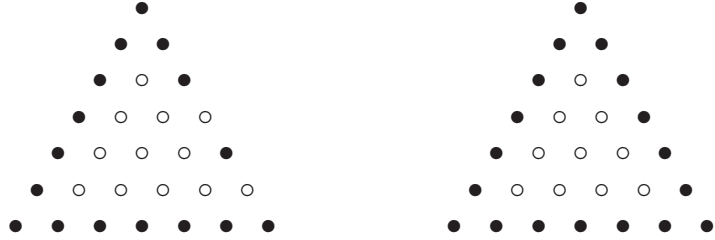
\begin{align*}
\xymatrix@C=\meiabase@R=\altura@!0{
&&&&&&{\bullet}\\
&&&&&{\bullet}&&{\bullet}\\
&&&&{\bullet}&&{\circ}&&{\bullet}\\
&&&{\bullet}&&{\circ}&&{\circ}&&{\circ}\\
&&{\bullet}&&{\circ}&&{\circ}&&{\circ}&&{\bullet}\\
&{\bullet}&&{\circ}&&{\circ}&&{\circ}&&{\circ}&&{\circ}\\
{\bullet}&&{\bullet}&&{\bullet}&&{\bullet}&&{\bullet}&&{\bullet}&&{\bullet}\\
}
&&
\xymatrix@C=\meiabase@R=\altura@!0{
&&&&&&{\bullet}\\
&&&&&{\bullet}&&{\bullet}\\
&&&&{\bullet}&&{\circ}&&{\bullet}\\
&&&{\bullet}&&{\circ}&&{\circ}&&{\bullet}\\
&&{\bullet}&&{\circ}&&{\circ}&&{\circ}&&{\bullet}\\
&{\bullet}&&{\circ}&&{\circ}&&{\circ}&&{\circ}&&{\bullet}\\
{\bullet}&&{\bullet}&&{\bullet}&&{\bullet}&&{\bullet}&&{\bullet}&&{\bullet}\\
}
\end{align*}

\begin{p}\label{exteriortriangle}
For any integers~$d$ and~$n$ such that ${d+2<n\le3d}$, except for ${(d,n)=(2,5)}$, we can obtain a family of~$n$ monomials in $K\left[X_0,X_1,X_2\right]$ of degree~$d$ such that their syzygy bundle is stable. For ${(d,n)=(2,5)}$, there is a family of~$5$ monomials such that their syzygy bundle is semistable.
\end{p}
\begin{dem}
Suppose that ${d\ge4}$, and consider the set
\begin{equation*}
\I_{2,d,n}':=\big\{{X_0}^d,\,{X_0}^{d-1}X_1,\ldots,
    X_0{X_1}^{d-1},\,{X_1}^d,\,{X_2}^d\big\}
\end{equation*}
and the sequence
\begin{align*}
\big(&X_0{X_2}^{d-1},\,{X_0}^2{X_2}^{d-2},\ldots,
        {X_0}^{d-2}{X_2}^2,\\
    &{X_1}^{d-2}{X_2}^2,\,{X_0}^{d-1}X_2,\\
    &X_1{X_2}^{d-1},\,{X_1}^2{X_2}^{d-2},\ldots,
        {X_1}^{d-3}{X_2}^3,\,{X_1}^{d-1}X_2\big).
\end{align*}
Let $1\le i\le3d-2$ be such that ${n=d+2+i}$. Let $\I_{2,d,n}''$ be the set of the first~$i$ monomials in this sequence, and let ${\I_{2,d,n}=\I_{2,d,n}'\cup \I_{2,d,n}''}$. In this way, the number of monomials in~$\I_{2,d,n}$ is~$n$.

For ${0<d_J<d}$, and again following strategy~\ref{X_0}, no monomial of degree~$d_J$ divides a greater number of monomials in~$\I_{2,d,n}$ than ${X_0}^{d_J}$.

If $i\le d-2$, the set of multiples of~${X_0}^{d_J}$ in~$\I_{2,d,n}$ is
\[\big\{{X_0}^d,\ldots,{X_0}^{d_J}X_1^{d-d_J},
    {X_0}^{d_J}{X_2}^{d-d_J},\ldots,{X_0}^e{X_2}^{d-e}\big\}\mbox{,}\]
where $e:=\max\{i,d_J-1\}$ and the list ${{X_0}^{d_J}{X_2}^{d-d_J},\ldots,{X_0}^e{X_2}^{d-e}}$ is understood to be empty if ${e=d_J-1}$. The number of monomials in this set is ${k=d-2d_J+e+2}$, and we get
\begin{align*}
(d-d_J)n+d_J-dk&=i(d-d_J)+dd_J-d_J-de>0.
\end{align*}
If $i=d-1$, the set of multiples of~${X_0}^{d_J}$ is
\[\big\{{X_0}^d,\ldots,{X_0}^{d_J}X_1^{d-d_J},
    {X_0}^{d_J}{X_2}^{d-d_J},\ldots,{X_0}^{d-2}{X_2}^2\big\}.\]
(The list ${X_0}^{d_J}{X_2}^{d-d_J},\ldots,{X_0}^{d-2}{X_2}^2$ is again understood to be empty whenever ${d_J=d-1}$.) The number of monomials in this set is ${k=2d-2d_J}$, and we get
\begin{align*}
(d-d_J)n+d_J-dk&=d>0.
\end{align*}
If $i\ge d$, the set of multiples of~${X_0}^{d_J}$ is
\[\big\{{X_0}^d,\ldots,{X_0}^{d_J}X_1^{d-d_J},
    {X_0}^{d_J}{X_2}^{d-d_J},\ldots,{X_0}^{d-1}X_2\big\}.\]
The number of monomials in this set is ${k=2d-2d_J+1}$, and we get
\begin{align*}
(d-d_J)n+d_J-dk&\ge d-d_J>0.
\end{align*}

In all cases inequality~(\ref{6.6}) is strictly satisfied, which makes the corresponding syzygy bundle stable.

\bigskip

Finally, for the cases left out, ${d=2}$ and ${d=3}$. If ${d=2}$, we get ${4<n\le6}$. For ${n=5}$, we can consider
\[\I_{2,2,5}:=\big\{{X_0}^2,\,{X_1}^2,\,{X_2}^2,\,X_0X_1,
    \,X_0X_2\big\}.\]

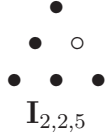
\[
\xymatrix@C=\meiabase@R=\altura@!0{
&&{\bullet}\\
&{\bullet}&&{\circ}\\
{\bullet}&&{\bullet}&&{\bullet}\\
&&\I_{2,2,5}
}
\]

The relevant case to verify is ${J:=\big\{{X_0}^2,\,X_0X_1,\,X_0X_2\big\}}$, and we see that
\begin{align*}
(d-d_J)n+d_J-dk&=5(2-1)+1-2\cdot3=0\mbox{,}
\end{align*}
which means that \[\syz\big({X_0}^2,\,{X_1}^2,\,{X_2}^2,\,X_0X_1,
    \,X_0X_2\big)\]
is strictly semistable.

For ${n=6}$, we can consider
\[\I_{2,2,6}:=\big\{{X_0}^2,\,{X_1}^2,\,{X_2}^2,\,X_0X_1,
    \,X_0X_2,\,X_1X_2\big\}.\]

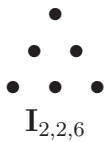
\[
\xymatrix@C=\meiabase@R=\altura@!0{
&&{\bullet}\\
&{\bullet}&&{\bullet}\\
{\bullet}&&{\bullet}&&{\bullet}\\
&&\I_{2,2,6}
}
\]

The relevant case to verify is again ${J:=\big\{{X_0}^2,\,X_0X_1,\,X_0X_2\big\}}$, but this time we get
\begin{align*}
(d-d_J)n+d_J-dk&=6(2-1)+1-2\cdot3=1>0\mbox{,}
\end{align*}
which means that \[\syz\big({X_0}^2,\,{X_1}^2,\,{X_2}^2,\,X_0X_1,
    \,X_0X_2,\,X_1X_2\big)\]
is stable.

If ${d=3}$, we get ${5<n\le9}$. For ${n=6}$, we can consider
\[\I_{2,3,6}:=\big\{{X_0}^3,\,{X_1}^3,\,{X_2}^3,\,{X_0}^2X_1,
    \,X_0{X_2}^2,\,{X_1}^2X_2\big\}.\]

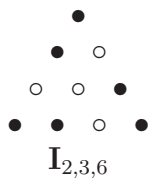
\[
\xymatrix@C=\meiabase@R=\altura@!0{
&&&{\bullet}\\
&&{\bullet}&&{\circ}\\
&{\circ}&&{\circ}&&{\bullet}\\
{\bullet}&&{\bullet}&&{\circ}&&{\bullet}\\
&&&\I_{2,3,6}
}
\]

Since strategy~\ref{X_0} applies to this set, we can look at multiples of~${X_0}^{d_J}$, for ${0<d_J<3}$. For ${d_J=1}$,
the relevant case to verify is \[J:=\big\{{X_0}^3,\,{X_0}^2X_1,\,X_0{X_2}^2\big\}\mbox{,}\]
and we see that
\begin{align*}
(d-d_J)n+d_J-dk&=6(3-1)+1-3\cdot3=4>0.
\end{align*}
For ${d_J=2}$,
the relevant case to verify is ${J:=\big\{{X_0}^3,\,{X_0}^2X_1\big\}}$, and we see that
\begin{align*}
(d-d_J)n+d_J-dk&=6(3-2)+2-3\cdot2=2>0\mbox{,}
\end{align*}
which means that \[\syz\big({X_0}^3,\,{X_1}^3,\,{X_2}^3,\,{X_0}^2X_1,
    \,X_0{X_2}^2,\,{X_1}^2X_2\big)\]
is stable.

For ${n=7}$, we will betray this section's title and add a monomial in the middle of the triangle. Consider
\[\I_{2,3,7}:=\big\{{X_0}^3,\,{X_1}^3,\,{X_2}^3,\,{X_0}^2X_1,
    \,X_0{X_2}^2,\,{X_1}^2X_2,\,X_0X_1X_2\big\}.\]

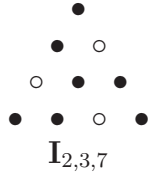
\[
\xymatrix@C=\meiabase@R=\altura@!0{
&&&{\bullet}\\
&&{\bullet}&&{\circ}\\
&{\circ}&&{\bullet}&&{\bullet}\\
{\bullet}&&{\bullet}&&{\circ}&&{\bullet}\\
&&&\I_{2,3,7}
}
\]

Since strategy~\ref{X_0} also applies to this set, we can again look at multiples of~${X_0}^{d_J}$, for ${0<d_J<3}$. For ${d_J=1}$,
the relevant case to verify is \[J:=\big\{{X_0}^3,\,{X_0}^2X_1,\,X_0{X_2}^2,\,X_0X_1X_2\big\}\mbox{,}\]
and we see that
\begin{align*}
(d-d_J)n+d_J-dk&=7(3-1)+1-3\cdot4=3>0.
\end{align*}
For ${d_J=2}$,
the relevant case to verify is ${J:=\big\{{X_0}^3,\,{X_0}^2X_1\big\}}$, and we see that
\begin{align*}
(d-d_J)n+d_J-dk&=7(3-2)+2-3\cdot2=3>0\mbox{,}
\end{align*}
which means that \[\syz\big({X_0}^3,\,{X_1}^3,\,{X_2}^3,\,{X_0}^2X_1,
    \,X_0{X_2}^2,\,{X_1}^2X_2,\,X_0X_1X_2\big)\]
is stable.

For ${n=8}$, we can consider
\[\I_{2,3,8}:=\big\{{X_0}^3,\,{X_1}^3,\,{X_2}^3,\,{X_0}^2X_1,
    \,X_0{X_1}^2,\,{X_0}^2X_2,\,X_0{X_2}^2,\,{X_1}^2X_2\big\}.\]

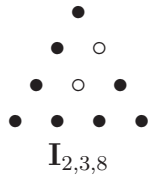
\[
\xymatrix@C=\meiabase@R=\altura@!0{
&&&{\bullet}\\
&&{\bullet}&&{\circ}\\
&{\bullet}&&{\circ}&&{\bullet}\\
{\bullet}&&{\bullet}&&{\bullet}&&{\bullet}\\
&&&\I_{2,3,8}
}
\]

Since strategy~\ref{X_0} applies again to this set, we can once more look at multiples of~${X_0}^{d_J}$, for ${0<d_J<3}$. For ${d_J=1}$,
the relevant case to verify is
\[J:=\big\{{X_0}^3,\,{X_0}^2X_1,\,X_0{X_1}^2,\,{X_0}^2X_2,
    \,X_0{X_2}^2\big\}\mbox{,}\]
and we see that
\begin{align*}
(d-d_J)n+d_J-dk&=8(3-1)+1-3\cdot5=2>0.
\end{align*}
For ${d_J=2}$,
the relevant case to verify is ${J:=\big\{{X_0}^3,\,{X_0}^2X_1,\,{X_0}^2X_2\big\}}$, and we see that
\begin{align*}
(d-d_J)n+d_J-dk&=8(3-2)+2-3\cdot3=1>0\mbox{,}
\end{align*}
which means that \[\syz\big({X_0}^3,\,{X_1}^3,\,{X_2}^3,\,{X_0}^2X_1,
    \,X_0{X_1}^2,\,{X_0}^2X_2,\,X_0{X_2}^2,\,{X_1}^2X_2\big)\]
is stable.

For ${n=9}$, we can consider
\[\I_{2,3,9}:=\big\{{X_0}^3,\,{X_1}^3,\,{X_2}^3,\,{X_0}^2X_1,
    \,X_0{X_1}^2,\,{X_0}^2X_2,\,X_0{X_2}^2,\,{X_1}^2X_2,\,X_1{X_2}^2\big\}.\]

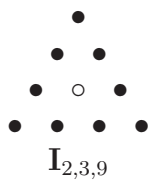
\[
\xymatrix@C=\meiabase@R=\altura@!0{
&&&{\bullet}\\
&&{\bullet}&&{\bullet}\\
&{\bullet}&&{\circ}&&{\bullet}\\
{\bullet}&&{\bullet}&&{\bullet}&&{\bullet}\\
&&&\I_{2,3,9}
}
\]

Since strategy~\ref{X_0} applies again to this set, and, with respect to the previous case, only monomial $X_1{X_2}^2$ was added, the fact that $X_0$ does not occur in this monomial implies that each of the relevant subsets ${J\subseteq\I_{2,3,9}}$ are the same, whereas the set $\I_{2,3,9}$ has one more element. Therefore the fact that sequence $(a_{d,j})_{j\le2}$ is monotonically increasing is enough to establish stability for the syzygy bundle
\[\syz\big({X_0}^3,\,{X_1}^3,\,{X_2}^3,\,{X_0}^2X_1,
    \,X_0{X_1}^2,\,{X_0}^2X_2,\,X_0{X_2}^2,\,{X_1}^2X_2,\,X_1{X_2}^2\big)\]
(see remark~\ref{adj}).
\end{dem}

\section{Triangle's interior}\label{t2interior}

In this section cases ${3d<n\le\tbinom{d+2}{2}}$ are solved, by taking a triangle with full sides and filling its interior starting by the monomials next to the sides and moving gradually to the centre.

\bigskip

We will start by the last case, ${n=\tbinom{d+2}{2}}$. This has already been solved, as mentioned in the previous chapter, in \cite{Fle84} (only semistability) and in \cite{Bal92} and \cite{Pao95} (stability). The proof here presented is again an illustration of Brenner's criterion.

\begin{p}\label{casemaximaln}
The syzygy bundle
\[\mathrm{Syz}\left(\big\{{X_0}^{i_0}{X_1}^{i_1}{X_2}^{i_2}:
    i_0+i_1+i_2=d\big\}\right)\]
is stable.
\end{p}
\begin{dem}
Let ${n=\tbinom{d+2}{2}}$ and let $\I_{2,d,n}:=\big\{{X_0}^{i_0}{X_1}^{i_1}{X_2}^{i_2}: i_0+i_1+i_2=d\big\}$. If~$g$ is the greatest common divisor of monomials in a subset $J\subseteq \I_{2,d,n}$, all monomials in~$J$ are of the form $gh$, with~$h$ a monomial of degree ${d-d_J}$, where~$d_J$ is the degree of~$g$. There are $\tbinom{d-d_J+2}{2}$ monomials of degree ${d-d_J}$, so
\[k=|J|\le\tbinom{d-d_J+2}{2}.\]
Now
\begin{align*}
(d-d_J)n+d_J-dk&=(d-d_J)\tbinom{d+2}{2}+d_J-d\tbinom{d-d_J+2}{2}
    =\tfrac{1}{2}dd_J(d-d_J)>0.
\end{align*}
Therefore inequality~(\ref{6.6}) holds.
\end{dem}

Having this case solved, we can proceed to the general case ${3d<n\le\tbinom{d+2}{2}}$.

\begin{p}
For any integers~$d$ and~$n$ such that ${3d<n\le\tbinom{d+2}{2}}$, there is a family~$\I_{2,d,n}$ of~$n$ monomials in $K\left[X_0,X_1,X_2\right]$ of degree~$d$ such that the corresponding syzygy bundle is stable.
\end{p}
\begin{dem}
We divide this proof in three cases. Let $j\ge1$ be such that $3j<d$ and suppose that
\[\tbinom{d+2}{2}-\tbinom{d+2-3j}{2}<n\le
    \tbinom{d+2}{2}-\tbinom{d+2-3(j+1)}{2}.\]
Note that as $j$ varies, we get all values of~$n$ mentioned, except $\tbinom{d+2}{2}$ when~$d$ is a multiple of~$3$. However, for this highest possible value of~$n$, we have already inequality~(\ref{6.6}) established, by proposition~\ref{casemaximaln}.

\bigskip

\noindent\textsc{Case~1.}
Suppose that
\[{n=\tbinom{d+2}{2}-\tbinom{d+2-3j}{2}+i=
    3dj-\tfrac{9j(j-1)}{2}+i}\mbox{,}\]
with ${1\le i\le d-3j+1}$ and consider the set
\[\I_{2,d,n}':=\big\{{X_0}^{i_0}{X_1}^{i_1}{X_2}^{i_2}:i_0+i_1+i_2=d
    \mbox{ and }(i_0<j\lor i_1<j\lor i_2<j)\big\}.\]
Note that if $j=1$, $\I_{2,d,n}'$ is a set of monomials for which inequality~(\ref{6.6}) holds, by proposition~\ref{exteriortriangle}. Consider the sequence
\begin{align*}
\big(&{X_0}^{d-2j}{X_1}^j{X_2}^j,\,
    {X_0}^{d-2j-1}{X_1}^{j+1}{X_2}^j,
    \ldots,{X_0}^j{X_1}^{d-2j}{X_2}^j\big).
\end{align*}
Let $\I_{2,d,n}''$ be the set of the first~$i$ monomials in this sequence and let us consider ${\I_{2,d,n}:=\I_{2,d,n}'\cup \I_{2,d,n}''}$. Then~$\I_{2,d,n}$ has~$n$ monomials and we will verify that it strictly satisfies inequality~(\ref{6.6}).

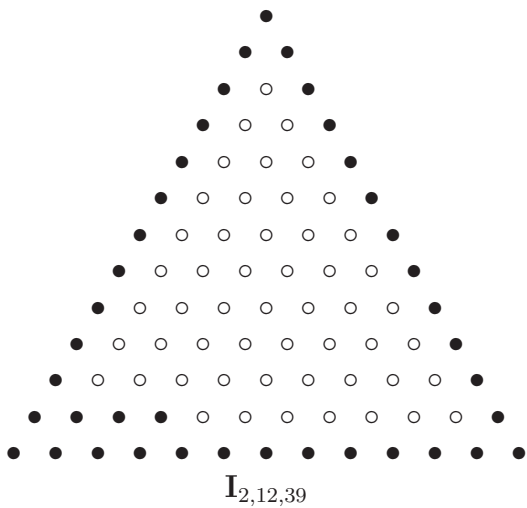
\[
\xymatrix@C=\meiabase@R=\altura@!0{
&&&&&&&&&&&&{\bullet}\\
&&&&&&&&&&&{\bullet}&&{\bullet}\\
&&&&&&&&&&{\bullet}&&{\circ}&&{\bullet}\\
&&&&&&&&&{\bullet}&&{\circ}&&{\circ}&&{\bullet}\\
&&&&&&&&{\bullet}&&{\circ}&&{\circ}&&{\circ}&&{\bullet}\\
&&&&&&&{\bullet}&&{\circ}&&{\circ}&&{\circ}&&{\circ}&&{\bullet}\\
&&&&&&{\bullet}&&{\circ}&&{\circ}&&{\circ}&&{\circ}&&{\circ}&&
    {\bullet}\\
&&&&&{\bullet}&&{\circ}&&{\circ}&&{\circ}&&{\circ}&&{\circ}&&
    {\circ}&&{\bullet}\\
&&&&{\bullet}&&{\circ}&&{\circ}&&{\circ}&&{\circ}&&{\circ}&&
    {\circ}&&{\circ}&&{\bullet}\\
&&&{\bullet}&&{\circ}&&{\circ}&&{\circ}&&{\circ}&&{\circ}&&
    {\circ}&&{\circ}&&{\circ}&&{\bullet}\\
&&{\bullet}&&{\circ}&&{\circ}&&{\circ}&&{\circ}&&
    {\circ}&&{\circ}&&{\circ}&&{\circ}&&{\circ}&&
    {\bullet}\\
&{\bullet}&&{\bullet}&&{\bullet}&&{\bullet}&&{\circ}&&
    {\circ}&&{\circ}&&{\circ}&&{\circ}&&{\circ}&&
    {\circ}&&{\bullet}\\
{\bullet}&&{\bullet}&&{\bullet}&&{\bullet}&&{\bullet}&&
    {\bullet}&&{\bullet}&&{\bullet}&&{\bullet}&&{\bullet}&&
    {\bullet}&&{\bullet}&&{\bullet}\\
&&&&&&&&&&&&\I_{2,12,39}
}
\]

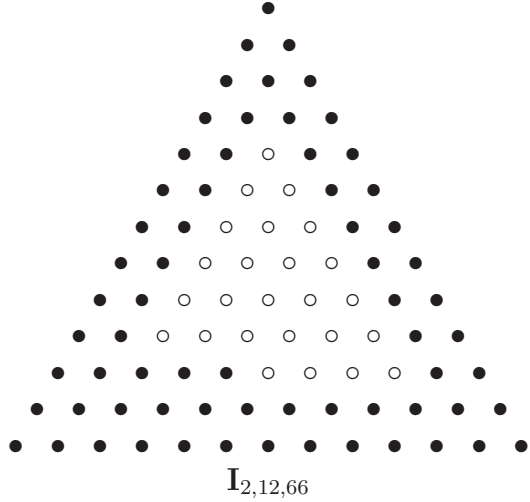
\[
\xymatrix@C=\meiabase@R=\altura@!0{
&&&&&&&&&&&&{\bullet}\\
&&&&&&&&&&&{\bullet}&&{\bullet}\\
&&&&&&&&&&{\bullet}&&{\bullet}&&{\bullet}\\
&&&&&&&&&{\bullet}&&{\bullet}&&{\bullet}&&{\bullet}\\
&&&&&&&&{\bullet}&&{\bullet}&&{\circ}&&{\bullet}&&{\bullet}\\
&&&&&&&{\bullet}&&{\bullet}&&{\circ}&&{\circ}&&{\bullet}&&{\bullet}\\
&&&&&&{\bullet}&&{\bullet}&&{\circ}&&{\circ}&&{\circ}&&{\bullet}&&
    {\bullet}\\
&&&&&{\bullet}&&{\bullet}&&{\circ}&&{\circ}&&{\circ}&&{\circ}&&
    {\bullet}&&{\bullet}\\
&&&&{\bullet}&&{\bullet}&&{\circ}&&{\circ}&&{\circ}&&{\circ}&&
    {\circ}&&{\bullet}&&{\bullet}\\
&&&{\bullet}&&{\bullet}&&{\circ}&&{\circ}&&{\circ}&&{\circ}&&
    {\circ}&&{\circ}&&{\bullet}&&{\bullet}\\
&&{\bullet}&&{\bullet}&&{\bullet}&&{\bullet}&&{\bullet}&&
    {\circ}&&{\circ}&&{\circ}&&{\circ}&&{\bullet}&&
    {\bullet}\\
&{\bullet}&&{\bullet}&&{\bullet}&&{\bullet}&&{\bullet}&&
    {\bullet}&&{\bullet}&&{\bullet}&&{\bullet}&&{\bullet}&&
    {\bullet}&&{\bullet}\\
{\bullet}&&{\bullet}&&{\bullet}&&{\bullet}&&{\bullet}&&
    {\bullet}&&{\bullet}&&{\bullet}&&{\bullet}&&{\bullet}&&
    {\bullet}&&{\bullet}&&{\bullet}\\
&&&&&&&&&&&&\I_{2,12,66}
}
\]

For ${0<d_J<d}$, no monomial of degree~$d_J$ divides a greater number of monomials in~$\I_{2,d,n}$ than~${X_0}^{d_J}$. Therefore all we have to do is count, in each case, the number of multiples of~${X_0}^{d_J}$ which are present in~$\I_{2,d,n}$.

For $d-2j\le d_J<d$, all monomials of degree~$d$ of type ${X_0}^{i_0}{X_1}^{i_1}{X_2}^{i_2}$, with $i_0\ge d_J$, are in~$\I_{2,d,n}$. Therefore the number of multiples of~${X_0}^{d_J}$ in~$\I_{2,d,n}$ is
\[k=\tbinom{d-d_J+2}{2}\]
and we get
\begin{align*}
(d-d_J)n+d_J-dk&=
    (d-d_J)\left(3dj-\tfrac{9j(j-1)}{2}+i\right)
        +d_J-d\tbinom{d-d_J+2}{2}.
\end{align*}
This expression can be rewritten in the two following ways:
\begin{align*}
&\tfrac{1}{2}d(d-d_J)(d_J+j-d)
        +\tfrac{5}{2}(d-3j)(d-d_J)(j-1)
        + 3(d-d_J)j(j-1) \\
    &\qquad +d(d-d_J) + (i-1)(d-d_J)\\
\intertext{and}
&\tfrac{1}{2}d(d-d_J-j)(d_J+2j-d)
        + \tfrac{3}{2}(d-d_J-j)^2(j-1)\\
    &\qquad+ \tfrac{3}{2}(d-d_J-j)d_J(j-1)
        + \tfrac{1}{2}(d_J+2j-d)j^2\\
    &\qquad+ \tfrac{5}{2}(d_J-j)j(j-1)
        + \tfrac{3}{2}(d-d_J-j)j + d_Jj
        + \tfrac{1}{2}j^2 + (i-1)(d-d_J).
\end{align*}
From the first one, we can see that the expression above is positive for
\[{d-j\le d_J<d}\mbox{,}\]
and the second shows us positivity for
\[{d-2j\le d_J<d-j}\]
(since ${3j<d}$, we get in this case ${j<d_J}$).

For $j\le d_J<d-2j$, the monomials in~$\I_{2,d,n}$ that are multiples of~${X_0}^{d_J}$ are the ones in the set
\begin{multline*}
J:=\big\{{X_0}^{i_0}{X_1}^{i_1}{X_2}^{i_2}\in \I_{2,d,n}':
    i_0\ge d_J\big\}\cup{}\\
    \cup\big\{{X_0}^{d-2j}{X_1}^j{X_2}^j,\,
    {X_0}^{d-2j-1}{X_1}^{j+1}{X_2}^j,
    \ldots,{X_0}^{d-2j-e}{X_1}^{j+e}{X_2}^j\big\}\mbox{,}
\end{multline*}
where $e:=\min(i-1,d-2j-d_J)$. Therefore their number is
\[{k=\tbinom{d-d_J+2}{2}-\tbinom{d-2j-d_J+2}{2}+e}.\]
If $i-1\leq d-2j-d_J$, we get
\begin{align*}
(d-d_J)n+d_J-dk&=
    (d-d_J)\left(3dj-\tfrac{9j(j-1)}{2}+i\right)+d_J\\
    &\quad-d\left(\tbinom{d-d_J+2}{2}
        -\tbinom{d-2j-d_J+2}{2}+i-1\right)\\
    &=(d-j-d_J)(d-2j-d_J)j + (d-2j-d_J)d_J(j-1) \\
    &\qquad + \tfrac{1}{2}dj^2 + \tfrac{7}{2}(d_J-j)j(j-1)
        + \tfrac{3}{2}j^2(j-1) \\
    &\qquad + \tfrac{3}{2}(d-2j-d_J)j + \tfrac{1}{2}d_Jj
        + j^2 + d \\
    &\qquad + (d-2j-d_J+1-i)d_J>0
\end{align*}
since $d-2j-d_J>0$ and $j\geq1$. If $i-1>d-2j-d_J$, we get
\begin{align*}
(d-d_J)n+d_J-dk&=
    (d-d_J)\left(3dj-\tfrac{9j(j-1)}{2}+i\right)+d_J\\
    &\quad-d\left(\tbinom{d-d_J+2}{2}
        -\tbinom{d-2j-d_J+2}{2}+d-2j-d_J\right)\\
    &=(d-d_J-2j)(d-j)(j-1) + \tfrac{1}{2}dj^2 +
        \tfrac{7}{2}(d_J-j)j(j-1) \\
    &\quad +\tfrac{3}{2}j^3 + (d-d_J-2j)^2 +
        \tfrac{5}{2}(d-d_J-2j)j
        + \tfrac{1}{2}(d_J-j)j \\
    &\quad + d + (i-1-d+2j+d_J)(d-d_J)>0
\end{align*}
since $d-2j-d_J>0$ and $j\geq1$.

For $0<d_J<j$, the number of monomials in~$\I_{2,d,n}$ that are multiples of~${X_0}^{d_J}$ is
\[{k=\tbinom{d-d_J+2}{2}-\tbinom{d-3j+2}{2}+i}\mbox{,}\]
and we get
\begin{align*}
(d-d_J)n+d_J-dk&=
    (d-d_J)\left(3dj-\tfrac{9j(j-1)}{2}+i\right)+d_J\\
    &\quad-d\left(\tbinom{d-d_J+2}{2}
        -\tbinom{d-3j+2}{2}+i\right)\\
    &=(d-3j)^2d_J + 3(d-3j)(j-d_J)d_J
        + \tfrac{5}{2}(d-3j){d_J}^2\\
    &\quad + \tfrac{9}{2}j(j-d_J)d_J + 3j{d_J}^2
        + \tfrac{1}{2}(d-3j)d_J\\
    &\quad + (d-3j+1-i)d_J>0.
\end{align*}
In all cases, inequality~(\ref{6.6}) is strictly satisfied, which makes the corresponding syzygy bundle stable.

\bigskip
\noindent\textsc{Case~2.}
Now suppose that
\[n=\tbinom{d+2}{2}-\tbinom{d+1-3j}{2}+i=
    3dj+d+1-\tfrac{3j(3j-1)}{2}+i\mbox{,}\]
with ${1\le i\le d-3j}$ and consider the set
\[\I_{2,d,n}':=\big\{{X_0}^{i_0}{X_1}^{i_1}{X_2}^{i_2}:
    i_0+i_1+i_2=d\mbox{ and }
    (i_0<j\lor i_1<j\lor i_2\le j)\big\}.\]
Consider the sequence
\begin{align*}
\big(&{X_0}^j{X_1}^j{X_2}^{d-2j},\,
    {X_0}^{j+1}{X_1}^j{X_2}^{d-2j-1},
    \ldots,{X_0}^{d-2j-1}{X_1}^j{X_2}^{j+1}\big).
\end{align*}
Let $\I_{2,d,n}''$ be the set of the first~$i$ monomials in this sequence and let us consider ${\I_{2,d,n}:=\I_{2,d,n}'\cup \I_{2,d,n}''}$. Then~$\I_{2,d,n}$ has~$n$ monomials and we will verify that it strictly satisfies inequality~(\ref{6.6}).

As in the previous case, we follow strategy~\ref{X_0}, and therefore for ${0<d_J<d}$, no monomial of degree~$d_J$ divides a greater number of monomials in~$\I_{2,d,n}$ than~${X_0}^{d_J}$. Therefore all we have to do is count, in each case, the number of multiples of~${X_0}^{d_J}$ which are present in~$\I_{2,d,n}$.

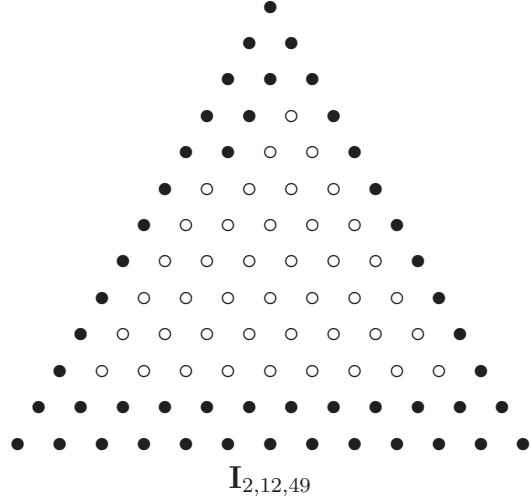
\[
\xymatrix@C=\meiabase@R=\altura@!0{
&&&&&&&&&&&&{\bullet}\\
&&&&&&&&&&&{\bullet}&&{\bullet}\\
&&&&&&&&&&{\bullet}&&{\bullet}&&{\bullet}\\
&&&&&&&&&{\bullet}&&{\bullet}&&{\circ}&&{\bullet}\\
&&&&&&&&{\bullet}&&{\bullet}&&{\circ}&&{\circ}&&{\bullet}\\
&&&&&&&{\bullet}&&{\circ}&&{\circ}&&{\circ}&&{\circ}&&{\bullet}\\
&&&&&&{\bullet}&&{\circ}&&{\circ}&&{\circ}&&{\circ}&&{\circ}&&
    {\bullet}\\
&&&&&{\bullet}&&{\circ}&&{\circ}&&{\circ}&&{\circ}&&{\circ}&&
    {\circ}&&{\bullet}\\
&&&&{\bullet}&&{\circ}&&{\circ}&&{\circ}&&{\circ}&&{\circ}&&
    {\circ}&&{\circ}&&{\bullet}\\
&&&{\bullet}&&{\circ}&&{\circ}&&{\circ}&&{\circ}&&{\circ}&&
    {\circ}&&{\circ}&&{\circ}&&{\bullet}\\
&&{\bullet}&&{\circ}&&{\circ}&&{\circ}&&{\circ}&&{\circ}&&
    {\circ}&&{\circ}&&{\circ}&&{\circ}&&{\bullet}\\
&{\bullet}&&{\bullet}&&{\bullet}&&{\bullet}&&{\bullet}&&{\bullet}&&
    {\bullet}&&{\bullet}&&{\bullet}&&{\bullet}&&{\bullet}&&
    {\bullet}\\
{\bullet}&&{\bullet}&&{\bullet}&&{\bullet}&&{\bullet}&&{\bullet}&&
    {\bullet}&&{\bullet}&&{\bullet}&&{\bullet}&&{\bullet}&&
    {\bullet}&&{\bullet}\\
&&&&&&&&&&&&\I_{2,12,49}
}
\]

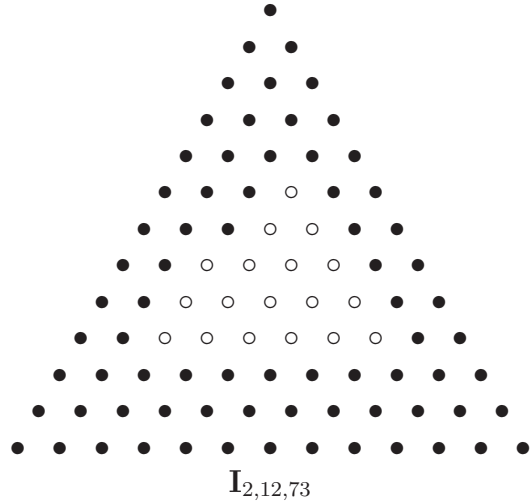
\[
\xymatrix@C=\meiabase@R=\altura@!0{
&&&&&&&&&&&&{\bullet}\\
&&&&&&&&&&&{\bullet}&&{\bullet}\\
&&&&&&&&&&{\bullet}&&{\bullet}&&{\bullet}\\
&&&&&&&&&{\bullet}&&{\bullet}&&{\bullet}&&{\bullet}\\
&&&&&&&&{\bullet}&&{\bullet}&&{\bullet}&&{\bullet}&&{\bullet}\\
&&&&&&&{\bullet}&&{\bullet}&&{\bullet}&&{\circ}&&{\bullet}&&{\bullet}\\
&&&&&&{\bullet}&&{\bullet}&&{\bullet}&&{\circ}&&{\circ}&&{\bullet}&&
    {\bullet}\\
&&&&&{\bullet}&&{\bullet}&&{\circ}&&{\circ}&&{\circ}&&{\circ}&&
    {\bullet}&&{\bullet}\\
&&&&{\bullet}&&{\bullet}&&{\circ}&&{\circ}&&{\circ}&&{\circ}&&
    {\circ}&&{\bullet}&&{\bullet}\\
&&&{\bullet}&&{\bullet}&&{\circ}&&{\circ}&&{\circ}&&{\circ}&&
    {\circ}&&{\circ}&&{\bullet}&&{\bullet}\\
&&{\bullet}&&{\bullet}&&{\bullet}&&{\bullet}&&{\bullet}&&{\bullet}&&
    {\bullet}&&{\bullet}&&{\bullet}&&{\bullet}&&{\bullet}\\
&{\bullet}&&{\bullet}&&{\bullet}&&{\bullet}&&{\bullet}&&{\bullet}&&
    {\bullet}&&{\bullet}&&{\bullet}&&{\bullet}&&{\bullet}&&
    {\bullet}\\
{\bullet}&&{\bullet}&&{\bullet}&&{\bullet}&&{\bullet}&&{\bullet}&&
    {\bullet}&&{\bullet}&&{\bullet}&&{\bullet}&&{\bullet}&&
    {\bullet}&&{\bullet}\\
&&&&&&&&&&&&\I_{2,12,73}
}
\]

For $d-2j\le d_J<d$, all monomials of degree~$d$ of type ${X_0}^{i_0}{X_1}^{i_1}{X_2}^{i_2}$, with $i_0\ge d_J$, are in~$\I_{2,d,n}$. Therefore the number of multiples of~${X_0}^{d_J}$ in~$\I_{2,d,n}$ is
\[k=\tbinom{d-d_J+2}{2}\mbox{,}\]
as it was in case~1, and we can claim that since all values are the same except for $n$, which is bigger, inequality~(\ref{6.6}) is strictly satisfied, due to the fact that sequence $(a_{d,j})_{j\le2}$ is monotonically increasing (see remark~\ref{adj}).

For $j\le d_J<d-2j$, the monomials in~$\I_{2,d,n}$ that are multiples of~${X_0}^{d_J}$ are the ones in the set
\begin{multline*}
J:=\big\{{X_0}^{i_0}{X_1}^{i_1}{X_2}^{i_2}\in \I_{2,d,n}':
    i_0\ge d_J\big\}\cup
    \big\{{X_0}^{d_J}{X_1}^j{X_2}^{d-j-d_J},\\
    {X_0}^{d_J+1}{X_1}^j{X_2}^{d-j-d_J-1},
    \ldots,{X_0}^{j+i-1}{X_1}^j{X_2}^{d-2j-i+1}\big\}\mbox{,}
\end{multline*}
where this last set is understood to be empty if ${j+i-1<d_J}$. Therefore their number is
\[k=\tbinom{d-d_J+2}{2}-\tbinom{d-2j-d_J+1}{2}+\max(0,j+i-d_J)
    \mbox{.}\]
If $j+i\le d_J$, we get (keeping in mind that $i\geq1$)
\begin{align*}
(d-d_J)n+d_J-dk&=
    (d-d_J)\left(3dj+d+1-\tfrac{3j(3j-1)}{2}+i\right)+d_J\\
    &\quad-d\left(\tbinom{d-d_J+2}{2}-\tbinom{d-2j-d_J+1}{2}
        \right)\\
    &=(d-j)(d-2j-d_J)j + \tfrac{1}{2}dj^2
        + \tfrac{7}{2}(d_J-j)j(j-1) + \tfrac{3}{2}j^3\\
    &\quad + \tfrac{1}{2}(d-2j-d_J)j + \tfrac{5}{2}(d_J-j)j
        + i(d-d_J)>0.
\end{align*}
If $j+i>d_J$, we get
\begin{align*}
(d-d_J)n+d_J-dk&=
    (d-d_J)\left(3dj+d+1-\tfrac{3j(3j-1)}{2}+i\right)+d_J\\
    &\quad-d\left(\tbinom{d-d_J+2}{2}-\tbinom{d-2j-d_J+1}{2}
        +j+i-d_J\right)\\
    &=(d-j)(d-2j-d_J)j + \tfrac{1}{2}(d-2j-d_J)j(j-1)\\
    &\quad + 4(d_J-j)j^2 + 3j^3+ (d_J-j)j+ (d-3j-i)d_J>0.
\end{align*}

For $0<d_J<j$, the number of monomials in~$I_{d,n}$ that are multiples of~${X_0}^{d_J}$ is
\[k=\tbinom{d-d_J+2}{2}-\tbinom{d-3j+1}{2}+i\mbox{,}\]
and we get
\begin{align*}
(d-d_J)n+d_J-dk&=
    (d-d_J)\left(3dj+d+1-\tfrac{3j(3j-1)}{2}+i\right)+d_J\\
    &\quad-d\left(\tbinom{d-d_J+2}{2}-\tbinom{d-3j+1}{2}+i
        \right)\\
    &=(d-j)(d-3j)d_J + d(j-d_J)d_J
        +\tfrac{1}{2}dd_J(d_J-1)\\
    &\quad + \tfrac{3}{2}j^2d_J
        + \tfrac{3}{2}jd_J+ (d-3j-i)d_J>0.
\end{align*}
Again in all cases, inequality~(\ref{6.6}) is strictly satisfied, which makes the corresponding syzygy bundle stable.

\bigskip
\noindent\textsc{Case~3.}
If $d=3j+1$, case~2 has exhausted all possible monic monomials of degree~$d$, and this proof is ended.

If ${d>3j+1}$, then suppose that
\[n=\tbinom{d+2}{2}-\tbinom{d-3j}{2}+i=
    3dj+2d+1-\tfrac{3j(3j+1)}{2}+i\mbox{,}\]
with ${1\le i\le d-3j-1}$ and consider the set
\[\I_{2,d,n}':=\big\{{X_0}^{i_0}{X_1}^{i_1}{X_2}^{i_2}:
    i_0+i_1+i_2=d\mbox{ and }
    (i_0<j\lor i_1\le j\lor i_2\le j)\big\}.\]
Consider the sequence
\begin{align*}
\big(&{X_0}^j{X_1}^{j+1}{X_2}^{d-2j-1},\,
    {X_0}^j{X_1}^{j+2}{X_2}^{d-2j-2},
    \ldots,{X_0}^j{X_1}^{d-2j-1}{X_2}^{j+1}\big).
\end{align*}
Let $\I_{2,d,n}''$ be the set of the first~$i$ monomials in this sequence and let us consider ${\I_{2,d,n}:=\I_{2,d,n}'\cup \I_{2,d,n}''}$. Then~$\I_{2,d,n}$ has~$n$ monomials and we will verify that it strictly satisfies inequality~(\ref{6.6}).

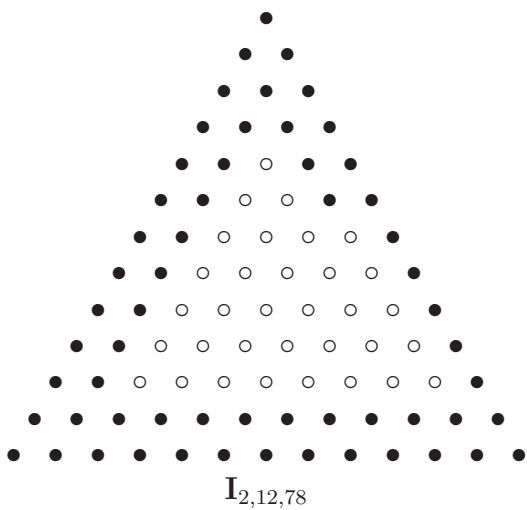
\[
\xymatrix@C=\meiabase@R=\altura@!0{
&&&&&&&&&&&&{\bullet}\\
&&&&&&&&&&&{\bullet}&&{\bullet}\\
&&&&&&&&&&{\bullet}&&{\bullet}&&{\bullet}\\
&&&&&&&&&{\bullet}&&{\bullet}&&{\bullet}&&{\bullet}\\
&&&&&&&&{\bullet}&&{\bullet}&&{\circ}&&{\bullet}&&{\bullet}\\
&&&&&&&{\bullet}&&{\bullet}&&{\circ}&&{\circ}&&{\bullet}&&{\bullet}\\
&&&&&&{\bullet}&&{\bullet}&&{\circ}&&{\circ}&&{\circ}&&{\circ}&&
    {\bullet}\\
&&&&&{\bullet}&&{\bullet}&&{\circ}&&{\circ}&&{\circ}&&{\circ}&&
    {\circ}&&{\bullet}\\
&&&&{\bullet}&&{\bullet}&&{\circ}&&{\circ}&&{\circ}&&{\circ}&&
    {\circ}&&{\circ}&&{\bullet}\\
&&&{\bullet}&&{\bullet}&&{\circ}&&{\circ}&&{\circ}&&{\circ}&&
    {\circ}&&{\circ}&&{\circ}&&{\bullet}\\
&&{\bullet}&&{\bullet}&&{\circ}&&{\circ}&&{\circ}&&{\circ}&&
    {\circ}&&{\circ}&&{\circ}&&{\circ}&&{\bullet}\\
&{\bullet}&&{\bullet}&&{\bullet}&&{\bullet}&&{\bullet}&&{\bullet}&&
    {\bullet}&&{\bullet}&&{\bullet}&&{\bullet}&&{\bullet}&&
    {\bullet}\\
{\bullet}&&{\bullet}&&{\bullet}&&{\bullet}&&{\bullet}&&{\bullet}&&
    {\bullet}&&{\bullet}&&{\bullet}&&{\bullet}&&{\bullet}&&
    {\bullet}&&{\bullet}\\
&&&&&&&&&&&&\I_{2,12,78}
}
\]

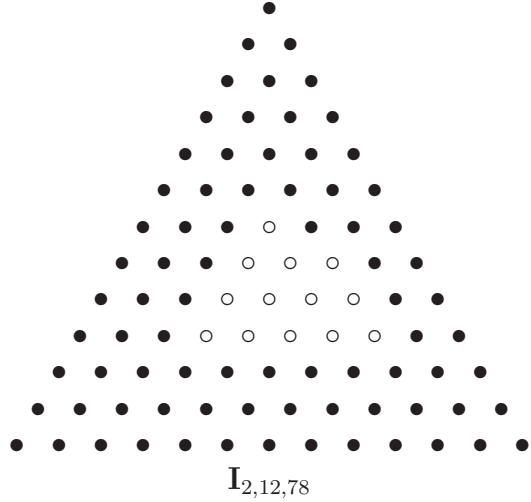
\[
\xymatrix@C=\meiabase@R=\altura@!0{
&&&&&&&&&&&&{\bullet}\\
&&&&&&&&&&&{\bullet}&&{\bullet}\\
&&&&&&&&&&{\bullet}&&{\bullet}&&{\bullet}\\
&&&&&&&&&{\bullet}&&{\bullet}&&{\bullet}&&{\bullet}\\
&&&&&&&&{\bullet}&&{\bullet}&&{\bullet}&&{\bullet}&&{\bullet}\\
&&&&&&&{\bullet}&&{\bullet}&&{\bullet}&&{\bullet}&&{\bullet}&&{\bullet}\\
&&&&&&{\bullet}&&{\bullet}&&{\bullet}&&{\circ}&&{\bullet}&&{\bullet}&&
    {\bullet}\\
&&&&&{\bullet}&&{\bullet}&&{\bullet}&&{\circ}&&{\circ}&&{\circ}&&
    {\bullet}&&{\bullet}\\
&&&&{\bullet}&&{\bullet}&&{\bullet}&&{\circ}&&{\circ}&&{\circ}&&
    {\circ}&&{\bullet}&&{\bullet}\\
&&&{\bullet}&&{\bullet}&&{\bullet}&&{\circ}&&{\circ}&&{\circ}&&
    {\circ}&&{\circ}&&{\bullet}&&{\bullet}\\
&&{\bullet}&&{\bullet}&&{\bullet}&&{\bullet}&&{\bullet}&&{\bullet}&&
    {\bullet}&&{\bullet}&&{\bullet}&&{\bullet}&&{\bullet}\\
&{\bullet}&&{\bullet}&&{\bullet}&&{\bullet}&&{\bullet}&&{\bullet}&&
    {\bullet}&&{\bullet}&&{\bullet}&&{\bullet}&&{\bullet}&&
    {\bullet}\\
{\bullet}&&{\bullet}&&{\bullet}&&{\bullet}&&{\bullet}&&{\bullet}&&
    {\bullet}&&{\bullet}&&{\bullet}&&{\bullet}&&{\bullet}&&
    {\bullet}&&{\bullet}\\
&&&&&&&&&&&&\I_{2,12,78}
}
\]

As in the previous cases, strategy\ref{X_0} is followed, and for ${0<d_J<d}$, no monomial of degree~$d_J$ divides a greater number of monomials in~$\I_{2,d,n}$ than ${X_0}^{d_J}$. Therefore all we have to do is count, in each case the number of multiples of~${X_0}^{d_J}$ which are present in~$\I_{2,d,n}$.

For $d-2j\le d_J<d$, all monomials of degree~$d$ of type ${X_0}^{i_0}{X_1}^{i_1}{X_2}^{i_2}$, with $i_0\ge d_J$, are in~$\I_{2,d,n}$. Therefore the number of multiples of~${X_0}^{d_J}$ in~$\I_{2,d,n}$ is
\[k=\tbinom{d-d_J+2}{2}\mbox{,}\]
as it was in cases~1 and~2, and we can claim that since all values are the same except for~$n$, which is bigger, inequality~(\ref{6.6}) is strictly satisfied, due to the fact that sequence $(a_{d,j})_{j\le2}$ is monotonically increasing (see remark~\ref{adj}).

For $j<d_J<d-2j$, an analogous argument based on calculations for step~2 allows us to claim that inequality~(\ref{6.6}) is strictly satisfied.

For $0<d_J\le j$, the number of monomials in~$\I_{2,d,n}$ that are multiples of~${X_0}^{d_J}$ is
\[k=\tbinom{d-d_J+2}{2}-\tbinom{d-3j}{2}+i\mbox{,}\]
and we get (keeping in mind that $i\leq d-3j-1$)
\begin{align*}
(d-d_J)n+d_J-dk&=
    (d-d_J)\left(3dj+2d+1-\tfrac{3j(3j+1)}{2}+i\right)+d_J\\
    &\quad-d\left(\tbinom{d-d_J+2}{2}-\tbinom{d-3j}{2}+i
        \right)\\
    &=(d-2j)(d-3j)d_J + 2(d-j)(j-d_J)d_J\\
    &\quad + \tfrac{3}{2}(d-j)d_J(d_J-1)
        + \tfrac{1}{2}j(j-d_J)d_J\\
    &\quad + 3jd_J + d_J + (d-3j-1-i)d_J>0\mbox{.}
\end{align*}
Again in all cases, inequality~(\ref{6.6}) is strictly satisfied, which makes the corresponding syzygy bundle stable and concludes the proof.
\end{dem}

\section{Main theorem}\label{t2main}

In this last section of the chapter, and adding up all the results in this chapter, we conclude the following.

\begin{te}\label{main2}
Let~$d$ and~$n$ be integers such that $3\le n\le\tbinom{d+2}{2}$, and $(d,n)\ne(2,5)$. Then there is a family of~$n$ monomials in $K\left[X_0,X_1,X_2\right]$ of degree~$d$ such that their syzygy bundle is stable. If $(d,n)=(2,5)$, there is a family of~$n$ monomials such that their syzygy bundle is semistable.
\end{te}

Problem~\ref{6.9}, presented by Brenner in \cite{Bre08a,Bre08b}, is thus completely solved for ${N=2}$. This will allow us to assert that for $3\le n\le\tbinom{d+2}{2}$, with $(d,n)\ne(2,5)$, a syzygy bundle associated to a family of~$n$ general homogeneous polynomials in $K\left[X_0,X_1,X_2\right]$ of degree~$d$ is stable, for the condition of stability in the moduli space is an open one.

\bigskip

We know there is no family of $5$ quadratic monomials in $K\left[X_0,X_1,X_2\right]$ such that their syzygy bundle is stable, but the following problem remains open:

\begin{pr}
Is there a family of $5$ quadratic homogeneous polynomials in $K\left[X_0,X_1,X_2\right]$ such that their syzygy bundle is stable?
\end{pr}

Given $5$ quadratic homogeneous polynomials ${f_1,\ldots,f_5}$ in $\kk\left[X_0,X_1,X_2\right]$, their syzygy bundle's presenting sequence is
\[\xymatrix{0\ar[r]&\syz(f_1,\ldots,f_5)\ar[r]&
    \OO(-2)^5\ar[r]&\OO\ar[r]&0\mbox{.}}\]
To check whether this bundle is stable, one may be tempted to use Hoppe's criterion (lemma~\ref{Hop2.6}). Unfortunately, it leads nowhere in this case. Let us see why. Since this bundle's rank is $4$, we have to compute the global sections of $\syz(f_1,\ldots,f_5)_{\mathrm{norm}}$, $\big(\Lambda^2\syz(f_1,\ldots,f_5)\big)_{\mathrm{norm}}$, and
$\big(\Lambda^3\syz(f_1,\ldots,f_5)\big)_{\mathrm{norm}}$.

Since the syzygy bundle's first Chern class is ${c_1\big(\syz(f_1,\ldots,f_5)\big)=-10}$ and its rank is $4$, its normalised is
\[\syz(f_1,\ldots,f_5)_{\mathrm{norm}}=
    \syz(f_1,\ldots,f_5)(2)\mbox{,}\]
for ${c_1\big(\syz(f_1,\ldots,f_5)(2)\big)=-10+4\cdot2=-2}$. Now the minimal resolution\linebreak[4]
of~${\tfrac{R}{I}}$, where ${R:=\kk\left[X_0,X_1,X_2\right]}$ and ${I:=(f_1,\ldots,f_5)}$, is
\[\xymatrix{0\ar[r]&R(-5)\ar[r]&R(-3)^5\ar[r]&R(-2)^5
    \ar[r]&R\ar[r]&\tfrac{R}{I}\ar[r]&0\mbox{.}}\]
If we consider the corresponding sequence of sheaves, we see that since the syzygy bundle is the kernel of \smash{$\xymatrix@1{\OO(-2)^5\ar[r]&\OO\mbox{,}}$} it is isomorphic to the cokernel of \smash{$\xymatrix@1{\OO(-5)\ar[r]&\OO(-3)^5\mbox{.}}$} We have therefore another exact sequence
\[\xymatrix{0\ar[r]&\OO(-5)\ar[r]&\OO(-3)^5\ar[r]&
    \syz(f_1,\ldots,f_5)\ar[r]&0\mbox{,}}\]
which gives us
\[\xymatrix{0\ar[r]&\OO(-3)\ar[r]&\OO(-1)^5\ar[r]&
    \syz(f_1,\ldots,f_5)_{\mathrm{norm}}\ar[r]&0\mbox{.}}\]
Therefore $\syz(f_1,\ldots,f_5)_{\mathrm{norm}}$ has no non\h{zero} global sections.

To compute the global sections of $\big(\Lambda^3\syz(f_1,\ldots,f_5)\big)_{\mathrm{norm}}$, observe that
\begin{align*}
\Lambda^3\syz(f_1,\ldots,f_5)&\cong
    \syz(f_1,\ldots,f_5)^\vee\otimes
        \Lambda^4\syz(f_1,\ldots,f_5)\\
    &\cong\syz(f_1,\ldots,f_5)^\vee\otimes\OO(-10)\\
    &=\syz(f_1,\ldots,f_5)^\vee(-10)\mbox{,}
\end{align*}
and the first Chern class of this bundle is
\[c_1\big(\syz(f_1,\ldots,f_5)^\vee(-10)\big)=
    -c_1\big(\syz(f_1,\ldots,f_5)\big)+4\cdot(-10)=-30\mbox{.}\]
Therefore to get its normalised, we must make the tensor product with ${\OO(7)}$, to get ${\syz(f_1,\ldots,f_5)^\vee(-3)}$, since
\[c_1\big(\syz(f_1,\ldots,f_5)^\vee(-3)\big)=-2\mbox{.}\]

Tensoring the dual of the syzygy bundle's presenting sequence by ${\OO(-3)}$, we get
\[\xymatrix{0\ar[r]&\OO(-3)\ar[r]&\OO(-1)^5\ar[r]&
    \syz(f_1,\ldots,f_5)^\vee(-3)\ar[r]&0\mbox{,}}\]
which tells us that $\big(\Lambda^3\syz(f_1,\ldots,f_5)\big)_{\mathrm{norm}}$ has no non\h{zero} global sections.

Finally, we have computed $H^0\big(\Lambda^2\syz(f_1,\ldots,f_5)\big)_{\mathrm{norm}}$ using Macaylay~2 \cite{GS}, with $5$ randomly chosen homogeneous polynomials ${f_1,\ldots,f_5}$ in\linebreak[4]
$\qq\left[X_0,X_1,X_2\right]$, and obtained non\h{zero} global sections. Therefore there seems to be no way of using Hoppe's criterion for this bundle. Since it is only a sufficient condition for stability, we can conclude nothing from here.

\chapter{Monomials in four or more variables}
\label{t3}

In this chapter an answer to problem~\ref{6.9} is presented for ${N\ge3}$. In each case there is a family of monomials whose corresponding syzygy bundle is stable.

In general monomials in $\kk[X_0,\ldots,X_N]$ of a given degree~$d$ can be represented in a hypertetrahedron, in an analogous manner to the triangles we have seen in case ${N=2}$ (recall we are assuming monomials to be monic). This hypertetrahedron is no more than the graph whose vertexes are all monomials of degree~$d$, and where two monomials are connected by an edge if and only if their greatest common divisor has degree ${d-1}$. We shall call the \emph{$i$th face} of this hypertetrahedron the set of monomials where the variable $X_i$ does not occur.

The chapter is divided into different sections, according to different values of~$n$. Recall that we have ${N+1\le n\le\tbinom{d+N}{N}}$. For the first cases, in section~\ref{t3first}, with
\[{N+1\le n\le\tbinom{d+N-1}{N-1}+1}\mbox{,}\]
each family of ${n-1}$ monomials in $\kk[X_0,\ldots,X_{N-1}]$ whose syzygy bundle over $\pp^{N-1}$ is stable generates a family of~$n$ monomials in $\kk[X_0,\ldots,X_N]$ whose syzygy bundle over $\pp^N$ is also stable. Cases \[{\tbinom{d+N-1}{N-1}+1<n\le\tbinom{d+N}{N}-\tbinom{d-1}{N}}\]
are solved in section~\ref{t3faces} by taking the $N$th face and the vertex ${X_N}^d$ of the hypertetrahedron, and adding monomials in the remaining faces. Taking the set of all the hypertetrahedron's faces and adding the monomials in its interior which are closest to the vertexes gives us a solution to cases \[{\tbinom{d+N}{N}-\tbinom{d-1}{N}<n\le\tbinom{d+N}{N}
    -\tbinom{d-1}{N}+N+1}\mbox{,}\]
treated in the beginning of section~\ref{t3interior}. The last cases, with \[{\tbinom{d+N}{N}-\tbinom{d-1}{N}+N+1<n\le\tbinom{d+N}{N}}\mbox{,}\] also in section~\ref{t3interior}, are solved by taking a family of monomials of lower degree whose syzygy bundle is stable, multiplying them by ${X_0\cdots X_N}$, and adding all monomials in the faces of the hypertetrahedron. This is a generalization of a lemma by Brenner, made for the case ${N=2}$ in his notes \cite{Bre}, which he kindly shared. Finally, gathering all information from both this chapter and the previous one, the main theorem in this work is stated in section~\ref{t3main}.

We conclude this chapter stating two open questions that arise naturally, now that problem~\ref{6.9} is solved. The first one is the generalisation of this problem one gets if we consider polynomials of arbitrary degree. The second one is to know whether syzygy bundles of higher order are stable, or semistable. This last question involves the open problem of finding the minimal free resolutions of general forms.

\section{First cases}\label{t3first}

In this section cases
\[N+1\le n\le \tbinom{d+N-1}{N-1}+1\]
will be explored. The first result provides a way of finding a solution, when a previous case, for syzygy bundles over $\pp^{N-1}$, is solved. In addition, a generalisation to proposition~\ref{casemaximaln} is presented, and as a consequence of these, case ${n=\tbinom{d+N-1}{N-1}+1}$ is solved.

The next lemma will give us a solution to the problem for the lowest values of~$n$ mentioned above, provided that the problem is solved for ${N-1}$. A \emph{well\h{behaved}} family of monomials in the $N$th face of the hypertetrahedron will be chosen and the vertex $X_N$ will be added. Well behaved means simply that this family of monomials in $\kk[X_0,\ldots,X_{N-1}]$  satisfies inequality~(\ref{6.6}).

\begin{lem}\label{part of a face and opposite vertex}
If $N\ge3$, $N+1\le n\le\tbinom{d+N-1}{N-1}+1$, and $\I_{N-1,d,n-1}$ is a family of ${n-1}$ monomials in $\kk\left[X_0,\ldots,X_{N-1}\right]$ of degree~$d$ such that their syzygy bundle is stable, then  ${\I_{N,d,n}:=\I_{N-1,d,n-1}\cup\big\{{X_N}^d\big\}}$ is a family of~$n$ monomials in $\kk\left[X_0,\ldots,X_N\right]$ of degree~$d$ whose associated syzygy bundle is stable.
\end{lem}
\begin{dem}
Note that the ideal generated by~$\I_{N,d,n}$ is primary. Let ${J\subseteq \I_{N,d,n}}$ be a subset with at least two monomials. If ${J\subseteq \I_{N-1,d,n-1}}$, then by hypotheses, inequality~(\ref{6.6}) holds. If not, then ${X_N}^d\in J$, and since~$J$ has at least another monomial, where the variable~$X_N$ does not occur, the greatest common divisor is~$1$, and $d_J=0$, so inequality~(\ref{6.6}) holds, for the sequence $\left(a_{d,j}\right)_{j\ge2}$ is monotonically increasing (see remark~\ref{adj}).
\end{dem}

A direct application of this lemma which will become handy for proofs to follow is to take the whole $N$th face and add ${X_N}^d$ to obtain a well\h{behaved} family again. To this end, let us state the case of the highest possible~$n$.

\begin{p}\label{casemaximalngeneral}
For any ${N\ge2}$, the syzygy bundle
\[\syz\left(\big\{{X_0}^{i_0}\cdots{X_N}^{i_N}:
    i_0+\cdots+i_N=d\big\}\right)\]
is stable on~$\pp^N$.
\end{p}
\begin{dem}
Let $\I_{N,d,n}:=\big\{{X_0}^{i_0}\cdots{X_N}^{i_N}: i_0+\cdots+i_N=d\big\}$. If~$g$ is the greatest common divisor of monomials in a subset $J\subseteq \I_{N,d,n}$, all monomials in~$J$ are of the form $gh$, with~$h$ a monomial of degree ${d-d_J}$, where~$d_J~$is the degree of~$g$. There are $\tbinom{N+d-d_J}{d-d_J}$ monomials of degree ${d-d_J}$, so
\[k=|J|\le\tbinom{N+d-d_J}{d-d_J}.\]
Now
\begin{align*}
(d-d_J)n+d_J-dk&=(d-d_J)\tbinom{d+N}{N}+d_J
    -d\tbinom{d-d_J+N}{N}>0\mbox{,}
\end{align*}
which can be proved by induction on~$N$. Therefore inequality~(\ref{6.6}) holds.
\end{dem}

Now using lemma~\ref{part of a face and opposite vertex} and this proposition, we get the following immediate consequence.

\begin{p}\label{whole face and opposite vertex}
For any ${N\ge3}$, the syzygy bundle associated to the family
\[\big\{{X_0}^{i_0}\cdots{X_{N-1}}^{i_{N-1}}:
    i_0+\cdots+i_{N-1}=d\big\}\cup\big\{{X_N}^d\big\}\]
is stable on~$\pp^N$.
\end{p}

\section{Faces of the hypertetrahedron}\label{t3faces}

In this section cases
\[\tbinom{d+N-1}{N-1}+1<n\le\tbinom{d+N}{N}-\tbinom{d-1}{N}\mbox{,}\]
are solved. These cases concern all faces of the hypertetrahedron but the $N$th, which was dealt with in the previous section.

Next proposition starts with the $N$th face and the vertex $X_N$. To these, a number of monomials belonging to face ${N-1}$ will be added, starting with the ones which are closer to the vertex $X_N$. When this face is complete, monomials in face ${N-2}$ are added, an so forth. Note that in case ${d\le N+1}$, this stops at face ${N+1-d}$, since the union of faces ${N+1-d}$ to~$N$ is already the whole hypertetrahedron. In particular, in case ${d=1}$, this process does not even begin. Note that this case is just the well\h{known} bundle of syzygies of ${N+1}$ linear forms over~$\pp^N$, which is nothing but the cotangent bundle~$\Omega^1_{\pp^N}$, which is well known to be stable.

With this method, we will see that we can give a positive answer to problem~\ref{6.9}, in cases ${N\ge3}$ and ${\tbinom{d+N-1}{N-1}+1<n\le\tbinom{d+N}{N}-\tbinom{d-1}{N}}$.

\begin{p}\label{faces}
Let ${N\ge3}$, ${d\ge2}$ and ${\tbinom{d+N-1}{N-1}+1<n\le\tbinom{d+N}{N}-\tbinom{d-1}{N}}$. Then there is a family of~$n$ monomials in $K\left[X_0,\ldots,X_N\right]$ of degree~$d$ whose associated syzygy bundle is stable.
\end{p}
\begin{dem}
Let ${1\leq r\leq\min(d-1,N)}$ and ${0\leq l\leq d-r-1}$ be such that
\[\tbinom{d+N}{N}-\tbinom{d-r+N}{N}+\tbinom{l+N-1}{N-1}<n
    \leq\tbinom{d+N}{N}-\tbinom{d-r+N}{N}+\tbinom{l+N}{N-1}\mbox{,}\]
and let
\[I'_r:=\big\{{X_0}^{i_0}\cdots{X_N}^{i_N}:
        i_0+\cdots+i_N=d\mbox{ and }
        i_{N-r+1}\cdots i_N=0\big\}.\]
This set has all monomials in faces  ${N-r+1}$ to~$N$, and its cardinality is
\[\left|I'_r\right|=\tbinom{d+N}{N}-\tbinom{d-r+N}{N}\]
(we are counting all possible monomials of degree~$d$, and subtracting the ones of type ${X_{N-r+1}\cdots X_Nf}$, where~$f$ is a monomial of degree ${d-r}$).

Now let
\begin{multline*}
I''_{r,l}:=\big\{{X_0}^{i_0}\cdots{X_N}^{i_N}:
        i_0+\cdots+i_N=d,\,i_{N-r+1}\cdots i_N\ne0,\\
    i_{N-r}=0\mbox{ and }i_N\geq d-r-l+1\big\}.
\end{multline*}
This is the set of the monomials in face ${N-r}$ with degree in $X_N$ greater than ${d-r-l}$ that do not belong to $I'_r$, i.e.\ the ones of type
\[{X_{N-r+1}\cdots X_{N-1}{X_N}^{d-r-l+1}f}\mbox{,}\]
where~$f$ is a monomial of degree ${l+1}$ where the variable $X_{N-r}$ does not occur. Therefore the cardinality of~$I''_{r,l}$ is
\[\left|I''_{r,l}\right|=\tbinom{l+N-1}{N-1}.\]
Let ${1\le i \le\tbinom{l+N-1}{N-2}}$ be such that
\[n=\tbinom{d+N}{N}-\tbinom{d-r+N}{N}+\tbinom{l+N-1}{N-1}+i\mbox{,}\]
and let~$I'''_{r,l}$ be a set of~$i$ monomials of degree~$d$ in $K\left[X_0,\ldots,X_N\right]$ of the form
\begin{equation}\label{monimialform}
{X_{N-r+1}\cdots X_{N-1}{X_N}^{d-r-l}f}\mbox{,}
\end{equation}
where~$f$ is a monomial of degree~${l+1}$, where variables $X_{N-r}$ and $X_N$ do not occur. Let us choose these monomials in such a way that the degree of~$X_0$ in~$f$ is the highest possible, i.e.\ such that if $g$ is a monomial as (\ref{monimialform}) that has higher degree in $X_0$ than another monomial in~$I'''_{r,l}$, then ${g\in I'''_{r,l}}$. Let ${\I_{N,d,n}:=I'_r\cup I''_{r,l}\cup I'''_{r,l}}$. Since~${I'_r\cup I''_{r,l}}$ is a set for which the claim in strategy~\ref{X_0} is true, the way we choose monomials for~$I'''_{r,l}$ guarantees that strategy~\ref{X_0} can be applied to~$\I_{N,d,n}$.

As always, it is enough to verify inequality~(\ref{6.6}) for ${0<d_J<d}$. We shall see the cases ${0<d_J\le l}$, ${d_J=l+1}$ and ${l+1<d_J<d}$ separately.

\bigskip
\noindent{\sc Case 1:} ${d_J\leq l}$.\\
In this case, if~$k$ is the number of multiples of ${X_0}^{d_J}$ in~$\I_{N,d,n}$, we have
\[k=\tbinom{d-d_J+N}{N}-\tbinom{d-d_J-r+N}{N}
    +\tbinom{l-d_J+N-1}{N-1}+\min\left[i,\tbinom{l-d_J+N-1}{N-2}\right].\]
Therefore
\begin{align*}
(d-d_J)n+d_J-dk&=(d-d_J)\left[\tbinom{d+N}{N}
        -\tbinom{d-r+N}{N}+\tbinom{l+N-1}{N-1}+i\right]+d_J\\
    &\qquad -d\left[\tbinom{d-d_J+N}{N}
        -\tbinom{d-d_J-r+N}{N}+\tbinom{l-d_J+N-1}{N-1}\right.\\
    &\qquad\qquad\quad\left.+\min\left[i,\tbinom{l-d_J+N-1}{N-2}
        \right]\right]\\
    &=(d-d_J)\left[\tbinom{d+N}{N}
        -\tbinom{d-r+N}{N}+\tbinom{l+N-1}{N-1}\right]+d_J\\
    &\qquad -d\left[\tbinom{d-d_J+N}{N}
        -\tbinom{d-d_J-r+N}{N}
        +\tbinom{l-d_J+N-1}{N-1}\right]\\
    &\qquad +(d-d_J)\left[i-\min\left[i,
        \tbinom{l-d_J+N-1}{N-2}\right]\right]\\
    &\qquad -d_J\min\left(i,\tbinom{l-d_J+N-1}{N-2}\right)\\
    &=(d-d_J)\left[\tbinom{d+N}{N}
        -\tbinom{d-r+N}{N}+\tbinom{l+N-1}{N-1}\right]+d_J\\
    &\qquad -d\left[\tbinom{d-d_J+N}{N}
        -\tbinom{d-d_J-r+N}{N}
        +\tbinom{l-d_J+N-1}{N-1}\right]\\
    &\qquad -d_J\tbinom{l-d_J+N-1}{N-2}.
\end{align*}
Concluding this case amounts to showing that this last expression is positive. To do this, let us give it a name. Let
\begin{align*}
T(N,d,d_J,r,l)&:=(d-d_J)\left[\tbinom{d+N}{N}
        -\tbinom{d-r+N}{N}+\tbinom{l+N-1}{N-1}\right]+d_J\\
    &\qquad -d\left[\tbinom{d-d_J+N}{N}
        -\tbinom{d-d_J-r+N}{N}
        +\tbinom{l-d_J+N-1}{N-1}\right]\\
    &\qquad -d_J\tbinom{l-d_J+N-1}{N-2}.
\end{align*}

Let us start by showing that~$T$ increases with~$r$:
\begin{align*}
T(&N,d,d_J,r+1,l)-T(N,d,d_J,r,l)=\\
    &=-(d-d_J)\left[\tbinom{d-r-1+N}{N}-\tbinom{d-r+N}{N}\right]
        +d\left[\tbinom{d-d_J-r-1+N}{N}
        -\tbinom{d-d_J-r+N}{N}\right]\displaybreak[0]\\
    &=(d-d_J)\tbinom{d-r+N-1}{N-1}-d\tbinom{d-d_J-r+N-1}{N-1}
\end{align*}
Note that if ${r>d-d_J}$, then \[{\tbinom{d-d_J-r+N-1}{N-1}=0\mbox{,}}\]
in which case this expression is clearly positive. Otherwise, we get
\begin{align*}
T(&N,d,d_J,r+1,l)-T(N,d,d_J,r,l)=\\
    &=\tfrac{1}{(N-1)!}\left[(d-d_J)\prod_{s=1}^{N-1}(d-r+s)
        -d\prod_{s=1}^{N-1}(d-d_J-r+s)\right]\displaybreak[0]\\
    &=\tfrac{1}{(N-1)!}\left[(d-d_J)(d-r+1)
        \prod_{s=2}^{N-1}(d-r+s)\right.\\
    &\qquad\qquad\qquad \left.-d(d-d_J-r+1)
        \prod_{s=2}^{N-1}(d-d_J-r+s)\right].
\end{align*}
This last expression is non\h{negative}, since
\[(d-d_J)(d-r+1)-d(d-d_J-r+1)=(r-1)d_J\ge0\mbox{,}\]
and ${d-r+s>d-d_J-r+s}$.

We can therefore look at the case ${r=1}$, since if~$T$ is positive in this case, it will always be positive.

Let us see now that~$T(N,d,d_J,1,l)$ increases with~$l$. Suppose ${d_J\le l\le d-3}$. We get
\begin{align*}
T(&N,d,d_J,1,l+1)-T(N,d,d_J,1,l)=\\
    &=(d-d_J)\left[\tbinom{l+N}{N-1}-\tbinom{l+N-1}{N-1}\right]
        -d\left[\tbinom{l-d_J+N}{N-1}
        -\tbinom{l-d_J+N-1}{N-1}\right]\\
    &\qquad -d_J\left[\tbinom{l-d_J+N}{N-2}
        -\tbinom{l-d_J+N-1}{N-2}\right]\\
    &=(d-d_J)\tbinom{l+N-1}{N-2}-d\tbinom{l-d_J+N-1}{N-2}
        -d_J\tbinom{l-d_J+N-1}{N-3}\\
    &=\tfrac{1}{(N-2)!}\left[(d-d_J)\prod_{s=2}^{N-1}(l+s)
        -d\prod_{s=2}^{N-1}(l-d_J+s)\right.\\
    &\qquad\qquad\qquad \left.-d_J(N-2)
        \prod_{s=3}^{N-1}(l-d_J+s)\right]\\
    &=\tfrac{1}{(N-2)!}\left[(d-d_J)\prod_{s=2}^{N-1}(l+s)\right.\\
    &\qquad\qquad\qquad \left.-\big[d(l-d_J+2)+d_J(N-2)\big]
        \prod_{s=2}^{N-2}(l-d_J+s+1)\right].
\end{align*}
This last expression is never negative, since for ${N=3}$ we have
\begin{align*}
T(&3,d,d_J,1,l+1)-T(3,d,d_J,1,l)=\\
    &=(d-d_J)\prod_{s=2}^{2}(l+s)-\big[d(l-d_J+2)+d_J\big]
        \prod_{s=2}^{1}(l-d_J+s+1)\\
    &=(d-l-3)d_J\ge0\mbox{,}
\end{align*}
and for ${N\ge4}$ we can write
\begin{align*}
T&(N,d,d_J,1,l+1)-T(N,d,d_J,1,l)=\\
    &=\tfrac{1}{(N-2)!}\left[(d-d_J)(l+2)(l+N-1)
        \prod_{s=3}^{N-2}(l+s)\right.\\
    &\qquad\qquad\quad \left.{}-\big[d(l-d_J+2)+d_J(N-2)\big]
        (l-d_J+3)\prod_{s=3}^{N-2}(l-d_J+s+1)\right].
\end{align*}
This is not negative, since
\begin{align*}
(d-d_J)&(l+N-1)(l+2)-\big[d(l-d_J+2)+d_J(N-2)\big](l-d_J+3)=\\
    &=(d-l-3)l(N-2)+(l-d_J)^2(N-3)+2(d-l)(N-3)\\
    &\qquad +5(l-d_J)(N-3)+2(d-l-1)(l-d_J)(d_J-1)\\
    &\qquad +(d-l-3){d_J}^2+l(l-d_J)(d_J-1)+3(d-l-2)(d_J-1)\\
    &\qquad +(l-d_J)(d_J-1)+(d-l-3)+3(d_J-1)(d_J+1)\mbox{,}
\end{align*}
which is non\h{negative}, and for ${3\le s\le N-2}$
\[l+s\ge l-d_J+s+1.\]

Therefore we can look at the case ${l=d_J}$, for if $T$ is positive in this case, it will always be positive.

We shall look at two cases separately now: ${d>2d_J}$ and ${d\le2d_J}$. In the former, we get
\begin{align*}
T(&N,d,d_J,1,d_J)=(d-d_J)\left[\tbinom{d+N-1}{N-1}
        +\tbinom{d_J+N-1}{N-1}\right]-d\tbinom{d-d_J+N-1}{N-1}\\
    &\qquad -d-d_J(N-2)=\\
    &=\tfrac{1}{(N-1)!}\left[(d-d_J)\prod_{s=1}^{N-1}(d+s)
        -d\prod_{s=1}^{N-1}(d-d_J+s)\right.\\
    &\qquad\qquad\qquad\left.{}+(d-2d_J)\left[
        \prod_{s=1}^{N-1}(d_J+s)-(N-1)!\right]\right.\\
    &\qquad\qquad\qquad\left.{}+d_J\left[\prod_{s=1}^{N-1}(d_J+s)
        -N!\right]\right]\displaybreak[0]\\
    &=\tfrac{1}{(N-1)!}\left[(d-d_J)(d+1)(d+2)
        \prod_{s=3}^{N-1}(d+s)\right.\\
    &\qquad\qquad\qquad\left.{}-d(d-d_J+1)(d-d_J+2)
        \prod_{s=3}^{N-1}(d-d_J+s)\right.\\
    &\qquad\qquad\qquad\left.{}+(d-2d_J)\left[
        \prod_{s=1}^{N-1}(d_J+s)-(N-1)!\right]\right.\\
    &\qquad\qquad\qquad\left.{}+d_J\left[\prod_{s=1}^{N-1}(d_J+s)
        -N!\right]\right].
\end{align*}
Now
\begin{multline*}
(d-d_J)(d+1)(d+2)-d(d-d_J+1)(d-d_J+2)=\\
    =(d-2d_J)dd_J+(d-2d_J){d_J}^2+2d_J\left({d_J}^2-1\right),
\end{multline*}
which is never negative, and for ${2\le s\le N-1}$ we have ${d+s>d-d_J+s}$. Furthermore
\begin{align*}
\prod_{s=1}^{N-1}(d_J+s)-N!&\ge\prod_{s=1}^{N-1}(1+s)-N!=0,
\end{align*}
and the term
\[(d-2d_J)\left[\prod_{s=1}^{N-1}(d_J+s)-(N-1)!\right]\]
is strictly positive. Therefore~$T$ is positive in this case.

In case ${d\le2d_J}$, since ${1\le d_J=l\le d-2}$, we get ${d\ge3}$, and therefore ${d_J\ge2}$. We can write
\begin{align*}
T(&N,d,d_J,1,d_J)=(d-d_J)\left[\tbinom{d+N-1}{N-1}
        +\tbinom{d_J+N-1}{N-1}\right]-d\tbinom{d-d_J+N-1}{N-1}\\
    &\qquad -d-d_J(N-2)\\
    &=\tfrac{1}{(N-1)!}\left[(d-d_J)\prod_{s=1}^{N-1}(d+s)
        -d\prod_{s=1}^{N-1}(d-d_J+s)\right.\\
    &\qquad\qquad\qquad\left.{}+(d-d_J)\prod_{s=1}^{N-1}(d_J+s)
        -d_JN!\right]+2d_J-d\displaybreak[0]\\
    &=\tfrac{1}{(N-1)!}\left[(d-d_J)(d+1)(d+2)
        \prod_{s=3}^{N-1}(d+s)\right.\\
    &\qquad\qquad\qquad\left.{}-d(d-d_J+1)(d-d_J+2)
        \prod_{s=3}^{N-1}(d-d_J+s)\right.\\
    &\qquad\qquad\qquad\left.{}+(d-d_J)(d_J+1)
        \prod_{s=3}^N(d_J-1+s)
        -2d_J\prod_{s=3}^Ns\right]+2d_J-d.
\end{align*}
Now if we observe that
\begin{multline*}
(d-d_J)(d+1)(d+2)-d(d-d_J+1)(d-d_J+2)=\\
    =d(d-d_J-1)d_J+(d-2)d_J>0\mbox{,}
\end{multline*}
we see that the difference
\[(d-d_J)(d+1)(d+2)\prod_{s=3}^{N-1}(d+s)-d(d-d_J+1)(d-d_J+2)
    \prod_{s=3}^{N-1}(d-d_J+s)\]
is strictly positive. Finally, since ${d-d_J\ge2}$, we get
\[(d-d_J)(d_J+1)-2d_J\ge2(d_J+1)-2d_J>0\mbox{,}\]
the difference
\[(d-d_J)(d_J+1)\prod_{s=3}^N(d_J-1+s)-2d_J\prod_{s=3}^Ns\]
is positive.

Therefore $T$ is positive also in this case, and hence it is always positive.

\bigskip
\noindent{\sc Case 2:} ${d_J=l+1}$.\\
To count the multiples of ${X_0}^{d_J}$, i.e.\ ${X_0}^{l+1}$, in $I'$, we count all possible multiples of ${X_0}^{l+1}$ and subtract the multiples of the form
\[{X_0}^{l+1}fX_{N-r+1}\cdots X_N\mbox{,}\]
where~$f$ is a monomial of degree ${d-l-1-r}$, which gives us
\[\tbinom{d-l-1+N}{N}-\tbinom{d-l-1-r+N}{N}.\]
In $I''$ there is only one multiple of ${X_0}^{d_J}$, namely
\[{X_0}^{l+1}X_{N-r+1}\cdots X_{N-1}{X_N}^{d-r-l}.\]
Therefore if~$k$ is the number of multiples of ${X_0}^{d_J}$ in~$\I_{N,d,n}$, we have
\[k:=\tbinom{d-l-1+N}{N}-\tbinom{d-l-1-r+N}{N}+1.\]
we get
\begin{align*}
(d-d_J)n+d_J-dk&=(d-l-1)\left[\tbinom{d+N}{N}
        -\tbinom{d-r+N}{N}+\tbinom{l+N-1}{N-1}+i\right]\\
    &\qquad +l+1-d\left[\tbinom{d-l-1+N}{N}
        -\tbinom{d-l-1-r+N}{N}+1\right]\\
    &\ge(d-l-1)\left[\tbinom{d+N}{N}
        -\tbinom{d-r+N}{N}+\tbinom{l+N-1}{N-1}+1\right]\\
    &\qquad +l+1-d\left[\tbinom{d-l-1+N}{N}
        -\tbinom{d-l-1-r+N}{N}+1\right]\\
    &=(d-l-1)\left[\tbinom{d+N}{N}
        -\tbinom{d-r+N}{N}+\tbinom{l+N-1}{N-1}\right]\\
    &\qquad -d\left[\tbinom{d-l-1+N}{N}
        -\tbinom{d-l-1-r+N}{N}\right]
\end{align*}
Again all we have to do is to prove that this expression is positive. Let
\begin{align*}
U(N,d,r,l)&:=(d-l-1)\left[\tbinom{d+N}{N}
        -\tbinom{d-r+N}{N}+\tbinom{l+N-1}{N-1}\right]\\
    &\qquad -d\left[\tbinom{d-l-1+N}{N}
        -\tbinom{d-l-1-r+N}{N}\right]\mbox{,}
\end{align*}
and, as we did in the previous case, let us see that this expression increases with~$r$.
\begin{align*}
U(N,d,r+1,l)-U(&N,d,r,l)=\\
    &=-(d-l-1)\left[\tbinom{d-r-1+N}{N}
        -\tbinom{d-r+N}{N}\right]\\
    &\qquad +d\left[\tbinom{d-l-2-r+N}{N}
        -\tbinom{d-l-1-r+N}{N}\right]\\
    &=(d-l-1)\tbinom{d-r-1+N}{N-1}-d\tbinom{d-l-2-r+N}{N-1}\\
    &=\tfrac{1}{(N-1)!}\left[(d-l-1)(d-r+1)
        \prod_{s=2}^{N-1}(d-r+s)\right.\\
    &\qquad\qquad\qquad\left.{}-d(d-l-r)
        \prod_{s=2}^{N-1}(d-l-1-r+s)\right].
\end{align*}
This is never negative, since
\[(d-l-1)(d-r+1)-d(d-l-r)=(l+1)(r-1)\ge0\mbox{,}\]
and for ${2\le s\le N-1}$, ${d-r+s\ge d-l-1-r+s}$. Therefore, to see that~$U$ is positive, it is enough to check the case ${r=1}$.
\begin{align*}
U(N,d,1,l)&=(d-l-1)\left[\tbinom{d+N}{N}
        -\tbinom{d-1+N}{N}+\tbinom{l+N-1}{N-1}\right]\\
    &\qquad -d\left[\tbinom{d-l-1+N}{N}
        -\tbinom{d-l-2+N}{N}\right]\displaybreak[0]\\
    &=(d-l-1)\left[\tbinom{d+N-1}{N-1}+\tbinom{l+N-1}{N-1}\right]
        -d\tbinom{d-l+N-2}{N-1}\displaybreak[0]\\
    &=\tfrac{1}{(N-1)!}\left[(d-l-1)(d+1)(d+2)
        \prod_{s=3}^{N-1}(d+s)\right.\\
    &\qquad\qquad\qquad\left.{}-d(d-l)(d-l+1)
        \prod_{s=3}^{N-1}(d-l-1+s)\right.\\
    &\qquad\qquad\qquad\left.{}+(d-l-1)\prod_{s=1}^{N-1}(l+s)\right].
\end{align*}
since
\begin{multline*}
(d-l-1)(d+1)(d+2)-d(d-l)(d-l+1)=\\
    =d(d-l-2)l+d(d-l-2)+(d-2)l+(d-2)\ge0\mbox{,}
\end{multline*}
and ${(d-l-1)\prod_{s=1}^{N-1}(l+s)>0}$, we get that~$U$ is positive.

\bigskip
\noindent{\sc Case 3:} ${d_J>l+1}$.\\
In this case, we shall use induction on~$r$. For ${r=1}$, we get that all multiples of ${X_0}^{d_J}$ in~$\I_{N,d,n}$ are in~$I'_1$, and since~${I'_1\cup I''_{1,0}}$ is the set mentioned in proposition~\ref{whole face and opposite vertex}, inequality~(\ref{6.6}) is satisfied for~${I'_1\cup I''_{1,0}}$. The fact that the sequence $\left(a_{d,j}\right)_{j\ge2}$ is monotonically increasing guarantees that inequality~(\ref{6.6}) is also satisfied for~$\I_{N,d,n}$.

For the induction step, suppose that for a given~$r$, inequality~(\ref{6.6}) is satisfied for~$I'_r\cup I''_{r,0}$. Note that if ${r>1}$ this set is just~$\I_{N,d,n}$, in case ${r-1}$ and ${l=d-r}$. Then again the fact that the sequence $\left(a_{d,j}\right)_{j\ge2}$ is monotonically increasing guarantees that inequality~(\ref{6.6}) is also satisfied for~$\I_{N,d,n}$ (see remark~\ref{adj}).
\end{dem}

\section{Hypertetrahedron's interior}\label{t3interior}

In this section a tool to solve cases
\[\tbinom{d+N}{N}-\tbinom{d-1}{N}<n\le\tbinom{d+N}{N}\]
is given. The first result adds up to ${N+1}$ monomials to a hypertetrahedron with all faces completed, and the second result gives us a way of filling the interior of the hypertetrahedron, once we have a suitable family of lower\h{degree} monomials whose syzygy bundle is stable.

Let $\F_{N,d}$ denote the family of monomials
\[\F_{N,d}:=\big\{{X_0}^{i_0}\cdots{X_N}^{i_N}:
    i_0+\cdots+i_N=d\mbox{ and }i_0\cdots i_N=0\big\}\mbox{,}\]
that is, the monomials living on the faces of the hypertetrahedron defined above. Counting all possible monomials of degree~$d$ and subtracting the ones which are multiple of ${X_0\cdots X_N}$, we see that the cardinality of~$\F_{N,d}$ is
\[\tbinom{d+N}{N}-\tbinom{d-1}{N}.\]

When ${d=N+1}$, the proposition above leaves out only one monomial, namely ${X_0\cdots X_N}$. Proposition~\ref{casemaximalngeneral} shows that for ${n=\tbinom{d+N}{N}=\tbinom{2N+1}{N}}$, the whole hypertetrahedron is a family whose associated syzygy bundle is stable. For ${d>N+1}$, the next proposition starts with $\F_{N,d}$ and adds up to ${N+1}$ monomials.

\begin{p}\label{facesanddots}
Let ${N\ge3}$, ${d>N+1}$ and \[\tbinom{d+N}{N}-\tbinom{d-1}{N}<n\le
    \tbinom{d+N}{N}-\tbinom{d-1}{N}+N+1.\]
Then there is a family of~$n$ monomials in $K\left[X_0,\ldots,X_N\right]$ of degree~$d$ whose associated syzygy bundle is stable.
\end{p}
\begin{dem}
Let ${1\le i\le N+1}$ be such that
\[n=\tbinom{d+N}{N}-\tbinom{d-1}{N}+i\mbox{,}\]
and let~$I'$ be the set of the first~$i$ monomials in the sequence
\[\big({X_0}^{d-N}X_1\cdots X_N,\,X_0{X_1}^{d-N}X_2\cdots X_N,
    \ldots,X_0\cdots X_{N-1}{X_N}^{d-N}\big).\]
Let ${\I_{N,d,n}:=\F_{N,d}\cup I'}$. Let us check that inequality~(\ref{6.6}) holds for ${0<d_J<d}$. Since~$\I_{N,d,n}$ satisfies strategy~\ref{X_0}, we can look only at multiples of~${{X_0}^{d_J}}$.

If $d-N<d_J<d$, all multiples of~${{X_0}^{d_J}}$ belong to~$\F_{N,d}$, and since, by the previous proposition, inequality~(\ref{6.6}) holds for~$\F_{N,d}$, the fact that the sequence $\left(a_{d,j}\right)_{j\ge2}$ is monotonically increasing guarantees that inequality~(\ref{6.6}) holds for~$\I_{N,d,n}$.

If $1<d_J\le d-N$ the number of multiples of~${{X_0}^{d_J}}$ is
\[k:=\tbinom{d-d_J+N}{N}-\tbinom{d-d_J}{N}+1\mbox{,}\]
and we get
\begin{align*}
(d-d_J)n+d_J-dk&=(d-d_J)\left[\tbinom{d+N}{N}-\tbinom{d-1}{N}+i
        \right]\\
    &\qquad +d_J-d\left[\tbinom{d-d_J+N}{N}
        -\tbinom{d-d_J}{N}+1\right]\\
    &\ge(d-d_J)\left[\tbinom{d+N}{N}-\tbinom{d-1}{N}+1\right]\\
    &\qquad +d_J-d\left[\tbinom{d-d_J+N}{N}
        -\tbinom{d-d_J}{N}+1\right]\\
    &=(d-d_J)\left[\tbinom{d+N}{N}-\tbinom{d-1}{N}\right]
        -d\left[\tbinom{d-d_J+N}{N}
        -\tbinom{d-d_J}{N}\right].
\end{align*}
Let
\[V(d,d_J,N):=(d-d_J)\left[\tbinom{d+N}{N}-\tbinom{d-1}{N}\right]
    -d\left[\tbinom{d-d_J+N}{N}-\tbinom{d-d_J}{N}\right].\]
If we look at~$V$ as a function on~$d_J$, its second derivative is
\[-\tfrac{2d}{N!}\cdot\sum_{1\le s<t\le N}\
    \left[\sideset{}{_{t\neq r\neq s}}\prod_{r=1}^N(d-d_J+r)
    -\sideset{}{_{t\neq r\neq s}}\prod_{r=1}^N
    (d-d_J-N+r)\right].\]
Since this is negative for ${d_J\in[1,d-N]}$, the function's minimum in this interval is at one of its ends. (Note that we are dealing with the case ${1<d_J\le d-N}$, and therefore the lowest value for $d_J$ is~$2$, and not~$1$, but for the sake of simplicity in calculations, we can look at ${1\le d_J\le d-N}$.) For ${d_J=1}$, we get
\begin{multline*}
V(d,1,N)=(d-1)\left[\tbinom{d+N}{N}-\tbinom{d-1}{N}\right]
        -d\left[\tbinom{d-1+N}{N}-\tbinom{d-1}{N}\right]\\
    =\tfrac{1}{N!}\left[(d-1)(d+1)(d+2)\prod_{s=2}^N(d+s)
        -d^2(d+1)\prod_{s=2}^N(d-1+s)\right]+\tbinom{d-1}{N}.
\end{multline*}
This is easily seen to be positive, since
\[(d-1)(d+1)(d+2)-d^2(d+1)=(d-2)(d+1)>0.\]
For ${d_J=d-N}$, we get
\begin{align*}
V(d,\,d-N,\,N)&=N\left[\tbinom{d+N}{N}-\tbinom{d-1}{N}\right]
        -d\left[\tbinom{2N}{N}-\tbinom{N}{N}\right]\\
    &=\tfrac{1}{N!}\left[N(d+1)(d+2)\prod_{s=3}^N(d+s)
        \right.\\
    &\qquad \left.{}-N(d-N)(d-N+1)\prod_{s=3}^N(d-N-1+s)
        \right.\\
    &\qquad \left.{}-d(N+1)(N+2)\prod_{s=3}^N(N+s)\right]+d.
\end{align*}
This is easily seen to be positive, since
\begin{multline*}
N(d+1)(d+2)-N(d-N)(d-N+1)-d(N+1)(N+2)=\\
    =(d-N)N(N-2)+(d-N)(N-2)>0,
\end{multline*}
which means that
\[N(d+1)(d+2)>N(d-N)(d-N+1)+d(N+1)(N+2)\mbox{,}\]
and for $3\le s\le N$, ${d+s>d-N-1+s}$ and ${d+s>N+s}$.

Finally, if $d_J=1$, the number of multiples of~${{X_0}^{d_J}}$ is
\[k:=\tbinom{d-1+N}{N}-\tbinom{d-1}{N}+i\mbox{,}\]
and we get
\begin{align*}
(d-d_J)n+d_J-dk&=(d-1)\left[\tbinom{d+N}{N}-\tbinom{d-1}{N}+i
        \right]\\
    &\qquad +1-d\left[\tbinom{d-1+N}{N}
        -\tbinom{d-1}{N}+i\right]\\
    &=(d-1)\left[\tbinom{d+N}{N}-\tbinom{d-1}{N}\right]\\
    &\qquad -d\left[\tbinom{d-1+N}{N}
        -\tbinom{d-1}{N}\right]+1-i\\
    &\ge(d-1)\left[\tbinom{d+N}{N}-\tbinom{d-1}{N}\right]
        -d\left[\tbinom{d-1+N}{N}-\tbinom{d-1}{N}\right]-N\\
    &=d\tbinom{d+N-1}{N-1}-\tbinom{d+N}{N}+\tbinom{d-1}{N}-N\\
    &=\left(d-\tfrac{d+N}{N}\right)\tbinom{d+N-1}{N-1}
        +\tbinom{d-1}{N}-N.
\end{align*}
This is positive, since
\[d-\tfrac{d+N}{N}=\tfrac{(d-2)(N-1)+N-2}{N}>0\mbox{,}\]
and
\[\tbinom{d-1}{N}>\tbinom{N+1}{N}>N.\]

Therefore inequality~(\ref{6.6}) strictly holds in all cases, so the syzygy bundle associated to~$\I_{N,d,n}$ is stable.
\end{dem}

\bigskip

For some computations in the next proof, we will need a result which is a simple consequence of the fact that for any numbers $a,b_1,\ldots,b_n$,
\[\prod_{s=1}^p(a+b_s)=\sum_{s=0}^p\
    \sum_{1\le t_1<\cdots<t_s\le p}
    a^{p-s}b_{t_1}\cdots b_{t_1}.\]
Note that if ${N\ge1}$ and ${d\ge0}$, we get
\[\tbinom{d+N}{N}-\tbinom{d-1}{N}=\tfrac{N+1}{(N-1)!}d^{N-1}+
    \mbox{ some positive terms}.\]
\begin{lem}\label{Brenner2}
Let ${N\ge1}$ and ${d\ge0}$. Then
\[\tbinom{d+N}{N}-\tbinom{d-1}{N}\ge\tfrac{N+1}{(N-1)!}d^{N-1}.\]
\end{lem}

The following lemma will be the key to prove the main result of this chapter. As mentioned in the beginning of this section, it is no more than a generalization to any~${N\ge3}$ of a lemma due to Brenner \cite{Bre}.

\begin{lem}\label{Brenner}
Let ${N\ge3}$, ${d>N+1}$, and ${\tbinom{d+N}{N}-\tbinom{d-1}{N}<n\le\tbinom{d+N}{N}}$. If~$I'$ is a family of
\[{n':=n-\left[\tbinom{d+N}{N}-\tbinom{d-1}{N}\right]}\]
monomials in $K\left[X_0,\ldots,X_N\right]$ of degree~${d':=d-N-1}$ such that their syzygy bundle is stable, then
\[\I_{N,d,n}:=\F_{N,d}\cup\big\{X_0\cdots X_Nf:f\in I'\big\}\]
is a family of $n$ monomials in $K\left[X_0,\ldots,X_N\right]$ of degree~$d$ whose associated syzygy bundle is stable.
\end{lem}
\begin{dem}
Let ${J\subset \I_{N,d,n}}$, with ${k:=|J|\ge2}$, and let~$d_J$ be the degree of the greatest common divisor of monomials in~$J$. To verify that~$\I_{N,d,n}$ satisfies inequality~(\ref{6.6}), we can assume ${0<d_J<d}$, since the fact that~$J$ has at least two elements makes ${d_J\ne d}$, and for ${d_J=0}$, the fact that the sequence $\left(a_{d,j}\right)_{j\ge2}$ is monotonically increasing is enough. We may also assume that~$J$ has all multiples of its greatest common divisor, since if a degree~$d_J$ is fixed, the higher~$k$ is, the harder it is to guarantee inequality~(\ref{6.6}).

Suppose~$J$ intersects exactly~$i$ faces of~$\F_{N,d}$, with ${0\le i\le N}$. Then exactly ${\tbinom{d-d_J+N}{N}-\tbinom{d-d_J+N-i}{N}}$ of its monomials are in~$\F_{N,d}$, and the remaining admit a greatest common divisor of degree ${d_J+i}$, and come from a subset~$J'\subseteq I'$, admitting a greatest common divisor of degree~${d_{J'}\ge d_J-N-1+i}$. Since ${N+1-i}$ faces of~$\F_{N,d}$ do not intersect~$J$, its greatest common divisor is multiple of the variables that are not present in those faces, and consequently ${d_J\ge N+1-i}$. Let
\[k':=\left|J'\right|=k-
    \left[\tbinom{d-d_J+N}{N}-\tbinom{d-d_J+N-i}{N}\right].\]
Observe that
\[\big[d'-(d_J-N-1+i)\big]n'+(d_J-N-1+i)-d'k'\ge
    \left(d'-d_{J'}\right)n'+d_{J'}-d'k'\mbox{,}\]
and the last expression is strictly positive, since~$I'$ satisfies inequality~(\ref{6.6}). We can see that
\begin{align*}
(d-d_J)n+d_J-dk&=\big[d'-(d_J-N-1+i)\big]n'+(d_J-N-1+i)-d'k'\\
    &\qquad +P(n',k',N,d,d_J,i)+Q(N,d,d_J,i)\mbox{,}
\end{align*}
where
\begin{align*}
P(n',k',N,d,d_J,i)&=i\left(n'-k'\right)+
        (N+1-i)\left[\tbinom{d-d_J+N-i}{N}-k'+1\right]\\
\intertext{and}
Q(N,d,d_J,i)&=(d-d_J)\left[\tbinom{d+N}{N}
        -\tbinom{d-1}{N}\right]-d\tbinom{d-d_J+N}{N}\\
    &\qquad +(d-N-1+i)\tbinom{d-d_J+N-i}{N}.
\end{align*}
Now~$P$ is clearly positive, since $\tbinom{d-d_J+N-i}{N}$ is the highest possible cardinality for~$J'$. If we can guarantee that~$Q$ is non\h{negative}, inequality~(\ref{6.6}) is strictly satisfied.

Suppose ${d_J=1}$. In this case, we have ${i=N}$, since ${d_J\ge N+1-i}$. Therefore we get
\begin{align*}
Q(N,d,1,N)&=(d-1)\left[\tbinom{d+N}{N}-\tbinom{d-1}{N}\right]
        -d\tbinom{d-1+N}{N}+(d-1)\tbinom{d-1}{N}\\
    &=\tfrac{(N-2)d+d-N}{d}\tbinom{d-1+N}{N}\mbox{,}
\end{align*}
and the last expression is positive for ${N\ge3}$.

Now suppose ${d_J\ge2}$. Since the last term in~$Q$ vanishes for ${i>d-d_J}$, we shall consider the cases ${i\le d-d_J}$ and ${i>d-d_J}$ separately.

\bigskip
\noindent{\sc Case 1:} ${i\le d-d_J}$.\\
Since
\[Q(N,d,d_J,i+1)-Q(N,d,d_J,i)=
    -\tfrac{(N-1)(d-N-2)+i(N+1)+(d_J-2)}{d-d_J-i}
    \tbinom{d-d_J+N-i-1}{N}\mbox{,}\]
and this is negative, we know that~$Q$ decreases as~$i$ gets higher. Therefore we should pay attention to its highest value $i=\min(N,d-d_J)$.

Let us start with the case ${d=N+2}$. In this case, ${i=N+2-d_J}$, and
\begin{align*}
Q(N,\,N+2,\,d_J,\,N+2-d_J)&=(N+2-d_J)\left[\tbinom{2N+2}{N}
        -\tbinom{N+1}{N}\right]\\
    &\qquad -(N+2)\tbinom{2N+2-d_J}{N}
        +(N+3-d_J)\tbinom{N}{N}\\
    &=(N+2-d_J)\tbinom{2N+2}{N}-(N+2)\tbinom{2N+2-d_J}{N}\\
    &\qquad -N^2-2N++Nd_J+1.
\end{align*}
If we look at this as a function on~$d_J$, its second derivative is
\[-2\cdot\tfrac{N+2}{N!}\cdot\sum_{1\le s<t\le N}\
    \prod_{\substack{1\le r\le N\\t\neq r\neq s}}
    (N+2-d_J+r).\]
Since this is negative for ${d_J\in[2,N+1]}$, the function's minimum in this interval is at one of its ends. For ${d_J=2}$, we have
\begin{align*}
Q(N,\,N+2,\,2,\,N)&=N\tbinom{2N+2}{N}-(N+2)\tbinom{2N}{N}-N^2+1\\
    &=\tfrac{3N^2-2N-4}{N+2}\tbinom{2N}{N}-N^2+1\mbox{,}
\end{align*}
and this is positive for ${N\ge3}$. For ${d_J=N+1}$, we have
\begin{align*}
Q(N,\,N+2,\,N+1,\,1)&=\tbinom{2N+2}{N}-(N+2)\tbinom{N+1}{N}-N+1\\
    &\ge\tbinom{2N+2}{3}-(N+2)(N+1)-N+1\\
    &=\tfrac{1}{3}\left(4N^3+3N^2-10N-3\right)\mbox{,}
\end{align*}
and this is also positive for ${N\ge3}$.

\bigskip

Now, for general ${d\ge N+2}$, we shall see two subcases separately, namely ${N\le d-d_J}$ and ${N>d-d_J}$.

\bigskip
\noindent{\sc Subcase 1.1:} ${N\le d-d_J}$.\\
Looking at~$Q(N,d,d_J,N)$ as a function on~$d_J$ again, its second derivative is
\[-\tfrac{1}{N!}\cdot\sum_{1\le s<t\le N}
    \left[d\cdot\sideset{}{_{t\neq r\neq s}}\prod_{r=1}^N
    (d-d_J+r)
    -(d-1)\cdot\sideset{}{_{t\neq r\neq s}}\prod_{r=1}^N
    (d-d_J-N+r)\right].\]
This is negative for ${d_J\in[2,d-N]}$, since ${d-d_J+r\ge d-d_J-N+r\ge1}$ for ${1\le r\le N}$. Therefore the function's minimum in this interval is at one of its ends. If ${d_J=2}$, we get
\begin{multline*}
Q(N,d,2,N)=\\
    =\tfrac{2d^2(N-2)+d(N-2)+2d(d-N-2)
        +3(d-N-2)(N-1)+N(N-3)+8(N-1)+2}{d(d-1)}
        \tbinom{d+N-2}{N}\\
    -\tfrac{(d-1)(d-N-2)}{d-N-1}\tbinom{d-2}{N}\ge0.
\end{multline*}
If ${d_J=d-N}$, we get
\[Q(N,\,d,\,d-N,\,N)=N\left[\tbinom{d+N}{N}-\tbinom{d-1}{N}\right]
    -d\tbinom{2N}{N}+d-1.\]
If we look at this expression now as a function on~$d$, its second derivative is
\[\tfrac{2}{(N-1)!}\cdot\sum_{1\le s<t\le N}
    \left[\sideset{}{_{t\neq r\neq s}}\prod_{r=1}^N(d+r)
    -\sideset{}{_{t\neq r\neq s}}\prod_{r=1}^N
    (d-N-1+r)\right]>0\mbox{,}\]
since ${d+r\ge d-N-1+r\ge1}$ for ${1\le r\le N}$. This means its first derivative increases with~$d$. When we evaluate this first derivative at ${d=N+2}$, we get
\begin{equation}\label{Brenner1}
\tfrac{1}{(N-1)!}\sum_{s=1}^N\left[\sideset{}{_{r\neq s}}
    \prod_{r=1}^N(N+2+r)-\sideset{}{_{r\neq s}}
    \prod_{r=1}^N(1+r)\right]-\tbinom{2N}{N}+1.
\end{equation}
Note that, except for the last term in the sum, the last factor of the first product is ${2(N+1)}$, and the last factor of the second is ${N+1}$. Splitting the last term in the sum and factoring out the last factor of both products, this becomes
\begin{multline*}
\tfrac{N+1}{(N-1)!}\sum_{s=1}^{N-1}\left[2\cdot
        \sideset{}{_{r\neq s}}\prod_{r=1}^{N-1}(N+2+r)
        -\sideset{}{_{r\neq s}}\prod_{r=1}^{N-1}(1+r)
        \right]\\
    +\tfrac{1}{(N-1)!}\left[\prod_{r=1}^{N-1}(N+2+r)
        -\prod_{r=1}^{N-1}(1+r)\right]-\tbinom{2N}{N}+1\mbox{,}
\end{multline*}
and this can be rearranged as
\begin{multline*}
\tfrac{N+1}{(N-1)!}\sum_{s=1}^{N-1}
    \left[\sideset{}{_{r\neq s}}\prod_{r=1}^{N-1}(N+2+r)
    -\sideset{}{_{r\neq s}}\prod_{r=1}^{N-1}(1+r)\right]\\
    +\tfrac{1}{(N-1)!}\left[\prod_{r=1}^{N-1}(N+2+r)
    -\prod_{r=1}^{N-1}(1+r)\right]\\
    +\tfrac{1}{N!}\left[\sum_{s=1}^{N-1}(N+1)N
    \sideset{}{_{r\neq s}}\prod_{r=1}^{N-1}(N+2+r)
    -\prod_{r=1}^{N}(N+r)\right]+1.
\end{multline*}
Now the first two terms are clearly positive and for the third one we can see that
\begin{multline*}
(N+1)N\sum_{s=1}^{N-1}\sideset{}{_{r\neq s}}
        \prod_{r=1}^{N-1}(N+2+r)-\prod_{r=1}^{N}(N+r)\\
    \ge(N+1)^2N\prod_{r=1}^{N-2}(N+2+r)-\prod_{r=1}^{N}(N+r)\\
    =(N+1)^2N\prod_{r=1}^{N-2}(N+2+r)-2N(2N-1)
        \prod_{r=1}^{N-2}(N+r)>0\mbox{,}
\end{multline*}
since ${(N+1)^2N>2N(2N-1)>1}$ for ${N\ge3}$, and ${N+2+r>N+r>1}$ for ${1\le r\le N-2}$. Therefore all the expression is positive, and we see that~$Q$ increases with~$d$. Since we have already seen that~$Q$ is positive for the first value ${d=N+2}$, we have that~$Q$ is positive in this case.

\bigskip
\noindent{\sc Subcase 1.2:} ${N>d-d_J}$.\\
Looking once more at~$Q(N,d,d_J,N)$ as a function on~$d_J$, its second derivative is
\[-\tfrac{2d}{N!}\cdot\sum_{1\le s<t\le N}
    \sideset{}{_{t\neq r\neq s}}\prod_{r=1}^N
    (d-d_J+r).\]
This is negative for ${d_J\in[d-N,d-1]}$, since ${d-d_J+r\ge1}$ for ${1\le r\le N}$.\linebreak[4] Therefore the function's minimum in this interval is at one of its ends. If ${d_J=d-N}$, we are exactly in the same situation as before, so we already know that~$Q$ is non\h{negative}. If ${d_J=d-1}$, we get
\[Q(N,\,d,\,d-1,\,N)=\tbinom{d+N}{N}-\tbinom{d-1}{N}-N(d+1).\]
Therefore, if we apply lemma~\ref{Brenner2}, we get
\begin{align*}
Q(N,\,d,\,d-1,\,N)&\geq\tfrac{N+1}{(N-1)!}d^{N-1}-N(d+1)\\
    &\ge\tfrac{(N+1)^{N-1}}{(N-1)!}d-N(d+1)\\
    &>(N+1)d-N(d+1)=d-N>0.
\end{align*}

\bigskip
\noindent{\sc Case 2:} ${i>d-d_J}$.\\
In this case, we get ${d_J>d-N}$, since ${i\le N}$. Let us start again with the case ${d=N+2}$. In this case, ${3\leq d_J\leq N+1}$, and
\[Q(N,\,N+2,\,d_J,\,i)=(N+2-d_J)\left[\tbinom{2N+2}{N}
    -\tbinom{N+1}{N}\right]-(N+2)\tbinom{2N+2-d_J}{N}.\]
If we look at this as a function on~$d_J$, its second derivative is
\[-2\cdot\tfrac{N+2}{N!}\cdot\sum_{1\le s<t\le N}
    \prod_{\substack{1\le r\le N\\t\neq r\neq s}}
    (N+2-d_J+r).\]
Since this is negative for ${d_J\in[3,N+1]}$, the function's minimum in this interval is at one of its ends. For ${d_J=3}$, we have
\begin{align*}
Q(N,\,N+2,\,3,\,i)&=(N-1)\left[\tbinom{2N+2}{N}
        -\tbinom{N+1}{N}\right]-(N+2)\tbinom{2N-1}{N}\\
    &=\tfrac{7N^2-8N-8}{N+2}\tbinom{2N-1}{N}-N^2+1\mbox{,}
\end{align*}
and this is positive for $N\ge3$. For ${d_J=N+1}$, we have
\begin{align*}
Q(N,\,N+2,\,N+1,\,i)&=\tbinom{2N+2}{N}-(N+3)(N+1)\\
    &\ge\tbinom{2N+2}{3}-(N+3)(N+1)\\
    &=\tfrac{1}{3}\left(4N^3+3N^2-4N+9\right)\mbox{,}
\end{align*}
and this is also positive for ${N\ge3}$.

In general, for ${d\ge N+2}$, we get
\[Q(N,d,d_J,i)=(d-d_J)\left[\tbinom{d+N}{N}
    -\tbinom{d-1}{N}\right]-d\tbinom{d-d_J+N}{N}.\]
Looking at this expression as a function on~$d_J$, its second derivative is
\[-\tfrac{2d}{N!}\cdot\sum_{1\le s<t\le N}
    \sideset{}{_{t\neq r\neq s}}\prod_{r=1}^N(d-d_J+r)\le0.\]
Since this is negative for ${d_J\in[d-N,d-1]}$, the function's minimum in this interval is again at one of its ends.

For ${d_J=d-N}$, we have
\[Q(N,\,d,\,d-N,\,i)=N\left[\tbinom{d+N}{N}
    -\tbinom{d-1}{N}\right]-d\tbinom{2N}{N}.\]
If we look at this expression now as a function on~$d$, its second derivative is
\[\tfrac{2}{(N-1)!}\cdot\sum_{1\le s<t\le N}
    \left[\sideset{}{_{t\neq r\neq s}}\prod_{r=1}^N(d+r)
    -\sideset{}{_{t\neq r\neq s}}\prod_{r=1}^N
    (d-N-1+r)\right]>0\mbox{,}\]
since ${d+r\ge d-N-1+r\ge1}$ for ${1\le r\le N}$. This means its first derivative increases with~$d$. When we evaluate this first derivative at ${d=N+2}$, we get
\[\tfrac{1}{(N-1)!}\sum_{s=1}^N\left[\sideset{}{_{r\neq s}}
    \prod_{r=1}^N(N+2+r)-\sideset{}{_{r\neq s}}
    \prod_{r=1}^N(1+r)\right]-\tbinom{2N}{N}.\]
this is the same expression as~(\ref{Brenner1}) above, except for the final term ($+1$). Since this term is ignored in the argument that follows, we arrive to the same conclusions here, and guarantee that~$Q$ is positive.

For ${d_J=d-1}$, we have
\[Q(N,\,d,\,d-1,\,i)=\tbinom{d+N}{N}-\tbinom{d-1}{N}-d(N+1).\]
Again we can apply lemma~\ref{Brenner2}, and get
\begin{align*}
Q(N,\,d,\,d-1,\,i)&\geq\tfrac{N+1}{(N-1)!}d^{N-1}-d(N+1)\\
    &>\tfrac{(N+1)^{N-1}}{(N-1)!}d-d(N+1)\\
    &>(N+1)d-d(N+1)=0.
\end{align*}

\bigskip

We have verified that~$Q$ is positive in all cases, so inequality~(\ref{6.6}) is strictly satisfied, and~$\I_{N,d,n}$ is a family of $n$ monomials whose associated syzygy bundle is stable.
\end{dem}

\section{Main theorem}\label{t3main}

The main theorem of this work is stated in this section.
The fact that for the case ${N=2}$ problem~\ref{6.9} is solved (except for the case ${d=2}$ and ${n=5}$), combined with the results in this chapter, will allow us to assert its main theorem. To get round that exception, we can see a particular case.

\begin{lem}\label{case326}
The syzygy bundle associated to
\[\I_{3,2,6}:=\big\{{X_0}^2,\,{X_1}^2,\,{X_2}^2,\,{X_3}^2,
    \,X_0X_1,\,X_2X_3\big\}\]
is stable.
\end{lem}
\begin{dem}
Note that the ideal generated by~$\I_{3,2,6}$ is primary, and that the relevant sets to verify inequality~(\ref{6.6}) have two elements and a linear greatest common divisor. Therefore
\[(d-d_J)n+d_J-dk=(2-1)\cdot6+1-2\cdot2=3>0\mbox{,}\]
so stability is guaranteed.
\end{dem}

\begin{te}\label{main3}
Let~$N$, $d$ and~$n$ be integers such that ${N\ge2}$, ${(N,d,n)\ne(2,2,5)}$, and ${N+1\le n\le\tbinom{d+N}{N}}$. Then there is a family of~$n$ monomials in $K\left[X_0,\ldots,X_N\right]$ of degree~$d$ such that their syzygy bundle is stable.

For ${(N,d,n)=(2,2,5)}$, there are~$5$ monomials of degree~$2$ in $K\left[X_0,X_1,X_2\right]$ such that their syzygy bundle is semistable.
\end{te}
\begin{dem}
Case ${N=2}$ was already stated in theorem~\ref{main2}.

For ${N\ge3}$, this can be done by induction on~$N$. For ${N=3}$, lemma~\ref{part of a face and opposite vertex} gives us an answer for
\[4\le n\le\tbinom{d+2}{2}+1\mbox{,}\]
except for the case where ${d=2}$ and ${n=6}$, for which we have lemma~\ref{case326}; proposition~\ref{faces} takes care of cases
\[\tbinom{d+2}{2}+1<n\le\tbinom{d+3}{3}-\tbinom{d-1}{3}\mbox{;}\]
proposition~\ref{facesanddots} gives an answer for
\[\tbinom{d+3}{3}-\tbinom{d-1}{3}<n\le\tbinom{d+3}{3}
    -\tbinom{d-1}{3}+4\mbox{;}\]
finally, lemma~\ref{Brenner} takes care of all other cases.

Now if we suppose the answer is positive for some~$N$, lemma~\ref{part of a face and opposite vertex} provides a positive answer for the first cases of ${N+1}$, and proposition~\ref{faces}, proposition~\ref{facesanddots}, and lemma~\ref{Brenner} take care of the rest.
\end{dem}

The theorem above is the main theorem in this thesis. It provides a complete answer to problem~\ref{6.9}, as presented by Brenner (see \cite{Bre08a} and \cite{Bre08b}).

As it generalises the results in the previous chapter, it will allow us to assert that for $N+1\le n\le\tbinom{d+N}{N}$, with $(N,d,n)\ne(2,2,5)$, a syzygy bundle associated to a family of~$n$ general homogeneous polynomials in $K\left[X_0,\ldots,X_N\right]$ of degree~$d$ is stable, for the condition of stability in the moduli space is an open one.

\section{Open questions}\label{t3open}

Once problem~\ref{6.9} is solved, two natural questions arise, both as a possible generalisation. In this section we state them and make some comments on their resolutions.

\begin{pr}\label{openseveraldegrees}
To determine integers $N$, $n$, and ${d_1\le\cdots\le d_n}$, with ${n\ge N+1}$, such that there is a family of~$n$ monomials ${f_1,\ldots,f_n}$ in $\kk\left[X_0,\ldots,X_N\right]$ of degrees ${d_1,\ldots,d_n}$, respectively, whose syzygy bundle is semistable, or stable.
\end{pr}

Brenner's criterion applies to families of monomials in any degree, so it can be of use for this question.

The other possible generalisation is related so syzygies of higher order. Given homogeneous polynomials ${f_1,\ldots,f_n}$, if
\[\xymatrix{0\ar[r]&\mathcal{F}_s\ar[r]^-{\varphi_s}&\cdots
    \ar[r]^-{\varphi_3}&\mathcal{F}_2
    \ar[r]^-{\varphi_2}&{\bigoplus\limits_{i=1}^n\OO(-d_i)}
    \ar[rr]^-{f_1,\ldots,f_n}&&\OO\ar[r]&0}\]
is a minimal free resolution, we say that the kernel of $\varphi_i$ is the \emph{$i$th syzygy sheaf} of ${f_1,\ldots,f_n}$.

\begin{pr}\label{openhighersyzygies}
Given a family ${f_1,\ldots,f_n}$ of homogeneous polynomials in\linebreak[4] ${\kk[X_0,\ldots,X_N]}$, of degrees ${d_1,\ldots,d_n}$, respectively, is the $i$th syzygy sheaf\linebreak[4] $\syz^i(f_1,\ldots,f_n)$ stable (or, at least, semistable)?
\end{pr}

This problem is closely related to the one of finding minimal free resolutions of $n$ general forms, which is far from being solved. In 2003, Migliore and Miró\h{Roig} \cite{MMR03} solved
some instances of the case ${n=N+1}$. In the paper's introduction, we can find an account of contributions made so far. Another recent contribution is a paper by Francisco \cite{Chr04}, where he proves the lex\h{plus}\h{powers} conjecture for almost complete intersections.

\chapter{Moduli spaces of syzygy bundles}
\label{t4}

In this chapter we present three results on syzygy bundles and their moduli space. The first one states that syzygy bundles given by a family of general homogeneous polynomials of a given degree in three or more variables are stable, except in the case of $5$~quadric homogeneous polynomials in three variables, where we can only guarantee semistability. The second gives us the dimension of the stratum of the moduli space corresponding to syzygy bundles. The third asserts that the irreducible component of the moduli space corresponding to the isomorphism class of a stable syzygy bundle over $\pp^N$, with ${N\ge2}$, is generically smooth, and also gives its dimension. As a consequence, we determine whether the closure of the stratum corresponding to syzygy bundles is the whole irreducible component, and in case it is not, we get its codimension.

\bigskip


Since stability of vector bundles is an open property, we get an immediate consequence from theorems \ref{main2} and \ref{main3} in the previous chapters.

\begin{te}\label{geral4}
Let $N\ge2$, $d\ge1$, and $N+1\le n\le\tbinom{d+N}{N}$, with ${(N,d,n)\ne(2,2,5)}$. If ${f_1,\ldots,f_n}$ are general forms of degree~$d$ in $\kk[X_0,\ldots,X_N]$, such that the ideal ${(f_1,\ldots,f_n)}$ is \mbox{$\mathfrak{m}$-primary}, then the syzygy bundle ${\syz(f_1,\ldots,f_n)}$ is stable.

If ${(N,d,n)=(2,2,5)}$ and ${f_1,\ldots,f_5}$ are general quadratic forms such that the ideal ${(f_1,\ldots,f_5)}$ is \mbox{$\mathfrak{m}$-primary}, then the syzygy bundle ${\syz(f_1,\ldots,f_5)}$ is semistable.
\end{te}

We denote by $M(r;c_1,\ldots,c_s)$ the moduli space of rank~$r$ stable vector bundles~$E$ on $\pp^N$ with fixed Chern classes $c_i(E)=c_i$, for ${1\le i\le s}$, where ${s:=min(r,N)}$. The existence of the moduli space $M(r;c_1,\ldots,c_s)$ was shown by Maruyama in 1977 (see \cite{Mar77} and \cite{Mar78}) and once the existence of the moduli space is established, the question arises as what can be said about its structure, both locally and globally. More precisely, what does the moduli space look like as an algebraic variety? Is it, for example, connected, irreducible, rational or smooth? What does it look like as a topological space? What is its geometry?  Until now, there is no general answer to these questions and the aim of this chapter is to study the points which parameterise stable syzygy bundles.

\bigskip

For each ${N\ge2}$, ${d\ge1}$, and ${N+1\le n\le\tbinom{d+N}{N}}$ such that ${(N,d,n)\ne(2,2,5)}$, let $E_{N,d,n}$ be a stable syzygy bundle of~$n$ forms of degree~$d$ in $\kk[X_0,\ldots,X_N]$. Let us denote by ${c_i^{N,d,n}:=c_i(E_{N,d,n})}$, for ${1\le i\le N}$, the $i$th Chern class of $E_{N,d,n}$. As we have seen in chapter~\ref{t1}, from the syzygy bundle's presenting sequence
\begin{equation}\label{presE_{N,d,n}}
\xymatrix{0\ar[r]&E_{N,d,n}\ar[r]&
    {\OO(-d)^n}\ar[r]&\OO\ar[r]&0}
\end{equation}
we get its Chern polynomial
\[c_t\left(E_{N,d,n}\right)=(1-dt)^n\in
    \tfrac{\zz[t]}{(t^{N+1})}.\]
Therefore the Chern classes of $E_{N,d,n}$ satisfy
\[\sum_{i=1}^Nc_i^{N,d,n}t^i=\sum _{i=1}^N\tbinom{n}{i}(-dt)^i\mbox{,}\]
i.e.\ we have, for ${1\le i\le N}$,
\[c_i^{N,d,n}=\tbinom{n}{i}(-d)^i.\]

Let
\[{M_{N,d,n}:=M\left(n-1;c_1^{N,d,n},\ldots,c_N^{N,d,n}\right)}\]
be the moduli space of rank ${n-1}$ stable vector bundles~$E$ on $\pp^N$ with fixed Chern classes ${c_i(E)=c_i^{N,d,n}}$, for ${1\le i\le N}$, and let ${S_{N,d,n}\subseteq M_{N,d,n}}$ be the stratum in $M_{N,d,n}$ consisting of isomorphism classes of rank ${n-1}$ stable syzygy bundles $E_{N,d,n}$ defined by the exact sequence (\ref{presE_{N,d,n}}). We shall call $S_{N,d,n}$ the \emph{syzygy locus} of $M_{N,d,n}$. In this chapter, we address the following problems:
\begin{enumerate}
  \item to determine the dimension of $S_{N,d,n}$ in terms of $N$, $d$ and $n$;
  \item to determine whether the closure of $S_{N,d,n}$ is an irreducible component of $M_{N,d,n}$, or to find its codimension in $M_{N,d,n}$, if that is not the case; and
  \item to determine whether $M_{N,d,n}$ is smooth along $S_{N,d,n}$.
\end{enumerate}

\section{The dimension of the syzygy locus}\label{t4locus}

The goal of this section is to compute the dimension of the stratum $S_{N,d,n}$ in $M_{N,d,n}$ parameterising rank ${n-1}$ stable syzygy bundles $E_{N,d,n}$ fitting into an exact sequence
\[\xymatrix{0\ar[r]&E_{N,d,n}\ar[r]&
    {\OO(-d)^n}\ar[r]&\OO\ar[r]&0\mbox{.}}\]
This is done in the following proposition.

\begin{p}\label{syzlocus} Fix integers $N$, $d$ and $n$ such that ${N\ge2}$, ${d\ge1}$,
\[{N+1\le n\le \tbinom{d+N}{N}}\mbox{,}\]
and ${(N,d,n)\ne (2,2,5)}$. Then
\[\dim S_{N,d,n}=n\tbinom{d+N}{N}-n^2.\]
\end{p}
\begin{dem}
First of all, we dualise the exact sequence
\[\xymatrix{0\ar[r]&E_{N,d,n}\ar[r]&
    {\OO(-d)^n}\ar[r]&\OO\ar[r]&0\mbox{,}}\]
and we get
\[\xymatrix{0\ar[r]&\OO\ar[r]&
    {\OO(d)^n}\ar[r]&{E_{N,d,n}}^\vee\ar[r]&0\mbox{.}}\]
Hence, ${E_{N,d,n}}^\vee$ is the cokernel of a morphism in ${H:=\Hom\big(\OO,\OO(d)^n\big)}$. However, different morphisms ${f,g\in H}$ may define the same cokernel. Indeed, ${\Aut(\OO)\times\Aut(\OO(d)^n)}$ acts on~$H$ as follows:
\[\xymatrix@R=0pt{{\big(\Aut(\OO)\times\Aut(\OO(d)^{n})\big)\times
        H}\ar[r]&H\\
    {\big((\varphi,\psi),f\big)}\ar@{|->}[r]&(\varphi,\psi)
        \circ f:=\psi^{-1}f\varphi}\]
and we get
\[\dim S_{N,d,n}=\dim H-\dim\Aut(\OO)-\dim\Aut(\OO(d)^n)
    +\dim I_f\mbox{,}\]
where
\[{I_f:=\{(\varphi,\psi)\in\Aut(\OO)\times\Aut(\OO(d)^n):
    (\varphi,\psi)\circ f=f\}}\]
is the isotropy group of a general morphism ${f\in H}$. Since
\begin{align*}
\dim\Aut(\OO)&=1\mbox{,}\\
\dim\Aut(\OO(d)^n)&=n^2\mbox{,}\\
\dim H&=n\tbinom{d+N}{N}\mbox{,}\\
\intertext{and}
\dim I_f&=1\mbox{,}
\end{align*}
we get
\[\dim S_{N,d,n}=n\tbinom{d+N}{N}-n^2\mbox{,}\]
as we wanted.
\end{dem}

\section{The moduli space of syzygy bundles}\label{t4moduli}

The goal of this section is to determine the unobstructedness of stable syzygy bundles $E_{N,d,n}$ on $\pp^N$ and whether the closure of the syzygy locus $S_{N,d,n}$ is an irreducible component of $M_{N,d,n}$.

\begin{te}\label{moduli}
Let $N$, $d$ and $n$ be integers such that ${N\ge2}$, ${d\ge1}$,
\[{N+1\le n\le\tbinom{d+N}{N}}\mbox{,}\]
and ${(N,d,n)\ne(2,2,5)}$. Then
\begin{enumerate}
  \item\label{moduli1} the syzygy bundle $E_{N,d,n}$ is unobstructed and its isomorphism class belongs to a generically smooth irreducible component of the moduli space~$M_{N,d,n}$, of dimension ${n\tbinom{d+N}{N}-n^2}$, if $N\ge3$, and ${n\tbinom{d+2}{2}+n\tbinom{d-1}{2}-n^2}$, if ${N=2}$;
  \item\label{moduli2} if $N\ge 3$, then the closure of the syzygy locus $S_{N,d,n}$ is an irreducible component of $M_{N,d,n}$; if ${N=2}$, the closure of $S_{2,d,n}$ has codimension $n\tbinom{d-1}{2}$ in $M_{2,d,n}$.
\end{enumerate}
\end{te}
\begin{dem}
To establish number \ref{moduli1}, observe that from deformation theory (see proposition~\ref{Har78,4.1}), we know that the Zariski tangent space of~$M_{N,d,n}$ at $\left[E_{N,d,n}\right]$ is canonically given by
\begin{align*}
T_{\left[E_{N,d,n}\right]}M_{N,d,n}&\cong
        \Ext^{1}(E_{N,d,n},E_{N,d,n})\\
    &\cong H^1\left(E_{N,d,n} \otimes
        {E_{N,d,n}}^{\vee}\right)\\
    &\cong H^1\big(\mathcal{E}nd(E_{N,d,n})\big)\mbox{,}
\end{align*}
and the obstruction space of the local ring $\OO_{M,[E_{N,d,n}]}$ is a subspace of
\begin{align*}
\Ext^2(E_{N,d,n},E_{N,d,n})&\cong
        H^2\left(E_{N,d,n}\otimes {E_{N,d,n}}^{\vee}\right)\\
    &\cong H^2\big(\mathcal{E}nd(E_{N,d,n})\big)\mbox{.}
\end{align*}
Thus, if ${H^2\big(\mathcal{E}nd(E_{N,d,n})\big)=0}$, then the moduli space $M_{N,d,n}$ is smooth at the point $\left[E_{N,d,n}\right]$ and in this last case
\[\dim_{K}\Ext^{1}(E_{N,d,n},E_{N,d,n})=
    \dim_{\left[E_{N,d,n}\right]}M_{N,d,n}\]
(see proposition~\ref{Har78,4.1}, but also \cite{Mar77} and \cite{Mar78}).

To compute $\Ext^{i}(E_{N,d,n},E_{N,d,n})$, we consider the exact sequence
\begin{equation}\label{suc1}
\xymatrix{
0\ar[r] &E_{N,d,n}\ar[r] &\OO(-d)^n\ar[r]
    &\OO\ar[r] &0
}
\end{equation}
and its dual
\begin{equation}\label{suc2}
\xymatrix{
0\ar[r] &\OO\ar[r] &\OO(d)^{n}\ar[r]
    &{E_{N,d,n}}^{\vee}\ar[r] &0.
}
\end{equation}

From the exact sequence (\ref{suc1}), we get the the cohomological exact sequence
\begin{equation*}
\begin{split}
\xymatrix@R=0pt{
0\ar[r] &H^0(E_{N,d,n})\ar[r] &H^0\big(\OO(-d)\big)^n
    \ar[r] &H^0\OO\ar[r] &\\
\ar[r] &H^1(E_{N,d,n})\ar[r] &H^1\big(\OO(-d)\big)^n\ar[r]
    &H^1\OO\ar[r] &\\
\ar[r] &H^2(E_{N,d,n})\ar[r] &H^2\big(\OO(-d)\big)^n\ar[r]
    &H^2\OO\ar[r] &\\
\ar[r] &H^3(E_{N,d,n})\ar[r] &H^3\big(\OO(-d)\big)^n\ar[r]
    &H^3\OO\ar[r] &\cdots.
}
\end{split}
\end{equation*}
Now, we know that $h^0\OO=1$, and for ${i>0}$, $h^i\OO=0$. We also know that
\[h^0\big(\OO(-d)\big)=h^1\big(\OO(-d)\big)=
    h^3\big(\OO(-d)\big)=0\mbox{,}\]
and that $h^2\big(\OO(-d)\big)=\tbinom{d-1}{2}$, if ${N=2}$, and $h^2\big(\OO(-d)\big)=0$, if ${N\ge4}$. From here we get
\begin{equation} \label{coho1}
\begin{array}{l}
h^0(E_{N,d,n})=0\mbox{;} \\
h^1(E_{N,d,n})=1\mbox{;} \\
h^2(E_{N,d,n})=
\left\{
    \begin{array}{ll}
        n\tbinom{d-1}{2}\mbox{,} & \mbox{if} \quad N=2\mbox{;}\\
        0\mbox{,} & \mbox{if} \quad N\ge3\mbox{;}
    \end{array}
\right.\\
h^3(E_{N,d,n})=
\left\{
    \begin{array}{ll}
        n\tbinom{d-1}{3}\mbox{,} & \mbox{if} \quad N=3\mbox{;}\\
        0\mbox{,} & \mbox{if} \quad N\ne3\mbox{.}
    \end{array}
\right.
\end{array}
\end{equation}

Denote by ${F:=E_{N,d,n}\otimes {E_{N,d,n}}^{\vee}}$. If we twist the sequence (\ref{suc2}) by $\otimes E_{N,d,n}$, we get
\[
\xymatrix{
0\ar[r] &E_{N,d,n}\ar[r] &E_{N,d,n}(d)^{n}\ar[r] &F\ar[r] &0.
}
\]
The corresponding cohomological exact sequence is
\begin{equation}\label{ext2}
\begin{split}
\xymatrix@R=0pt{
0\ar[r] &H^0(E_{N,d,n})\ar[r] &H^0\big(E_{N,d,n}(d)\big)^n
    \ar[r] &H^0(F)\ar[r] &\\
\ar[r] &H^1(E_{N,d,n})\ar[r] &H^1\big(E_{N,d,n}(d)\big)^n\ar[r]
    &H^1(F)\ar[r] &\\
\ar[r] &H^2(E_{N,d,n})\ar[r] &H^2\big(E_{N,d,n}(d)\big)^n\ar[r]
    &H^2(F)\ar[r] &\\
\ar[r] &H^3(E_{N,d,n})\ar[r] &H^3\big(E_{N,d,n}(d)\big)^n\ar[r]
    &H^3(F)\ar[r] &\cdots.
}
\end{split}
\end{equation}
Since $E_{N,d,n}$ is stable, it is simple, i.e.\ $H^0(F)=\kk$. Thus, from the exact sequence (\ref{ext2}), and the fact that by (\ref{coho1}), ${h^0(E_{N,d,n})=0}$, we get ${h^0\big(E_{N,d,n}(d)\big)=0}$.

Twisting by $\OO(d)$ the exact sequence (\ref{suc1}), and taking cohomology, we deduce
\begin{equation} \label{coho2}
\begin{split}
    &h^1\big(E_{N,d,n}(d)\big)=\tbinom{N+d}{d}-n\mbox{,} \\
    &h^2\big(E_{N,d,n}(d)\big)=0\mbox{,} \\
    &h^3\big(E_{N,d,n}(d)\big)=0\mbox{.}
\end{split}
\end{equation}

In particular, from (\ref{ext2}) we get
\[H^2(F)\cong\Ext^2(E_{N,d,n},E_{N,d,n})=
    \left\{
    \begin{array}{ll}
        n\tbinom{d-1}{3}\mbox{,} & \mbox{if} \quad N=3\mbox{;}\\
        0\mbox{,} & \mbox{if} \quad N\ne3\mbox{.}
    \end{array}
    \right.
    \]
This automatically gives us that if ${N\ne3}$, then the syzygy bundle $E_{N,d,n}$ is unobstructed.

From the exact sequence
\[
\xymatrix@C-.5em{
0\ar[r] &\kk\ar[r] &\kk\ar[r] &H^1\big(E_{N,d,n}(d)\big)^n\ar[r]
    &H^1(F)\ar[r] &H^2(E_{N,d,n})\ar[r] &0\mbox{,}
}
\]
we conclude
\[ h^1(F) = ext^1(E_{N,d,n},E_{N,d,n})= \left\{
    \begin{array}{ll}
        n\tbinom{d+2}{2}+n\tbinom{d-1}{2}-n^2\mbox{,} & \mbox{if}
            \quad N=2\mbox{,}\\
        n\tbinom{N+d}{d}-n^2\mbox{,} & \mbox{if} \quad N \ge4\mbox{.}
   \end{array}
\right.
\]
Now, for ${N=3}$, we have
\begin{align*}
n\tbinom{d+3}{3}-n^2&=\dim S_{3,d,n}\\
    &\le\dim_{[E_{3,d,n}]}M_{3,d,n}\\
    &\le\dim T_{[E_{3,d,n}]}M_{3,d,n}\\
    &=n\tbinom{d+3}{3}-n^2\mbox{,}
\end{align*}
and therefore
\[\dim S_{3,d,n}=\dim T_{[E_{3,d,n}]}M_{3,d,n}.\]
So we conclude that also in the case ${N=3}$ the syzygy bundle $E_{3,d,n}$ is unobstructed.

\smallskip

Finally, number~\ref{moduli2} easily follows by comparing the dimension of $S_{N,d,n}$ to the dimension of the irreducible component of $M_{N,d,n}$ passing through the isomorphism class of a syzygy bundle $E_{N,d,n}$. Note that $M_{2,d,n}$ is irreducible and therefore we can talk about the codimension of $S_{2,d,n}$ in $M_{2,d,n}$ (see \cite{Ell83}, or \cite{Bar77}, \cite{Hul80}, \cite{Hul79}, and \cite{ES81} for particular cases).
\end{dem}


\renewcommand{\starchapter}{Resum en català}
\chapter*{\starchapter\markboth{\starchapter}{\starchapter}}
\addcontentsline{toc}{chapter}{\starchapter}

\selectlanguage{catalan}

Determinar si un fibrat de sizígies sobre $\pp^N$ és estable, o semiestable, és un problema amb una llarga història en geometria algebraica. Està estretament relacionat amb el problema de trobar la resolució lliure minimal de l'anell de coordenades de la varietat definida per una família de polinomis homogenis genèrics ${f_1,\ldots,f_n}$ en $\kk[X_0,\ldots,X_N]$. Aquest problema data almenys dels anys vuitanta, quan Fröberg l'estudia al seu article \cite{Fro85} i troba una estimació per a un minorant de la sèrie de Hilbert d'aquell anell en termes dels graus dels polinomis ${f_1,\ldots,f_n}$.

Un fibrat de sizígies és, per definició, el nucli d'un epimorfisme
\[\xymatrix{{\bigoplus\limits_{i=1}^n\OO_{\pp^N}(-d_i)}
    \ar[rr]^-{f_1,\ldots,f_n}&&\OO_{\pp^N}\mbox{,}}\]
donat per ${(g_1,\ldots,g_n)\mapsto f_1g_1+\cdots+f_ng_n}$, on ${f_1,\ldots,f_n}$ són polinomis homogenis en $\kk[X_0,\ldots,X_N]$ de graus ${d_1,\ldots,d_n}$, respectivament, tals que l'ideal  ${(f_1,\ldots,f_n)}$ és \mbox{$\mathfrak{m}$-primari}, amb
${\mathfrak{m}=(X_0,\ldots,X_N)}$.

En aquesta tesi, considerem el cas de fibrats de sizígies definits per formes genèriques ${f_1,\ldots,f_n}$ d'un mateix grau~$d$, i demostrem la seva estabilitat i no obstrucció per a ${N\ge2}$, excepte en el cas ${(N,d,n)=(2,2,5)}$, on només la semiestabilitat està garantida. Per dur a terme aquesta tasca, ens restringirem primer al cas de monomis i en traurem conseqüències per al cas de formes genèriques. Per això, l'objectiu principal d'aquesta tesi és donar una resposta completa al problema següent:
\begin{prcat*}[\ref{6.9}]
Existeix per a cada~$d$ i cada $n\le\tbinom{d+N}{N}$ una família de~$n$ monomis en $\kk\left[X_0,\ldots,X_N\right]$ de grau~$d$ tal que el seu fibrat de sizígies és semiestable?
\end{prcat*}
L'estabilitat és una propietat oberta per a la topologia de Zariski, com la semiestabilitat. Per tant, si aquest problema té una resposta positiva, o millor encara, si trobem una família de monomis d'un grau fix~$d$ tal que el seu fibrat de sizígies és estable, llavors sabrem que una família de formes genèriques ${f_1,\ldots,f_n}$ de grau~$d$ també dóna lloc a un fibrat de sizígies estable.

El problema~\ref{6.9} ha estat presentat per Brenner a \cite{Bre08b}, article on podem trobar una interpretació de la \textit{tight closure} en termes de fibrats vectorials. La \textit{tight closure} és una tècnica de la teoria d'anells en característica positiva introduïda per Hochster i Huneke a finals dels vuitanta (vegeu \cite{HH88} i \cite{HH90}). És una operació sobre ideals que contenen cossos, o submòduls d'un mòdul donat. En el cas d'ideals de paràmetres, l'existència d'àlgebres Cohen\h{Macaulay} grans és inherent a l'estudi de les seves \textit{tight closures} \cite{HH92}. També ens donen informació sobre singularitats racionals \cite{Smi97a}, el teorema d'anu{\lgem}ació de Kodaira \cite{HS97} i el teorema de Briançon\h{Skoda} (vegeu \cite{HH90}, \cite{HH94} i \cite{Smi97b}). Una descripció d'aquests resultats es pot trobar a \cite{Hun98}. El 1994 Hochster va introduir la \textit{solid closure}, una generalització de la \textit{tight closure}, i va caracteritzar la \textit{tight closure} d'un ideal ${I=(f_1,\ldots,f_n)}$ en un anell~$R$ en termes de la cohomologia local de les \textit{forcing algebras}
\[\frac{R[T_1,\ldots,T_n]}{(f_1T_1+\cdots+f_nT_n+f)}\mbox{.}\]
Aquesta caracterització ha estat descrita per Brenner a \cite[teorema~2.3]{Bre08b}, i és una conseqüència  del coro{\lgem}ari~2.4, de la proposició~5.3 i dels teoremes~8.5 i~8.6 a \cite{Hoc94}. Brenner va associar a aquestes \textit{forcing algebras} un fibrat vectorial, el fibrat de sizígies dels polinomis ${f_1,\ldots,f_n}$, denotat per ${\syz(f_1,\ldots,f_n)}$, i que és la restricció de
\[\mathrm{Spec}\
    \frac{R[T_1,\ldots,T_n]}{(f_1T_1+\cdots+f_nT_n)}\]
a l'obert ${U:=D(f_1,\ldots,f_n)}$, dels punts que no són zeros comuns de ${f_1,\ldots,f_n}$ (vegeu \cite[definició~2.12]{Bre08b}). Aquest fibrat és a la successió exacta
\[\xymatrix{0\ar[r]&\syz(f_1,\ldots,f_n)\ar[r]&
    {{\OO_U}^n}\ar[rr]^{f_1,\ldots,f_n}&&\OO_U\ar[r]&0\mbox{,}}\]
a la qual diem la seva successió de presentació. Per a anells graduats de dimensió dos, alguns problemes oberts sobre la \textit{tight closure} han estat resolts fent servir aquests fibrats, com ara el fet de que la \textit{tight closure} és el mateix que la \textit{plus closure} si el cos de base és finit \cite{Bre06b}, i que la multiplicitat de Hilbert\h{Kunz} és un nombre racional \cite{Bre06a}.

\bigskip

Hi ha d'altres raons importants per estudiar fibrats vectorials sobre varietats algebraiques. Entre elles, que donen informació sobre les respectives varietats i els espais de moduli a què donen lloc, els quals esdevenen exemples interessants de varietats de dimensions elevades.

Els espais de moduli de fibrats vectorials són esquemes que parametritzen famílies de classes d'isomorfia de fibrats vectorials. La seva existència depèn de la forma com està escollida cadascuna d'aquestes famílies. El conjunt de totes les classes d'isomorfia de fibrats vectorials sobre una varietat~$X$ és en general massa gran perquè es pugui parametritzar. Encara que imposem restriccions tals fixar el rang donat i les classes de Chern, no hi ha esperança de trobar un esquema de tipus finit que parametritzi aquesta família de classes d'isomorfia de fibrats vectorials. Tanmateix, si imposem una restricció addicional, l'estabilitat, hi ha una forma natural de fer aquesta parametrització, com va demostrar Maruyama en \cite{Mar76}. Mumford va introduir la noció d'estabilitat per als fibrats vectorials sobre corbes \cite{Mum62} precisament per resoldre aquest problema. Takemoto la va generalitzar després a fibrats vectorials sobre superfícies en \cite{Tak72} i \cite{Tak73}. La noció d'estabilitat coneguda avui com estabilitat de Mumford\h{Takemoto} o \mbox{$\mu$-estabilitat} és una generalització d'aquesta per als feixos lliures de torsió sobre varietats de qualsevol dimensió.

\bigskip

Com hem dit abans, aquest treball està dedicat a l'estudi dels fibrats de sizígies sobre l'espai projectiu. Aquests fibrats són en particular feixos de sizígies, els quals estan definits de forma similar com al nucli d'un epimorfisme
\[\xymatrix{{\bigoplus\limits_{i=1}^n\OO_{\pp^N}(-d_i)}
    \ar[rr]^-{f_1,\ldots,f_n}&&\OO_{\pp^N}\mbox{,}}\]
donat per ${(g_1,\ldots,g_n)\mapsto f_1g_1+\cdots+f_ng_n}$, on ${f_1,\ldots,f_n}$ són polinomis homogenis en $\kk[X_0,\ldots,X_N]$ de graus ${d_1,\ldots,d_n}$, respectivament. La diferència aquí és que per a la definició de feixos de sizígies, l'ideal ${(f_1,\ldots,f_n)}$ no ha de ser \mbox{$\mathfrak{m}$-primari}. Un feix de sizígies és localment lliure sobre l'obert ${D(f_1,\ldots,f_n)}$, i aquesta és la raó per la qual si volem que aquests feixos siguin fibrats vectorials sobre tot l'espai projectiu, els polinomis ${f_1,\ldots,f_n}$ no poden tenir haver zeros comuns, el que és el mateix que dir que l'ideal ${(f_1,\ldots,f_n)}$ ha de ser \mbox{$\mathfrak{m}$-primari}. Aquesta última condició implica que ${n\ge N+1}$, ja que un ideal \mbox{$\mathfrak{m}$-primari} no pot tenir un nombre de generadors inferior al nombre de variables \cite[teorema~11.14]{AM69}.

Per veure un exemple molt senzill i ben conegut d'un fibrat de sizígies, considereu el fibrat cotangent $\Omega_{\pp^N}^1$, ja que és el nucli de
\[\xymatrix{\OO_{\pp^N}(-1)^{N+1}\ar[rr]^-{X_0,\ldots,X_N}&&
    \OO_{\pp^N}\mbox{.}}\]

Fins ara, molt poc era conegut sobre l'estabilitat dels fibrats de sizígies. Se sap que el fibrat cotangent és estable; en característica zero, si en lloc de formes lineals prenem ${N+1}$ polinomis homogenis d'un grau donat~$d$, també sabem que el fibrat de sizígies corresponent és estable, gràcies a un resultat de \mbox{Bohnhorst} i Spindler \cite{BS92}. Encara en el context de característica zero, un altre fibrat de sizígies estable és el que podem obtenir a partir d'una família linealment independent de $\tbinom{d+N}{N}$ polinomis homogenis de grau~$d$. Flenner en va demostrar la semiestabilitat el 1984 \cite{Fle84}, i Ballico en va demostrar el 1992 l'estabilitat \cite{Bal92}. El 2008 Hein va demostrar la semiestabilitat d'un fibrat de sizígies definit per una família de~$n$ polinomis homogenis genèrics, amb ${2\le n\le d(N+1)}$, en qualsevol característica
\cite[apèndix~A]{Bre08a}.

També el 2008 Brenner va presentar la condició suficient següent perquè un fibrat de sizígies sigui estable, o semiestable (vegeu \cite{Bre08a} o \cite{Bre08b}).

\begin{cocat*}[\ref{cor6.6}]
Siguin $f_i$, amb $i\in I$, monomis de l'anell ${\kk[X_0,\ldots,X_N]}$ de graus $d_i$, tals que l'ideal ${(f_i,i\in I)}$ sigui \mbox{$\mathfrak{m}$-primari}. Suposem que per a cada subconjunt ${J\subseteq I}$, amb ${|J|\ge2}$, la desigualtat
\[\frac{d_J-\sum_{i\in J}d_i}{|J|-1}
    \le\frac{-\sum_{i\in I}d_i}{|I|-1}\]
se satisfà, on $d_J$ és el grau del màxim comú divisor de la subfamília $\{f_i:i\in J\}$. Llavors el fibrat de sizígies ${\syz(f_i,i\in I)}$ és semiestable (i estable si la desigualtat és sempre estricta).
\end{cocat*}

La desigualtat que s'ha de satisfer en aquest coro{\lgem}ari depèn només de les cardinalitats dels conjunts esmentats i dels graus dels monomis. Això redueix el problema de decidir si un fibrat de sizígies donat és estable (o semiestable) a un conjunt finit de càlculs. És clar que, com que el nombre de subconjunts d'una família creix exponencialment amb la seva cardinalitat, el nombre de càlculs a fer pot tornar-se bastant feixuc. Per tant, qualsevol estratègia d'elecció de la família que pugui reduir el nombre de càlculs pot ser útil.

Quan va presentar aquest resultat, Brenner va proposar alhora el problema~\ref{6.9}, esmentat abans. Els capítols \ref{t2} i \ref{t3} d'aquesta tesi estan dedicats a respondre aquesta qüestió per a ${N\ge2}$. Hi podem veure que la resposta és afirmativa en tots els casos. De fet, els resultats obtinguts són més forts. Amb l'excepció del cas ${(N,d,n)=(2,5,5)}$, hem trobat en cada cas una família de $n$ monomis tal que el seu fibrat de sizígies és estable.

Per fer\h{se} una idea de què hem fet per resoldre aquest problema, noteu en primer lloc que com que només considerem monomis d'un mateix grau~$d$, la desigualtat del coro{\lgem}ari~\ref{cor6.6} esdevé simplement
\[\frac{d_J-kd}{k-1}\le\frac{-nd}{n-1}\mbox{,}\]
on ${n:=|I|}$ i ${k:=|J|}$. Si ${d_J=0}$, aquesta desigualtat és certa, ja que tenim sempre que ${k\le n}$ i la successió $\left(a_{d,j}\right)_{j\ge2}$, amb $a_{d,j}:=-\tfrac{jd}{j-1}$, és creixent. Per tant, si volem trobar un conjunt $I$ que satisfaci aquesta desigualtat, hem de garantir que per a cada subconjunt~${J\subseteq I}$, el grau~$d_J$ del màxim comú divisor dels seus elements sigui prou baix.

Com es pot veure amb detall al capítol~\ref{t2}, el conjunt dels monomis en tres variables d'un grau donat pot ser representat convenientment en forma de triangle, com en l'exemple següent, per a grau tres.
{\small%
\[
\xymatrix@C=\meiabasex@R=\meiaalturax@!0{
&&&{X_2}^3\\
\\
&&{X_0}{X_2}^2&&{X_1}{X_2}^2\\
\\
&{X_0}^2X_2&&{X_0}X_1X_2&&{X_1}^2X_2\\
\\
{X_0}^3&&{X_0}^2X_1&&{X_0}{X_1}^2&&{X_1}^3\\
}
\]
}%
D'aquesta manera, quant més propers són dos monomis, més alt és el grau del seu màxim comú divisor. D'aquí podem veure que per trobar una família de monomis tal que el seu fibrat de sizígies sigui estable, els monomis han de ser escollits de forma que estiguin prou apartats entre ells en aquesta representació.

Diversos mètodes han estat adoptats amb aquest propòsit. Per als valors més baixos de~$n$ (fins a $18$), s'ha trobat una solució particular per a cada cas. Després, per a ${18<n\le d+2}$, hem buscat el nombre triangular~$T$ més alt que no fos més gran que~$n$, i hem escollit $T$ monomis disposats de la forma més uniforme possible, posant els $n-T$ monomis que faltaven en llocs apropiats als costats del triangle. Per a ${d+2<n\le3d}$, tots els monomis han estat escollits als costats del triangle, i de ${3d+1}$ endavant, el seu interior s'ha omplert, començant per les capes més properes dels costats.

En quatre o més variables, són possibles representacions anàlogues, amb la forma d'un tetraedre, o d'un hipertetraedre, respectivament. Aquests casos han estat resolts d'una manera més senzilla. Informalment, aquest problema esdevé una qüestió de trobar espai suficient per distribuir monomis en aquests triangles, tetraedres i hipertetraedres. El fet que en dimensions superiors hi hagi més espai no déu ser una gran sorpresa.

Per als primers valors de~$n$, la solució ha estat considerar una família de monomis en ${\kk[X_0,\ldots,X_{N-1}]}$ tal que el seu fibrat de sizígies sobre $\pp^{N-1}$ fos estable i afegir\h{hi} el monomi ${X_N}^d$. Això correspon a escollir monomis convenients en una cara de l'hipertetraedre i afegir\h{hi} el vèrtex oposat. Després, anàlogament al cas de tres variables, les cares de l'hipertetraedre han estat omplertes. Finalment, seguint una idea de Brenner, l'interior de l'hipertetraedre ha estat omplert amb una família de monomis, obtinguda a partir d'un conjunt de monomis de grau inferior que donaven lloc a un fibrat de sizígies estable.

D'aquesta forma, hem aconseguit donar una resposta completa al problema posat per Brenner. A més, el fet que l'estabilitat sigui una condició oberta ens ha permès concloure que un fibrat de sizígies donat per una família de~$n$ polinomis homogenis genèrics en ${\kk[X_0,\ldots,X_N]}$ d'un grau~$d$ donat és estable (excepte en el cas ${(N,d,n)=(2,2,5)}$, on només la semiestabilitat està garantida).

\bigskip

L'estabilitat d'un fibrat de sizígies associat a una família de polinomis homogenis genèrics d'un grau donat ens dóna informació sobre l'espai de moduli que correspon a la seva classe d'isomorfia. Se sap molt poc sobre els espais de moduli de fibrats vectorials estables en general. Un cop tenim la seva existència establerta, és natural preguntar\h{se} què es pot dir sobre la seva estructura, tant localment com globalment. En particular, seria interessant saber si és connex, irreductible, racional o llis, i com es caracteritza com a espai topològic.

Pel que respecta als espais de moduli de fibrats de sizígies, hem demostrat que un fibrat de sizígies estable és no obstruït i la seva classe d'isomorfisme pertany a una component irreductible genèricament llisa de l'espai de moduli, i hem donat la dimensió d'aquesta component. També hem determinat els casos en els quals la clausura de l'estrat de l'espai de moduli corresponent a fibrats de sizígies estables és tota la component irreductible en la que està contingut, i hem calculat la seva codimensió en el cas que no ho sigui.

\bigskip

Farem ara una descripció amb més detall de l'estructura d'aquesta tesi, i en destacarem els resultats principals.

\smallskip

En el \textbf{capítol~\ref{t1}} donem les definicions dels objectes estudiats en aquest treball, i algunes eines que seran utilitzades en els capítols següents. Comencem per definir fibrats i feixos de sizígies en la secció~\ref{t1syzygy}. En la secció~\ref{t1stability}, definim la noció d'estabilitat per als feixos coherents i discutim algunes de les seves propietats. Després fem una compilació de les contribucions fetes fins ara a l'estudi de l'estabilitat dels fibrats de sizígies. En la secció~\ref{t1segment}, donem una solució al problema~\ref{6.9} per al cas de monomis en dues variables. Finalment, en la secció~\ref{t1moduli}, presentem una definició formal d'espai de moduli i un resultat que dóna condicions suficients perquè sigui no singular en un punt donat.

\smallskip

En el \textbf{capítol~\ref{t2}} responem al problema~\ref{6.9} per a ${N=2}$. L'estabilitat està garantida en tots els casos excepte per a  ${(N,d,n)=(2,2,5)}$, on només s'obté la semiestabilitat. El capítol està dividit en diferents seccions, d'acord amb els diferents valors de~$n$. Recordeu que, per a ${N=2}$, tenim ${3\le n\le\tbinom{d+2}{2}}$. Per als primers casos, amb ${3\le n\le18}$, com hem dit abans, hem trobat una solució particular per a cada valor de~$n$; aquestes solucions estan descrites en la secció~\ref{t2first}. Com hem vist, els monomis en tres variables d'un grau donat poden ser representats en un triangle d'una forma convenient per abordar aquest problema. Per això, la clau per resoldre els casos  ${18<n\le d+2}$, encara a la secció~\ref{t2first}, ha estat  escollir una bona posició per a~$T$ monomis en aquest triangle, on $T$ és el nombre triangular més gran no superior a~$n$. Els casos ${d+2<n\le3d}$ han estat resolts a la secció~\ref{t2sides} prenent un triangle amb un costat ple i omplint els altres dos costats. Els últims casos, a la secció~\ref{t2interior}, per a ${3d<n\le\tbinom{d+2}{2}}$, han estat resolts prenent un triangle amb tots els costats plens i omplint el seu interior de forma convenient. Finalment, a la secció~\ref{t2main}, ajuntant tots els resultats d'aquest capítol, es presenta el teorema següent.

\begin{tecat*}[\ref{main2}]
Siguin~$d$ i~$n$ enters tals que $3\le n\le\tbinom{d+2}{2}$ i $(d,n)\ne(2,5)$. Aleshores hi ha una família de~$n$ monomis en $K\left[X_0,X_1,X_2\right]$ de grau~$d$ tal que el seu fibrat de sizígies és estable. Si $(d,n)=(2,5)$, hi ha una família de~$5$ monomis tal que el seu fibrat de sizígies és semiestable.
\end{tecat*}

\smallskip

En el \textbf{capítol~\ref{t3}} es presenta una resposta al problema~\ref{6.9} per a ${N\ge3}$. En cada cas, hi ha una família de monomis tal que el seu fibrat de sizígies és estable.

En general, els monomis en $\kk[X_0,\ldots,X_N]$ d'un grau~$d$ donat poden ser representats en un hipertetraedre, de forma anàloga als triangles que hem vist en el cas ${N=2}$. Diem \emph{cara $i$} d'aquest hipertetraedre al conjunt de monomis on la variable $X_i$ no apareix.

Aquest capítol està dividit en diferents seccions, tal com l'anterior, d'acord amb els diferents valors de~$n$. Recordeu que tenim ${N+1\le n\le\tbinom{d+N}{N}}$. Per als primers casos, a la secció~\ref{t3first}, amb
\[{N+1\le n\le\tbinom{d+N-1}{N-1}+1}\mbox{,}\]
cada família de ${n-1}$ monomis en $\kk[X_0,\ldots,X_{N-1}]$ tal que el seu fibrat de sizígies sobre $\pp^{N-1}$ és estable genera a una família de~$n$ monomis en $\kk[X_0,\ldots,X_N]$ tal que el seu fibrat de sizígies sobre $\pp^N$ també ho és. Els casos
\[{\tbinom{d+N-1}{N-1}+1<n\le\tbinom{d+N}{N}-\tbinom{d-1}{N}}\]
han estat resolts a la secció~\ref{t3faces} prenent la cara $N$ i el vèrtex ${X_N}^d$ de l'hipertetraedre, i afegint\h{hi} monomis en les altres cares. Prendre el conjunt de tots les cares de l'hipertetraedre i afegir\h{hi} els monomis del seu interior que són més propers als vèrtexs ens ha donat una solució per als casos
\[{\tbinom{d+N}{N}-\tbinom{d-1}{N}<n\le\tbinom{d+N}{N}
    -\tbinom{d-1}{N}+N+1}\mbox{,}\]
que estan tractats a l'inici de la secció~\ref{t3interior}. Els últims casos, amb
\[{\tbinom{d+N}{N}-\tbinom{d-1}{N}+N+1<n\le\tbinom{d+N}{N}}
    \mbox{,}\]
també a la secció~\ref{t3interior}, han estat resolts prenent una família de monomis de grau inferior tal que el seu fibrat de sizígies és estable, multiplicant-los per ${X_0\cdots X_N}$, i afegint\h{hi} tots els monomis de les cares de l'hipertetraedre. Això és una generalització d'un lema de Brenner, fet per al cas ${N=2}$ en les seves notes \cite{Bre}, que amablement ha compartit. Finalment, ajuntant la informació d'aquest capítol i de l'anterior, es presenta en la secció~\ref{t3main} el teorema següent, que és el principal resultat d'aquest treball.
\begin{tecat*}[\ref{main3}]
Siguin~$N$, $d$ i~$n$ enters tals que ${N\ge2}$, ${N+1\le n\le\tbinom{d+N}{N}}$ i\linebreak[4] ${(N,d,n)\ne(2,2,5)}$. Llavors hi ha una família de~$n$ monomis en $\kk\left[X_0,\ldots,X_N\right]$ de grau~$d$ tal que el seu fibrat de sizígies és estable.

Per a  ${(N,d,n)=(2,2,5)}$, hi ha una família de~$5$ monomis de grau~$2$ en l'anell $\kk\left[X_0,X_1,X_2\right]$ tal que el seu fibrat de sizígies és semiestable.
\end{tecat*}
Concloem aquest capítol posant dues qüestions obertes que sorgeixen de forma natural, un cop el problema~\ref{6.9} està resolt. La primera és la generalització d'aquest problema considerant polinomis homogenis de grau arbitrari. La segona és la de saber si els fibrats de sizígies d'ordre superior són estables, o semiestables. Aquesta última qüestió està relacionada amb el problema, també obert, de trobar resolucions lliures minimals de formes genèriques.

\smallskip

En el \textbf{capítol~\ref{t4}} presentem tres resultats sobre fibrats de sizígies i els seus espais de moduli. El primer és una conseqüència del fet que les condicions d'estabilitat i semiestabilitat siguin obertes a l'espai de moduli.
\begin{tecat*}[\ref{geral4}]
Siguin $N\ge2$, $d\ge1$ tals que
\[{N+1\le n\le\tbinom{d+N}{N}}\mbox{,}\]
amb ${(N,d,n)\ne(2,2,5)}$. Si ${f_1,\ldots,f_n}$ són formes genèriques de grau~$d$ en l'anell $\kk[X_0,\ldots,X_N]$, tals que l'ideal ${(f_1,\ldots,f_n)}$ és \mbox{$\mathfrak{m}$-primari}, aleshores el fibrat de sizígies ${\syz(f_1,\ldots,f_n)}$ és estable.

Si ${(N,d,n)=(2,2,5)}$ i ${f_1,\ldots,f_5}$ són formes quadràtiques genèriques tals que l'ideal ${(f_1,\ldots,f_5)}$ és \mbox{$\mathfrak{m}$-primari}, llavors el fibrat de sizígies ${\syz(f_1,\ldots,f_5)}$ és semiestable.
\end{tecat*}

Siguin $N$, $d$ i $n$ enters tals que ${N\ge2}$, ${d\ge1}$ i ${N+1\le n\le\tbinom{d+N}{N}}$, amb ${(N,d,n)\ne(2,2,5)}$. Denotem per $M_{N,d,n}$ l'espai de moduli dels fibrats vectorials estables sobre $\pp^N$ de rang ${n-1}$ amb classes de Chern ${c_i=\tbinom{n}{i}(-d)^i}$, per a ${1\le i\le N}$, i per $S_{N,d,n}$ l'estrat d'aquest espai de moduli que correspon als fibrats de sizígies.

El segon resultat d'aquest capítol ens dóna la dimensió de l'estrat $S_{N,d,n}$, en termes de $N$, $d$ i $n$.

\begin{pcat*}[\ref{syzlocus}]
Fixem enters $N$, $d$ i $n$ tals que ${N\ge2}$, ${d\ge1}$, ${N+1\le n\le \tbinom{d+N}{N}}$, i ${(N,d,n)\ne (2,2,5)}$. Aleshores
\[\dim S_{N,d,n}=n\tbinom{d+N}{N}-n^2.\]
\end{pcat*}

El tercer resultat afirma que la component irreductible de l'espai de moduli que correspon a la classe d'isomorfia d'un fibrat de sizígies estable $E_{N,d,n}$ sobre $\pp^N$, amb ${N\ge2}$, és genèricament llisa, i ens dóna la seva dimensió.

\begin{tecat*}[\ref{moduli}]
Siguin $N$, $d$ i $n$ enters tals que ${N\ge2}$, ${d\ge1}$, ${N+1\le n\le\tbinom{d+N}{N}}$, i ${(N,d,n)\ne(2,2,5)}$. Aleshores
\begin{enumerate}
  \item el fibrat de sizígies $E_{N,d,n}$ és no obstruït i la seva classe d'isomorfisme pertany a una component irreductible de l'espai de moduli~$M_{N,d,n}$ genèricament llisa, de dimensió ${n\tbinom{d+N}{N}-n^2}$, si $N\ge3$, i ${n\tbinom{d+2}{2}+n\tbinom{d-1}{2}-n^2}$, si ${N=2}$;
  \item si $N\ge 3$, llavors la clausura del lloc geomètric de sizígies $S_{N,d,n}$ és una component irreductible de $M_{N,d,n}$; si ${N=2}$, la clausura de $S_{N,d,n}$ té codimensió $n\tbinom{d-1}{2}$ en $M_{N,d,n}$.
\end{enumerate}
\end{tecat*}

\selectlanguage{english}

\bibliographystyle{amsalpha}
\cleardoublepage
\addcontentsline{toc}{chapter}{Bibliography}
\bibliography{tese}

\end{document}